\documentclass[11pt,letter]{article}

\usepackage{color}
\usepackage{appendix}
\usepackage{titletoc}

\usepackage{setspace}
\usepackage{bbm}
\usepackage{natbib}
\usepackage{amstext}
\usepackage{amsthm}
\usepackage{epsfig}
\usepackage{pdfpages}
\usepackage{amsmath,amssymb,amsthm}
\usepackage{graphicx}
\usepackage{fancyhdr}
\usepackage{calc}
\usepackage{dsfont}
\usepackage{subfigure}
\usepackage{rotating}
\usepackage{multirow}
\usepackage{amsmath}
\usepackage{lscape}
\usepackage{tabularx}
\usepackage{bm}
\usepackage{caption}
\usepackage{color}
\usepackage{enumerate,array}
\usepackage{verbatim}
\usepackage{fancyhdr}
\usepackage{stackengine}
\usepackage{enumitem}
\usepackage{algorithm}
\usepackage{algpseudocode}
\usepackage{multirow}

\makeatletter
\DeclareRobustCommand\check[1]{{\mathpalette\@widecheck{#1}}}
\def\@widecheck#1#2{%
    \setbox\z@\hbox{\m@th$#1#2$}%
    \setbox\tw@\hbox{\m@th$#1%
       \widehat{%
          \vrule\@width\z@\@height\ht\z@
          \vrule\@height\z@\@width\wd\z@}$}%
    \dp\tw@-\ht\z@
    \@tempdima\ht\z@ \advance\@tempdima2\ht\tw@ \divide\@tempdima\thr@@
    \setbox\tw@\hbox{%
       \raise\@tempdima\hbox{\scalebox{1}[-1]{\lower\@tempdima\box
\tw@}}}%
    {\ooalign{\box\tw@ \cr \box\z@}}}
\makeatother

\usepackage{hyperref}
\definecolor{cornellred}{rgb}{0.7, 0.11, 0.11}
\hypersetup{    
	colorlinks=true,
	linkcolor=blue,
	citecolor=blue,
	urlcolor=blue
}

\renewcommand{\baselinestretch}{1.5}

\usepackage{xcolor}

\numberwithin{equation}{section}

\allowdisplaybreaks

\newtheorem{theorem}{Theorem}[section]
\newtheorem{lemma}[theorem]{Lemma}
\newtheorem{proposition}[theorem]{Proposition}

\newtheorem{corollary}[theorem]{Corollary}

\theoremstyle{definition}

\newtheorem{remark}{Remark}[section]

\newtheorem{example}{Example}[section]

\newtheorem{assumption}{Assumption}

\newcommand{\cov}{{\rm cov}}

\newcommand{\E}{{\mathbb E}}

\newcommand{\I}{\mathcal{I}}

\renewcommand{\H}{\mathcal{H}}

\newcommand{\one}{\mathbbm{1}}

\def\d{{\rm d}}

\def\c{{\rm c}}

\def\trans{^{\mkern-1mu\mathsf{T}\mkern-1mu}}

\def\oracle{{\rm or}}
\def\G{\mathbb G}
\def\D{\mathcal D}
\def\B{\mathcal B}
\def\M{\mathbb M}

\DeclareMathOperator{\sgn}{sgn}
\newcommand{\var}{{\rm var}}
\newcommand{\Var}{{\rm Var}}

\newcommand{\e}{\epsilon}

\newcommand{\EE}{\mathcal{E}}

\newcommand{\C}{\mathcal{C}}

\newcommand{\X}{\mathcal{X}}

\newcommand{\Z}{\mathcal{Z}}
\newcommand{\N}{\mathcal{N}}

\newcommand{\s}{\mathfrak s}
\newcommand{\converged}{\overset{d.}{\longrightarrow}}
\newcommand{\weakconverge}{\rightsquigarrow}

\DeclareMathOperator*{\argmin}{arg\,min}

\newcommand{\diag}{{\rm diag}}

\renewcommand{\l}{\langle}
\renewcommand{\r}{\rangle}

\renewcommand{\phi}{\varphi}

\renewcommand{\tilde}{\widetilde}
\renewcommand{\hat}{\widehat}
\renewcommand{\epsilon}{\varepsilon}
\renewcommand{\P}{{\mathbb P}}

\def\boxit#1{\vbox{\hrule\hbox{\vrule\kern6pt  \vbox{\kern6pt#1\kern6pt}\kern6pt\vrule}\hrule}}

\def\boxit#1{\vbox{\hrule\hbox{\vrule\kern6pt
          \vbox{\kern6pt#1\kern6pt}\kern6pt\vrule}\hrule}}

\newtheorem{innercustomgeneric}{\customgenericname}
\providecommand{\customgenericname}{}
\newcommand{\newcustomtheorem}[2]{%
  \newenvironment{#1}[1]
  {%
   \renewcommand\customgenericname{#2}%
   \renewcommand\theinnercustomgeneric{##1}%
   \innercustomgeneric
  }
  {\endinnercustomgeneric}
}
\newcustomtheorem{customthm}{Assumption}
\newcustomtheorem{customcor}{Corollary}

\renewcommand{\baselinestretch}{1.4}

\makeatletter
\def\singlespace{\deltaf\baselinestretch{1}\@normalsize}

  \renewenvironment{thebibliography}[1]{%
    \begin{oldthebibliography}{#1}%
      \setlength{\parskip}{0.3ex}%
      \setlength{\itemsep}{0ex}%
  }%
  {%
    \end{oldthebibliography}%
  }

\begin{document}


\baselineskip=17pt

\begin{center}
{\bf \Large Simultaneous 
semiparametric inference for single-index models} \\[1.5em]

Jiajun Tang, Holger Dette\\[1.5em]
Department of Statistics, Harvard University, Cambridge, USA\\
Fakult\"at f\"ur Mathematik, Ruhr-Universit\"at Bochum, Bochum, Germany
\end{center}

\vspace{.5cm}

\begin{quotation}

\noindent {\bf Abstract:} 
In the common  partially linear single-index model we establish a Bahadur representation for a smoothing spline  estimator of  all model parameters and use this result  to prove the joint weak convergence of the   estimator of  the  index link function at a given point, together with the estimators  of the parametric regression coefficients. We obtain the surprising result that, despite of the nature of single-index models where the link function is evaluated at a linear combination of the index-coefficients, the estimator of the 
link function and   the estimator of the  index-coefficients are asymptotically independent.
Our  approach leverages a delicate analysis based on reproducing kernel Hilbert space and empirical process theory.

We  show that the smoothing spline estimator achieves the minimax optimal rate with respect to  the $L^2$-risk and 
consider several statistical applications where joint inference on all model parameters is of interest. In particular, we develop a simultaneous confidence band for the link function and propose inference tools  to investigate if  the maximum absolute  deviation between the (unknown)   link function and a given function exceeds a given  threshold. We also  construct tests for joint hypotheses regarding model parameters which involve both the nonparametric and parametric components and propose novel multiplier bootstrap procedures to avoid the estimation of unknown asymptotic quantities.

\medskip

\medskip


\medskip

\medskip


\end{quotation}

\clearpage



\section{Introduction}

The 
\textit{partially linear single-index model} 
\begin{align}\label{model}
Y=g_0(X\trans\beta_0)+Z\trans\gamma_0+\e,
\end{align}
with link function $g_0$, $p$-dimensional parameter $\beta_0$ and $q$-dimensional parameter $\gamma_0$
has found  considerable attention in the literature of statistics, econometrics, biometrics, and many other scientific fields due to its interpretability of the linear component and its flexibility due to the nonparametric linkage. Estimating the index direction $\beta_0$ (after postulating an appropriate identifiability condition)
and the index link function $g_0$ are of particular interest. In the broad scale of regression analysis, estimating $\beta_0$ falls within the scope of  dimension reduction 
and the methods for estimating   the index link function $g_0$ nonparametrically can be  divided in two classes.
 The first category consists in local approaches, including kernel smoothing \citep{ichimura1993semiparametric,hardle1993,carroll1997,xia1999single,huh2002likelihood}, the estimating equation approach \citep{cui2011}, profile likelihood \citep{liang2010estimation}, average derivative \citep{hardle1989investigating,stoker1986consistent}, minimum average variance estimation \citep{xia2007constructive,wang2010}, empirical likelihood \citep{zhu2006empirical} among many others. The second category is a series expansion approach, where, among others, polynomial spline \citep{wang2009spline}, orthogonal series expansion \citep{dong2015,dong2016}, penalized RKHS regression \citep{yu2002,kuchibhotla2020}, Baysian estimation \citep{antoniadis2004bayesian}, and profile-likelihood based approaches \citep{yu2017penalised} were proposed.

Statistical inference regarding the  parametric components was developed by \cite{gao1997statistical,liang2010estimation,dong2015,dong2016, cui2011,wang2010,kuchibhotla2020} among others. Semiparametric efficient scores for the parametric components $\beta_0$ or $\gamma_0$ were established via the nuisance tangent space by  treating the true index link function $g_0$ as a nuisance parameter \cite[see, for example,][among others]{newey1993efficiency,ma2013}. 
 For statistical inference regarding the  index link function $g_0$,  \cite{carroll1997,xia1999} and \cite{cui2011} derived point-wise asymptotic normality  of  estimators based on local approaches, while 
 \cite{hardle1998testing,hardle2001bootstrap} developed
 bootstrap procedures  for specification testing of $g_0$. 
 The asymptotic properties of smoothing spline estimators  in  the single-index model were investigated by \cite{wang2009spline} and \cite{kuchibhotla2020} who also established  the asymptotic efficiency of the  corresponding   index direction estimator.

A common feature in all these references  is its focus  on \textit{marginal} (asymptotic) properties of estimators 
either referring to the 
 estimators of the parametric components $(\beta_0,\gamma_0)$ or
 to a  nonparametric estimator of the link function $g_0$  at a given point. 
 On the other hand, joint asymptotic  properties of estimators 
 have  been investigated by  \cite{cheng2015} for a partially-linear generalized semi-parametric regression model $\E [Y|X,Z] =F( g_0(X)+Z\trans\gamma_0) $ with a univariate predictor $X$ and a known link function $F$. The fact that the regression function in \eqref{model} contains a further unknown non-linear parameter $\beta_0$  poses substantial  challenges in the theoretical analysis of estimators in  single-index models and to  the best of our knowledge results about joint asymptotic properties and  corresponding inference tools for  all model parameters $(g_0,\beta_0,\gamma_0)$ in   model \eqref{model} do not exist. 
\smallskip

\textbf{Main contributions.} 
In this article we leverage the theory of reproducing kernel Hilbert spaces  (RKHS) to establish the joint weak convergence of  a smoothing spline estimator, say $(\hat g_n,\hat\beta_n,\hat\gamma_n)$
for all parameters in the single-index model \eqref{model}. 
We establish a \textit{joint Bahadur representation} for this estimator, which to our knowledge is the first one in the literature. This useful tool has several theoretical and statistical applications. First, it enables us to derive the joint asymptotic normality of the estimator  $\big[\hat g_n(s),\hat\beta_n\trans,\hat\gamma_n\trans\big]\trans$ as a vector in $\mathbb R^{p+q+1}$ (after appropriate rescaling) and to prove that  $\hat g_n(s)$ and $(\hat\beta_n,\hat\gamma_n)$ are asymptotically independent.  We use this result  to construct consistent tests for the  hypotheses that  the model parameters (including the index function $g_0$) satisfy a structural equation.
Second, we establish that the index link function estimator $\hat g_n$  is (asymptotically)  minimax optimal   with respect to the $L^2$-risk. Third, the joint Bahadur representation enables us to derive the marginal Bahadur representation of $\hat g_n$ and $(\hat\beta_n,\hat\gamma_n)$, respectively.
We use these results to  show that  $\hat\beta_n$ and $\hat\gamma_n$ are $O_p(n^{-1/2})$-consistent and, more importantly,  to construct simultaneous confidence bands for the index link function $g_0$. 

Additionally, we address the important problem  of testing the form of the index link function $g_0$ from a different perspective. In contrast to most of the literature which test for the   exact equality  $H_0: g_0  \equiv g_*$ for a given function  $g_*$
 \citep{hardle1998testing,hardle2001bootstrap} we develop tests for  the hypothesis that the maximum absolute  deviation between $g_0$ and $g_*$ does not exceed a (small) prespecified threshold $ \Delta >0$, that is 
$H_0:   \sup_{s } |g_0(s) -  g_*(s) | \leq   \Delta$. We argue that these hypotheses better reflect the practical requirements, because in applications it is rarely believed that $g_0 $ and $ g_*$ exactly coincide on their complete 
domain, but one would work with   the function $g_*$ in  the single index  \eqref{model} if the deviation between 
 $g_0 $ and $ g_*$ is small.  For example, if $g_* \equiv 0$ this point of view  allows  to investigate if the impact of  
 the single-index $X\trans\beta_0$ on the response is  negligible instead of testing  that there is no
 effect $X\trans\beta_0$ on $Y$.  
 This formulation  of the hypotheses reflects  a point of view mentioned by   \cite{tukey1991}, who  argues in the context of multiple comparisons of means that  
\textit{``All we know about the world teaches us that the effects of A and B are always different -- in some decimal place -- for any A and B. Thus asking `Are the effects different?' is foolish.''}

Furthermore, due to the fact that the joint asymptotic distributions entail unknown quantities which are difficult to estimate in practice, we propose novel multiplier bootstrap procedures for these statistical inference problems and prove their validity.
\smallskip

\textbf{Organization.} In Section~\ref{sec2}, we derive the joint Bahadur representation and joint weak convergence of the smoothing spline estimator.
In Section~\ref{sec:stat} we apply these results to several statistical inference problems. To be specific, Section~\ref{sec:minimax} establishes that the index link estimator is minimax optimal in the $L^2$ sense. Simultaneous confidence bands, inference tools for relevant hypotheses and joint hypotheses are derived in Sections~\ref{sec:band}, \ref{sec:rele}, and \ref{sec:joint}, respectively. Section~\ref{sec:finite} investigates the finite sample properties of the proposed estimators and tests. Finally all proofs are deferred to an online supplement.

\section{Joint asymptotic properties}
\label{sec2}

Suppose that  $(X_1,Y_1,Z_1,\e_1),\ldots,(X_n,Y_n,Z_n,\e_n)$ is a sample of  i.i.d.~copies of $(X,Y,Z,\e)$ generated by the partially linear single-index model \eqref{model}, 
where $\E(\e|X,Z)=\E(\e|X\trans\beta_0,Z)=0$. In model \eqref{model}, $g_0$ is the unknown index link function defined on a
compact interval $\I \subset \mathbb{R} $, and $\beta_0\in\mathbb R^p,\gamma_0\in\mathbb R^q$ are unknown coefficient vectors, where $\beta_0$ is the index direction consisting of the single-index $X\trans\beta_0$. The covariates $X$ and $Z $  in   model \eqref{model}  take values in  compact sets  $\X  \subset \mathbb R^p $ and $\Z  \subset \mathbb R^q $, respectively. 
Note that these  assumptions  are standard  in the literature of nonparametric and semiparametric statistics; see, for example, \cite{hardle1993,mammen1997,kuchibhotla2020,shang2013} among many others.

The identifiability condition for the single-index model \eqref{model} is well understood in the literature: it is assumed that the index link function $g_0$ is non-constant and 
\begin{align*}
\beta_0\in \B:=\big\{\beta\in\mathbb{R}^p:\|\beta\|_2=1,\,\text{the first nonzero element of }\beta~\text{is positive}\big\}
\end{align*}
 \citep[see][among many others]{carroll1997,yu2002,cui2011,kuchibhotla2020}. Note that for  $Z\equiv0$, model \eqref{model} becomes the classical single-index model \citep[][]{ichimura1993semiparametric,hardle1993}.

Following \cite{yu2002,kuchibhotla2020} we consider a penalized least square smoothing spline estimator. To be precise, let $$
\H_m=\big\{g:\I\to\mathbb{R}\,\big|\,g^{(k)}\text{ is absolutely continuous},\,0\leq k\leq m-1\,;g^{(m)}\in L^2(\I)\big\}
$$
denote the Sobolev space on the interval $\I$ of order $m$. The smoothing spline estimator    is given  by
\begin{align}\label{hatwn0}
(\hat g_n,\hat\beta_n,\hat\gamma_n)=\argmin_{(g,\beta ,\gamma) \in\H_m\times\B\times\mathbb R^q}\,\bigg[\frac{1}{2n}\sum_{i=1}^n\big\{Y_i-g(X_i\trans \beta)-Z_i\trans\gamma\big\}^2+\frac{\lambda}{2}\,J(g,g)\bigg],
\end{align}
where $\lambda>0$ is a  regularization parameter,
and $J$ is a smoothness penalty functional defined by
\begin{align}\label{J}
J(g_1,g_2)=\int_\I g_1^{(m)}(s)\, g_2^{(m)}(s)\,\d s,\qquad g_1,g_2\in\H_m.
\end{align}
Here, and throughout this paper  we do not reflect the dependence of the estimator $(\hat g_n,\hat\beta_n,\hat\gamma_n)$ on $\lambda$ in our notations.  In the following discussion we will derive the asymptotic properties 
of this estimator. In particular we will prove joint asymptotic normality and show that the nonparametric component $\hat g_n$ and the parametric estimators $(\hat\beta_n,\hat\gamma_n)$ are asympotically independent.  For this purpose we introduce several useful notations.

\smallskip 

\textbf{Notations.} 
Let ``$\converged$" denote convergence in distribution in $\mathbb{R}^k$. For a square matrix $M\in\mathbb R^{k \times k}$, let $\|M\|_2$ denote its spectral norm, ${\rm tr}(M)$ its trace, and $\lambda_{\min}(M)$ its smallest eigenvalue; for a vector $v\in\mathbb R^p$, let $\|v\|_2$ be its $\ell_2$-norm. For two real sequences  $\{a_n\},\{b_n\}$, we write $a_n\asymp b_n$ and $a_n\lesssim b_n$ to indicate that there exist constants $c_1,c_2>0$ such that $\lim_{n\to\infty}a_n/b_n=c_1$, and $a_n/b_n\leq c_2$ for all $n\in\mathbb N$, respectively. Let $\partial \mathcal{G}$ denote the boundary of the set $\mathcal{G}$. Let $0_p\in\mathbb R^p ,0_{p,q}\in\mathbb R^{p\times q}$ be the vector and matrix whose entries are all zero, respectively, and let $I_p$ denote the $p$-dimensional identity matrix.  Let $\mathbf{1}\{\cdot\}$ denote the indicator function. Define 
\begin{align}\label{Theta}
&\Theta=\H_m\times\mathbb R^{p-1}\times\mathbb R^q=\{(g,\theta,\gamma):g\in\H_m,\,\theta\in\mathbb R^{p-1},\,\gamma\in\mathbb{R}^q\}
\end{align} 
as the parameter space and 
\begin{equation}\label{sigma2}
\sigma^2(X,Z) =\E(\e^2|X,Z),\qquad \sigma_0^2(s)=\E\{\sigma^2(X,Z)|X\trans\beta_0=s\},\qquad s\in\I
\end{equation}
(in particular, we assume the existence of these expectations).

\subsection{Defining a norm on the parameter space $\Theta = \H_m\times\mathbb R^{p-1}\times\mathbb R^q$}\label{sec:definenorm}

Let $w_0=(g_0,0_{p-1},\gamma_0)\in\Theta$.
For $w=(g,\theta,\gamma)\in\Theta$, define 
\begin{align}\label{norm}
\|w\|^2=\|(g,\theta,\gamma)\|^2=\E\big[ \sigma^2(X,Z)\{g(X\trans\beta_0)+g_0'(X\trans\beta_0)X\trans Q_{\beta_0}\theta+Z\trans\gamma\}^2\big]+\lambda J(g,g),
\end{align}
where $\sigma^2(X,Z)$ is defined in \eqref{sigma2}, $J$ is the penalty functional defined in \eqref{J}   and   $Q_{\beta_0}$ is a $p\times(p-1)$ matrix such that $[\beta_0,Q_{\beta_0}]$ is a $p\times p$ orthogonal matrix, that is   $Q_{\beta_0}\trans\beta_0=0_{p-1}$, $Q_{\beta_0}\trans Q_{\beta_0}=I_{p-1}$, and $Q_{\beta_0}Q_{\beta_0}\trans+\beta_0\beta_0\trans =I_p$. 
The intuition of the definition of $\|\cdot\|$ in \eqref{norm} is driven by the   approximation 
\begin{align}\label{ad}
g(x\trans  \beta)
\approx g(x\trans \beta_0)+g_0'(x\trans \beta_0)x\trans  Q_{\beta_0}\theta_{\beta},\qquad\mbox{where }\theta_{\beta}:=Q_{\beta_0}\trans\beta\in\mathbb R^{p-1},
\end{align}
for $(g,\beta)\in\H_m\times\B$ in a neighborhood of the ground truth $(g_0,\beta_0)$, which will be discussed in detail in Lemma~\ref{lem:unionbound2} of the online supplement.
Note that when $(g,\beta)=(g_0,\beta_0)$, the two quantities on the both sides of \eqref{ad} are exactly the same.

In order to  show that the function $\|\cdot\|$ defined in \eqref{norm} is a well-defined norm on the parameter space $\Theta=\H_m\times\mathbb R^{p-1}\times\mathbb R^q$, we    introduce some further notations.
We define  $(p+q) \times (p+q)$ and $(p+q-1) \times (p+q-1)$ matrices $\Omega$ and $\Omega_{\beta_0}$  by
\begin{align}\label{Omega}
\Omega=\Biggl[
\begin{matrix}
\Omega_{1}&\Omega_{3}\trans\\
\Omega_{3}&\Omega_{2}
\end{matrix}
\Biggl],
\qquad\Omega_{\beta_0}=\Biggl[
\begin{matrix}
Q_{\beta_0}\trans&{0}_{p,q}\\
{0}_{q,p-1}&I_q
\end{matrix}
\Biggl]\Omega\Biggl[
\begin{matrix}
Q_{\beta_0}\trans&{0}_{p,q}\\
{0}_{q,p-1}&I_q
\end{matrix}
\Biggl]\trans,
\end{align}
where the matrices 
$\Omega_1\in\mathbb R^{p
\times  p},\Omega_2\in\mathbb R^{q \times q},\Omega_3\in\mathbb R^{q \times p}$ are given by
\begin{equation}
 \label{Omega123}   
\begin{split}
\Omega_1&=\E\big[\sigma^2(X,Z)\{g_0'(X\trans\beta_0)\}^2\{X-R_X(X\trans\beta_0)\}\{X-R_X(X\trans\beta_0)\}\trans\big],\\
\Omega_2&=\E\big[\sigma^2(X,Z)\{Z-R_Z(X\trans\beta_0)\}\{Z-R_Z(X\trans\beta_0)\}\trans\big],\\
\Omega_3&=\E\big[\sigma^2(X,Z)g_0'(X\trans\beta_0)\{Z-R_Z(X\trans\beta_0)\}\{X-R_X(X\trans\beta_0)\}\trans\big],
\end{split}
\end{equation}
respectively. Here, $R_X:\I\to\mathbb R^p$ and $R_Z:\I\to\mathbb R^q$ are vector-valued functions respectively defined by
\begin{align}\label{vxz}
R_X(s)=\frac{\E\{\sigma^2(X,Z)X|X\trans\beta_0=s\}}{\E\{\sigma^2(X,Z)|X\trans\beta_0=s\}},\qquad R_Z(s)=\frac{\E\{\sigma^2(X,Z)Z|X\trans\beta_0=s\}}{\E\{\sigma^2(X,Z)|X\trans\beta_0=s\}},\qquad s\in\I.
\end{align}
We make the following mild assumption on the ranks of the matrices $\Omega_1,\Omega_2$ and $\Omega_3$, and the probability density function $f_{X\trans\beta_0}$ and $\sigma_0$ in \eqref{sigma2}.

\begin{assumption}~\label{a:rank}

\begin{enumerate}[label={\rm(\ref*{a:rank}.\arabic*)},series=conditionA,leftmargin=1.5cm,nolistsep]

\item\label{a:rank1} The matrix $\Omega_2$ in \eqref{Omega123} is invertible and ${\rm rank}(\Omega_1-\Omega_3\Omega_2^{-1}\Omega_3\trans)=p-1$.

\item\label{a:rank2} The density $f_{X\trans\beta_0}$ and the conditional variance $\sigma_0^2$ in \eqref{sigma2} are bounded away from zero and infinity on $\I$.
\end{enumerate}

\end{assumption}

Note that the vector $\beta_0$ is in the null space of the matrices  $\Omega_1$ and $\Omega_3$ in \eqref{Omega123}, which implies that,  under Assumption~\ref{a:rank1}, the matrix $Q_{\beta_0}\trans(\Omega_1-\Omega_3\Omega_2^{-1}\Omega_3\trans)Q_{\beta_0}$ is invertible and the Schur complement of the block $\Omega_2$ in the matrix   $\Omega_{\beta_0}$ in \eqref{Omega}. 
As a consequence Assumption~\ref{a:rank1} guarantees that $\Omega_{\beta_0}$ is a positive definite matrix, so $\lambda_{\min}(\Omega_{\beta_0})>0$. Assumption~\ref{a:rank2} is a standard condition in the literature of single-index models \citep[see, for example,][among many others]{hardle1993,mammen1997,shang2013,kuchibhotla2020}. The following proposition is proved in Section~\ref{app:proof:prop2.1} of the supplement. 

\begin{proposition}\label{prop:2.1}
Under Assumption~\ref{a:rank}, the matrix  $\Omega_{\beta_0}$ in \eqref{Omega} is positive definite, so the function $\|\cdot\|$ defined in \eqref{norm} is a well-defined norm on $\Theta=\H_m\times\mathbb R^{p-1}\times\mathbb R^{q}$ with the corresponding inner product 
\begin{align*}
\l w_1,w_2\r&=\E\big[ \sigma^2(X,Z)\{g_1(X\trans\beta_0)+g_0'(X\trans\beta_0)X\trans Q_{\beta_0}\theta_1+Z\trans\gamma_1\}\\
&\hspace{1cm}\times\{g_2(X\trans\beta_0)+g_0'(X\trans\beta_0)X\trans Q_{\beta_0}\theta_2+Z\trans\gamma_2\}\big]+\lambda J(g_1,g_2)
\end{align*}
for  $w_1=(g_1,\theta_1,\gamma_1),w_2=(g_2,\theta_2,\gamma_2)\in\Theta$.
Moreover, for any $(g,\theta,\gamma)\in\Theta$, it holds that $$\|\theta\|_2^2+\|\gamma\|_2^2\leq \{\lambda_{\min}(\Omega_{\beta_0})\}^{-1}\Vert (g,\theta,\gamma)\Vert^2.$$

\end{proposition}

In order to derive  the convergence rate of the smoothing spline estimator $(\hat g_n,\hat\beta_n,\hat\gamma_n)$ in \eqref{hatwn0} we make the following assumptions on the error distribution and the rate of the regularization parameter $\lambda$.

\begin{assumption}~\label{a:rate}
\begin{enumerate}[label={\rm(\ref*{a:rate}.\arabic*)},series=conditionA,leftmargin=1.5cm,nolistsep]

\item \label{a:02} There exists a constant $c_0>0$ such that $\E\{\exp(c_0\e^2)|X,Z\}<\infty$ almost surely.

\item\label{a:5.0} $\lambda=o(n^{-1/m})$ and $\lambda^{-1}=O(n^{2m/(2m+1)})$ as $n\to\infty$.

\end{enumerate}

\end{assumption}

 Our first main result  gives  the convergence rate of $\hat w_n=(\hat g_n,Q_{\beta_0}\trans\hat\beta_n,\hat\gamma_n)$ as an estimator of $w_0=(g_0,0_{p-1},\gamma_0)$   and is proved in Section~\ref{app:thm:rate} of the online supplement.
\begin{theorem}[Convergence rate]\label{thm:rate}

Under Assumptions~\ref{a:rank} and \ref{a:rate}, it is true that, as $n\to\infty$,
\begin{align*}
\|\hat w_n-w_0\|=\big\|(\hat g_n-g_0,\,Q_{\beta_0}\trans\hat\beta_n,\,\hat\gamma_n-\gamma_0)\big\|=O_p(r_n),
\end{align*}
where $r_n=\sqrt\lambda+n^{-1/2}\lambda^{-1/(4m)}$.
\end{theorem}

\begin{remark}\label{rem:rate}
A direct consequence of Theorem~\ref{thm:rate} is the observation that, in terms of estimation accuracy of $\hat w_n$ with respect to~the $\|\cdot\|$-norm, the optimal choice for the regularization parameter is $\lambda\asymp n^{-2m/(2m+1)}$. This leads to the typical nonparametric convergence rate $n^{-m/(2m+1)}$; see, for example \cite{tsybakov2009nonparametric}. In Section~\ref{sec:minimax} we will show that this rate is in fact minimax optimal.

\end{remark}

\subsection{Reproducing kernel Hilbert space}\label{sec:rkhs}


In order to determine  the joint estimating equation for  the smoothing spline estimator $(\hat g_n,\hat\beta_n,\hat\gamma_n)$ defined in \eqref{hatwn0}, we first derive a necessary condition. The following Lemma is proved in Section~\ref{app:prop:jointee} of the online supplement.

\begin{lemma}\label{prop:jointee}
The smoothing spline estimator $(\hat g_n,\hat\beta_n,\hat\gamma_n)$ satisfies, for any $(g,\theta,\gamma)\in\Theta$,
\begin{align}\label{scoren}
F_n(g,\theta,\gamma)&:=-\frac{1}{n}\sum_{i=1}^n\big\{Y_i-\hat g_n(X_i\trans \hat\beta_n)-Z_i\trans\hat\gamma_n\big\}\notag\\
&\qquad\times\big\{g(X_i\trans \hat\beta_n)+\hat g_n'(X_i\trans \hat\beta_n)X_i\trans Q_{\hat\beta_n}\theta+Z_i\trans\gamma\big\}+\lambda J(\hat g_n,g)\equiv0,
\end{align}
where $Q_{\hat\beta_n}\in\mathbb R^{p\times(p-1)}$ is such that $[\hat\beta_n,Q_{\hat\beta_n}]\in\mathbb R^{p\times p}$ is an orthogonal matrix.

\end{lemma}

We aim to leverage the reproducing kernel Hilbert space theory to identify a data-dependent functional $\mathcal R_n:\H_m\times\mathcal B\times\mathbb R^q\to\Theta$, such that
\begin{align}\label{approx}
F_n(g,\theta,\gamma)=\big\l \mathcal R_n(\hat g_n,\hat\beta_n,\hat\gamma_n),(g,\theta,\gamma)\big\r+\texttt{rem}_n,\qquad \mbox{for all }(g,\theta,\gamma)\in\Theta,
\end{align}
where the remainder term $\texttt{rem}_n$ is uniformly small  for all $(g,\theta,\gamma)\in\Theta$. Then, by Lemma~\ref{prop:jointee}, 
$$\mathcal R_n( g,\beta,\gamma)=0$$ 
is the {\em approximate joint estimating equation} for the smoothing spline estimator $(\hat g_n,\hat\beta_n,\hat\gamma_n)$. To achieve this, we define  the bilinear functionals on $\H_m$
\begin{align}\label{V}
\begin{split}
{V}(g_1,g_2)&=\E\{\sigma_0^2(X\trans\beta_0)\,g_1(X\trans\beta_0)\,g_2(X\trans\beta_0)\},
\\  
\l g_1,g_2\r_K &= {V}(g_1,g_2)+\lambda J(g_1,g_2) ,
\end{split}  
\end{align}
$(g_1,g_2\in\H_m)$ and note that under Assumption~\ref{a:rank2}, the Sobolev space $\H_m$ is a reproducing kernel Hilbert space (RKHS) equipped with the inner product $\l\cdot,\cdot\r_K$ defined in \eqref{V}; see, for example, \cite{cox1990asymptotic}. Let $\|\cdot\|_K$ denote the corresponding norm, $K$ the reproducing kernel, and define $K_s=K(s,\cdot)$, for $s\in\I$.

Theory for ordinary differential equations (ODE, \citealp[e.g.,][]{coddington1956theory}) enables the construction of a sequence of eigenfunctions $\{\phi_j\}_{j\geq1}$ for the inner product $\l\cdot,\cdot\r_K$ that simultaneously diagonalize the operators $V$ and $J$ defined in \eqref{V} and \eqref{J}, respectively, that is,
\begin{align}\label{diag}
{V}(\phi_{j},\phi_{j'})=\delta_{jj'}\,,\qquad \int_{\I} \phi_{j}^{(m)}(s)\,\phi_{j'}^{(m)}(s)\,\d s=\rho_{j}\delta_{jj'}\,,\qquad j,j'\geq1,
\end{align}
where the real-value sequence $\{\rho_j\}_{j\geq1}$ consists of the eigenvalues, and $\delta_{jj'}$ denotes the Kronecker delta.

\begin{proposition}
\label{prop:eigen}

Suppose that the functions $f_{X\trans\beta_0}$ and $\sigma_0$ are continuously differentiable on the inverval $\I$ up to order $(2m-1)$ and satisfy Assumption~\ref{a:rank2}.
Let $\{(\phi_j,\rho_j)\}_{j\geq1}$ be the solutions to the following ordinary differential equation with boundary conditions:
\begin{align}\label{id}
\begin{cases}
(-1)^{m}\phi_{j}^{(2m)}(s)=\rho_{j} \,\phi_{j}(s)f_{X\trans\beta_0}(s)\,\sigma_0^2(s),\\
\phi_{j}^{(\ell)}(s)=0\,,\qquad\qquad s\in\partial\I\,,\,m\leq\ell\leq 2m-1,
\end{cases}
\end{align}
normalized such that $V(\phi_j,\phi_j)=1$. Then, $\{(\rho_j,\phi_j)\}_{j\geq1}$ satisfies \eqref{diag} with $\sup_{j\geq1}\|\phi_j\|_\infty<\infty$ and $\rho_{j}\asymp j^{2m}$. Furthermore, any $g\in\H_{m}$ admits the expansion $g=\sum_{j\geq1}{V}(g,\phi_{j})\phi_{j}$ with convergence in $\H_{m}$ with respect to the norm $\|\cdot\|_K$.
\end{proposition}

A proof of Proposition~2.2 can be found in  \cite{shang2013}.
Denoting $\psi_j=\phi_j^{(m)}$, it follows by a  direct calculation that 
\begin{equation}\label{ode0}
D_m\psi_j=\rho_j\psi_j, \qquad\text{subject to  }\psi_j^{(\ell)}(0)=\psi_j^{(\ell)}(1)=0,\quad \ell=0,1,\ldots,m-1,
\end{equation}
where the differential operator $D_m$ is defined by
\begin{align*}
D_m=(-1)^m\sum_{\ell=0}^m{m\choose\ell}\frac{\d ^\ell}{\d s^\ell}\bigg(\frac{1}{\sigma_0^2(s)f_{X\trans\beta_0}(s)}\bigg)\,\d^{2m-\ell}.
\end{align*}
Then, $\rho_j$ is the eigenvalue of the ordinary differential equation (ODE) with boundary conditions in \eqref{ode0}, whose diverging rate can be obtained by leveraging the eigenvalue theory of differential operators 
\cite[see, for example,][]{cesari2012asymptotic}.
In principle, equation \eqref{ode0} (and as a consequence, equation \eqref{id}) can be solved numerically using the open-source Matlab software \texttt{Chebfun} \citep{driscoll2014chebfun}. 

\begin{example}\label{rem:eigen}
For illustration, consider an example where $\sigma_0^2(s)\equiv\sigma_0^2$ and $X\trans\beta_0$ follows a uniform distribution on $\I=[0,1]$, so $\sigma_0^2(s)f_{X\trans\beta_0}(s)\equiv\sigma_0^2>0$. In this case the eigenfunctions in Proposition~\ref{prop:eigen} are defined by 
\[
\phi_1\equiv\sigma_0^{-1},\qquad \phi_{2j}(s)=\sigma_0^{-1}\sqrt 2\cos(2\pi js),\qquad \phi_{2j+1}(s)=\sigma_0^{-1}\sqrt 2\sin(2\pi js),\qquad j\geq1,
\]
with corresponding eigenvalues $\rho_1=0$, and $\rho_{2j}=\rho_{2j+1}=\sigma_0^{-2}(2\pi j)^{2m}$, for $j\geq1$.
\end{example}

A direct consequence of Proposition~\ref{prop:eigen} is that $\l\phi_j,\phi_k\r_K=\delta_{jk}(1+\lambda\rho_j)$, so for any $g=\sum_{j\geq1}{V}(g,\phi_j)\phi_j\in\H_m$, it holds that 
\begin{align}\label{expansion}
\l g,\phi_j\r_K=\sum_{\ell\geq1}{V}(g,\phi_{\ell})\l\phi_{\ell},\phi_j\r_K=(1+\lambda\rho_j){V}(g,\phi_{j});\qquad g=\sum_{j\geq1}\frac{\l g,\phi_j\r_K}{1+\lambda\rho_j}\phi_{j}.
\end{align}
As a consequence, for the reproducing kernel $K$, we have $\phi_{j}(s)=\l K_{s},\phi_{j}\r_K$, and \eqref{expansion} implies
\begin{align}\label{kst2}
K_{s}=\sum_{j\geq1}\frac{\l K_s,\phi_j\r_K}{1+\lambda\rho_j}\phi_{j}=\sum_{j\geq1}\frac{\phi_{j}(s)}{1+\lambda\rho_{j}}\phi_{j};\qquad  K(s,s')=\sum_{j\geq1}\frac{\phi_{j}(s)\phi_{j}(s')}{1+\lambda\rho_{j}},\qquad s,s\in\I.
\end{align}

Observe that the functional defined by $\lambda J(g,\,\cdot\,)$ is a bounded linear functional on $\H_m$. By the Riesz representation theorem \citep[see, for example, Theorem~1.1 in][]{berlinet2011reproducing}, there exists a unique element $M_\lambda(g)\in\H_m$ such that $\l M_\lambda(g),\,\cdot\,\r_K=\lambda J( g,\,\cdot\,)$.
By definition,  we have $\l M_\lambda(\phi_{j}),\phi_{j'}\r_K=\lambda J(\phi_{j},\phi_{j'})=\lambda\rho_{j}\delta_{jj'}$ for the eigenfunctions $\{\phi_{j}\}_{j\geq 1}$ in Proposition~\ref{prop:eigen}, so \eqref{expansion} implies, for $g\in\H_m$,
\begin{align}\label{mlambda}
M_\lambda(\phi_{j})=\sum_{j,j'\geq1}\frac{\l M_\lambda(\phi_{j}),\phi_{j'}\r_K}{1+\lambda\rho_{j'}}\phi_{j'}=\frac{\lambda\rho_{j}\phi_{j}}{1+\lambda\rho_{j}},\qquad M_\lambda(g)=\lambda\sum_{j\geq1}{V}(g,\phi_j)\frac{\rho_{j}\phi_{j}}{1+\lambda\rho_{j}}.
\end{align}

In the sequel, for vector-valued functions $G=[ g_{1},\ldots, g_{1}]\trans\in\H_m^p$ and $\tilde G=[\tilde g_{1},\ldots,\tilde g_{q}]\trans\in\H_m^{q}$, denote the $p\times q$ matrices ${V}( G,\tilde G)$ and $\l G,\tilde G\r_K$ whose $(i,j)$-entry are ${V}( g_{i},\tilde g_{j})$ and $\l g_{i},\tilde g_{j} \r_K$,
 respectively. 
Recalling the definition of $R_X$ and $R_Z$ in \eqref{vxz}  we make the following assumption on the Fourier coefficients of the vector-valued functions $g_0' \cdot R_X$ and $R_Z$ with respect to the eigenfunctions $\{\phi_j\}_{j\geq1}$ in Proposition~\ref{prop:eigen}.

\begin{assumption}\label{a:finite}
The Fourier coefficient vectors $V(g_0'\cdot R_X,\phi_j)$ and $V(R_Z,\phi_j)$ satisfy
\begin{align*}
\sum_{j\geq1}\big\{\|{V}(g_0'\cdot R_X,\phi_j)\|_2^2+\|{V}(R_Z,\phi_j)\|_2^2\big\}< \infty\,.
\end{align*}
\end{assumption}
Assumption~\ref{a:finite} requires that $\|V(g_0'\cdot R_X,\phi_j)\|_2$ and $\|{V}(R_Z,\phi_j)\|_2$ decay at least at a polynomial rate $j^{-1/2}$. Note that this is a rather mild condition, which is satisfied provided that each of the diagonal entries of the matrices $V(g_0'\cdot R_X,g_0'\cdot R_X)$ and $V(R_Z,R_Z)$ are finite, since
\begin{align*}
&\sum_{j\geq1}\|{V}(g_0'\cdot R_X,\phi_j)\|_2^2={\rm tr}\{V(g_0'\cdot R_X,g_0'\cdot R_X)\},\qquad~~~ \sum_{j\geq1}\|{V}(R_Z,\phi_j)\|_2^2={\rm tr}\{V(R_Z,R_Z)\}.
\end{align*}
A direct consequence of Assumption~\ref{a:finite} is that the functionals defined by $V(g_0'\cdot R_X,\cdot)$ and $V(R_Z,\cdot)$ are bounded linear functionals, that is
\begin{align*}
\|V(g_0'\cdot R_X,g)\|_2^2\leq\|g\|_K^2 \sum_{j\geq1}\|V(g_0'\cdot R_X,\phi_j)\|_2^2\,,\qquad \|V(R_Z,g)\|_2^2\leq\|g\|_K^2 \sum_{j\geq1}\|V(R_Z,\phi_j)\|_2^2\,
\end{align*}
for $g\in\H_m$. 
It then follows from the Riesz representation theorem 
that there exist unique elements $A_X\in\H_m^{p}$ and $A_Z\in\H_m^q$ such that, for $g\in\H_m$,
\begin{align}\label{ab}
\l A_X,g\r_K&={V}(g_0'\cdot R_X,g)=\E \big\{\sigma_0^2(X\trans\beta_0)\,g_0'(X\trans\beta_0)\,R_X(X\trans\beta_0)\,g(X\trans\beta_0)\big\};\notag\\
\l A_Z,g\r_K&={V}(R_Z,g)=\E\big\{\sigma_0^2(X\trans\beta_0)\,R_Z(X\trans\beta_0)\,g(X\trans\beta_0)\big\}.
\end{align}
The series expansion in equation \eqref{expansion} then enables us to deduce that
\begin{align}\label{abx}
&A_X=\sum_{j\geq1}\frac{{V}(g_0'\cdot R_{X},\phi_j)}{1+\lambda\rho_j}\phi_j\,;\qquad
A_Z=\sum_{j\geq1}\frac{{V}(R_Z,\phi_j)}{1+\lambda\rho_j}\phi_j\,.
\end{align}

\subsection{Approximate joint estimating equation}\label{sec:jeq}

We will now leverage the RKHS theory developed in Section~\ref{sec:rkhs}  to establish a data-dependent functional $\mathcal R_n$ such  that  the approximation in \eqref{approx} is valid. For this purpose, we aim to identify two functionals, say $S_{x,z}$ and $S_\lambda^\circ$, such that for any $(g,\theta,\gamma)\in\Theta$,
\begin{equation}
\begin{gathered}
\big\l S_{x,z}(\hat g_n,\hat\beta_n),(g,\theta,\gamma)\big\r =g(x\trans\hat\beta_n)+\hat g_n'(x\trans\hat\beta_n)x\trans  Q_{\beta_0}\theta+z\trans\gamma,~~\quad
\\
\big\l S_\lambda^\circ(\hat g_n),(g,\theta,\gamma) \big\r=\lambda J(\hat g_n, g),\label{wlambdag}
\end{gathered}
\end{equation}
where $(x,z)\in\X\times\Z$. This is achieved by applying the following proposition, which gives the explicit expression for $S_{x,z}$ and $S_\lambda^\circ$ and is proved in Section~\ref{app:proof:ru} of the online supplement, respectively.

\begin{proposition}\label{cor:sxz}
Suppose Assumptions~\ref{a:rank}--\ref{a:finite} are valid. Equation~\eqref{wlambdag} is satisfied by taking $$S_{x,z}(g,\beta)=\big(H_{x,z}(g,\beta),N_{x,z}(g,\beta),T_{x,z}(g,\beta)\big),\ \, S_\lambda^\circ(g)=\big(H^\circ_{\lambda}(g),N^\circ_{\lambda}(g),T^\circ_{\lambda}(g)\big)\in\H_m\times\mathbb R^{p-1}\times\mathbb R^q,$$ where their components are respectively defined by
\begin{align*}
&\left[\begin{matrix}
H_{x,z}(g,\beta) & H^\circ_{\lambda}(g)\\ N_{x,z}(g,\beta) & N^\circ_{\lambda}(g)\\
 T_{x,z}(g,\beta) & T^\circ_{\lambda}(g)
\end{matrix}\right]=\left[
\begin{matrix}
K_{x\trans\beta} & M_\lambda (g)\\
0_{p-1} &0_{p-1}\\
0_q &0_q
\end{matrix}
\right]\\
&\hspace{1cm}-\left[
\begin{matrix}
A_X\trans Q_{\beta_0} & A_Z\trans\\
I_{p-1}&0_{p-1,q}\\
0_{q,p-1}&I_q
\end{matrix}
\right]{(\Omega_{\beta_0}+\Sigma_{\beta_0,\lambda})^{-1}}\Biggl[\begin{matrix}
-Q_{\beta_0}\trans  \{g'(x\trans\beta)x-A_X(x\trans\beta)\} & Q_{\beta_0}\trans{V}(g_0'\cdot R_X,M_\lambda g)\\
-\{z-A_Z(x\trans\beta)\} &
{V}(R_Z,M_\lambda g)
\end{matrix}\Biggl].
\end{align*}
Here, $\Omega_{\beta_0}$, $M_{\lambda}$, and $A_X,A_Z$ are given by \eqref{Omega}, \eqref{mlambda}, and \eqref{abx}, respectively; $\Sigma_{\beta_0,\lambda}$ is a $(p+q-1)$-square matrix such that $\lim_{\lambda\to0}\|\Sigma_{\beta_0,\lambda}\|_2=0$; see equation~\eqref{Sigmalambda} in the online supplement for its exact definition.

\end{proposition}

Using Proposition~\ref{cor:sxz}, we now define the functional $\mathcal R_n$ in equation \eqref{approx} by
\begin{align}\label{scoren2}
\mathcal R_n( g,\beta,\gamma)=-\frac{1}{n}\sum_{i=1}^n\big\{Y_i- g(X_i\trans \beta)-Z_i\trans\gamma\big\}S_{X_i,Z_i}( g,\beta)+S_\lambda^\circ(g).
\end{align}
The following proposition states that $\mathcal R_n(g,\beta,\gamma)=0$ is the approximate joint estimating equation for the smoothing spline estimator. Its proof leverages Proposition~\ref{prop:2.1} and Lemma~\ref{prop:jointee}, and is given in Section~\ref{app:thm:score} of the supplement.

\begin{proposition}
\label{thm:score}
Suppose Assumptions~\ref{a:rank}--\ref{a:finite} are valid. Then, for the functional $\mathcal R_n$ in \eqref{scoren2}, the smoothing spline estimator $(\hat g_n,\hat\beta_n,\hat\gamma_n)$ defined by \eqref{hatwn0} satisfies 
$$\|\mathcal R_n(\hat g_n,\hat\beta_n,\hat\gamma_n)\|=o_p(n^{-1/2}).
$$

\end{proposition}

\subsection{Joint Bahadur representation and weak convergence}\label{sec:jointbahadur}

Now, we are ready to state the main result of this article.
Let
\begin{align}\label{proxy}
\tilde{\mathcal R}_n(w):=-\frac{1}{n}\sum_{i=1}^n\big\{Y_i-g(X_i\trans\beta_0)-g_0'(X_i\trans\beta_0)Q_{\beta_0}\theta-Z_i\trans\gamma\big\}S_{X_i,Z_i}(g_0,\beta_0)+S_\lambda^\circ(g)
\end{align}
be the proxy linear functional of $\mathcal{R}_n$ in \eqref{scoren2}, where $w=(g,\theta,\gamma)\in\Theta$. Note 
 that the Fr\'echet derivative of $\tilde{\mathcal R}_n$ satisfies $\big\l\mathcal D\tilde{\mathcal R}_n(w)w_1,w_2\big\r=\l w_1,w_2\r$, for $w,w_1,w_2\in\Theta$.
We now make the following assumption on the rate of the regularization parameter $\lambda$.

\begin{assumption}\label{a:reg}
For the sequence $r_n=\sqrt\lambda+n^{-1/2}\lambda^{-1/(4m)}$ in Theorem~\ref{thm:rate}, it holds that $n\lambda^{1/m}\to\infty$, $\sqrt n\lambda^{-1/(4m)} (\lambda^{-1/(4m-4)}r_n^2+r_n^{3/2})=o(1)$ and $a_n:=\lambda^{-(4m-1)/(8m^2)}r_n\sqrt{\log\log(n)}=o(1)$ as $n\to\infty$.

\end{assumption}

The conditions on the rate of $\lambda$ of the regularizing parameter in Assumption~\ref{a:reg} are used in deriving the joint Bahadur representation where we have to bound certain nonlinear remainder terms, which are explicitly stated in Lemmas~\ref{lem:n-or} and \ref{lem:bounds} in the online supplement. First, taking the optimal regularization parameter  $\lambda\asymp n^{-2m/(2m+1)}$ in view of Theorem~\ref{thm:rate} yields $n\lambda^{1/m}\asymp n^{(2m-1)/(2m+1)}\to\infty$, so $n\lambda^{1/m}\to\infty$ is a mild condition. Second, the condition that $\sqrt n\lambda^{-1/(4m)} (\lambda^{-1/(4m-4)}r_n^2+r_n^{3/2})=o(1)$ is essentially due to the error bound in the theoretical derivation regarding the index link derivative $g_0'$, which in principle requires smoothness condition of higher order compared to estimating $g_0$; see, for example, \cite{vandegeer2000,ruppert2003semiparametric,tsybakov2009nonparametric}.
Third, it is true that under Assumption~\ref{a:5.0}, $a_n\lesssim\big(\lambda^{(2m-1)^2/(8m^2)}+n^{-\frac{(2m-1)^2}{4m(2m+1)}}\big)\sqrt{\log\log(n)}$, so the assumption that $a_n=o(1)$ is a very mild condition as well.

Proposition~\ref{thm:score}, together with Lemmas~\ref{lem:n-or} and \ref{lem:bounds} in the online supplement enable us to derive the following theorem on the joint Bahadur representation of the smoothing spline estimator $(\hat g_n,\hat\beta_n,\hat\gamma_n)$ defined in \eqref{hatwn0}. 
We emphasize that this result is the first one in the literature of semiparametric single-index models.

\begin{theorem}\label{thm:bahadur}
Under Assumptions~\ref{a:rank}--\ref{a:reg}, we have for $\hat w_n=(\hat g_n,Q_{\beta_0}\trans\hat\beta_n,\hat\gamma_n)$ and $w_0=(g_0,0_{p-1},\gamma_0)$,
\begin{align*}
\|\hat w_n-w_0+\tilde{\mathcal R}_n(w_0)\|=\bigg\|\hat w_n-w_0-\frac{1}{n}\sum_{i=1}^n\epsilon_iS_{X_i,Z_i}(g_0,\beta_0)+S_\lambda^\circ(g_0)\bigg\|=o_p(n^{-1/2}),
\end{align*}
where $S_{x,z}(g,\beta)$ and $S_\lambda^\circ$ are defined in Proposition~\ref{cor:sxz}.
\end{theorem}

By leveraging the joint Bahadur representation in Theorem~\ref{thm:bahadur}, we are able to derive the joint weak convergence of the smoothing spline estimator $(\hat g_n,\hat\beta_n,\hat\gamma_n)$ defined in \eqref{hatwn0}.
The following theorem states that the joint vector of both the point-wise evaluation of the nonparametric index link function $\hat g_n(s)$ for a given $s\in\I$, and the parametric coefficients $\hat\beta_n$ and $\hat\gamma_n$ are asymptotically normally distributed. It is proved in Section~\ref{app:thm:asymp} of the online supplement.

\begin{theorem}\label{thm:asymp}

Let $s\in\I$ be arbitrary. Suppose Assumptions~\ref{a:rank}--\ref{a:reg} are valid and that the limit
\begin{align}\label{sigmas2}
\sigma_{(s)}^2:=\lim_{\lambda\to0}\lambda^{1/(2m)}{V}(K_{s},K_{s})
\end{align}
exists. In addition, assume that there exists a constant $\mu\in((2m)^{-1},1]$ such that $n \lambda^{\mu+1}=o(1)$ and
\begin{align}\label{newcond}
\sum_{j\geq1}\rho_j^{\mu}\big\{\|{V}(g_0'\cdot R_X,\phi_j)\|_2^2+\|{V}(R_Z,\phi_j)\|_2^2\big\}<\infty\,.
\end{align}
Then, it holds that
\begin{align}\label{main}
\left[\begin{matrix}
\sqrt{n}\lambda^{1/(4m)}\{\hat g_n(s)-g_0(s)+M_\lambda g_0(s)\}\\
\sqrt{n}(\hat\beta_n-\beta_0)\\
\sqrt{n}(\hat\gamma_n-\gamma_0)\\
\end{matrix}\right]\converged N\left(\Bigg[\begin{matrix}
0 \\ 0_{p+q}
\end{matrix}\Bigg],\Bigg[\begin{matrix}
\sigma_{(s)}^2 & 0_{p+q}\trans\\
0_{p+q} & \Omega^+
\end{matrix}\Bigg]\right)\,,
\end{align}
where $M_\lambda$ is defined in \eqref{mlambda} and $\Omega^+$ is the Moore–Penrose inverse of the matrix $\Omega$ in \eqref{Omega} defined by
\begin{align}\label{Omegastar}
\Omega^+&=\Biggl[
\begin{matrix}
Q_{\beta_0}& 0_{p,q}\\
 0_{q,p-1}&I_q
\end{matrix}
\Biggl]\Omega_{\beta_0}^{-1}\Biggl[
\begin{matrix}
Q_{\beta_0}& 0_{p,q}\\
 0_{q,p-1}&I_q
\end{matrix}
\Biggl]\trans\,.
\end{align}

\end{theorem}

The joint asymptotic normality established in Theorem~\ref{thm:asymp} shows that $\hat g_n(s)$ and $(\hat\beta_n,\hat\gamma_n)$ are asymptotically independent. This result is somewhat surprising, considering that the index link function $g_0$ and index direction $\beta_0$ have nonlinear relations in the single-index model \eqref{model}. Note that the joint weak convergence enables us to construct consistent tests for the hypothesis that  the model parameters (including the index link function $g_0$) satisfy a structural equation; see Section~\ref{sec:joint} for details.

\begin{example}
 We retrieve Example \ref{rem:eigen} where $\sigma_0^2(s)\equiv\sigma_0^2$ and $f_{X\trans\beta_0}\equiv1$ on $\I=[0,1]$. In this case we have
\begin{gather*}
{V}(K_{s},K_{s})=\sigma_0^{-2}
+2\sigma_0^{-2}\sum_{j\geq1}\frac{\sin^2\{(2j-1)\pi s\}+\cos^2(2j\pi s)}{\big\{1+\sigma_0^{-2}\lambda(2j\pi)^{2m}\big\}^2},
\\
\sigma_{(s)}^2\asymp\sigma_0^{-2}+c_m\pi^{-1}\sigma_0^{-2+1/m},\qquad\text{where }c_m=\int_0^{+\infty}\frac{1}{(1+x^{2m})^2}\d x>0.
\end{gather*}
To be more specific, the point-wise variances at points $s=0,1/4,1/2,3/4$ can be explicitly computed as
\begin{align*}
&\sigma^2_{(1/2)}=\sigma_0^{-2}+2c_m\pi^{-1}\sigma_0^{-2+1/m},\qquad \sigma_{(0)}^2=\sigma_{(1/4)}^2=\sigma_{(3/4)}^2=\sigma_0^{-2}+c_m\pi^{-1}\sigma_0^{-2+1/m}.
\end{align*}
\end{example}

\begin{remark} ~~ 
\begin{itemize}
    \item[(a)]
Note that the  (marginal) asymptotic distribution of the smoothing spline estimator $\hat\beta_n$ in Theorem~\ref{thm:asymp} is the same as the one derived for the estimators in \cite{hardle1993,chang2010asymptotically,kuchibhotla2020} for single-index models with no partially linear structure (i.e., $Z\equiv0$). Moreover, for partially linear single-index models, \cite{carroll1997,yu2002,wang2010,ma2013} among others showed the marginal weak convergence of estimators for the parametric components $(\beta_0,\gamma_0)$, and our marginal limit distribution is the same with theirs. Finally, we emphasize that, the root-$n$ weak convergence rate of the parametric components $(\hat\beta_n,\hat\gamma_n)$ implies the efficiency of the smoothing spline estimator defined in \eqref{hatwn0}.
    \item[(b)]
The condition in \eqref{newcond} essentially means that the entries of $g_0'\cdot R_X$ and $R_Z$ reside in the Sobolev space $\H_{m\mu}(\I)$ on $\I$ of order $m\mu$, which is satisfied when $\|V(g_0'\cdot R_X,\phi_j)\|_2$ and $\|V(R_Z,\phi_j)\|_2$ decay at a polynomial rate that is faster than $j^{-1/2}\rho_j^{-\mu/2}\asymp j^{-m\mu-1/2}$. Note that this polynomial-decay condition on Fourier coefficients of the type in \eqref{newcond} is standard in the literature of nonparametric and semiparametric statistics; we refer to, for example Condition~(2.3) in \cite{cox1990asymptotic} and the discussion on the smoothness conditions in \cite{mammen1997}.
For the sake of completeness, we provide in Theorem~\ref{thm:asymp:prep} of the 
online supplement a result on the  weak convergence without condition~\eqref{newcond}.
    \item[(c)]
If $\sqrt n\lambda^{1/(4m)}M_\lambda g_0(s)=o(1)$ as $n\to\infty$, the point-wise index link estimator $\hat g_n(s)$ is asymptotically unbiased. In either case, we propose to use a term $M_\lambda \hat g_n(s)$ for bias adjustment; see  Lemma~\ref{prop:bias} for a theoretical justification.
\end{itemize}
\end{remark}

\subsection{Marginal Bahadur representation}\label{sec:marginalbahadur}

The final result of this section gives a marginal Bahadur representation of the index link estimator $\hat g_n$, together with that of the $(p+q)$-vector $[\hat\beta_n\trans,\hat\gamma_n\trans]\trans$ consisting of the smoothing spline estimators of the parametric components in \eqref{model}. Its proof is given in Section~\ref{app:marginal} of the online supplement.

\begin{theorem}\label{thm:marginal}
If Assumptions~\ref{a:rank}--\ref{a:reg} and Condition~\eqref{newcond} are  satisfied, we have
\begin{gather}
\sup_{s\in\I}\bigg|\hat g_n(s)-g_0(s)+M_\lambda  g_0(s)- \frac{1}{n}\sum_{i=1}^n\e_i K_{X_i\trans\beta_0}(s)\bigg|=o_p(n^{-1/2}\lambda^{-1/(4m)}),\notag
\\
\left\|\Bigg[\begin{matrix}
\hat\beta_n-\beta_0\\
\hat\gamma_n-\gamma_0
\end{matrix}\Bigg]-\frac{1}{n}\sum_{i=1}^n\e_i\Omega^+\Bigg[\begin{matrix} g_0'(X_i\trans\beta_0)\{X_i-R_X(X_i\trans\beta_0)\}\\
Z_i-R_Z(X_i\trans\beta_0)
\end{matrix}
\Bigg]\right\|_2=o_p(n^{-1/2}),\label{m22}
\end{gather}
where $M_\lambda$, $\Omega^+$, and $R_X,R_Z$ are defined in \eqref{mlambda}, \eqref{Omegastar} and \eqref{vxz}, respectively.

\end{theorem}

The corresponding marginal Bahadur representation for $\hat g_n$ established in Theorem~\ref{thm:marginal} is useful for deriving inference tools regarding the index link function $g_0$ (see Sections~\ref{sec:band} and \ref{sec:rele} for more details). In addition, a direct consequence of the marginal Bahadur representation for $(\hat\beta_n,\hat\gamma_n)$ is that the estimators $\hat\beta_n$ and $\hat\gamma_n$ are $O_p(n^{-1/2})$-consistent, which is formally stated in the following proposition. We refer to Section~\ref{app:marginal} of the supplement online for its proof.

\begin{proposition}\label{prop:rootn}
Let Assumptions~\ref{a:rank}--\ref{a:reg}, and Condition \eqref{newcond} be satisfied. Then, it holds that $\|\hat\beta_n-\beta_0\|_2+\|\hat\gamma_n-\gamma_0\|_2=O_p(n^{-1/2}).$

\end{proposition}

\section{Statistical consequences}\label{sec:stat}

In this section, we use the theory developed in Section \ref{sec2} to develop several statistical inference tools for the single-index model in \eqref{model}. We first show in Section~\ref{sec:minimax} that the estimator $\hat g _n$ in \eqref{hatwn0} achieves the minimax optimal convergence rate. In Section~\ref{sec:band}, we develop multiplier bootstrap simultaneous confidence bands for the index link function $ g _0$.
In Section~\ref{sec:rele}, we propose a test for the hypothesis that the maximum absolute deviation between the index link function $g_0$ and a given function $g_*$ does not exceed a given threshold. Finally, we construct a multiplier bootstrap test for the joint hypotheses testing problem in \eqref{hp}.

\subsection{Minimax optimality}\label{sec:minimax}

Let $(\tilde g_n,\tilde\beta_n,\tilde\gamma_n)\in\H_m\times\B\times\mathbb R^q$ be a generic estimator of $(g_0,\beta_0,\gamma_0)$ from  the i.i.d.~data $\{(X_i,Y_i,Z_i)\}_{i=1}^n$ defined by  model \eqref{model}. We consider the $L^2$-risk of $(\tilde g_n,\tilde\beta_n,\tilde\gamma_n)$ for estimating $g_0$ and for the joint estimation of $(g_0,\beta_0,\gamma_0)$, respectively, defined by
\begin{align}\label{l2risk}
\mathfrak R^2(\tilde g_n,\tilde\beta_n,\tilde\gamma_n)=\int_{\I}\{\tilde g_n(s)-g_0(s)\}^2\d s+\|\tilde\beta_n-\beta_0\|_2^2+\|\tilde\gamma_n-\gamma_0\|_2^2.
\end{align}
The convergence rate of $\hat w_n=(\hat g_n,Q_{\beta_0}\trans\hat\beta_n,\hat\gamma_n)$ w.r.t.~the $\|\cdot\|$-norm established in Theorem~\ref{thm:rate} enables us to derive the upper bound of the convergence rate of the smoothing spline estimator $(\hat g_n,\hat\beta_n,\hat\gamma_n)$ with respect to the $L^2$-risk $\mathfrak R$ in \eqref{l2risk}. Moreover, we establish that this rates is of the same order as the lower bound for the estimation of $(g_0,\beta_0,\gamma_0)$, which implies that the smoothing spline estimator $(\hat g _n,\hat\beta_n,\hat\gamma_n)$  achieves minimax optimal rate under the $L^2$-risk $\mathfrak R$. To be precise, let $\mathcal G_n\subset\H_m\times\B\times\mathbb R^q$ denote the set consisting of all measureable functions of the data $\{(X_i,Y_i,Z_i)\}_{i=1}^n$, and let $\mathcal F$ denote the collection of the joint distribution functions $F_{X,Y,Z}$ of $(X,Y,Z)$ that satisfy Assumptions~\ref{a:rank} and \ref{a:rate}, according to the partially linear single-index model in \eqref{model}. 
The following theorem establishes the optimal convergence rate of the smoothing spline estimator \eqref{hatwn0}
and is proved in  Section~\ref{app:thm:minimax} of the online supplement.

\begin{theorem} ~~\\
\label{thm:minimax}
(i) If the regularization parameter satisfies $\lambda\asymp n^{-2m/(2m+1)}$, the smoothing spline estimator $(\hat g_n,\hat\beta_n,\hat\gamma_n)$ in \eqref{hatwn0} satisfies
\begin{align*}
\lim_{c\to \infty}\,\limsup_{n\to\infty}\sup_{\substack{F_{X,Y,Z}\in\mathcal{F}\\ (g _0,\beta_0,\gamma_0)\in\mathcal{H}_m\times\B\times\mathbb R^q}}\P\big\{\mathfrak R(\hat g_n,\hat\beta_n,\hat\gamma_n)\geq cn^{-m/(2m+1)}\big\}=0\,.
\end{align*}
(ii) There exists a constant $c_0>0$ such that
\begin{align*}
\liminf_{n\to\infty}\inf_{(\tilde g _n,\tilde\beta_n,\tilde\gamma_n)\in\mathcal G_n}\sup_{\substack{F_{X,Y,Z}\in\mathcal{F}\\ (g _0,\beta_0,\gamma_0)\in\mathcal{H}_m\times\B\times\mathbb R^q}}\P\big\{\mathfrak R(\tilde g_n,\tilde\beta_n,\tilde\gamma_n)\geq c_0n^{-m/(2m+1)}\big\}>0\,.
\end{align*}
\end{theorem}

\subsection{Simultaneous confidence bands}\label{sec:band}

In principle, the joint weak convergence established in Theorem~\ref{thm:asymp} provides  a point-wise asymptotic $(1-\alpha)$-confidence interval of $g_0(s)$, for given $s\in\I$, defined by
\begin{align*}
\big[\hat g_{n}(s)+M_\lambda g_0(s)-n^{-1/2}\lambda^{-1/(4m)}\mathcal Q_{1-\alpha/2}\,\sigma_{(s)}\,,\hat g_{n}(s)+M_\lambda g_0(s)+n^{-1/2}\lambda^{-1/(4m)}\mathcal Q_{1-\alpha/2}\,\sigma_{(s)}\big],
\end{align*}
where $\sigma_{(s)}$ is defined in \eqref{sigmas2}, and $\mathcal Q_{1-\alpha/2}$ denotes the $(1-\alpha/2)$-quantile of the standard normal distribution. However, this interval cannot be directly implemented and 
several issues have to be resolved. First, we propose a bias adjustment via a term $M_\lambda(\hat g_n)$. Second, the asymptotic variance $\sigma_{(s)}$ in \eqref{sigmas2} is rarely available in practice and in general difficult to estimate. 
We address this problem as a by-product, when we are constructing simultaneous confidence bands for the index link function $g_0$. This   is a even more complex problem and the main focus here. Point-wise confidence intervals, which avoid the estimation  of the nuisance parameters, can be obtained by  similar arguments and are briefly discussed  in  Remark~\ref{rem:pointwise} at the end of this section.
 We make the following assumptions, where we observe  that, by Assumption~\ref{a:rank2},  the bilinear functional $V$ in \eqref{V} defines a norm on $\H_m$, which we denote by $\|\cdot\|_V^2$ in the following discussion.

\begin{assumption}\label{a:b}~

\begin{enumerate}[label={\rm(\ref*{a:b}.\arabic*)},series=conditionA,leftmargin=1.5cm,nolistsep]

\item\label{a:bias} Assume as $n\to\infty$, (i) $\sqrt n\lambda^{(2m+1)/(4m)}=O(1)$; (ii) $\sup_{s\in\I}\|M_\lambda (K_s)\|_K=o(1)$.

\item\label{holder} There exist positive constants $c_0,\nu,\vartheta$ such that
\begin{align*}
\limsup_{\lambda\to0}\,\lambda^{\nu}\|K_{s}-K_{t}\|_V^2\leq c_0|s-t|^{2\vartheta},\qquad s,t\in\I,
\end{align*}
where the constant $\nu$ satisfies either one of the following two conditions: (a) $\nu< 1/(2m)$, $\vartheta\geq0$; (b) $\nu=1/(2m)$, $\vartheta>1$.

\end{enumerate}
\end{assumption}

\begin{remark}
Assumption~\ref{a:bias} is a mild condition which essentially justifies approximating the bias term $M_\lambda(g_0)$ in the margianl Bahadur representation for $g_0$ in \eqref{m22} by its natural estimator $M_\lambda(\hat g_n)$. Condition (i) is satisfied if $\lambda=O(n^{-2m/(2m+1)})$ (see Remark~\ref{rem:rate} and Theorem~\ref{thm:minimax}), which is implied by the optimal regularization parameter $\lambda\asymp n^{-2m/(2m+1)}$. Condition (ii) in Assumption~\ref{a:bias} requires that the operator $M_\lambda$ in equation~\eqref{mlambda} satisfies
\begin{align*}
\sup_{s\in\I}\|M_\lambda(K_s)\|_K^2=\lambda^2\sup_{s\in\I}\sum_{j\geq1}\frac{\rho_j^2\phi_j^2(s)}{(1+\lambda\rho_j)^3}=o(1).
\end{align*}
Assumption~\ref{holder} is a H\"older-type condition on the reproducing kernel $K$ with respect to the norm $\| \cdot\|_V$ in \eqref{V},
which is also mild and useful for deriving simultaneous confidence bands for the index link function $g_0$.
Note that a similar assumption has been made by \cite{shang2013} in the context of nonparametric regression. Compared with the equivalent kernel condition and  exponential envelop condition in equations~(5.1)--(5.4) in this reference, we observe that the H\"older-type condition~\ref{holder} is weaker. In particular, \cite{shang2013} assumed a uniform upper bound on the partial derivative of the reproducing kernel $K$ is additionally.
\end{remark}

The following lemma provides a theoretical justification of the bias adjustment by approximating $M_\lambda (g_0)$ by $M_\lambda(\hat g_n)$. Its proof is given in Section~\ref{app:prop:bias} of the supplement.

\begin{lemma} \label{prop:bias}
If Assumptions~\ref{a:rank}--\ref{a:reg} and \ref{a:bias} are valid.
Then, it holds
\begin{align*}
\sup_{s\in\I}|M_\lambda\hat g_n(s)-M_\lambda g_0(s)|=o_p(n^{-1/2}\lambda^{-1/(4m)}).
\end{align*}
As a consequence, it is true that
\begin{align*}
\sup_{s\in\I}\bigg|\hat g_n(s)-g_0(s)+M_\lambda\hat g_n(s)- \frac{1}{n}\sum_{i=1}^n\e_i K_{X_i\trans\beta_0}(s)\bigg|=o_p(n^{-1/2}\lambda^{-1/(4m)}).
\end{align*}
\end{lemma}

Now, we are ready to construct simultaneous confidence bands for the index link function $g_0$ based on the smoothing spline estimator in \eqref{hatwn0}. Define the process
\begin{align}\label{gns}
\mathbb G_{n}(s)=\sqrt{n}\lambda^{1/(4m)}\{\hat g_{n}(s)+M_\lambda\hat g_n(s)-g_{0}(s)\},\qquad T_{n} =\sup_{s\in\I}|\mathbb G_{n}(s)|.
\end{align}
Using Lemma \ref{prop:bias} it can be shown that 
the statistic $T_{n} =\sup_{s\in\I}|\mathbb G_{n}(s)|$ converges weakly with a non-degenerate limit distribution. Consequently, if  $c_{1-\alpha}>0$ denotes  the $(1-\alpha)$-quantile of this limit distribution, then 
\begin{align}\label{cv}
\mathfrak{C}_n(\alpha)=\Big\{g:\I\to\mathbb R\, \Big |\,    |  g(s) - \hat g_{n}(s)-M_\lambda\hat g_n(s) 
 |
\leq n^{-1/2}\lambda^{-1/(4m)} c_{1-\alpha} ,
~ \forall  s\in\I  \Big\}
\end{align}
 defines a  simultaneous asymptotic $(1-\alpha)$-confidence region for $g_0$, that is, $\lim_{n\to\infty}\P\{g_0\in\mathfrak{C}_n(\alpha)\}=1-\alpha$.
However, the critical value $c_{1-\alpha}$ in \eqref{cv} depends in an intricate way on the unknown limiting distribution of $T_n$ and is difficult to estimate in practice (see Section~\ref{app:thm:bootstrap} of the online supplement for more theoretical details).
To circumvent this difficulty, we propose to apply the multiplier bootstrap to approximate the critical value $c_{1-\alpha}$. The procedure is described in Algorithm~\ref{algo:band}, and is theoretically justified by the following result, which is  proved in Section~\ref{app:thm:bootstrap} of the online supplement. 

\begin{algorithm}[tb!]
\caption{Simultaneous confidence band for the index link function\label{algo:band}}
\medskip
\noindent
{\bf Input:}  Data $\{(Y_i,X_i,Z_i)\}$;  regularization parameter $\lambda$; nominal level $\alpha$; bootstrap sample size $B$.

\begin{enumerate}

\item Generate i.i.d.~bootstrap weights $\{W_{i,b}\}_{1\leq i\leq n,1\leq b\leq B}$ independent of the data $\{(X_i,Y_i,Z_i)\}_{i=1}^n$ from a two-point distribution: taking the value $1-1/\sqrt{2}$ with probability $2/3$, and $1+\sqrt{2}$ with probability $1/3$, such that $\E(W_{i,b})=\Var(W_{i,b})=1$.

\item Compute  $(\hat g_{n},\hat\beta_n,\hat\gamma_n)$ in \eqref{hatwn0}. For each $1\leq b\leq B$, compute the multiplier bootstrap estimator
\begin{align}\label{hatbetastar}
(\hat g_{n,b}^*,\hat\beta_{n,b}^*,\hat\gamma_{n,b}^*)&=\argmin_{(g,\beta,\gamma)\in\H_m\times\B\times\mathbb R^q}\Bigg\{\frac{1}{2n}\sum_{i=1}^nW_{i,b}\big[Y_i-g(X_i\trans \beta)-Z_i\trans\gamma\big]^2+\frac{\lambda}{2}\,J(g,g)\Bigg\}.
\end{align}

\item Compute the multiplier bootstrap estimates: for $1\leq b\leq B$,
\begin{align}\label{gnqstar}
\mathbb G_{n,b}^*(s)=\sqrt{n}\lambda^{1/(4m)}\{\hat g_{n,b}^*(s)+M_\lambda \hat g_{n,b}^*(s)-\hat g_n(s)-M_\lambda\hat g_n(s)\},\quad \hat T_{n,b}^* =\sup_{s\in\I}|\mathbb G_{n,b}^*(s)|.
\end{align}

\item Compute the empirical $(1-\alpha)$-quantile of the  multiplier bootstrap estimates $\{\hat T_{n,b}^*\}_{b=1}^B$:
\begin{align*}
\mathcal Q_{1-\alpha}(\hat T_{n,B}^*)=\inf\Big\{\tau\in\mathbb R:\frac{1}{B}\sum_{b=1}^B\mathbf{1}\{\hat T_{n,b}^*\leq \tau\}\geq 1-\alpha\Big\}.
\end{align*}

\item Let $\hat g_{n,B}^{*^\pm}(s)=\hat g_{n}(s)+M_\lambda\hat g_n(s)\pm n^{-1/2}\lambda^{-1/(4m)}\mathcal Q_{1-\alpha}(\hat T^*_{n,B})$, $s\in\I$.

\end{enumerate}

{\bf Output:} Simultaneous $(1-\alpha)$-confidence band $\mathfrak{C}_{n,B}^*(\alpha)=\big\{g\,|\, \hat g_{n,B}^{*^-}(s)\leq g(s)\leq \hat g_{n,B}^{*^+}(s);s\in\I\big\}$.
 
\end{algorithm}

\begin{theorem}\label{thm:cb:boot}

Suppose Assumptions~\ref{a:rank}--\ref{a:b} and \eqref{newcond} are valid. Then, the set $\mathfrak{C}_{n,B}^*(\alpha)$ obtained by Algorithm~\ref{algo:band} defines a simultaneous asymptotic $(1-\alpha)$ confidence band for the index link function $g_0$ in model \eqref{model}, that is
\begin{align}
\label{consetboot}
\lim_{n,B\to\infty}\P\{g_0\in\mathfrak{C}^*_{n,B}(\alpha)\}=1-\alpha\,.
\end{align}
\end{theorem}

\begin{remark}\label{rem:pointwise}

Using similar ideas, we are able to construct a point-wise confidence interval for the index link function. To be specific, let $\mathcal Q_{1-\alpha}\{\mathbb G_{n,B}^*(s)\}$ be the empirical $(1-\alpha)$-quantile of the multiplier bootstrap sample $\{\mathbb G_{n,b}^*(s)\}_{b=1}^B$ defined by \eqref{gnqstar} in Algorithm \ref{algo:band}. Then, the interval with endpoints 
\begin{align*}
\hat g_{n}(s)+M_\lambda\hat g_n(s)\pm n^{-1/2}\lambda^{-1/(4m)}\mathcal Q_{1-\alpha}\{\mathbb G_{n,B}^*(s)\}
\end{align*}
defines an asymptotic $(1-\alpha)$-confidence interval for $g_0(s)$. The proof is similar to that of Theorem~\ref{thm:cb:boot} and omitted for brevity.

\end{remark}

\subsection{Relevant hypotheses}\label{sec:rele}

We aim to test whether the maximum deviation between  a given function $ g _*$
and the unknown true index link function $g_0$ exceeds a given value $\Delta>0$:
\begin{align}\label{rele}
H_0:\,d_\infty:=\sup_{s\in\I}|g_0(s)-g_*(s)|\leq\Delta\qquad&\text{versus}\qquad H_1:\, d_\infty>\Delta.
\end{align}
A simple ad-hoc test for these hypotheses can be constructed  observing that the 
condition 
$$
g_0\in \mathfrak B_*:=\big\{g:\I\to\mathbb R\,|\,g_*(s)-\Delta\leq g(s)\leq g_*(s)+\Delta,\text{ for all }s\in\I\big\}
$$
is a sufficient condition for $H_0$ in \eqref{rele}.
Due to the duality between  hypotheses testing and confidence regions, one can construct an (asymptpotic)  level-$\alpha$ test for the relevant hypotheses \eqref{rele}  on the basis of  the simultaneous confidence bands derived in Section~\ref{sec:band}, that is 
\begin{align}\label{rele:band}
\text{reject the null hypotheis  $H_0$ in \eqref{rele} whenever }~ \mathfrak B_*\not\subset \mathfrak{C}_{n,B}^*(\alpha),
\end{align}
where $\mathfrak{C}_{n,B}^*(\alpha)$ is the simultaneous confidence band defined by Algorithm~\ref{algo:band}. However,  tests construct by this general principle are often very conservative and might have low  power in practice (we will confirm this for test \eqref{rele:band} in Section~\ref{sec:finite}, where the study the finite finite-sample properties). This motivates us to construct an exact test based on a consistent estimator of the maximum deviation $d_\infty$ in \eqref{rele} which will be the main focus of the discussion of this section.

In view of the bias-adjusted index link estimator in Lemma~\ref{prop:bias}, a reasonable decision rule rejects the null hypothesis in \eqref{rele} for large values of the statistic
\begin{equation} \label{estsup}
    \hat d_\infty=\sup_{s\in\I}|\hat g_{n}(s)+M_\lambda\hat g_n(s)-g_*(s)|.
\end{equation}
Compared with the simultaneous confidence bands approach, the construction of an asymptotic level $\alpha$ test based on $ \hat d_\infty$
 is substantially more difficult due to additional nuisance components in the asymptotic distribution of the difference 
$ \hat d_\infty -   d_\infty$ under the null hypothesis in \eqref{rele}.

For a precise description of a test for the relevant hypotheses \eqref{rele} based on $\hat d_\infty$, let $$\EE^{\pm}=\big\{s\in\I: g _0(s)- g _*(s)=\pm d_\infty\big\}$$
denote the set of {\it  extremal points} where the function   $ g _0- g _*$ attains it sup-norm (the set  $\EE^{+}$)
or its negative sup-norm  (the set $\EE^{-}$).
Here we take the convention that $\EE^{\pm}=\I$ if $d_\infty=0$. 
It turns out that the asymptotic distribution of $\hat d_\infty$ (after normalization) depends on the extremal sets $\EE^\pm$. Therefore, 
we propose to estimate these sets  by
\begin{align}\label{extremal}
&\hat{\EE}_n^{\pm}=\Big\{s\in\I:\pm\{\hat g _{n}(s)+M_\lambda\hat g_n(s)- g _*(s)\}\geq(1-n^{-1/2})\hat d_\infty\Big\}.
\end{align}
Combining a multiplier  bootstrap  similar to the one introduced in Section~\ref{sec:band}, for the estimates $\G_{n,b}^*$ in \eqref{gnqstar}, we define the bootstrap statistics 
\begin{align}\label{hattnq}
\hat T_{\EE,n,b}^*=\max\Big\{\sup_{s\in\hat{\EE}_n^+}\G_{n,b}^*(s)\,,\sup_{s\in\hat{\EE}_n^-}\{-\G_{n,b}^*(s)\}\Big\}\,,\qquad 1\leq b\leq B ~.
\end{align}
Then, an asymptotic level-$\alpha$ test is defined by rejecting $H_0$ in \eqref{rele}, 
whenever 
\begin{align}
\label{hd11}
\hat d_\infty=\sup_{s\in\I}|\hat g _{n}(s)+M_\lambda\hat g _{n}(s)- g _*(s)|>\Delta+\frac{\mathcal Q_{1-\alpha}(\hat T_{\EE,n,B}^*)}{\sqrt{n}\lambda^{1/(4m)}}\,,
\end{align}
where $\mathcal Q_{1-\alpha}(\hat T_{\EE,n,B}^*)$ denotes the empirical $(1-\alpha)$ quantile of the bootstrap sample $\{\hat T_{\EE,n,b}^*\}_{b=1}^B$. 
We summarize this procedure  in Algorithm~\ref{algo:rt}. Its consistency is stated in the following proposition, which is proved in Section~\ref{app:thm:extreme} of the online supplement.

\begin{algorithm}[tb!]
\caption{Multiplier bootstrap  test for relevant hypotheses \label{algo:rt}}
\medskip

\begin{enumerate}[nolistsep,leftmargin=3.7em]
\item[\textbf{Input:}] Data $\{(Y_i,X_i,Z_i)\}_{i=1}^n$; target index link function $g_*$; threshold $\Delta$; regularization parameter $\lambda$; nominal level $\alpha$; bootstrap sample size $B$.
\end{enumerate}

\begin{enumerate}

\item Compute the smoothing spline estimator $(\hat g_n,\hat\beta_n,\hat\gamma_n)$ in \eqref{hatwn0}.

\item Compute $\hat d_\infty=\sup_{s\in\I}|\hat g_{n}(s)+M_\lambda\hat g_n(s)-g_*(s)|$ and the empirical extremal sets
\begin{align*}
\hat\EE^{\pm}_n=\Big\{s\in\I:\pm\{\hat g _{n}(s)+M_\lambda\hat g_n(s)- g _*(s)\}\geq(1-n^{-1/2})\hat d_\infty\Big\}.
\end{align*}

\item Generate i.i.d.~bootstrap weights $\{W_{i,b}\}_{1\leq i\leq n,1\leq b\leq B}$ independent of the data as in Algorithm~\ref{algo:band}. Compute the multiplier bootstrap estimators $\G_{n,b}^*$ in \eqref{gnqstar}.

\item Compute the multiplier bootstrap quantile estimates
\begin{align*}
\mathcal Q_{1-\alpha}(\hat T_{\EE,n,B}^*)=\inf\bigg\{\tau\in\mathbb R:\frac{1}{B}\sum_{b=1}^B\mathbf{1}\Big\{\sup_{s\in\hat{\EE}^+_n}\G_{n,b}^*(s)\leq \tau,\inf_{s\in\hat\EE_n^-}\mathbb G_{n,b}^*(s)\geq-\tau\Big\}\geq 1-\alpha\bigg\}.
\end{align*}

\end{enumerate}

{\bf Output:} Rejection of the null hypothesis in \eqref{rele} at nominal level $\alpha$ if
\begin{align}\label{d2}
\hat d_\infty>\Delta+\frac{\mathcal Q_{1-\alpha}(\hat T_{\EE,n,B}^*)}{\sqrt{n}\lambda^{1/(4m)}}\,.
\end{align} 
\end{algorithm}

\begin{proposition}\label{thm:rt:boot}
Suppose Assumptions~\ref{a:rank}--\ref{a:b} and \eqref{newcond} hold. Then, the decision rule \eqref{d2} in Algorithm~\ref{algo:rt} defines 
a consistent asymptotic level $\alpha$ test for the hypotheses \eqref{rele}. More precisely,
\begin{align}\label{eq:rele:boot}
\lim_{n,B\to\infty}\P\bigg\{\hat d_\infty>\Delta+\frac{\mathcal Q_{1-\alpha}(\hat T_{\EE,n,B}^*)}{\sqrt{n}\lambda^{1/(4m)}}\bigg\}=
\left\{\begin{array}{ll}
\vspace{-0.5em}0 & \quad \text{if}\ d_\infty<\Delta\\
\vspace{-0.5em}\alpha &\quad \text{if}\ d_\infty=\Delta\\
1 &\quad \text{if}\ d_\infty>\Delta
\end{array}\right. ~.
\end{align}
\end{proposition}

\begin{remark} \label{remthresh}
A crucial step in testing relevant hypotheses regarding the index link function is the specification of the threshold $\Delta$, which is application-specific. As this choice is often difficult we note that the hypotheses  in \eqref{rele} are nested. It is clear that  the function $\Delta \to  \hat d_\infty   - \Delta$  is  decreasing. 
  Consequently,  rejecting $H_0$ 
  by the test \eqref{eq:rele:boot}
  for $\Delta= \Delta_0$ also yields (asymptotically) rejection of $H_{0}$ 
  for all  $\Delta\leq \Delta_0$.
  We can therefore  simultaneously test the  hypotheses  in \eqref{rele}  for different thresholds $\Delta \geq 0$, 
  starting at $\Delta  = 0$ and 
   increasing  $\Delta $ to 
   find the minimum value of $\Delta $, say 
   $
 \hat \Delta_\alpha$, 
   for which 
   $H_0$ is not rejected with a controlled type I error $\alpha$.
 In particular, $\hat \Delta_\alpha $ could be interpreted as a measure of evidence against the null hypothesis in \eqref{rele}. We will illustrate this approach in Section \ref{sec:realdata}.
\end{remark}

\subsection{Joint hypothesis testing}\label{sec:joint}

We will use the  joint  weak convergence established in Theorem~\ref{thm:asymp} 
to construct a test for  the  joint hypotheses 
\begin{align}\label{hp}
H_0: g_0(x\trans\beta_0)+z\trans\gamma_0=y,\qquad\text{versus}\qquad H_1: g_0(x\trans\beta_0)+z\trans\gamma_0\neq y
\end{align}
regarding all model parameters  $(g_0,\beta_0,\gamma_0)$ in the single-index model \eqref{model}, where the vector $(x, y, z)$ is  given.
In view of the bias adjustment in Lemma~\ref{prop:bias}, a reasonable test for \eqref{hp} is based on the statistic defined by
\begin{align}\label{Tn}
\hat T_n=\sqrt n\lambda^{1/(4m)}\big\{\hat g_n(x\trans\hat\beta_n)+M_\lambda\hat g_n(x\trans\hat\beta_n)+z\trans\hat\gamma_n-y\big\}.
\end{align}
The following proposition states the asymptotic distribution of the test statistic defined in \eqref{Tn}, and is proved in Section~\ref{app:prop:tn} of the supplement.

\begin{proposition}\label{prop:tn}
Suppose Assumptions~\ref{a:rank}--\ref{a:b} and \eqref{newcond} are satisfied. It holds that
\begin{align} \nonumber 
\left\{\sigma_{(x\trans\beta_0)}^2+\lambda^{1/(2m)}\Bigg[\begin{matrix}
g_0'(x\trans\beta_0) Q_{\beta_0}\trans x\\
z
\end{matrix}\Bigg]\trans\Omega_{\beta_0}^{-1}\Bigg[\begin{matrix}
g_0'(x\trans\beta_0) Q_{\beta_0}\trans x\\
z
\end{matrix}\Bigg]\right\}^{-1/2}\hat T_n\converged N(0,1),
\end{align}
where $\sigma^2_{(\cdot)}$ and $\Omega_{\beta_0}$ are defined in \eqref{sigmas2} and \eqref{Omega}, respectively.
\end{proposition}

By Proposition~\ref{prop:tn}, a consistent decision rule could be defined by rejecting the null hypothesis at nominal level $\alpha$ if
\begin{equation}
\label{tn}
\big(\sigma_{(x\trans\beta_0)}^2+\lambda^{1/(2m)}v_{x,z}\trans\Omega_{\beta_0}^{-1}v_{x,z}\big)^{-1/2}|\hat T_n|>\mathcal Q_{1-\alpha/2},
\end{equation}
where $v_{x,z}:=\big[g_0'(x\trans\beta_0)x\trans Q_{\beta_0};z\trans\big]\trans\in\mathbb R^{p+q-1}$ and $\mathcal Q_{1-\alpha/2}$ denotes the $(1-\alpha/2)$-quantile of the standard normal distribution.
However, the normalizing quantity in \eqref{tn} involves unknown quantities such as the derivative of the true unknown index link function $g_0'$ and the unknown functions $R_X,R_Z$ defined in \eqref{vxz} (see the definitions of $\Omega_{1}$ and $\Omega_3$ in equation \eqref{Omega123}), and the unknown asymptotic variance $\sigma_{(\cdot)}^2$ in Theorem~\ref{thm:asymp}, which in practice can be either intractable or difficult to estimate. We therefore propose to combine the multiplier bootstrap approach similar to the ones developed in Sections~\ref{sec:band} and \ref{sec:rele} to avoid the estimation of these nuisance parameters. This procedure is described in Algorithm~\ref{algo:joint}. The validity  of the test is stated in the following proposition. Note that its proof follows similar arguments as that of Propositions~\ref{thm:cb:boot} and \ref{thm:rt:boot}, and is therefore omitted for brevity.

\begin{proposition}
Suppose Assumptions~\ref{a:rank}--\ref{a:b} and \eqref{newcond} are satisfied, 
then, the decision rule in Algorithm~\ref{algo:joint} defines a consistent asymptotic level-$\alpha$ test for the hypotheses \eqref{hp}.
\end{proposition}

\begin{algorithm}[tb!]
\caption{Multiplier bootstrap joint hypotheses tests\label{algo:joint}}
\medskip

\begin{enumerate}[nolistsep,leftmargin=3.7em]
\item[\textbf{Input:}] Data $\{(Y_i,X_i,Z_i)\}_{i=1}^n$; constants $(x,y,z)$ in \eqref{hp}; nominal level $\alpha$; regularization parameter $\lambda$; bootstrap sample size $B$.
\end{enumerate}

\begin{enumerate}

\item Generate i.i.d.~bootstrap weights $\{W_{i,b}\}_{1\leq i\leq n,1\leq b\leq B}$ as in Algorithm~\ref{algo:band}.

\item Compute the smoothing spline estimator $(\hat g_n,\hat\beta_n,\hat\gamma_n)$ defined in \eqref{hatwn0} and the multiplier bootstrap estimators $(\hat g_{n,b}^*,\hat\beta_{n,b}^*,\hat\gamma_{n,b}^*)$, for $1\leq b\leq B$, defined in \eqref{hatbetastar}.

\item Compute $\hat T_n=\sqrt n\lambda^{1/(4m)}\big\{\hat g_n(x\trans\hat\beta_n)+M_\lambda\hat g_n(x\trans\hat\beta_n)+z\trans\hat\gamma_n-y\big\}$. 

\item Compute the multiplier bootstrap test statistics
\begin{align}\label{tnb}
\hat T_{n,b}^*& =\sqrt n\lambda^{1/(4m)}\big\{\hat g_{n,b}^*(x\trans\hat\beta_{n,b}^*)+M_\lambda\hat g_{n,b}^*(x\trans\hat\beta_{n,b}^*)+z\trans\hat\gamma_{n,b}^*-y\big\}.
\end{align}

\item Compute the multiplier bootstrap quantile estimates
\begin{align*}
\mathcal Q_{1-\alpha/2}(\hat T_{n,B}^*)=\inf\Bigg\{\tau\in\mathbb R:\frac{1}{B}\sum_{b=1}^B\mathbf{1}\{\hat T_{n,b}^*\leq \tau\}\geq 1-\alpha/2\Bigg\}.
\end{align*}

\end{enumerate}

{\bf Output:} Rejection of the null hypothesis in \eqref{hp} at nominal level $\alpha$ if $|\hat T_n|>\mathcal Q_{1-\alpha/2}(\hat T_{n,B}^*).$

\end{algorithm}

\section{Finite sample properties}\label{sec:finite}

\subsection{Simulation study}

Our simulation study consists of two parts. We first investigate the finite-sample properties of the  smoothing spline estimator 
\eqref{hatwn0} in terms of estimation accuracy. 
Then, we present simulation results for the three  statistical inference problems in Section~\ref{sec:stat}. All results presented in this section are based on 500 Monte-Carlo samples, and for the  multiplier bootstrap, we took the bootstrap sample size $B=200$.

We use an open-source Matlab software \texttt{Chebfun} \citep[][available at \nolinkurl{https://www.chebfun.org/}]{driscoll2014chebfun}  to obtain the eigensystem in Proposition~\ref{prop:eigen}.
For choosing the regularization parameter $\lambda$ in \eqref{hatwn0}, we propose to adopt the well-known generalized cross-validation approach \cite[GCV, see, for example,][]{whaba1990,green1993nonparametric,ruppert2003semiparametric}, by minimizing the GCV score
\begin{align}\label{gcv}
{\rm GCV}(\lambda)=\frac{1}{n-{\rm tr}\{H(\lambda)\}}\sum_{i=1}^n\{Y_i-\hat g_{n,\lambda}(X_i\trans\hat\beta_{n,\lambda})-M_\lambda\hat g_{n,\lambda}(X_i\trans\hat\beta_{n,\lambda})-Z_i\trans\hat\gamma_{n,\lambda}\}^2,
\end{align}
where $H(\lambda)$ is known as the ``hat matrix'' \cite[see, for example, Section~3.2 in][]{green1993nonparametric} and  $(\hat g_{n,\lambda},\hat\beta_{n,\lambda},\hat\gamma_{n,\lambda})$ is the smoothing spline estimator with the regularization parameter $\lambda$. We refer to \cite{ruppert2003semiparametric} for a comprehensive discussion on tuning parameter selection for single-index models. We refer to Section~\ref{app:numerical} of the online supplement for a detailed description of our implementation and the description of our algorithm.

\medskip

\begin{figure}
\centering
\includegraphics[width=0.32\textwidth]{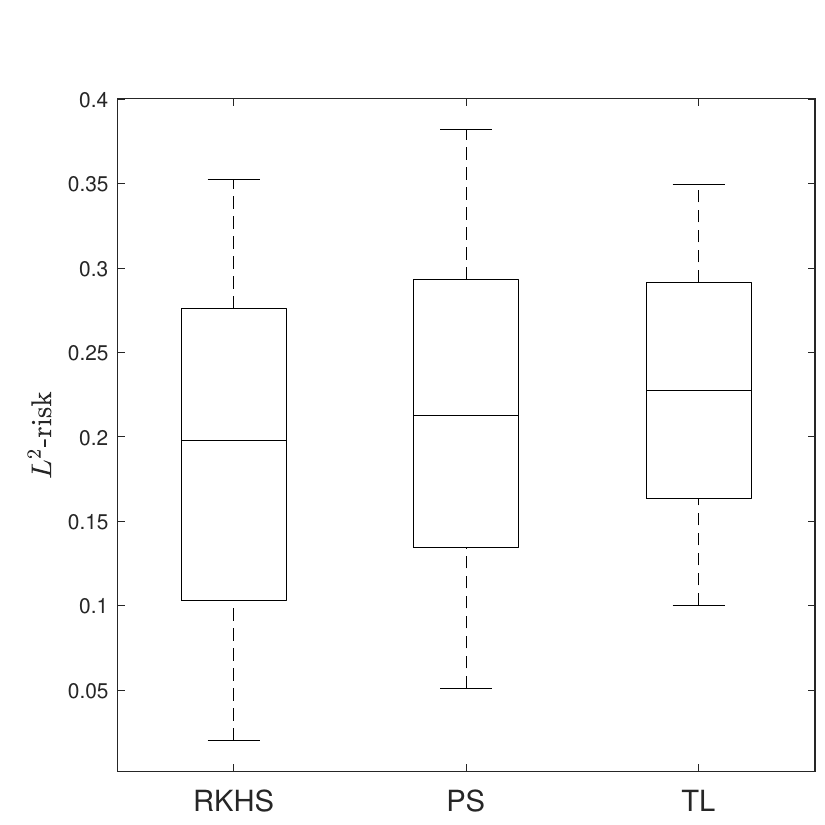}
\includegraphics[width=0.32\textwidth]{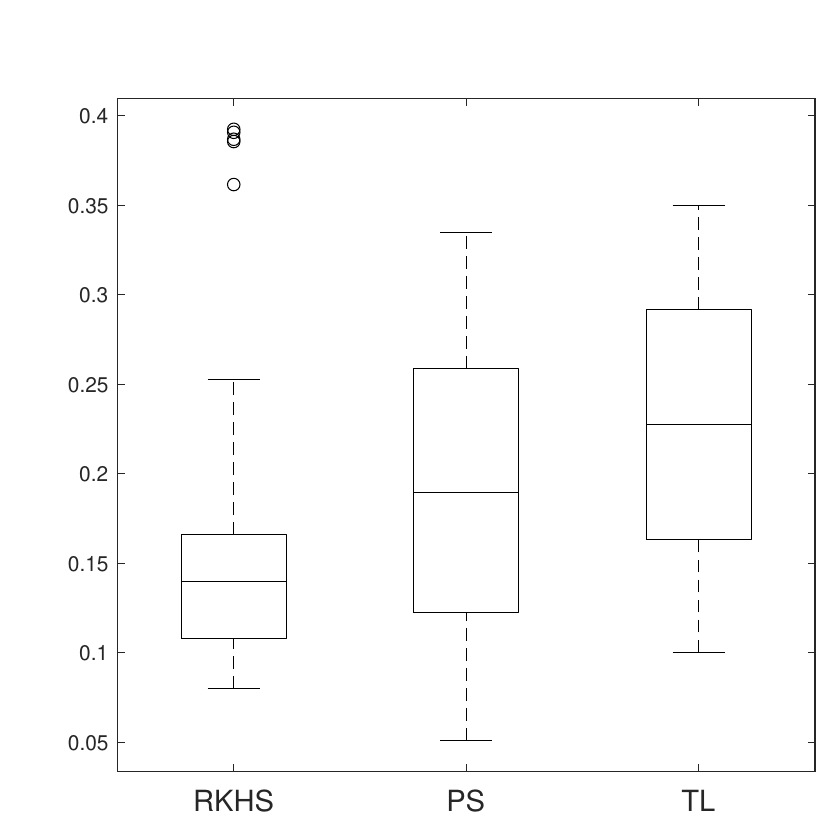}
\includegraphics[width=0.32\textwidth]{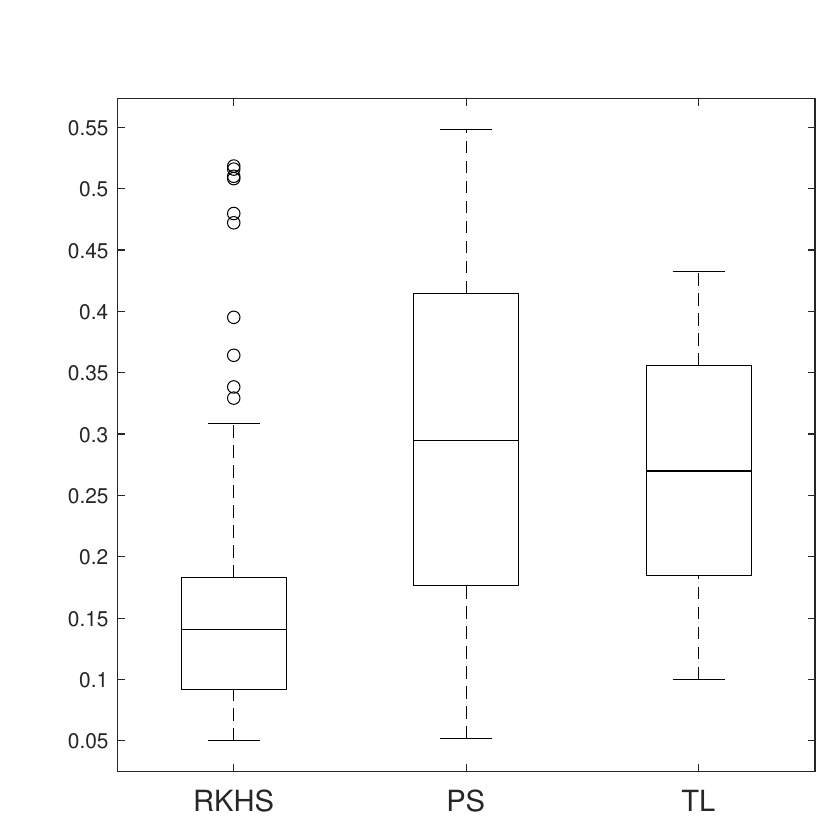}
\caption{\em 
Empirical  $L^2$-risk $\mathfrak R$ of different estimators under the different error settings (1) -- (3) in \eqref{simmod1} (left--right).
\label{fig:2}}
\end{figure}

\noindent
\textbf{Estimation accuracy.}
We compare the estimation accuracy of our approach (denoted by RKHS) with two competing methods for partially linear single-index models in the literature: the P-spline approach in \cite{yu2002} (denoted by PS), and the two-stage local smoothing approach in \cite{wang2010} (denoted by TL). We evaluate the estimation accuracy of the different estimators via the $L^2$-criteria  defined in \eqref{l2risk} considering  a similar setting as in \cite{wang2010} and \cite{kuchibhotla2020}.
To be specific, we consider model \eqref{model} with $p=6$ and $q=1$
and index link function 
$g_0(s)=s^2$. The predictor
$X=[X_1, \ldots, X_6]\trans$ is generated according to the following setting: 
\begin{equation}
    \label{simmod1}
    \begin{split}
& X_1,X_2 \text{ independent}  \sim  \text{Uniform}[-1,1]  ,\\
&X_3=0.2 X_1+0.2(X_2+2)^2+0.2 U_1, \\ &X_4=0.1+0.1(X_1+X_2)+0.3(X_1+1.5)^2+0.2 U_2, \\ 
& X_5 \sim \operatorname{Bernoulli}\big(e^{X_1} /(1+e^{X_1})\big), \\
& X_6 \sim \operatorname{Bernoulli}\big(e^{X_2} /(1+e^{X_2})\big),
\end{split}  
\end{equation}
where
$U_1,U_2$ are independent ${\rm Uniform}[-1,1]$-distributed random variables. In addition, $\beta_0=\tilde\beta_0/\|\tilde\beta_0\|_2$, where $\tilde\beta_0=[1.3,-1.3,1,-0.5,-0.5,-0.5]\trans$; $\gamma_0=1$ and $Z=2\tilde Z-1$, where $\tilde Z\sim{\rm Bernoulli}\big(e^{X\trans\beta_0} /(1+e^{X\trans\beta_0})\big)$. In this case we have $\I=[-1.04,1.00]$.
For  the error in model \eqref{model}
we consider three cases
\begin{equation}
\label{simmod2}
    \begin{split}
 (1) & ~~~ \epsilon \sim N(0,1); \\
 (2) & ~~~ \epsilon | X,Z \sim N\big(0, \log \{2+(X\trans \beta_0)^2+Z\gamma_0\}\big); \\ 
(3) & ~~~  \epsilon | \xi \sim(-1)^{\xi} \operatorname{Beta}(2,3)      ; 
    \end{split}      
\end{equation}
where $\xi \sim \operatorname{Bernoulli}\big(e^{X\trans\beta_0+Z\gamma_0} /(1+e^{X\trans\beta_0+Z\gamma_0})\big)$ and Beta$(a,b)$ denotes a Beta distribution with parameters $a,b>0$.  In 
Figure~\ref{fig:2} we present   boxplots of the simulated  $L^2$-risk $\mathfrak R$ for the different estimators which demonstrate some advantages of the RKHS approach proposed in this paper  over the  two other methods.

\medskip 

\noindent
\textbf{Simultaneous confidence bands and relevant hypotheses.} We investigate the finite sample properties  of the multiplier bootstrap simultaneous confidence band in Algorithm~\ref{algo:band}, where we use the same setting as in \eqref{simmod1} and \eqref{simmod2}. The empirical coverage probabilities of the  $90\%$ and $95\%$ simultaneous  confidence bands are displayed in  Table~\ref{table:cp}
under the three error settings
\eqref{simmod2} for  sample sizes $n=200,1000$.

\begin{figure}
\centering
\includegraphics[width=0.32\textwidth]{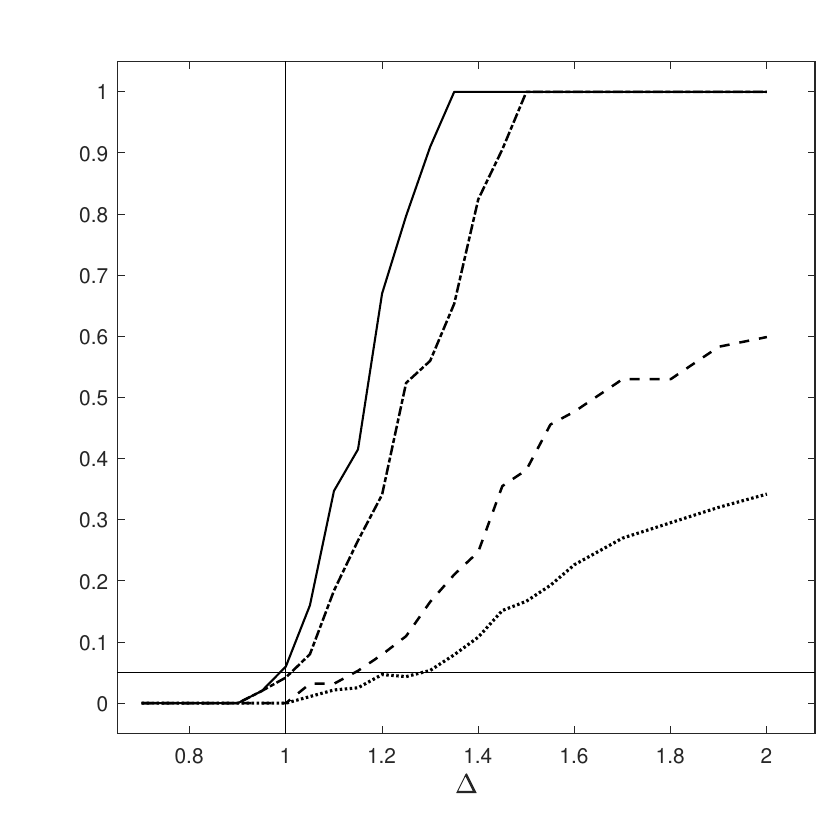}
\includegraphics[width=0.32\textwidth]{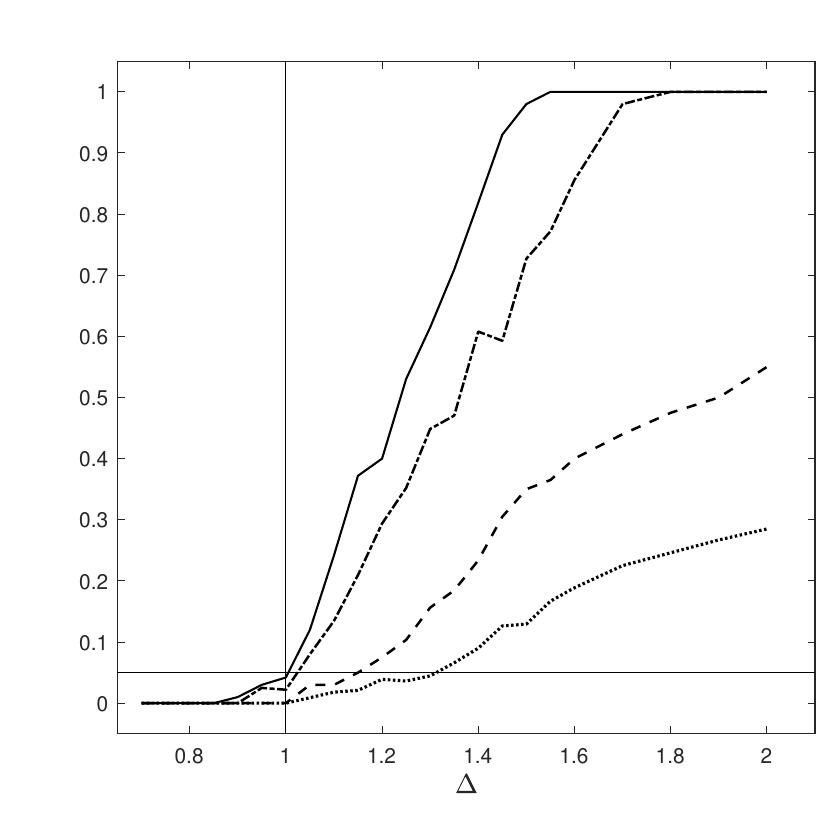}
\includegraphics[width=0.32\textwidth]{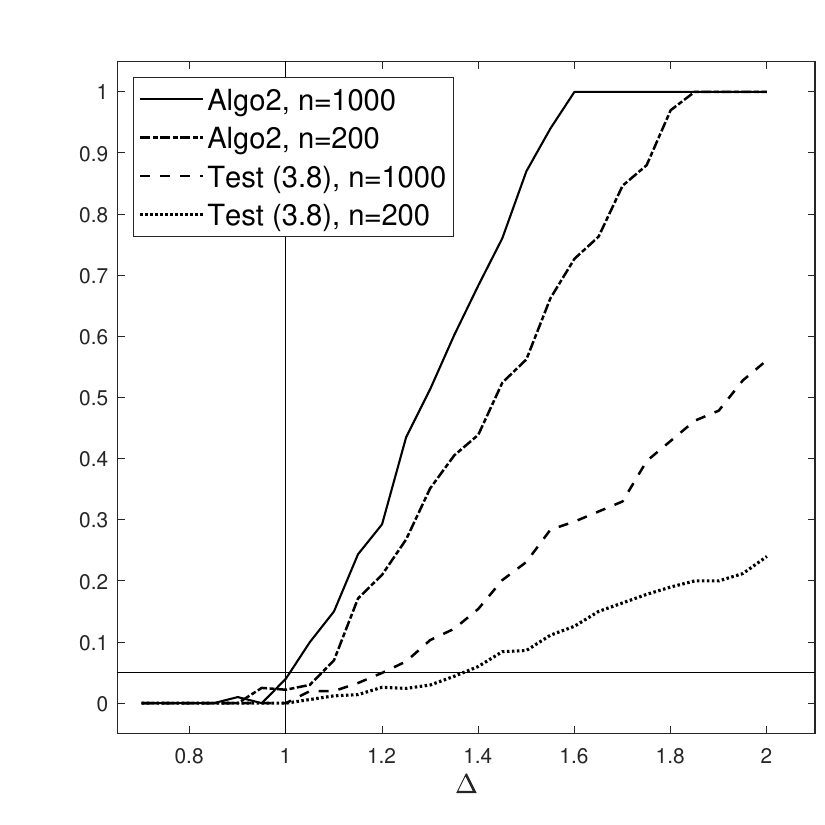}
\caption{\em 
Empirical rejection probabilities  of the test \eqref{rele:band} and the test  in Algorithm~\ref{algo:rt} for the relevant hypotheses \eqref{rele} with sample sizes $n=200,1000$.
\label{fig:3}}
\end{figure}

Next we investigate  the finite sample properties of the 
 test \eqref{rele:band}  and the bootstrap  test in Algorithm~\ref{algo:rt} for the relevant  hypotheses \eqref{rele} with $g_*\equiv 0$.  
 The empirical rejection probabilities of both tests are 
displayed in Figure~\ref{fig:3} for  various choices of the threshold $\Delta$ 
in \eqref{rele}. The  horizontal line corresponds  the nominal level $\alpha=0.05$ while  the vertical line corresponds to $\Delta=\|g_0\|_\infty$. The simulation results 
for the test  in Algorithm~\ref{algo:rt}
are in accordance with the theoretical findings in  Propositions~\ref{thm:cb:boot} and \ref{thm:rt:boot}. 
More precisely, in the ``interior'' of the null hypothesis ($d_\infty < \Delta$) the rejection probabilities quickly approach $0$, the simulated  level is close to $\alpha$ at the ``boundary''  ($d_\infty =\Delta$) of the hypotheses, and under the alternative the test has good power.
Moreover, we observe that this test  
outperforms  the   test \eqref{rele:band} which is based on the duality between confidence regions and hypotheses testing.

\medskip 

\noindent
\textbf{Testing of joint hypotheses.} For the setting 
\eqref{simmod1} and \eqref{simmod2} we applied the multiplier bootstrap  tests in Algorithm~\ref{algo:joint}  for the joint hypotheses in \eqref{hp}
with $x=[0,0,1,0,1,1]\trans$, $z=1$, and took $y=1$ and $0$, corresponding to the case of true $H_0$ and $H_1$, respectively. 
The empirical rejection probabilities of the test 
in Algorithm \ref{algo:joint}
under the various settings are displayed in  Table~\ref{table:joint}. We observe   a  good approximation of the nominal level under the null hypothesis and reasonable empirical power under the alternative hypothesis.

\begin{table}
\centering
\begin{small}
\begin{tabular}{c|cc|cc|cc} 
\cline{1-7}
Setting&\multicolumn{2}{c|}{(1)}&\multicolumn{2}{c|}{(2)}&\multicolumn{2}{c}{(3)}\\
\cline{1-7}
$\alpha$&0.10&0.05&0.10&0.05&0.10&0.05\\
\cline{1-7}
$n=200$&0.844&0.982&0.834&0.938&0.882&0.974\\
$n=1000$&0.872&0.932&0.918&0.974&0.888&0.940\\
\cline{1-7}
\end{tabular} \caption{\it Empirical coverage probabilities  of the simultaneous confidence band in Algorithm~\ref{algo:band}.
\label{table:cp}}
\end{small}
\end{table}

\begin{table}
\centering
\begin{small}
\begin{tabular}{c|cc|cc|cc|cc|cc|cc} 
\cline{1-13}
$y$&\multicolumn{6}{c|}{1 (true $H_0$)}&\multicolumn{6}{c}{0 (true $H_1$)}\\
\cline{1-13}
Setting&\multicolumn{2}{c|}{1}&\multicolumn{2}{c|}{2}&\multicolumn{2}{c|}{3}&\multicolumn{2}{c|}{1}&\multicolumn{2}{c|}{2}&\multicolumn{2}{c}{3}\\
\cline{1-13}
$\alpha$&0.10&0.05&0.10&0.05&0.10&0.05 &0.10&0.05&0.10&0.05&0.10&0.05\\
\cline{1-13}
$n=200$&0.062&0.054&0.128&0.032& 0.094&0.060 &0.900 &0.820 & 0.888&  0.804 & 0.832& 0.786\\
$n=1000$&0.084&0.062& 0.082&0.056&0.116&0.028 &0.942  & 0.864 & 0.920&0.854 &0.914 &0.802 \\
\cline{1-13}
\end{tabular} \caption{\it Empirical rejection probabilities  of the  test in Algorithm~\ref{algo:joint} for the joint hypotheses \eqref{hp} .
\label{table:joint}}
\end{small}
\end{table}

\subsection{Real data example}\label{sec:realdata}

\begin{figure}
\begin{minipage}{0.5\linewidth}
\centering
\includegraphics[width=1\textwidth]{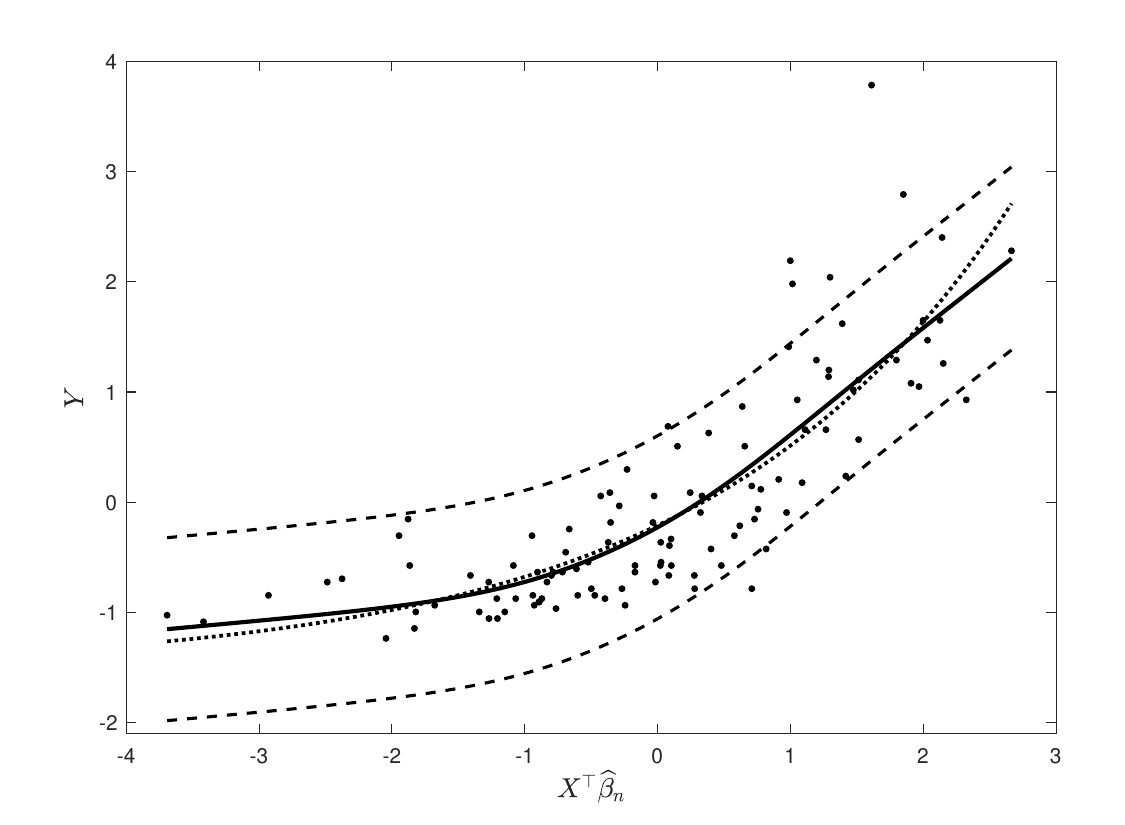}
\end{minipage}
~~~~~~~~~~~~~~~~~~
\begin{minipage}{0.5\linewidth}
\def\arraystretch{1.5}
\begin{tabular}{c|ccc}
\cline{1-4}
$\Delta\backslash\alpha$ &$0.01$ &$0.05$ &$0.10$ \\\cline{1-4}
0.38 &R &R &R \\ 
0.39 &-- &R &R \\\cline{1-4}
0.53 &-- &R &R \\ 
0.54 &-- &-- &R \\\cline{1-4}
0.67 &-- &-- &R \\ 
0.68 &-- &-- &-- \\
\cline{1-4}
\end{tabular}
\end{minipage}
\caption{\em 
Left: the estimated index link function (solid line) and its 95\%-simultaneous confidence band (dashed line); dotted line: fitted index link function $g_*(s)=1.31\exp(0.44s)-1.52$ for the parametric model \eqref{para}. Right: decision boundaries for relevant hypotheses \eqref{rele}.
\label{fig:4}}
\end{figure}

We illustrate the finite sample properties of the  new smoothing spline estimator and the new inference tools via a simulation study on the NYC air pollution data consisting of measurements of three covariates, namely solar radiation, temperature and wind speed, for the response ozone concentration for $n=111$ days, available in the R package \texttt{EnvStats} \citep{millard2013package}. We follow the analysis in \cite{yu2002,kuchibhotla2020} and  use a conventional single-index model where the single-index $X\trans\beta_0$ consists of the above three covariates, and $Z\equiv0$ in \eqref{model}; the data are considered independent, and we standardized all  variables such that they have mean $0$  and  variance $1$. We applied the smoothing spline estimator to the dataset to obtain the estimated index direction estimator $\hat\beta_n=[0.39, 0.67, -0.63]\trans$, and the index link estimator together with its 95\%-simultaneous confidence band are displayed in Figure~\ref{fig:4}. 
The width of the simultaneous confidence band is $1.66$. 

We also fitted  a parametric model 
\begin{equation}\label{para}
g(s)=c_1\exp(c_2s) +c_3,  ~~(c_1,c_2,c_3\in\mathbb R)  
\end{equation}
for the index link function $g_0$ to the data.  The resulting  estimators are $\hat c_1=1.31$, $\hat c_2=0.44$, and $\hat c_3=-1.52$, and the corresponding index direction estimator is given by $[0.40, 0.66, -0.63]\trans$. The average sum of squared residues for the nonparametric single-index model and the parametric model in \eqref{para} are 0.29 and 0.31, respectively. 

For the relevant hypotheses \eqref{rele}, we set $g_*$ to be the fitted parametric model $g_*(s)=1.31\exp(0.44s)-1.52$. We applied the multiplier bootstrap test in Algorithm~\ref{algo:rt} to obtain the decision boundaries 
$\hat \Delta_\alpha$ for nominal levels $\alpha=0.01,0.05,0.10$ (see Remark \ref{remthresh}). The results are shown in the right part of Figure~\ref{fig:4}, where ``R'' and ``--'' stand for rejection and no rejection, respectively. The values  in this table are interpreted as follows. First, there is no evidence to deduce $\|g_0-g_*\|_\infty>0.68$ under $\alpha=0.10$, but strong evidence for $\|g_0-g_*\|_\infty>0.38$ ($\alpha=0.01$). In addition, the relevant uniform difference $\|g_0-g_*\|_\infty\leq\Delta$ should be rejected at level $\alpha\geq0.05$ for $\Delta\in[0.39,0.53]$. On the other hand, for $\Delta\in[0.54,0.67]$, it should only be rejected at level $\alpha\geq0.10$, which indicates weaker evidence  of the alternative hypothesis that $\|g_0-g_*\|_\infty>\Delta$.

\baselineskip=11pt
\vspace{1cm}

\begin{center}
{\large\textbf{Acknowledgements}}
\end{center}
\noindent 
This work was partially supported by the DFG Research unit 5381 \textit{Mathematical Statistics in the Information Age}, project number 460867398.

\baselineskip=11pt
\vspace{1cm}

\begingroup
\renewcommand{\section}[2]{\subsection#1{#2}}
{\centering
\bibliographystyle{apalike}
\bibliography{reference.bib}

\begin{thebibliography}{}

\bibitem[Agmon, 2010]{agmon2010}
Agmon, S. (2010).
\newblock {\em Lectures on elliptic boundary value problems (Vol. 369)}.
\newblock American Mathematical Soc.

\bibitem[Antoniadis et~al., 2004]{antoniadis2004bayesian}
Antoniadis, A., Gr{\'e}goire, G., and McKeague, I.~W. (2004).
\newblock Bayesian estimation in single-index models.
\newblock {\em Statistica Sinica}, 14:1147--1164.

\bibitem[Berlinet and Thomas-Agnan, 2011]{berlinet2011reproducing}
Berlinet, A. and Thomas-Agnan, C. (2011).
\newblock {\em Reproducing kernel Hilbert spaces in probability and
  statistics}.
\newblock Springer Science \& Business Media.

\bibitem[B{\"u}cher and Kojadinovic, 2019]{bucher2019note}
B{\"u}cher, A. and Kojadinovic, I. (2019).
\newblock A note on conditional versus joint unconditional weak convergence in
  bootstrap consistency results.
\newblock {\em Journal of Theoretical Probability}, 32(3):1145--1165.

\bibitem[Carroll et~al., 1997]{carroll1997}
Carroll, R.~J., Fan, J., Gijbels, I., and Wand, M.~P. (1997).
\newblock Generalized partially linear single-index models.
\newblock {\em Journal of the American Statistical Association}, 92:477--489.

\bibitem[Cesari, 2012]{cesari2012asymptotic}
Cesari, L. (2012).
\newblock {\em Asymptotic behavior and stability problems in ordinary
  differential equations}, volume~16.
\newblock Springer Science \& Business Media.

\bibitem[Chang et~al., 2010]{chang2010asymptotically}
Chang, Z., Xue, L., and Zhu, L. (2010).
\newblock On an asymptotically more efficient estimation of the single-index
  model.
\newblock {\em Journal of Multivariate Analysis}, 101(8):1898--1901.

\bibitem[Cheng and Shang, 2015]{cheng2015}
Cheng, G. and Shang, Z. (2015).
\newblock Joint asymptotics for semi-nonparametric regression models with
  partially linear structure.
\newblock {\em Annals of Statistics}, 43:1351--1390.

\bibitem[Coddington et~al., 1956]{coddington1956theory}
Coddington, E.~A., Levinson, N., and Teichmann, T. (1956).
\newblock {\em Theory of ordinary differential equations}.
\newblock American Institute of Physics.

\bibitem[Cox and O'Sullivan, 1990]{cox1990asymptotic}
Cox, D.~D. and O'Sullivan, F. (1990).
\newblock Asymptotic analysis of penalized likelihood and related estimators.
\newblock {\em The Annals of Statistics}, pages 1676--1695.

\bibitem[Cui et~al., 2011]{cui2011}
Cui, X., H\"ardle, W.~K., and Zhu, L. (2011).
\newblock The {EFM} approach for single-index models.
\newblock {\em The Annals of Statistics}, 39:1658--1688.

\bibitem[Dette and Kokot, 2021]{dette2021bio}
Dette, H. and Kokot, K. (2021).
\newblock Bio-equivalence tests in functional data by maximum deviation.
\newblock {\em Biometrika}, 108(4):895--913.

\bibitem[Dette et~al., 2020]{dette2020functional}
Dette, H., Kokot, K., and Aue, A. (2020).
\newblock Functional data analysis in the banach space of continuous functions.
\newblock 48(2):1168--1192.

\bibitem[Dong et~al., 2015]{dong2015}
Dong, C., Gao, J., and Peng, B. (2015).
\newblock Semiparametric single-index panel data models with cross-sectional
  dependence.
\newblock {\em Journal of Econometrics}, 188:301--312.

\bibitem[Dong et~al., 2016]{dong2016}
Dong, C., Gao, J., and Tj{\o}stheim, D. (2016).
\newblock Estimation for single-index and partially linear single-index
  integrated models.
\newblock {\em The Annals of Statistics}, 44:425--453.

\bibitem[Driscoll et~al., 2014]{driscoll2014chebfun}
Driscoll, T.~A., Hale, N., and Trefethen, L.~N. (2014).
\newblock {\em Chebfun guide}.
\newblock Pafnuty Publications, Oxford.

\bibitem[Ferreira and Menegatto, 2013]{ferreira2013positive}
Ferreira, J. and Menegatto, V.~A. (2013).
\newblock Positive definiteness, reproducing kernel {Hilbert} spaces and
  beyond.
\newblock {\em Annals of Functional Analysis}, 4(1):64--88.

\bibitem[Gao and Liang, 1997]{gao1997statistical}
Gao, J. and Liang, H. (1997).
\newblock Statistical inference in single-index and partially nonlinear models.
\newblock {\em Annals of the Institute of Statistical Mathematics},
  49:493--517.

\bibitem[Green and Silverman, 1993]{green1993nonparametric}
Green, P.~J. and Silverman, B.~W. (1993).
\newblock {\em Nonparametric regression and generalized linear models: a
  roughness penalty approach}.
\newblock CRC Press.

\bibitem[H\"ardle et~al., 1993]{hardle1993}
H\"ardle, W., Hall, P., and Ichimura, H. (1993).
\newblock Optimal smoothing in single-index models.
\newblock {\em The Annals of Statistics}, 21:157--178.

\bibitem[H{\"a}rdle et~al., 1998]{hardle1998testing}
H{\"a}rdle, W., Mammen, E., and M{\"u}ller, M. (1998).
\newblock Testing parametric versus semiparametric modeling in generalized
  linear models.
\newblock {\em Journal of the American Statistical Association},
  93(444):1461--1474.

\bibitem[H{\"a}rdle et~al., 2001]{hardle2001bootstrap}
H{\"a}rdle, W., Mammen, E., and Proen{\c{c}}a, I. (2001).
\newblock A bootstrap test for single index models.
\newblock {\em Statistics}, 35(4):427--451.

\bibitem[H{\"a}rdle and Stoker, 1989]{hardle1989investigating}
H{\"a}rdle, W. and Stoker, T.~M. (1989).
\newblock Investigating smooth multiple regression by the method of average
  derivatives.
\newblock {\em Journal of the American statistical Association},
  84(408):986--995.

\bibitem[Huh and Park, 2002]{huh2002likelihood}
Huh, J. and Park, B. (2002).
\newblock Likelihood-based local polynomial fitting for single-index models.
\newblock {\em Journal of Multivariate Analysis}, 80(2):302--321.

\bibitem[Ichimura, 1993]{ichimura1993semiparametric}
Ichimura, H. (1993).
\newblock Semiparametric least squares {(SLS)} and weighted sls estimation of
  single-index models.
\newblock {\em Journal of Econometrics}, 58(1-2):71--120.

\bibitem[Kley et~al., 2016]{kley2016quantile}
Kley, T., Volgushev, S., Dette, H., and Hallin, M. (2016).
\newblock Quantile spectral processes: Asymptotic analysis and inference.
\newblock {\em Bernoulli}, 22(3):1770--1807.

\bibitem[Kuchibhotla and Patra, 2020]{kuchibhotla2020}
Kuchibhotla, A.~K. and Patra, R.~K. (2020).
\newblock Efficient estimation in single index models through smoothing
  splines.
\newblock {\em Bernoulli}, 26:1587--1618.

\bibitem[Kuchibhotla et~al., 2021]{kuchibhotla2021}
Kuchibhotla, A.~K., Patra, R.~K., and Sen, B. (2021).
\newblock Semiparametric efficiency in convexity constrained single-index
  model.
\newblock {\em Journal of the American Statistical Association}, 118:272--286.

\bibitem[Liang et~al., 2010]{liang2010estimation}
Liang, H., Liu, X., Li, R., and Tsai, C.-L. (2010).
\newblock Estimation and testing for partially linear single-index models.
\newblock {\em Annals of statistics}, 38(6):3811.

\bibitem[Ma and Zhu, 2013]{ma2013}
Ma, Y. and Zhu, L. (2013).
\newblock Doubly robust and efficient estimators for heteroscedastic partially
  linear single-index models allowing high dimensional covariates.
\newblock {\em Journal of the Royal Statistical Society: Series B (Statistical
  Methodology)}, 75:305--322.

\bibitem[Mammen and van~de Geer, 1997]{mammen1997}
Mammen, E. and van~de Geer, S. (1997).
\newblock Penalized quasi-likelihood estimation in partial linear models.
\newblock {\em The Annals of Statistics}, 25:1014--1035.

\bibitem[Millard, 2013]{millard2013package}
Millard, S.~P. (2013).
\newblock Package for environmental statistics, including us epa guidance.

\bibitem[Newey and Stoker, 1993]{newey1993efficiency}
Newey, W.~K. and Stoker, T.~M. (1993).
\newblock Efficiency of weighted average derivative estimators and index
  models.
\newblock {\em Econometrica}, 61:1199--1223.

\bibitem[Oden and Reddy, 2012]{oden2012}
Oden, J.~T. and Reddy, J.~N. (2012).
\newblock {\em An introduction to the mathematical theory of finite elements}.
\newblock Courier Corporation.

\bibitem[Pollard, 1990]{pollard1990empirical}
Pollard, D. (1990).
\newblock {\em Empirical processes: theory and applications}.
\newblock IMS.

\bibitem[Ruppert et~al., 2003]{ruppert2003semiparametric}
Ruppert, D., Wand, M.~P., and Carroll, R.~J. (2003).
\newblock {\em Semiparametric regression}.
\newblock Number~12. Cambridge university press.

\bibitem[Shang and Cheng, 2013]{shang2013}
Shang, Z. and Cheng, G. (2013).
\newblock Local and global asymptotic inference in smoothing spline models.
\newblock {\em The Annals of Statistics}, 41:2608--2638.

\bibitem[Stoker, 1986]{stoker1986consistent}
Stoker, T.~M. (1986).
\newblock Consistent estimation of scaled coefficients.
\newblock {\em Econometrica}, 54:1461--1481.

\bibitem[Tsybakov, 2009]{tsybakov2009nonparametric}
Tsybakov, A.~B. (2009).
\newblock {\em Introduction to Nonparametric Estimation}.
\newblock Springer.

\bibitem[Tukey, 1991]{tukey1991}
Tukey, J.~W. (1991).
\newblock The philosophy of multiple comparisons.
\newblock {\em Statistical Science}, 6(1):100--116.

\bibitem[van~de Geer, 2000]{vandegeer2000}
van~de Geer, S.~A. (2000).
\newblock {\em Applications of empirical process theory}.
\newblock Cambridge University Press, Cambridge.

\bibitem[van~der Vaart, 2000]{van2000asymptotic}
van~der Vaart, A.~W. (2000).
\newblock {\em Asymptotic statistics}, volume~3.
\newblock Cambridge university press.

\bibitem[van~der Vaart and Wellner, 1996]{vaart1996}
van~der Vaart, A.~W. and Wellner, J.~A. (1996).
\newblock {\em Weak convergence and empirical processes}.
\newblock Springer.

\bibitem[Wahba, 1990]{whaba1990}
Wahba, G. (1990).
\newblock {\em Spline models for observational data}.
\newblock Society for Industrial and Applied Mathematics.

\bibitem[Wang et~al., 2010]{wang2010}
Wang, J.~L., Xue, L., Zhu, L., and Chong, Y.~S. (2010).
\newblock Estimation for a partial-linear single-index model.
\newblock {\em The Annals of Statistics}, 38:246--274.

\bibitem[Wang and Yang, 2009]{wang2009spline}
Wang, L. and Yang, L. (2009).
\newblock Spline estimation of single-index models.
\newblock {\em Statistica Sinica}, 19:765--783.

\bibitem[Wen and Yin, 2013]{wen2013feasible}
Wen, Z. and Yin, W. (2013).
\newblock A feasible method for optimization with orthogonality constraints.
\newblock {\em Mathematical Programming}, 142(1-2):397--434.

\bibitem[Xia, 2007]{xia2007constructive}
Xia, Y. (2007).
\newblock A constructive approach to the estimation of dimension reduction
  directions.
\newblock {\em The Annals of Statistics}, 35:2654--2690.

\bibitem[Xia and Li, 1999]{xia1999single}
Xia, Y. and Li, W.~K. (1999).
\newblock On single-index coefficient regression models.
\newblock {\em Journal of the American Statistical Association},
  94(448):1275--1285.

\bibitem[Xia et~al., 1999]{xia1999}
Xia, Y., Tong, H., and Li, W.~K. (1999).
\newblock On extended partially linear single-index models.
\newblock {\em Biometrika}, 86:831--842.

\bibitem[Yu and Ruppert, 2002]{yu2002}
Yu, Y. and Ruppert, D. (2002).
\newblock Penalized spline estimation for partially linear single-index models.
\newblock {\em Journal of the American Statistical Association}, 97:1042--1054.

\bibitem[Yu et~al., 2017]{yu2017penalised}
Yu, Y., Wu, C., and Zhang, Y. (2017).
\newblock Penalised spline estimation for generalised partially linear
  single-index models.
\newblock {\em Statistics and Computing}, 27:571--582.

\bibitem[Zhu and Xue, 2006]{zhu2006empirical}
Zhu, L. and Xue, L. (2006).
\newblock Empirical likelihood confidence regions in a partially linear
  single-index model.
\newblock {\em Journal of the Royal Statistical Society Series B: Statistical
  Methodology}, 68(3):549--570.

\end{thebibliography}
}
\endgroup

\clearpage

\appendix

\baselineskip=15pt

\begin{center}
{\bf \Large Supplement for ``Simultaneous 
semiparametric inference for single-index models''} \\[1.5em]

Jiajun Tang, Holger Dette 
\end{center}


\baselineskip=15pt

\section*{Organization and notations}
The supplementary material is organized as follows. Section~\ref{app:numerical} contains auxiliary numerical details. Proof of theoretical results in Sections~\ref{sec:definenorm}--\ref{sec:jeq} of the main paper are included in Section~\ref{app:proof}. In Section~\ref{app:C} we give the proofs of the theoretical results regarding Bahadur representations and and joint weak convergence in Sections~\ref{sec:jointbahadur} and \ref{sec:marginalbahadur} of the main paper. Proof of theoretical results in Section~\ref{sec:stat} of the main paper are given in Section~\ref{app:sec:stat}. Finally, Section~\ref{app:aux} provide auxiliary lemmas which are used in Sections~\ref{app:proof}--\ref{app:sec:stat} along with their proofs.

For $v\in\mathbb R^{p\times q}$, write $v^{\otimes2}=vv\trans\in\mathbb R^{p\times p}$. In the sequel, for a vector $v\in\mathbb R^p$, let $\|v\|_\infty$ denote its maximum norm. Let ${\rm span}\{\beta_0\}$ and ${\rm span}\{\beta_0\}^\perp$ denote the linear span of $\beta_0$ and its orthogonal complement.
In the sequel, we use ``$c$'' to denote a generic positive constant that might differ from line to line. Given the data $\{(X_i,Y_i,Z_i)\}_{i=1}^n$, define
\begin{align}\label{elln}\tag{$\star$}
\ell_n(g,\beta,\gamma)=\frac{1}{n}\sum_{i=1}^n\big\{Y_i-g(X_i\trans \beta)-Z_i\trans\gamma\big\}^2,\qquad\quad (g,\beta,\gamma)\in\H_m\times\mathcal B\times\mathbb R^q.
\end{align}

\section{Additional numerical details}\label{app:numerical}

We shall only present our approach  for computing  the bootstrap estimator $(\hat g_{n,b}^*,\hat\beta_{n,b}^*,\hat\gamma_{n,b}^*)$ in \eqref{hatbetastar}, since  the original estimator \eqref{hatwn0} can be viewed as a special case of \eqref{hatbetastar} by taking $W_{i,b}\equiv1$ for any $1\leq i\leq n$ and $1\leq b\leq B$. For the $b$-th multiplier bootstrap estimator, for the moment suppose the regularization parameter $\lambda_b$ is determined, which we will soon introduce how to obtain. Our approach for computing the smoothing spline estimator for the single-index model \eqref{model} consists of the following points:

\begin{itemize}

\item Given an estimate $\tilde\beta$ of $\beta_{0}$, compute the eigen-system in Proposition~\ref{prop:eigen}.

\item 
Based on the estimator $\tilde\beta$, compute the estimator $(\tilde g,\tilde\gamma)$ of $(g_0,\gamma_0)$.

\item  Given the estimator $(\tilde g,\tilde\gamma)$, compute an updated  estimator $\tilde\beta$ of $\beta_0$.

\end{itemize}

Consider the first point. Suppose $\tilde\beta$ is an estimator of $\beta_0$. In view of Proposition~\ref{prop:eigen}, we consider the case of $m=3$ and the ordinary differential equation (ODE) in equation \eqref{id}.
The probability density function $f_{X\trans\beta_0}$ in \eqref{id} is replaced via a kernel density estimator $\hat f_{X\trans\beta_0}$ based on the input data $\{X_i\trans\tilde\beta\}_{i=1}^n$. Then, we find the normalized eigenfunctions $\hat \phi_{j}$ and the eigenvalues $\hat\rho_j$, using the Matlab package \texttt{Chebfun} \citep{driscoll2014chebfun}.

Consider the second point. For a parameter $v$, we propose to approximate $\H_m$ by the linear space  spanned by the leading $v$ eigenvectors $\{\hat\phi_j\}_{j=1}^v$, so for an element $g$ in this space it admits a series expansion $g=\sum_{j=1}^va_j\hat\phi_j$, where $a_j=V(g,\hat\phi_j)$ is the corresponding Fourier coefficient. We propose to choose the parameter $v$ via five-fold cross validation. 
By Proposition~\ref{prop:eigen}, for $g=\sum_{j\geq1}a_j\hat\phi_j$, it is true that $J(g,g)=\sum_{j\geq1}a_j^2\rho_j$.
For $1\leq b\leq B$, for the $b$-th bootstrap estimator $(\hat g_{n,b}^*,\hat\beta_{n,b}^*,\hat\gamma_{n,b}^*)$ and for the bootstrap weights $\{W_{i,b}\}_{i=1}^n$ in Algorithms~\ref{algo:band}--\ref{algo:joint}, let $W_b={\rm diag}(W_{1,b},\ldots,W_{n,b})$ be an $n\times n$ diagonal matrix. We find the $\hat a_{1}^{(b)},\ldots,\hat a_v^{(b)},\hat\gamma_{n,b}^*$ by solving the following minimization problem
\begin{align}\label{lambda}
(\hat a_{1}^{(b)},\ldots,\hat a_v^{(b)},\hat\gamma_{n,b}^*)&=\underset{a_1,\ldots,a_v\in\mathbb R;\gamma\in\mathbb R^q}{\arg\min}\left\{\frac{1}{2n}\sum_{i=1}^nW_{i,b}\bigg[Y_i-\sum_{k=1}^{v}a_{k}\hat\phi_k(X_i\trans\tilde\beta)-Z_i\trans\gamma\bigg]^2+\frac{\lambda_b }{2}  \sum_{k=1}^{v}a_{k}^2\hat \rho_{k}\right\}\notag\\
&=\underset{a\in\mathbb R^v,\gamma\in\mathbb R^q}{\arg\min}\Big\{\big( Y-\hat\Phi a-Z\trans\gamma\big)\trans{W}_b\big( Y-\hat\Phi a-Z\trans\gamma\big)+n\lambda_b a\trans\,\hat\Lambda\, a\Big\}\,,
\end{align}
where we write $\hat a=[\hat a_1,\ldots,\hat a_v]\trans\in\mathbb{R}^v$, $Y=[Y_1,\ldots,Y_n]\trans\in\mathbb R^n$, $Z=[Z_1,\ldots,Z_n]\trans\in\mathbb R^{n\times b}$, $\hat\Phi=\big[\hat\phi_k(X_i\trans\tilde\beta)\big]_{i,k}\in\mathbb R^{n\times v}$, and $\hat\Lambda=\diag(\hat\rho_1,\ldots,\hat\rho_v)\in\mathbb R^{v\times v}$. By direct calculations, for $1\leq\ell\leq v$, we have 
\begin{align}\label{hatb}
\Bigg[\begin{matrix}
\hat a^{(b)}\\\hat \gamma^{(b)}
\end{matrix}\Bigg]=\big\{[\hat\Phi,Z]\trans{W}_b[\hat\Phi,Z]+n\lambda_b\diag(\hat\Lambda,0_{q,q})\big\}^{-1}[\hat\Phi,Z]\trans{W}_bY\,,
\end{align}
so that the bias-adjusted estimator of $g_0$ in view of Lemma~\ref{prop:bias} is defined by
\begin{align}\label{key}
\hat g_{n,b}^*+M_\lambda(\hat g_{n,b}^*)=\sum_{k=1}^v\hat a^{(b)}_k\hat\phi_k+\lambda_b\sum_{k=1}^v\frac{\hat a^{(b)}_k\hat\rho_k\hat\phi_k}{1+\lambda_b\hat\rho_k}.
\end{align}

Consider the third point. Suppose $(\tilde g,\tilde\gamma)$ are estimators of $(g_0,\gamma_0)$. To estimate $\beta_0$ given $(\tilde g,\tilde\gamma)$, we propose to follow the optimization in \cite{wen2013feasible} \citep[see also][]{kuchibhotla2021}, which we refer to for a more detailed derivation and describe now for completion. Recalling the definition of $\ell_n$ in equation~\eqref{elln}, let $\beta$ be an initial value of the optimization, define
$$
\mathfrak b:=\nabla_{\beta} \ell_n(\tilde g,\beta,\tilde\gamma)=\frac{2}{n}\sum_{i=1}^ng'(X_i\trans\beta)\big\{Y_i-g(X_i\trans \beta)-Z_i\trans\gamma\big\}X_i\in \mathbb{R}^p,
$$
where $\nabla_\beta$ denotes the gradient w.r.t.~$\beta$. Define the path $\tau \mapsto \beta_\tau$, where
\begin{align}\label{betatau}
\beta_\tau:=\frac{4+\tau^2\{(\beta\trans \mathfrak b)^2-\|\mathfrak b\|_2^2\}+4\tau \beta\trans \mathfrak b}{4-\tau^2(\beta\trans \mathfrak b)^2+\tau^2\|\mathfrak b\|_2^2} \beta-\frac{4\tau}{4-\tau^2(\beta\trans \mathfrak b)^2+\tau^2\|\mathfrak b\|_2^2}\mathfrak b,\qquad\quad \tau\in\mathbb R.
\end{align}
Find $\tau\in(\tau^-,\tau^+)$ through a line search along the path $\tau\mapsto\beta_\tau$ that minimizes $\ell_n(\tilde g,\beta_\tau,\tilde\gamma)$, where 
\begin{align}\label{taupm}
\tau^\pm:= \frac{2\{\beta\trans \mathfrak b-\mathfrak b(1) / \beta(1)\} \pm 2\sqrt{\{\beta\trans \mathfrak b-\mathfrak b(1) / \beta(1)\}^2+\|\mathfrak b\|_2^2-(\beta\trans \mathfrak b)^2}}{\|\mathfrak b\|_2^2-(\beta\trans \mathfrak b)^2} .
\end{align}
Note that the above three points can be applied recursively, and the entire procedure for computing the smoothing spline estimator for the single-index model is summarized in Algorithm~\ref{algo1}. 

Our approach requires an initial index direction estimator, which we denote by $\hat\beta_n^{\rm (ini)}$, using the original sample $\{(X_i,Y_i,Z_i)\}_{i=1}^n$, which we now introduce how to obtain. For some constant $C$, we first sample from the uniform distribution on the unit sphere $\mathcal S_{p-1}$ random variables $\kappa_1,\ldots,\kappa_C$. Note that this can be achieved by sampling i.i.d.~standard Gaussian random vectors in $\mathbb R^p$ and divide their corresponding $\ell^2$-norm. Then, for each $1\leq j\leq C$, pretend for the moment that $\kappa_j$ is the true index direction and treat the estimation task as a partially linear nonparametric regression problem:
\begin{align}\label{newn}
Y_i=g(X_i\trans\kappa_j)+Z_i\trans\gamma+e_i,\qquad 1\leq i\leq n,
\end{align}
where $(g,\gamma)$ are estimated through a polynomial spline estimator $(\tilde g_j,\tilde\gamma_j)$.
We then propose to take $\hat\beta_n^{\rm(ini)}$ to the value that minimizes the median of the absolute residual, that is,
\begin{align*}
\hat\beta_n^{\rm(ini)}=\kappa_{j_*},\qquad \mbox{where }j_*=\min_{1\leq j\leq C}R_j.
\end{align*}
where $R_j$ denotes the median of absolute residuals $\big\{|Y_i-\tilde g_j(X_i\trans\kappa_j)-Z_i\trans\tilde\gamma_j|\big\}_{i=1}^n$. Note that various selection criteria may be applied, but our simulation results show that our final estimators is robust to the choice of initial index direction estimator.

Generalized cross-validation (GCV, see, for example, \citealp{whaba1990}) is applied to choose the smoothing parameter $\lambda_b$ in \eqref{lambda}. For the $b$-th bootstrap estimator $(\hat g_{n,b}^*,\hat\beta_{n,b}^*,\hat\gamma_{n,b}^*)$, in view of \eqref{gcv}, we choose $\lambda_b$ that minimizes the GCV score
\begin{align*}
\text{GCV}(\lambda_b)=\frac{n^{-1}\Vert Y-\hat Y_b\Vert_{2}^2}{\{1-{\rm{tr}}(H_b)/n\}^2}\,,
\end{align*}
where the predicted response $\hat Y_b$ and the hat matrix $H_b$ are such that
\begin{align*}
&\hat Y_{b}=[\hat\Phi,Z]\big\{[\hat\Phi,Z]\trans{W}_b[\hat\Phi,Z]+n\lambda_b\diag(\hat\Lambda,0_{q,q})\big\}^{-1}[\hat\Phi,Z]\trans{W}_bY,
\\
&{\rm tr}(H_{b})={\rm tr}\Big([\hat\Phi,Z]\big\{[\hat\Phi,Z]\trans{W}_b[\hat\Phi,Z]+n\lambda_b\diag(\hat\Lambda,0_{q,q})\big\}^{-1}[\hat\Phi,Z]\trans{W}_b\Big).
\end{align*}

\begin{algorithm}[tb!]
\caption{Smoothing spline estimator for single-index model \eqref{model}\label{algo1}}
\medskip
\noindent

\begin{enumerate}[nolistsep,leftmargin=3.7em]
\item[\textbf{Input:}] Data $\{(Y_i,X_i,Z_i)\}_{i=1}^n$; regularization parameter $\lambda_b$; parameter $v$; multiplier bootstrap weights matrix $ W_b=\diag(W_{1,b},\ldots,W_{n,b})$.
\end{enumerate}

\begin{enumerate}

\item Compute the initial estimate $\hat\beta_n^{\rm (ini)}$ of $\beta_0$:

\begin{enumerate}[label={(1\alph*)}]
\item Generate i.i.d.~standard Gaussian random variables $\N_1,\ldots,\N_C$. For each $1\leq j\leq C$, compute $\kappa_j=\N_j\sgn\{\N_j(1)\}/\|\N_j\|_2$.

\item Compute the polynomial spline estimator $(\tilde g_j,\tilde\gamma_j)$ of the partially-linear nonparametric regression model in \eqref{newn} based on the input $\{(X_i\trans\kappa_j,Z_i)\}_{i=1}^n$.


\item Take $\hat\beta_n^{\rm(ini)}=\kappa_{j_*}$, where $\displaystyle j_*=\min_{1\leq j\leq C}\mathop{{\rm median}}_{1\leq i\leq n}|Y_i-\tilde g_j(X_i\trans\kappa_j)-Z_i\trans\tilde\gamma_j|$.

\end{enumerate}

\item Eigen-system construction:

\begin{enumerate}[label={(2\alph*)}]
\item  Estimate $f_{X\trans\beta_0}$ by $\hat f_{X\trans\beta_0}$, a kernel density estimator based on the input $\{X_i\trans\hat\beta_n^{{\rm (ini)}}\}_{i=1}^n$.


\item Compute the normalized eigenvalues and eigenfunctions $\{(\hat\phi_j,\hat\rho_j)\}_{j=1}^v$ of the following ordinary differential equation:
\begin{align*}
-\phi_{j}^{(2m)}(s)=\rho_{j}\sigma_0^2(s)\phi_{j}(s)\hat f_{X\trans\beta_0}(s),\qquad s\in\I,\,j\geq1,
\end{align*}
with boundary conditions given in $\phi_{j}^{(k)}(s)=0$, for $s\in\partial\I$ and $k=m,\ldots,2m-1$.

\end{enumerate}

\item Compute the estimators for $g_0$ and $\gamma_0$:

\begin{enumerate}[label={(3\alph*)}]
\item Compute $\big[
\hat a^{(b)\trans},\hat \gamma^{*\trans}_{n,b}
\big]\trans=\big\{[\hat\Phi,Z]\trans{W}_b[\hat\Phi,Z]+n\lambda_b\diag(\hat\Lambda,0_{q,q})\big\}^{-1}[\hat\Phi,Z]\trans{W}_bY$ in \eqref{hatb}.

\item Compute the bias-adjusted index link estimator $\hat g_{n,b}^*+M_\lambda(\hat g_{n,b}^*)$.

\end{enumerate}

\item Estimate $\beta_0$: compute via line search 
\begin{align*}
\tau^*=\argmin_{\tau\in(\tau^-,\tau^+)}\ell_n\big(\hat g_{n,b}^*+M_\lambda\hat g_{n,b}^*,\beta_\tau,\hat\gamma_{n,b}^*\big),
\end{align*}
for $\beta_\tau,\tau^\pm$ defined in \eqref{betatau} and \eqref{taupm}. Take $\hat\beta_{n,b}^*=\beta_{\tau^*}/\|\beta_{\tau^*}\|_2$.

\item Repeat Steps 2--4 until convergence.

\end{enumerate}
{\bf Output:} Index link estimator $\hat g_{n,b}^*+M_\lambda(\hat g_{n,b}^*)$ for $g_0$; estimators $(\hat\beta_{n,b}^*,\hat\gamma_{n,b}^*)$ for $(\beta_0,\gamma_0)$.

\end{algorithm}

\clearpage

\section{Proofs of theoretical results in Sections~\ref{sec:definenorm}--\ref{sec:jeq}}\label{app:proof}

\subsection{Proof of Proposition~\ref{prop:2.1}}\label{app:proof:prop2.1}

It is clear that $\l\cdot,\cdot\r$ is bilinear and symmetric. Suppose $g\in\H_m,\theta\in\mathbb R^{p-1},\gamma\in\mathbb R^q$. For the $R_X$ and $R_Z$ defined in \eqref{vxz}, observe that
\begin{align}\label{identity}
&\E\{\sigma^2(X,Z)g_0'(X\trans\beta_0) X|X\trans\beta_0\}=g_0'(X\trans\beta_0)\E\{\sigma^2(X,Z) X|X\trans\beta_0\}\notag\\
&\qquad=g_0'(X\trans\beta_0)\sigma_0^2(X\trans\beta_0) R_X(X\trans\beta_0)=\E\big[\sigma^2(X,Z)g_0'(X\trans\beta_0) R_X(X\trans\beta_0)|X\trans\beta_0\big],\notag\\
&\E\{\sigma^2(X,Z)Z|X\trans\beta_0\}=\sigma_0^2(X\trans\beta_0)R_Z(X\trans\beta_0)=\E\big[\sigma^2(X,Z)R_Z(X\trans\beta_0)|X\trans\beta_0\big].
\end{align}
Therefore, by the law of total expectation, we find
\begin{align*}
&\E\Big(\sigma^2(X,Z)\big[\{g'_0(X\trans \beta_0)X\trans Q_{\beta_0}-g_0'(X\trans\beta_0)R_X(X\trans\beta_0)\trans Q_{\beta_0}\}\theta+\{Z\trans-R_Z(X\trans\beta_0)\trans\} \gamma\big]\\
&\qquad\times\big[g(X\trans \beta_0)+ g_0'(X\trans\beta_0)R_X(X\trans\beta_0)\trans\theta+R_Z(X\trans\beta_0)\trans\gamma\big]\Big)\\
&=\E\Big\{\E\Big(\sigma^2(X,Z)\big[\{g'_0(X\trans \beta_0)X\trans Q_{\beta_0}-g'_0(X\trans \beta_0)R_X(X\trans\beta_0)\trans Q_{\beta_0}\}\theta+\{Z\trans-R_Z(X\trans\beta_0)\trans\} \gamma\big]\big|X\trans\beta_0\Big)\\
&\qquad\times\big[g(X\trans \beta_0)+ g_0'(X\trans\beta_0)R_X(X\trans\beta_0)\trans\theta+R_Z(X\trans\beta_0)\trans\gamma\big]\Big\}=0
\end{align*}
Assume $\|(g,\theta,\gamma)\|^2=0$. We then obtain from the above equation that
\begin{align}\label{zero1}
0&=\E\big[\sigma^2(X,Z)\big\{g(X\trans \beta_0)+g'_0(X\trans \beta_0)X\trans Q_{\beta_0}\theta+Z\trans \gamma\big\}^2\big]+\lambda J(g,g)\notag\\
&=\E\Big(\sigma^2(X,Z)\big[g(X\trans \beta_0)+ g_0'(X\trans\beta_0)R_X(X\trans\beta_0)\trans Q_{\beta_0}\theta+R_Z(X\trans\beta_0)\trans\gamma\notag\\
&\qquad+\{g'_0(X\trans \beta_0)X\trans Q_{\beta_0}-g'_0(X\trans \beta_0)R_X(X\trans\beta_0)\trans\}Q_{\beta_0}\theta+\{Z\trans-R_Z(X\trans\beta_0)\trans\} \gamma\big]^2\Big)+\lambda J(g,g)\notag\\
&=\E\Big(\sigma^2(X,Z)\big[g(X\trans \beta_0)+ g_0'(X\trans\beta_0)R_X(X\trans\beta_0)\trans Q_{\beta_0}\theta+R_Z(X\trans\beta_0)\trans\gamma\big]^2\Big)\notag\\
&\qquad+\E\Big(\sigma^2(X,Z)\big[g'_0(X\trans \beta_0)\{X-R_X(X\trans\beta_0)\}\trans Q_{\beta_0}\theta+\{Z-R_Z(X\trans\beta_0)\}\trans \gamma\big]^2\Big)+\lambda J(g,g)\notag\\
&\geq\E\Big(\sigma^2(X,Z)\big[g'_0(X\trans \beta_0)\{X-R_X(X\trans\beta_0)\}\trans Q_{\beta_0}\theta+\{Z-R_Z(X\trans\beta_0)\}\trans \gamma\big]^2\Big)\notag\\
&=\big[\theta\trans,\,\gamma\trans\big]\,\Omega_{\beta_0}\,\big[\theta\trans,\,\gamma\trans\big]\trans,
\end{align}
where $\Omega_{\beta_0}$ is defined in \eqref{Omega}. Note that it holds almost surely that
\begin{align}\label{beta0null}
\bigg[X-\frac{\E\{\sigma^2(X,Z)X|X\trans\beta_0\}}{\E\{\sigma^2(X,Z)|X\trans\beta_0\}}\bigg]\trans\beta_0=0.
\end{align}
Note that $Q_{\beta_0}\trans\Omega_1Q_{\beta_0}-Q_{\beta_0}\trans\Omega_3\Omega_2^{-1}\Omega_3\trans Q_{\beta_0}=Q_{\beta_0}\trans(\Omega_1-\Omega_3\Omega_2^{-1}\Omega_3\trans) Q_{\beta_0}$ is the Schur complement of the block $Q_{\beta_0}\trans\Omega_2Q_{\beta_0}$ in the matrix $\Omega_{\beta_0}$. In addition, it has full rank due to the above equation and Assumption~\ref{a:rank}. By the Guttman rank additivity, $\Omega_{\beta_0}$ is invertible, so is positive definite. Hence, $\lambda_{\min}(\Omega_{\beta_0})\geq c$ for some constant $c>0$. Therefore, from \eqref{zero1} we deduce that $\theta=0_{p-1}$ and $\gamma=0_q$. Furthermore, we have
\begin{align*}
0&=\|(g,\theta,\gamma)\|^2=\E\big\{\sigma^2(X,Z)\,g^2(X\trans\beta_0)\}+\lambda J(g,g)\geq\E\big\{\sigma^2(X,Z)\,g^2(X\trans\beta_0)\}\\
&=\E\big\{\sigma_0^2(X\trans\beta_0)\,g^2(X\trans\beta_0)\}=\int_\I \sigma_0^2(s) g^2(s)f_{X\trans\beta_0}(s)dx,
\end{align*}
which implies that $g\equiv0$ by assumption $\sigma_0^2f_{X\trans\beta_0}>0$ on $\I$. Therefore, $(g,\theta,\gamma)=0$, so $\|\cdot\|$ is a well-defined norm. Moreover, from the calculations in \eqref{zero1}, we obtain that, for $w=(g,Q_{\beta_0}\trans\beta,\gamma)\in\Theta$,
\begin{align*}
\|Q_{\beta_0}\trans\beta\|_2^2+\|\gamma\|_2^2&\leq[\lambda_{\min}(\Omega_{\beta_0})]^{-1}\times\big[(Q_{\beta_0}\trans\beta)\trans,\,\gamma\trans\big]\,\Omega\,\big[(Q_{\beta_0}\trans\beta)\trans,\,\gamma\trans\big]\trans\leq[\lambda_{\min}(\Omega_{\beta_0})]^{-1}\|w\|^2.
\end{align*}
This concludes the proof.

\subsection{Proof of Theorem~\ref{thm:rate}}\label{app:thm:rate}

We first introduce some notations.
For $(g,\beta)\in\H_m\times\B$, let $g\circ\beta$ denote the function on $\X$, defined by $g\circ\beta(x)=g(x\trans\beta)$, for $x\in\X$. Let $\varpi=(g\circ\beta,\gamma)$, and define
\begin{align}\label{norms}
&\|\varpi\|_n^2=\frac{1}{n}\sum_{i=1}^n\big\{g(X_i\trans\beta)+Z_i\trans\gamma\big\}^2,\qquad\quad\|\varpi\|_*^2=\E\big\{g(X\trans\beta)+Z\trans\gamma\big\}^2.
\end{align}
In addition, denote $\hat\varpi_n=(\hat g_n\circ\hat\beta_n,\hat\gamma_n)$ and $\varpi_0=(g_0\circ\beta_0,\gamma_0)$. Observe that
\begin{align}
&\|\hat\varpi_n-\varpi_0\|_n^2=\frac{1}{n}\sum_{i=1}^n\big\{\hat g_n(X_i\trans\hat\beta_n)-g_0(X_i\trans\beta_0)+Z_i\trans(\hat\gamma_n-\gamma_0)\big\}^2.
\end{align}

\subsubsection{A preparatory theorem (Theorem~\ref{prop:rate}) and its proof}

We first prove the following theorem, which is useful for proving Theorem~\ref{thm:rate}. We defer the proof of Theorem~\ref{thm:rate} to Section~\ref{app:sec:thm:rate}.

\begin{theorem}\label{prop:rate}

Under the assumptions of Theorem~\ref{thm:rate}, it holds that
\begin{align*}
\|\hat\varpi_n-\varpi_0\|_n=O_p(r_n);\quad
J(\hat g_n,\hat g_n)+\|\hat\gamma_n\|_2^2=O_p(1);\quad
\|\hat \varpi_n-\varpi_0\|_*=O_p(r_n);\quad
\|\hat g_n\|_\infty=O_p(1),
\end{align*}
where $r_n=\lambda^{1/2}+n^{-1/2}\lambda^{-1/(4m)}$.







\end{theorem}

\begin{proof}[\underline{Proof of Theorem~\ref{prop:rate}}]

Recall the definition of $\ell_n(g,\beta,\gamma)$ in \eqref{elln}.
Since $(\hat g_{n}, \hat\beta_n,\hat\gamma_n )$ minimizes the objective functional $\ell_n(g,\beta,\gamma)+\lambda J(g,g)$, we obtain
\begin{align}\label{rate1}
\ell_n(\hat g_n,\hat\beta_n,\hat\gamma_n)+\lambda J(\hat g_n,\hat g_n) \leq \ell_n(g_0, \beta_0,\gamma_0)+\lambda J(g_0,g_0),
\end{align}
which directly implies that
\begin{align}\label{rate2}
\ell_n(g_0, \beta_0,\gamma_0)-\ell_n(\hat g_{n}, \hat\beta_n,\hat\gamma_n ) \geq \lambda\{J(\hat g_{n},\hat g_n)-J(g_0,g_0)\} \geq-\lambda J(g_0,g_0) = o_p(1).
\end{align}
In addition, \eqref{rate1} implies that
\begin{align*}
\frac{1}{n}\sum_{i=1}^n\big\{Y_i-\hat g_n(X_i\trans \hat\beta_n)-Z_i\trans\hat\gamma_n\big\}^2+\lambda J(\hat g_n,\hat g_n) \leq \frac{1}{n}\sum_{i=1}^n\big\{Y_i-g_0(X_i\trans \beta_0)-Z_i\trans\gamma_0\big\}^2+\lambda J(g_0,g_0) .
\end{align*}
We therefore deduce from the above equation that
\begin{align*}
\|\hat \varpi_n-\varpi_0\|_n^2+\lambda J(\hat g_n,\hat g_n)\leq\frac{2}{n} \sum_{i=1}^n \epsilon_i\big\{\hat g_n(X_i\trans\hat\beta_n)+Z_i\trans\hat\gamma_n-g_0(X_i\trans\beta_0)-Z_i\trans\gamma_0\big\}+\lambda J (g_0,g_0) .
\end{align*}
On the other hand, we deduce by the Cauchy-Schwarz inequality that
\begin{align}\label{lndiff}
\ell_n(g_0, \beta_0,\gamma_0)-\ell_n(\hat g_{n}, \hat\beta_n ,\hat\gamma_n)&=\frac{2}{n} \sum_{i=1}^n \epsilon_i\big\{\hat g_n(X_i\trans\hat\beta_n)+Z_i\trans\hat\gamma_n-g_0(X_i\trans\beta_0)-Z_i\trans\gamma_0\big\}-\|\hat \varpi_n-\varpi_0\|_n^2 \notag\\
&\leq\bigg(\frac{4}{n} \sum_{i=1}^n \epsilon_i^2\bigg)^{1 / 2}\|\hat \varpi_n-\varpi_0\|_n-\|\hat \varpi_n-\varpi_0\|_n^2\notag \\
&\leq O(1)\times\|\hat \varpi_n-\varpi_0\|_n-\|\hat \varpi_n-\varpi_0\|_n^2,
\end{align}
where the last step follows from the fact that $n^{-1}\sum_{i=1}^n \epsilon_i^2=O(1)$ almost surely in view of Assumption~\ref{a:02}.

Combining the above equation with \eqref{rate2} yields
\begin{align*}
\|\hat \varpi_n-\varpi_0\|_n^2 &\leq O_p(1)\times\|\hat \varpi_n-\varpi_0\|_n  -\{\ell_n(g_0, \beta_0,\gamma_0)-\ell_n(\hat g_{n}, \hat\beta_n,\hat\gamma_n )\}
\\
&\leq O_p(1)\times\|\hat \varpi_n-\varpi_0\|_n+o_p(1) .
\end{align*}
The above equation implies that $\|\hat \varpi_n-\varpi_0\|_n^2=O_p(1)$, which  further yields
\begin{align}\label{rate3}
\|\hat \varpi_n\|_n^2=\frac{1}{n}\sum_{i=1}^n\big\{\hat g_n(X_i\trans\hat\beta_n)+Z_i\trans\hat\gamma_n\big\}^2=O_p(1)
\end{align}
due to the fact that $\|\varpi_0\|_n^2=n^{-1}\sum_{i=1}^n\big\{g_0(X_i\trans\beta_0)+Z_i\trans\gamma_0\big\}^2=O(1)$.

Next, we bound $\|\hat g_{n}\|_{\infty}$ using $J(\hat g_{n},\hat g_{n})$, by applying the Sobolev embedding theorem (see Lemma~\ref{lem:embedding} for its statement), which enables us to decompose
\begin{align}\label{rate5}
\hat g_n(s)=\hat g_{n,1}(s)+\hat g_{n,2}(s),\qquad s\in\I,
\end{align}
where $\hat g_{n,1}$ is a polynomial of degree $(m-1)$, defined by
\begin{align}\label{rate6}
\hat g_{n,1}(s)=\sum_{j=1}^m \hat b_{n,j}\, s^{j-1},
\end{align}
and $\|\hat g_{n,2}\|_\infty^2 \leq cJ(\hat g_n,\hat g_n)$, by the Sobolev embedding theorem in Lemma~\ref{lem:embedding}, where $c>0$ only depends on $\I$. Define $\|\hat\gamma_n\|_n^2=n^{-1}\sum_{i=1}^n(Z_i\trans\hat\gamma_n)^2$ and
\begin{align*}
\|\hat \varpi_{n,\ell}\|_n^2&=\frac{1}{n}\sum_{i=1}^n\big\{\hat g_{n,\ell}(X_i\trans\hat\beta_n)+Z_i\trans\hat\gamma_n\big\}^2,\qquad \ell=0,1.
\end{align*}
Observe that $\|\hat\gamma_n\|_n^2\leq c_0^2\|\hat\gamma_n\|_2^2$, where $c_0=\sup_{z\in\Z}\|z\|_2$, so that $\|\hat \varpi_{n,2}\|_n^2\leq\|\hat g_{n,2}\|_\infty^2+c_0^2\|\hat\gamma_n\|^2_2$. Therefore, we deduce from \eqref{rate3} that
\begin{align}\label{rate4}
&\frac{\|\hat \varpi_{n,1}\|_n}{1+\sqrt{J(g_0,g_0)}+\sqrt{J(\hat g_{n},\hat g_n})+c_0(\|\hat\gamma_n\|_2+\|\gamma_0\|_2)}\notag
\\
& \leq \frac{\|\hat \varpi_n \|_n+\|\hat \varpi_{n,2}\|_n}{1+\sqrt{J(g_0,g_0)}+\sqrt{J(\hat g_{n},\hat g_n})+c_0(\|\hat\gamma_n\|_2+\|\gamma_0\|_2)}=O_p(1) .
\end{align}

Note that
\begin{align*}
\|\hat \varpi_{n,1}\|_n^2&=\frac{1}{n}\sum_{i=1}^n\big\{\hat g_{n,1}(X_i\trans\hat\beta_n)+Z_i\trans\hat\gamma_n\big\}^2=\frac{1}{n}\sum_{i=1}^n\Bigg[\bigg\{\sum_{j=1}^m \hat b_{n,j}\, (X_i\trans\hat\beta_n)^{j-1}\bigg\}^2+Z_i\trans\hat\gamma_n\Bigg]^2.
\end{align*}
Define
\begin{align*}
A_n(\beta)&=\frac{1}{n} \sum_{i=1}^n [1, \beta\trans X_i,\ldots,(\beta\trans X_i)^{m-1},Z_i\trans]\trans[1, \beta\trans X_i,\ldots,(\beta\trans X_i)^{m-1},Z_i\trans] ,
\\
A(\beta)&=\E\{A_n(\beta)\}.
\end{align*}
By the Glivenko-Cantelli theorem, we have
\begin{align*}
\sup _{\beta\in\B}|\lambda_{\min}\{A_n(\beta)\}-\lambda_{\min}\{A(\beta)\}|=o_p(1),\qquad n\to\infty .
\end{align*}
Recall that $\lambda_{\min}\{A(\beta)\}$ is the minimum eigenvalue of $A(\beta)$.
Let $\tilde\lambda_{\min}(A)=\min_{\beta\in\B}\lambda_{\min}\{A(\beta)\}$. By Assumption~\ref{a:rank} and the fact that $\|\beta\|_2=1$, together with the identifiability condition for single-index models that the component of $X$ is not linear-dependent \citep[see for example][]{ichimura1993semiparametric,kuchibhotla2020} we obtain that $\inf_{\beta\in\B}|A(\beta)|>0$, so $\inf_{\beta\in\B}\tilde\lambda_{\min}\{A(\beta)\}>0$. Then it follows from the above equation that
\begin{align*}
\|\hat w_{n,1}\|_n^2& =\big[\hat{b}_{n,1},\ldots, \hat{b}_{n,m},\hat\gamma_n\trans\big] A_n(\hat\beta_n)\big[\hat{b}_{n,1},\ldots, \hat{b}_{n,m},\hat\gamma_n\trans\big]\trans \\
& \geq \lambda_{\min}\{A_n(\hat\beta_n)\}\bigg(\sum_{j=1}^m\hat b_{n,j}^2+\|\hat\gamma_n\|_2^2\bigg) \\
& =[\lambda_{\min}\{A_n(\hat\beta_n)\}-\lambda_{\min}\{A(\hat\beta_n)\}]\bigg(\sum_{j=1}^m\hat b_{n,j}^2+\sum_{\ell=1}^q\hat \gamma_{n,\ell}^2\bigg)+\lambda_{\min}\{A(\hat\beta_n)\}\bigg(\sum_{j=1}^m\hat b_{n,j}^2+\sum_{\ell=1}^q\hat \gamma_{n,\ell}^2\bigg) \\
& \geq o_p\bigg(\sum_{j=1}^m\hat b_{n,j}^2+\sum_{\ell=1}^q\hat \gamma_{n,\ell}^2\bigg)+\inf_{\beta\in\B}\tilde\lambda_{\min}\{A(\beta)\}\times \max \{\hat{b}_{n,1}^2,\ldots, \hat{b}_{n,m}^2,\hat\gamma_{n,1}^2,\ldots,\hat\gamma_{n,q}^2\} .
\end{align*}
Combining the above equation and  and \eqref{rate4} yields
\begin{align}\label{rate7}
\frac{\max (|\hat{b}_{n,1}|,\ldots, |\hat{b}_{n,m}|,|\hat\gamma_{n,1}|,\ldots,|\hat\gamma_{n,q}|)}{1+\sqrt{J(g_0,g_0)}+\sqrt{J(\hat g_{n},\hat g_{n})}+c_0(\|\hat\gamma_n\|_2+\|\gamma_0\|_2)}=O_p(1) .
\end{align}
The above equation and \eqref{rate6} imply that
\begin{align*}
\frac{\|\hat g_{n,1}\|_{\infty}+\|\hat\gamma_n\|_2} {1+\sqrt{J(g_0,g_0)}+\sqrt{J(\hat g_{n},\hat g_{n})}+c_0(\|\hat\gamma_n\|_2+\|\gamma_0\|_2)}=O_p(1).
\end{align*}
Combining the above equation with \eqref{rate5} yields
\begin{align}\label{rate44}
&\frac{\|\hat g_{n}\|_{\infty}+\|\hat\gamma_n\|_2}{1+\sqrt{J(g_0,g_0)}+\sqrt{J(\hat g_{n},\hat g_{n})}+c_0(\|\hat\gamma_n\|_2+\|\gamma_0\|_2)} \notag\\
&\qquad\leq \frac{\|\hat g_{n,1}\|_{\infty}+\|\hat g_{n,2}\|_\infty+\|\hat\gamma_n\|_2}{1+\sqrt{J(g_0,g_0)}+\sqrt{J(\hat g_{n},\hat g_{n})}+c_0(\|\hat\gamma_n\|_2+\|\gamma_0\|_2)}=O_p(1) .
\end{align}

Now define the function class
\begin{align}\label{gc}
\mathcal{G}_C&:=\bigg\{f_{g,\beta,\gamma}(x,z):=\frac{g(x\trans\beta)-g_0(x\trans\beta_0)+v\trans(\gamma-\gamma_0)}{1+\sqrt{J(g_0,g_0)}+\sqrt{J(g,g)}+c_0(\|\gamma\|_2+\|\gamma_0\|_2)}: g \in \H_m, \beta \in \B, \gamma\in\mathbb R^q,\notag\\
&\hspace{4cm} \frac{\|g\|_{\infty}+\|\gamma\|_2}{1+\sqrt{J(g_0,g_0)}+\sqrt{J(g,g)}+c_0(\|\gamma\|_2+\|\gamma_0\|_2)} \leq C\bigg\} .
\end{align}
Equation~\eqref{rate7} implies that for any $e>0$, there exists $C_{e}>0$ such that for any $n$,
\begin{align*}
\mathbb{P}\big(f_{\hat g_n,\hat\beta_n,\hat\gamma_n} \in \mathcal{G}_{C_{e}}\big) \geq 1-e.
\end{align*}
The following lemma states the bracketing number of $\mathcal G_C$ w.r.t.~the $\|\cdot\|_\infty$-norm, and its proof is deferred to Section~\ref{app:rate8} at the end of this Section.

\begin{lemma}\label{rate8}
For every fixed constant $C>0$, the bracketing number of $\mathcal{G}_C$ w.r.t.~the infinity norm satisfies
\begin{align*}
\log N(\eta, \mathcal{G}_C,\|\cdot\|_{\infty}) \lesssim \delta^{-1 / m}.
\end{align*}
\end{lemma}

Applying Lemma~\ref{rate8} and Lemma~8.4 of \cite{vandegeer2000} (see Lemma~\ref{lem:geer} in Section~\ref{app:aux} for its statement) yields
\begin{align*}
\frac{\big|n^{-1} \sum_{i=1}^n \epsilon_i\{\hat g_{n}(X_i\trans\hat\beta_n )-g_0(X_i\trans\beta_0)+Z_i\trans(\hat\gamma_n-\gamma_0)\}\big|}{\|\hat \varpi_n-\varpi_0\|_n^{1-1/(2m)}\big\{1+\sqrt{J(g_0,g_0)}+\sqrt{J(\hat g_{n},\hat g_{n})}+c_0(\|\hat\gamma_n\|_2+\|\gamma_0\|_2)\big\}^{1 / (2m)}}=O_p(n^{-1 / 2}) .
\end{align*}
The above equation and \eqref{lndiff} imply
\begin{small}
\begin{align}\label{new3}
& \lambda\{J(\hat g_{n},\hat g_{n})-J(g_0,g_0)\}\leq \ell_n(g_0, \beta_0,\gamma_0)-\ell_n(\hat g_{n}, \hat\beta_n,\hat\gamma_n ) \\
& \leq\|\hat \varpi_n-\varpi_0\|_n^{1-1/(2m)}\big\{1+\sqrt{J(g_0,g_0)}+\sqrt{J(\hat g_{n},\hat g_{n})}+c_0(\|\hat\gamma_n\|_2+\|\gamma_0\|_2)\big\}^{1 /(2m)} O_p(n^{-1 / 2})-\|\hat \varpi_n-\varpi_0\|_n^{2} .\notag
\end{align}
\end{small}
Consider the following two cases.

\textit{Case (i)}: $\sqrt{J(\hat g_{n},\hat g_{n})}+c_0\|\hat\gamma_n\|_2>1+\sqrt{J(g_0,g_0)}+c_0\|\gamma_0\|_2$. 

By the proof of Theorem 10.2 of \cite{vandegeer2000} where we take $v=2$ and $\alpha=1 /m$, we deduce that
\begin{align*}
\sqrt{J(\hat g_{n},\hat g_{n})}+c_0\|\hat\gamma_n\|_2=O_p(n^{-1 / 2}\lambda^{-1/2-1/(4m)}),  \quad \quad\|\hat\varpi_n-\varpi_0\|_n=O_p(n^{-1 / 2}\lambda^{-1 / (4m)})  .
\end{align*}
Therefore, in view of Assumption~\ref{a:5.0}, in this case it is true that
\begin{align}\label{new4}
\sqrt{J(\hat g_{n},\hat g_{n})}+c_0\|\hat\gamma_n\|_2=O_p(1),  \quad \quad\|\hat\varpi_n-\varpi_0\|_n=O_p(n^{-1 / 2}\lambda^{-1 / (4m)})  .
\end{align}

\textit{Case (ii)}: $\sqrt{J(\hat g_{n},\hat g_{n})}+c_0\|\hat\gamma_n\|_2\leq 1+\sqrt{J(g_0,g_0)}+c_0\|\gamma_0\|_2$.

In this case, \eqref{new3} implies,
\begin{align*}
\|\hat\varpi_n-\varpi_0\|_n^2 \leq\|\hat\varpi_n-\varpi_0\|_n^{1-1/(2m)}\{1+\sqrt{J(g_0,g_0)}\}^{1 /(2m)} O_p(n^{-1 / 2})+\lambda J(g_0,g_0) .
\end{align*}
Therefore, by Assumption~\ref{a:5.0}, it follows that
\begin{align*}
\text{either}\quad&\|\hat\varpi_n-\varpi_0\|_n \leq\{1+\sqrt{J(g_0,g_0)}\}^{2/(2m+1)} O_p(n^{-m / (2m+1)})=O_p(\lambda^{1/2})
\\
\text{or}\quad&\|\hat\varpi_n-\varpi_0\|_n \leq O_p(1) \{\lambda J(g_0,g_0)\}^{1/2}=O_p(\lambda^{1/2}) .
\end{align*}

In conclusion, combining the above equation with \eqref{new4} yields that 
\begin{align*}
&\|\hat w_n-w_0\|_n=O_p(\lambda^{1/2}+n^{-1 / 2}\lambda^{-1 / (4m)}),\qquad J(\hat g_{n},\hat g_{n})+\|\hat\gamma_n\|_2^2=O_p(1).
\end{align*}
In addition, combining the above equation and \eqref{rate44} implies $\|\hat g_{n}\|_{\infty}=O_p(1)$.
Moreover, the proof of the convergence rate $\|\hat\varpi_n-\varpi_0\|_*$ follows from the bound of $\|\hat\varpi_n-\varpi_0\|_n$ and Lemma 5.16 of \cite{vandegeer2000}; see Lemma~\ref{lem:geer2} for its statement. This concludes the proof of Theorem~\ref{prop:rate}.

\end{proof}

The proof of Theorem~\ref{thm:rate} leverages the following lemma, which is proved by applying Lemma~\ref{lem:g'infinity} in Section~\ref{app:aux} and a slight modification of the proof of Theorem~3 in \cite{kuchibhotla2020}; we omit its proof for brevity.
\begin{lemma}\label{jgn2}
Under the assumptions of Theorem~\ref{thm:rate}, it holds that
\begin{align*}
    \|\hat g_n'-g_0'\|_\infty=o_p(1),\qquad \|\hat\beta_n-\beta_0\|_2=o_p(1).
\end{align*}
\end{lemma}

\subsubsection{Proof of Theorem~\ref{thm:rate}}\label{app:sec:thm:rate}
Now, we prove Theorem~\ref{thm:rate} using Theorem~\ref{prop:rate}. Recall the definition of $r_x(g,\beta)$ in \eqref{rwx}, and recall that $w=(g,Q_{\beta_0}\trans\beta,\gamma)$, $\hat w_n=(\hat g_n,Q_{\beta_0}\trans\hat\beta_n,\hat\gamma_n)$, and $w_0=(g,0_{p-1},\gamma)$. 
Define
\begin{align*}
    &\|w\|_\ddagger^2=\E\big\{g(X\trans\beta_0)+g_0'(X\trans\beta_0)X\trans Q_{\beta_0}Q_{\beta_0}\trans\beta+Z\trans\gamma\big\}^2,\notag\\
&\|w\|_\Diamond^2=\E\Big[\sigma^2(X,Z)\big\{g(X\trans\beta_0)+g_0'(X\trans\beta_0)X\trans Q_{\beta_0}Q_{\beta_0}\trans\beta+Z\trans\gamma\big\}^2\Big].
\end{align*}
Observe that $\|w\|^2=\|w\|_\Diamond^2+\lambda J(g,g)$.
From the proof of Lemma~\ref{lem:unionbound2} in Section~\ref{app:lem:suprdiff}, we obtain
\begin{align*}
\sup_{x\in\X}|r_x(\hat g_n,\hat\beta_n)|\leq c\|\hat\beta_n-\beta_0\|_2^{3/2}+c\|\hat g'_n-g_0'\|_\infty\times\|\hat\beta_n-\beta_0\|_2+c\|\hat\beta_n-\beta_0\|_2^2=o_p(1)\times\|\hat\beta_n-\beta_0\|_2,
\end{align*}
due to the fact that $\|\hat g_n'-g_0'\|_\infty=o_p(1)$ and $\|\hat\beta_n-\beta_0\|_2=o_p(1)$ by Lemma~\ref{jgn2}. Observe that
\begin{align*}
\|w\|_\ddagger^2&=\E\big\{g(X\trans\beta)+Z\trans\gamma-r_X(g,\beta)\big\}^2\leq2\|\varpi\|_*^2+2\E|r_X(g,\beta)|^2\leq 2\|\varpi\|_*^2+2\sup_{x\in\X}|r_x(g,\beta)|^2.
\end{align*}
Note that $r_x(g_0,\beta_0)\equiv0$ for any $x\in\X$. The above two equations therefore yield
\begin{align}\label{rate9}
\|\hat w_n-w_0\|_\ddagger^2&\leq2\|\hat \varpi_n-\varpi_0\|_*^2+2\sup_{x\in\X}|r_x(\hat g_n,\hat\beta_n)|^2\leq 2\|\hat \varpi_n-\varpi_0\|_*^2+o_p(1)\times\|\hat\beta_n-\beta_0\|_2^2.
\end{align}
In addition, it follows from Lemma~\ref{lem:qbeta} and Proposition~\ref{prop:2.1} that
\begin{align*}
\|\hat\beta_n-\beta_0\|_2^2&=\|Q_{\beta_0}Q_{\beta_0}\trans(\hat\beta_n-\beta_0)\|_2^2+\|\beta_0\beta_0\trans(\hat\beta_n-\beta_0)\|_2^2\\
&=\|Q_{\beta_0}\trans\hat\beta_n\|_2^2+|\beta_0\trans\hat\beta_n-1|^2\leq 2\|Q_{\beta_0}\trans\hat\beta_n\|_2^2\leq2\{\lambda_{\min}(\Omega_{\beta_0})\}^{-1}\|\hat w_n-w_0\|^2.
\end{align*}
Note that the function $\sigma^2(x,z)$ is upper bounded on $\X\times\Z$ due to Assumption~\ref{a:rank1}. Combining the above equation with \eqref{rate9} yields
\begin{align*}
&\|\hat w_n-w_0\|^2_\Diamond\leq\|\sigma^2\|_\infty\times\|\hat w_n-w_0\|_\ddagger^2\leq2\|\sigma^2\|_\infty\times\|\hat \varpi_n-\varpi_0\|_*^2+o_p(1)\times\|\hat w_n-w_0\|^2
\end{align*}
Combining the above equation with the fact that $J(\hat g_n,\hat g_n)=O_p(1)$ in view of Theorem~\ref{prop:rate} yields
\begin{align*}
\|\hat w_n-w_0\|^2=\|\hat w_n-w_0\|^2_\Diamond+\lambda J(\hat g_n,\hat g_n)\leq c\| \hat\varpi_n-\varpi_0\|_*^2+O_p(\lambda)+o_p(1)\times\|\hat w_n-w_0\|^2.
\end{align*}
Since $\| \hat\varpi_n-\varpi_0\|_*^2+O_p(\lambda)=O_p(r_n^2)$ due to Theorem~\ref{prop:rate}, we deduce from the above equation that $\|\hat w_n-w_0\|^2=O_p(r_n^2)$.
This completes the proof of Theorem~\ref{thm:rate}.

\subsubsection{Proof of Lemma~\ref{rate8}}\label{app:rate8}

Recall that we aim to bound the bracketing entropy of the space $\mathcal G_C$ defined in \eqref{gc}.
For constants $c_1,c_2,c_3>0$, define
\begin{align*}
&\mathcal{G}_{c_1, c_2}=\left\{g \in \mathcal{H}_m:\|g\|_{\infty} \leq c_1, J(g,g) \leq c_2\right\},\qquad
\mathcal C_{c_3}=\{\gamma\in\mathbb R^q,\|\gamma\|_2\leq c_3\}.
\end{align*}
Note that w.r.t.~$\|\cdot\|_\infty$-norm, the covering number and the bracketing number are equivalent up to a constant.
By Theorem~2.4 in \cite{vandegeer2000} and Lemma~4.1 of \cite{pollard1990empirical}, we obtain that the covering numbers of $\mathcal G_{c_1,c_2}$, $\B$ and $\C_{c_3}$ satisfy
\begin{align*}
&N\big(\eta, \mathcal{G}_{c_1, c_2},\|\cdot\|_{\infty}\big) \leq c\exp\{(c_1+c_2)\eta^{-1/m}\},\quad N\big(\eta, \B,\|\cdot\|_2\big) \leq c\delta^{-p+1},\quad N\big(\eta, \mathcal C_{c_3},\|\cdot\|_2\big) \leq c \delta^{-q}.
\end{align*}
Define
\begin{align*}
\mathcal{M}_{c_1, c_2,c_3}:=\left\{f_{g,\beta,\gamma}:f_{g,\beta,\gamma}(x,z)=g(x\trans\beta)+z\trans\gamma, \beta\in \B, g \in\mathcal G_{c_1,c_2},\|\gamma\|_2\leq c_3\right\} .
\end{align*}
We claim that
\begin{align}\label{new2}
\log N_{}(\eta, \mathcal{M}_{c_1, c_2,c_3},\|\cdot\|_{\infty}) \leq c \eta^{-1/m}.
\end{align}
To show \eqref{new2}, let $\{g_1, g_2, \ldots, g_{r_1}\}$ be a $\eta$-net for $\mathcal{G}_{c_1, c_2}$ w.r.t.~the $\|\cdot\|_\infty$-norm, let $\{\beta_1, \beta_2, \ldots, \beta_{r_2}\}$ be a $\eta$-net of $\B$ w.r.t.~$\|\cdot\|_2$-norm, and let $\{\gamma_1,\gamma_2,\ldots,\gamma_{r_3}\}$ be a $\eta$-net of $\mathcal C_c$ w.r.t.~the $\|\cdot\|_2$-norm, where $r_1\leq c\exp(c\eta^{-1/m})$, $r_2\leq c\eta^{-p+1}$ and $r_3\leq c\eta^{-q}$. Then, it holds that
\begin{align*}
\tilde{\mathcal{M}}:=\Big\{\tilde f_{g_i,\beta_j,\gamma_k}:\tilde f_{g_i,\beta_j,\gamma_k}(x,z)=g_i(x\trans\beta_j)+z\trans\gamma_k,1\leq i\leq r_1,1\leq j\leq r_2,1\leq k\leq r_3\Big\}\subset\mathcal M_{c_1,c_2,c_3}
\end{align*}
forms a $\eta$-bracket of $\mathcal M_{c_1,c_2,c_3}$ w.r.t.~the $\|\cdot\|_\infty$-norm. To see this, by applying Lemma~\ref{lem:g'infinity}, we obtain that, for any $f_{g,\beta,\gamma}\in\mathcal M_{c_1,c_2,c_3}$, there exists $\tilde f_{g_i,\beta_j,\gamma_k}\in \tilde{\mathcal M}$, for some $i$, $j$ and $k$, defined by $\tilde f_{g_i,\beta_j,\gamma_k}(x,z)=g_i( x\trans\beta_j)+z\trans\gamma_k$, such that
\begin{align*}
&\sup_{x\in\X,z\in\Z}|f_{g,\beta,\gamma}(x,z)-\tilde f_{g_i,\beta_j,\gamma_k}(x,z)|=\sup_{x\in\X,z\in\Z}| g(x\trans\beta)  -g_i( x\trans\beta_j) +z\trans(\gamma-\gamma_k)|
\\
 &\qquad\leq\sup_{x\in\X}|g(x\trans\beta)-g(x\trans\beta_j)|+\sup_{x\in\X}|g( x\trans\beta_j)-g_i(x\trans\beta_j)| +\sup_{z\in\Z}|z\trans(\gamma-\gamma_k)|\\
& \qquad\leq c\|g'\|_{\infty}\times\|\beta-\beta_j\|_2+\| g-g_i \|_{\infty} +c\|\gamma-\gamma_k\|_2\leq c\eta ,
\end{align*}
which implies that $\log N_{}(\eta, \mathcal{M}_{c_1, c_2,c_3},\|\cdot\|_{\infty})\leq \log(r_1r_2) \leq c \eta^{-1/m}$  and proves \eqref{new2}.

Now we prove Lemma~\ref{rate8}. Define
\begin{align*}
&\mathcal{F}_{c}:=\bigg\{h_{g,\beta,\gamma}: h_{g,\beta,\gamma}(x,z)=\frac{g(x\trans\beta)+z\trans\gamma}{1+J(g_0,g_0)+J(g,g)+c_0(\|\gamma\|_2^2+\|\gamma_0\|_2^2)},\\
&\hspace{3cm} \beta \in \B, g \in \H_m, \gamma\in\mathbb R^q, \frac{\|g\|_{\infty}+\|\gamma\|_2}{1+J(g_0,g_0)+J(g,g)+c_0(\|\gamma\|_2^2+\|\gamma_0\|_2^2)} \leq c\bigg\},
\\
&\mathcal{J}:=\left\{h_{g,\gamma}: h_{g,\gamma}(x,z)=\frac{g_0(x\trans\beta_0)+z\trans\gamma_0}{1+J(g_0,g_0)+J(g,g)+c_0(\|\gamma\|_2^2+\|\gamma_0\|_2^2)},g \in \mathcal{H}_m,\gamma\in\mathbb R^q\right\}.
\end{align*}
Since $\mathcal{F}_{c} \subset \mathcal{M}_{c, 1,1}$, we can choose $\eta$-brackets $\left[g_{1,1}, g_{1,2}\right], \ldots,\left[g_{r_4, 1}, g_{r_4, 2}\right]$ of $\mathcal{F}_{c}$ such that for every $h_{g,\beta,\gamma} \in \mathcal{F}_{c}$ there exists $i$ such that $g_{i, 1} \leq h_{g,\beta,\gamma} \leq g_{i, 2}$. Here, $r_4$ is such that $\log(r_4)\lesssim \eta^{-1/m}$.
Note that since $\Z$ is compact, it holds that $N(\eta,\mathcal C_c,\|\cdot\|_\infty)\leq c\eta^{-q}$. Observe that $\mathcal{J} \subset \mathcal{G}_{c_0, 1}+\mathcal C_1$, where $c_0=\|g_0\|_{\infty} / J(g_0,g_0)$. Thus we can choose $\eta$-brackets $\left[b_{1,1}, b_{1,2}\right], \ldots,\left[b_{r_5, 1}, b_{r_5, 2}\right]$ of $\mathcal{J}$ such that for every $h \in \mathcal{J}$ there exists a $j$ such that $b_{j, 1}(x,z) \leq h(x,z) \leq b_{j, 2}( x,z)$. Here, $r_5$ is such that $\log(r_5)\lesssim \eta^{-1/m}$.
Therefore, we deduce that
\begin{align*}
&g_{i, 1}(x,z)-b_{j, 2}(x,z) 
\\
&\leq \frac{g(x\trans\beta)+z\trans\gamma}{1+J(g_0,g_0)+J(g,g)+c_0(\|\gamma\|_2^2+\|\gamma_0\|_2^2)}-\frac{g_0(x\trans\beta_0)+z\trans\gamma_0}{1+J(g_0,g_0)+J(g,g)+c_0(\|\gamma\|_2^2+\|\gamma_0\|_2^2)}
\\
& \leq g_{i, 2}(x,z)-b_{j, 1}(x,z),
\end{align*}
where $i$ depends on $(g, \beta,\gamma)$ and $j$ on $(g,\gamma)$.
Therefore, the following brackets
\begin{align*}
\big\{[g_{i, 1}-b_{j, 2}, g_{i, 2}-b_{j, 1}]\big\}_{1\leq i\leq r_4,1\leq j\leq r_5}
\end{align*}
cover the space $\mathcal G_C$ defined in \eqref{gc}. Hence, the bracketing entropy satisfies
\begin{align*}
\log N(\eta, \mathcal{G}_C,\|\cdot\|_{\infty}) \lesssim \log(r_4r_5)\lesssim\eta^{-1/m}.
\end{align*}
This completes the proof of Lemma~\ref{rate8}.

\subsection{Proof of Lemma~\ref{prop:jointee}}\label{app:prop:jointee}


For arbitrary $(g,\beta,\gamma)\in\Theta$, consider the path $\{(G_t,B_t,\Gamma_t)\}_{t\in[-1,1]}$ on $\Theta$ defined by, for $t\in[-1,1]$ and $u\in\I$,
\begin{align*}
&G_t(u)=tg(u)+(1-t)\hat g_n(u),\qquad B_t=\frac{t\beta+(1-t)\hat\beta_n}{\|t\beta+(1-t)\hat\beta_n\|_2},\qquad \Gamma_t=t\gamma+(1-t)\hat\gamma_n,
\end{align*}
with $(G_0,B_0,\Gamma_0)=(\hat g_n,\hat\beta_n,\hat\gamma_n)$ and $(G_1,B_1,\Gamma_1)=(g,\beta,\gamma)$.
Direct calculations yield
\begin{align*}
&\partial_tG_t(u)=g(u)-\hat g_n(u);\qquad \partial_uG_t(u)=tg'(u)+(1-t)\hat g_n'(u);\qquad \Gamma_t'\equiv\gamma-\hat\gamma_n;\\
&B'_t=\frac{\beta-\hat\beta_n}{\|t\beta+(1-t)\hat\beta_n\|_2}-\frac{(\beta-\hat\beta_n)\trans\{t\beta+(1-t)\hat\beta_n\}}{\|t\beta+(1-t)\hat\beta_n\|_2^3}\{t\beta+(1-t)\hat\beta_n\}.
\end{align*}
Since $\hat\beta_n\trans\hat\beta_n=\|\hat\beta_n\|_2^2=1$ and $I_p=Q_{\hat\beta_n}Q_{\hat\beta_n}\trans+\hat\beta_n\hat\beta_n\trans$, we obtain
\begin{align*}
&B'_0=(\beta-\hat\beta_n)-\hat\beta_n\hat\beta_n\trans(\beta-\hat\beta_n)=(I_p-\hat\beta_n\hat\beta_n\trans)(\beta-\hat\beta_n)=Q_{\hat\beta_n}Q_{\hat\beta_n}\trans\beta\,.
\end{align*}
For the objective functional in \eqref{hatwn0} defined by
\begin{align*}
L_n(g,\beta,\gamma)=\frac{1}{2n}\sum_{i=1}^n\big\{Y_i-g(X_i\trans \beta)-Z_i\trans\gamma\big\}^2+\frac{\lambda}{2}\,J(g,g),
\end{align*}
consider its constraint along the path $\{(G_z,B_z)\}_{z\in[-1,1]}$:
\begin{align*}
L_n(G_t,B_t,\Gamma_t)=\frac{1}{2n}\sum_{i=1}^n\big[Y_i-G_t(X_i\trans B_t)-Z_i\trans\Gamma_t\big]^2+\frac{\lambda}{2} J(G_t,G_t).
\end{align*}
Observe that
\begin{align*}
\frac{d}{dt}L_n(G_t,B_t,\Gamma_t)&=-\frac{1}{n}\sum_{i=1}^n\big\{Y_i-G_t(X_i\trans B_t)-Z_i\trans\Gamma_t\big\}\\
&\qquad\qquad\times\big\{\partial_tG_t(X_i\trans B_t)+\partial_u G_t(X_i\trans B_t)X_i\trans B'_t+Z_i\trans\Gamma_t'\big\}+\lambda J(G'_t,G_t).
\end{align*}
Since $t=0$ minimizes $L_n(G_t,B_t,\Gamma_t)$, we deduce that
\begin{align}\label{12}
0&=\frac{d}{dt}L_n(G_t,B_t,\Gamma_t)\bigg|_{t=0}\notag
\\
&=-\frac{1}{n}\sum_{i=1}^n\big\{Y_i-G_0(X_i\trans B_0)-Z_i\trans\Gamma_0\big\}\big\{\partial_tG_0(X_i\trans B_0)+\partial_uG_0(X_i\trans B_0)X_i\trans B'_0+Z_i\trans\Gamma_0'\big\}+\lambda J(G'_0,G_0)\notag\\
&=-\frac{1}{n}\sum_{i=1}^n\big\{Y_i-\hat g_n(X_i\trans \hat\beta_n)-Z_i\trans\hat\gamma_n\big\}\notag\\
&\quad\times\Big\{g(X_i\trans \hat\beta_n)-\hat g_n(X_i\trans \hat\beta_n)+\hat g_n'(X_i\trans \hat\beta_n)X_i\trans Q_{\hat\beta_n}Q_{\hat\beta_n}\trans\beta+Z_i\trans\gamma-Z_i\trans\hat\gamma_n\Big\}+\lambda J(\hat g_n,g-\hat g_n),
\end{align}
where $(g,\beta,\gamma)\in\Theta$ is arbitrary. Define
\begin{align*}
&\tilde F_n(g,\beta,\gamma):=-\frac{1}{n}\sum_{i=1}^n\big\{Y_i-\hat g_n(X_i\trans \hat\beta_n)-Z_i\trans\hat\gamma_n\big\}\Big\{g(X_i\trans \hat\beta_n)+\hat g_n'(X_i\trans \hat\beta_n)X_i\trans Q_{\hat\beta_n}Q_{\hat\beta_n}\trans\beta+Z_i\trans\gamma\Big\}+\lambda J(\hat g_n,g),
\end{align*}
so that \eqref{12} implies that $F(g,\beta,\gamma)\equiv F(\hat g_n,\hat\beta_n,\hat\gamma_n)$ for any $(g,\beta,\gamma)\in\Theta$. Due to the arbitrarity of $(g,\beta,\gamma)\in\Theta$, since $\|\beta+\hat\beta_n\|_2\neq0$, we obtain that
\begin{align*}
\tilde F_n(g,\beta,\gamma)=\tilde F_n\bigg(\frac{g}{\|\beta+\hat\beta_n\|_2},\frac{\beta+\hat\beta_n}{\|\beta+\hat\beta_n\|_2},\frac{\gamma}{\|\beta+\hat\beta_n\|_2}\bigg).
\end{align*}
The above equation together with the fact that $\|\beta+\hat\beta_n\|_2\neq1$ yields $$\tilde F_n(g,\beta,\gamma)\equiv0,\qquad(g,\beta,\gamma)\in\H_m\times\B\times\mathbb R^q.$$ Recall from \eqref{scoren} the definition of $F_n$. Due to the arbitrarity of $\beta\in\B$, we deduce that from the above equation that $F_n(g,\theta,\gamma)\equiv0$, where $(g,\theta,\gamma)\in\H_m\times\mathbb R^{p-1}\times\mathbb R^q$. The proof is therefore complete.

\subsection{Proof of Proposition~\ref{cor:sxz}}\label{app:proof:ru}

This section is organized as follows.  The proof regarding $S_{x,z}$ in Proposition~\ref{cor:sxz} is given in Section~\ref{app:prep1}. The proof regarding $S_\lambda^\circ$ in Proposition~\ref{cor:sxz} is given in Section~\ref{app:proof:plambda}. Finally, Section~\ref{app:prep3} states a corollary (Corollary~\ref{n}) and its proof, which is useful for proving the joint Bahadur representation in Theorem~\ref{thm:bahadur}.

\subsubsection{Proof regarding $S_{x,z}$ in Proposition~\ref{cor:sxz}
}\label{app:prep1}

We start by stating and proving a preparatory proposition, which is useful for computing $S_{x,z}$ in Proposition~\ref{cor:sxz}. Let $f_{(g,\theta,\gamma)}$ be the linear functional of $(g,\theta,\gamma)\in\Theta$, defined by
\begin{align}\label{f}
f_{(g,\theta,\gamma)}(s,u,v)=g(s)+u\trans \theta+v\trans \gamma,\qquad~~~ (s,u,v)\in\I\times\mathbb R^{p-1}\times\mathbb R^q.
\end{align}
The following proposition identifies the Riesz representation of $f_{(g,\theta,\gamma)}$ in $\Theta$, w.r.t.~the inner product $\l\cdot,\cdot\r$ defined in Proposition~\ref{prop:2.1}.

\begin{proposition}\label{prop:ru}
It is true that $f_{(g,\theta,\gamma)}(s,u,v)=\big\l S^\dagger_{s,u,v},(g,\theta,\gamma)\big\r$, for any $(s,u,v)\in\I\times\mathbb R^{p-1}\times\mathbb R^q$ and $(g,\theta,\gamma)\in\Theta$, where $S^\dagger_{s,u,v}=(H^\dagger_{s,u,v},N^\dagger_{s,u,v},T^\dagger_{s,u,v})\in \H_m\times\mathbb R^{p-1}\times\mathbb R^q$ is defined by
\begin{align*}
\Biggl[\begin{matrix}
N^\dagger_{s,u,v}\\
T^\dagger_{s,u,v}
\end{matrix}\Biggl]&=(\Omega_{\beta_0}+\Sigma_{\beta_0,\lambda})^{-1}\Biggl[\begin{matrix}
v-Q_{\beta_0}\trans A_X(s)\\
z-A_Z(s)
\end{matrix}\Biggl],\\
 H^\dagger_{s,u,v}&=K_{s}-A_X\trans Q_{\beta_0} N^\dagger_{s,u,v}-A_Z\trans T^\dagger_{s,u,v}.
\end{align*}
Here, $\Omega_{\beta_0}$ is defined in \eqref{Omega} and $\Sigma_{\beta_0,\lambda}$ is a $(p+q-1)$-square matrix defined by
\begin{align}\label{Sigmalambda}
\Sigma_{\beta_0,\lambda}&={V}\left\{\Biggl[
\begin{matrix}
Q_{\beta_0}\trans g_0'\cdot R_X\\
R_Z
\end{matrix}
\Biggl],\Biggl[
\begin{matrix}
Q_{\beta_0}\trans M_\lambda(g_0'\cdot R_X)\\
M_\lambda(R_Z)
\end{matrix}
\Biggl]\trans\right\}
\end{align}
and is such that $\lim_{\lambda\to0}\Sigma_{\beta_0,\lambda}=0_{p+q-1,p+q-1}$.

\end{proposition}

\begin{proof}[\underline{Proof of Proposition~\ref{prop:ru}}]

Denote the tuple $\xi=(s,u,v)\in\I\times\mathbb R^{p-1}\times\mathbb R^q$ for brevity. Let
\begin{align}\label{gu}
{G}=\Biggl[
\begin{matrix}
g_0'(X\trans\beta_0)\,Q_{\beta_0}\trans X\\
Z
\end{matrix}
\Biggl]\,.
\end{align}
First, observing \eqref{identity} and the definitions of the functions $A_X$ and $A_Z$ in \eqref{ab}, we have, for $g\in\H_m$,
\begin{align}\label{eab}
\E\big\{\sigma^2(X,Z)\,g(X\trans\beta_0){G}\big\}&=\E\Biggl[
\begin{matrix}
\sigma_0^2(X\trans\beta_0)\,g(X\trans\beta_0)\,g_0'(X\trans\beta_0)\,Q_{\beta_0}\trans R_X(X\trans\beta_0)\\
\sigma_0^2(X\trans\beta_0)\,g(X\trans\beta_0)\,R_Z(X\trans\beta_0)
\end{matrix}
\Biggl]=\Biggl[
\begin{matrix}
Q_{\beta_0}\trans\l A_X,g\r_K\\
\l A_Z,g\r_K
\end{matrix}
\Biggl]\,.
\end{align}

Recalling the definition of the inner product $\l\cdot,\cdot\r$ in Proposition~\ref{prop:2.1}, we have
\begin{align*}
&\l S_{\xi}^\dagger,(g,\theta,\gamma)\r=\l (H_{\xi}^\dagger,N_{\xi}^\dagger,T_{\xi}^\dagger),(g,\theta,\gamma)\r\notag\\
&=\E\Big[\sigma^2(X,Z)\big\{H_{\xi}^\dagger(X\trans \beta_0)+g'_0(X\trans \beta_0)X\trans Q_{\beta_0} N_{\xi}^\dagger+Z\trans  T_{\xi}^\dagger\big\}\notag\\
&\hspace{1cm}\times\big\{g(X\trans \beta_0)+g'_0(X\trans \beta_0)X\trans Q_{\beta_0}\theta+Z\trans \gamma\big\}\Big]+\lambda J(H_{\xi}^\dagger,g)\notag\\
&=\E\Big[\sigma^2(X,Z)g(X\trans \beta_0)\big\{g'_0(X\trans \beta_0)X\trans Q_{\beta_0} N_{\xi}^\dagger+Z\trans  T_{\xi}^\dagger\big\}\Big]\notag\\
&\hspace{1cm}+\E\Big[\sigma^2(X,Z)\big\{H_{\xi}^\dagger(X\trans \beta_0)+g'_0(X\trans \beta_0)X\trans Q_{\beta_0} N_{\xi}^\dagger+Z\trans  T_{\xi}^\dagger\big\}\big\{g'_0(X\trans \beta_0)X\trans Q_{\beta_0}\theta+Z\trans \gamma\big\}\Big]\notag\\
&\hspace{1cm}+\E\big\{\sigma^2(X,Z)H_{\xi}^\dagger(X\trans \beta_0)\,g(X\trans \beta_0)\big\}+\lambda J(H_{\xi}^\dagger,g)\notag\\
&=\Biggl[\begin{matrix}
\theta\\
\gamma
\end{matrix}
\Biggl]\trans\E\left\{\sigma^2(X,Z)H_{\xi}^\dagger(X\trans \beta_0)\Biggl[
\begin{matrix}
g_0'(X\trans\beta_0)\,Q_{\beta_0}\trans X\notag\\
Z
\end{matrix}
\Biggl]+\sigma^2(X,Z)\Biggl[
\begin{matrix}
g_0'(X\trans\beta_0)\,Q_{\beta_0}\trans X\\
Z
\end{matrix}
\Biggl]\Biggl[
\begin{matrix}
g_0'(X\trans\beta_0)\,Q_{\beta_0}\trans X\\
Z
\end{matrix}
\Biggl]\trans\right\}\Biggl[\begin{matrix}
N_{\xi}^\dagger\\
T_{\xi}^\dagger
\end{matrix}
\Biggl]\\
&\hspace{1cm}+\E\left\{\sigma^2(X,Z)g(X\trans\beta_0)\Biggl[
\begin{matrix}
g_0'(X\trans\beta_0)\,Q_{\beta_0}\trans X\notag\\
Z
\end{matrix}
\Biggl]\trans\right\}\Biggl[\begin{matrix}
N_{\xi}^\dagger\\
T_{\xi}^\dagger
\end{matrix}
\Biggl]+\l H_{\xi}^\dagger,g\r_K\,.
\end{align*}
In view of ${G}$ defined in \eqref{gu}, we obtain from the above equation and \eqref{eab} that
\begin{align}\label{ru1}
&\l S_{\xi}^\dagger,(g,\theta,\gamma)\r=\l (H_{\xi}^\dagger,N_{\xi}^\dagger,T_{\xi}^\dagger),(g,\theta,\gamma)\r\notag\\
&=\Biggl[\begin{matrix}
\theta\notag\\
\gamma
\end{matrix}
\Biggl]\trans\E\big\{\sigma^2(X,Z)H_{\xi}^\dagger(X\trans\beta_0)\,{G}+\sigma^2(X,Z){G}{G}\trans\big\}\Biggl[\begin{matrix}
N_{\xi}^\dagger\\
T_{\xi}^\dagger
\end{matrix}
\Biggl]+\E\big\{\sigma^2(X,Z)g(X\trans\beta_0)\,{G}\trans\big\}\Biggl[\begin{matrix}
N_{\xi}^\dagger\notag\\
T_{\xi}^\dagger
\end{matrix}
\Biggl]+\l H_{\xi}^\dagger,g\r_K\notag\\
%
%
&=\Biggl[\begin{matrix}
\theta\\
\gamma
\end{matrix}
\Biggl]\trans\E\big\{\sigma^2(X,Z)H_{\xi}^\dagger(X\trans\beta_0)\,{G}+\sigma^2(X,Z){G}{G}\trans\big\}\Biggl[\begin{matrix}
N_{\xi}^\dagger\\
T_{\xi}^\dagger
\end{matrix}
\Biggl]+\l H_{\xi}^\dagger+A_X\trans Q_{\beta_0} N_{\xi}^\dagger+A_Z\trans T_{\xi}^\dagger,g\r_K\,.
\end{align}
Observe the fact that $\l K_s,g\r_K=g(s)$. By setting 
\begin{align*}
\l S_{\xi}^\dagger,(g,\theta,\gamma)\r=g(s)+u\trans\theta+v\trans\gamma=\l K_s,g\r_K+u\trans\theta+v\trans\gamma,
\end{align*}
we deduce from \eqref{ru1} that
\begin{align}\label{equ}
\left\{
\begin{array}{l}
\Biggl[
\begin{matrix}
u\\
v
\end{matrix}
\Biggl]=\E\big\{\sigma^2(X,Z)H_{\xi}^\dagger(X\trans\beta_0)\,{G}\big\}+\E\{\sigma^2(X,Z){G} {G}\trans\}\Biggl[
\begin{matrix}
N_{\xi}^\dagger\\
T_{\xi}^\dagger
\end{matrix}
\Biggl]\\
K_{s}=H_{\xi}^\dagger+A_X\trans Q_{\beta_0} N_{\xi}^\dagger+A_Z\trans T_{\xi}^\dagger
\end{array}
\right.
\end{align}
Substituting $H_{\xi}^\dagger$ by $K_{s}-A_X\trans Q_{\beta_0} N_{\xi}^\dagger-A_Z\trans T_{\xi}^\dagger$, we obtain from the first equation in \eqref{equ} that
\begin{align}\label{equ2}
\Biggl[
\begin{matrix}
u\\
v
\end{matrix}
\Biggl]&=\E\{\sigma^2(X,Z){G} {G}\trans\}\Biggl[
\begin{matrix}
N_{\xi}^\dagger\\
T_{\xi}^\dagger
\end{matrix}
\Biggl]+\E\Big[\sigma^2(X,Z)\big\{K_{s}(X\trans\beta_0)-A_X(X\trans\beta_0)\trans Q_{\beta_0} N_{\xi}^\dagger-A_Z(X\trans\beta_0)\trans T_{\xi}^\dagger\big\}\,{G}\Big]\notag\\
&=\Biggl[
\begin{matrix}
Q_{\beta_0}\trans A_X(s)\\
A_Z(s)
\end{matrix}
\Biggl]+\left[\E\{\sigma^2(X,Z){G} {G}\trans\}-\E\left\{\sigma^2(X,Z){G}\Biggl[
\begin{matrix}
Q_{\beta_0}\trans A_X(X\trans\beta_0)\\
A_Z(X\trans\beta_0)
\end{matrix}
\Biggl]\trans\right\}\,\right]\Biggl[
\begin{matrix}
N_{\xi}^\dagger\\
T_{\xi}^\dagger
\end{matrix}
\Biggl]\\
&=\Biggl[
\begin{matrix}
Q_{\beta_0}\trans A_X(s)\\
A_Z(s)
\end{matrix}
\Biggl]+\E\left\{\sigma^2(X,Z)\Biggl[
\begin{matrix}
g_0'(X\trans \beta_0) Q_{\beta_0}\trans X\\
Z
\end{matrix}
\Biggl]\Biggl[
\begin{matrix}
Q_{\beta_0}\trans \{g_0'(X\trans\beta_0)X-A_X(X\trans\beta_0)\} \\
Z-A_Z(X\trans\beta_0)
\end{matrix}
\Biggl]\trans\right\}\Biggl[
\begin{matrix}
N_{\xi}^\dagger\\
T_{\xi}^\dagger
\end{matrix}
\Biggl]\,,\notag
\end{align}
where, in view of \eqref{eab}, we have used the fact that
\begin{align}\label{to}
\E\big\{\sigma^2(X,Z)K_{s}(X\trans\beta_0)\,{G}\big\}
=\Biggl[
\begin{matrix}
Q_{\beta_0}\trans\l A_X,K_{s}\r_K\\
\l A_Z,K_{s}\r_K
\end{matrix}
\Biggl]=\Biggl[
\begin{matrix}
Q_{\beta_0}\trans A_X(s)\\
A_Z(s)
\end{matrix}
\Biggl]\,.
\end{align}
Therefore, we deduce from \eqref{equ2} that
\begin{align}\label{ntu}
\Biggl[
\begin{matrix}
N_{\xi}^\dagger\\
T_{\xi}^\dagger
\end{matrix}
\Biggl]&=\left(\E\left\{\sigma^2(X,Z)G\Biggl[
\begin{matrix}
Q_{\beta_0}\trans \{g_0'(X\trans\beta_0) X -A_X(X\trans\beta_0)\} \\
Z-A_Z(X\trans\beta_0)
\end{matrix}
\Biggl]\trans\right\}\right)^{-1}\Biggl[\begin{matrix}
u-Q_{\beta_0}\trans A_X(s)\\
v-A_Z(s)
\end{matrix}\Biggl]\,.
\end{align}

Observe that \eqref{identity} implies that
\begin{align*}
\E\left\{\sigma^2(X,Z)\Bigg[
\begin{matrix}
g_0'(X\trans\beta_0)X\\
Z
\end{matrix}
\Bigg]\Bigg[
\begin{matrix}
g_0'(X\trans\beta_0)R_X(X\trans\beta_0)\\
R_Z(X\trans\beta_0)
\end{matrix}
\Bigg]\trans\right\}=\E\left\{\sigma_0^2(X\trans\beta_0)\Bigg[
\begin{matrix}
g_0'(X\trans\beta_0)R_X(X\trans\beta_0)\\
R_Z(X\trans\beta_0)
\end{matrix}
\Bigg]^{\otimes2}\right\}.
\end{align*}
Therefore,
\begin{align*}
\Omega=\E\{\sigma^2(X,Z)GG\trans\}-\E\left\{\sigma_0^2(X\trans\beta_0)\Bigg[
\begin{matrix}
g_0'(X\trans\beta_0)R_X(X\trans\beta_0)\\
R_Z(X\trans\beta_0)
\end{matrix}
\Bigg]^{\otimes2}\right\}.
\end{align*}
Therefore, combining the above equation with \eqref{ntu} yields 
\begin{align*}
\Sigma_{\beta_0,\lambda}&=\E\left\{\sigma_0^2(X\trans\beta_0)\Bigg[
\begin{matrix}
g_0'(X\trans\beta_0)Q_{\beta_0}\trans R_X(X\trans\beta_0)\\
R_Z(X\trans\beta_0)
\end{matrix}
\Bigg]^{\otimes2}\right\}\\
&\qquad-\E\left\{\sigma^2(X,Z)\Bigg[
\begin{matrix}
g_0'(X\trans\beta_0)Q_{\beta_0}\trans R_X(X\trans\beta_0)\\
R_Z(X\trans\beta_0)
\end{matrix}
\Bigg]\Biggl[
\begin{matrix}
Q_{\beta_0}\trans A_X(X\trans\beta_0) \\
A_Z(X\trans\beta_0)
\end{matrix}
\Biggl]\trans\right\}\\
&=\E\left\{\sigma_0^2(X\trans\beta_0)\Biggl[
\begin{matrix}
g_0'(X\trans\beta_0)Q_{\beta_0}\trans R_X(X\trans\beta_0)\\
R_Z(X\trans\beta_0)
\end{matrix}
\Biggl]\Biggl[
\begin{matrix}
Q_{\beta_0}\trans \{g_0'(X\trans\beta_0)R_X(X\trans\beta_0)- A_X(X\trans\beta_0)\}\\
R_Z(X\trans\beta_0)-A_Z(X\trans\beta_0)
\end{matrix}
\Biggl]\trans\right\}\\
&={V}\left\{\Biggl[
\begin{matrix}
Q_{\beta_0}\trans g_0'\cdot R_X\\
R_Z
\end{matrix}
\Biggl],\Biggl[
\begin{matrix}
Q_{\beta_0}\trans (g_0'\cdot R_X-A_X)\\
R_Z-A_Z
\end{matrix}
\Biggl]\trans\right\}={V}\left\{\Biggl[
\begin{matrix}
Q_{\beta_0}\trans g_0'\cdot R_X\\
R_Z
\end{matrix}
\Biggl],\Biggl[
\begin{matrix}
Q_{\beta_0}\trans M_\lambda(g_0'\cdot R_X)\\
M_\lambda(R_Z)
\end{matrix}
\Biggl]\trans\right\}.
\end{align*}

To show $\lim_{\lambda\to0}\Sigma_{\beta_0,\lambda}=0$, observing \eqref{expansion} and \eqref{abx}, we obtain
\begin{align*}
&g_0'\cdot R_X-A_X=M_\lambda(g_0'\cdot R_X)=\lambda\sum_{j\geq1}\frac{{V}(g_0'\cdot R_{X},\phi_j)}{1+\lambda\rho_j}\rho_j\phi_j;\\
&R_Z-A_Z=M_\lambda(R_Z)=\lambda\sum_{j\geq1}\frac{{V}(R_Z,\phi_j)}{1+\lambda\rho_j}\rho_j\phi_j.
\end{align*}
Therefore, we obtain
\begin{align*}
\Sigma_{\beta_0,\lambda}&=\lambda {V}\left\{\Biggl[
\begin{matrix}
Q_{\beta_0}\trans\sum_{j}{V}(g_0'\cdot R_{X},\phi_j)\phi_j\\
\sum_{j}{V}(R_Z,\phi_j)\phi_j
\end{matrix}
\Biggl],\Biggl[
\begin{matrix}
Q_{\beta_0}\trans\sum_{j}\frac{{V}(g_0'\cdot R_{X},\phi_j)}{1+\lambda\rho_j}\rho_j\phi_j\\
\sum_{j}\frac{{V}(R_Z,\phi_j)}{1+\lambda\rho_j}\rho_j\phi_j
\end{matrix}
\Biggl]\trans\right\}\\
&=\lambda\left[
\begin{matrix}
\sum_{j}\frac{Q_{\beta_0}\trans\{{V}(g_0'\cdot R_{X},\phi_j)\}^{\otimes2}Q_{\beta_0}}{1+\lambda\rho_j}\rho_j & \sum_{j}\frac{Q_{\beta_0}\trans{V}(g_0'\cdot R_{X},\phi_j){V}(R_Z,\phi_j)\trans }{1+\lambda\rho_j}\rho_j\\
\sum_{j}\frac{Q_{\beta_0}\trans{V}(R_Z,\phi_j){V}(g_0'\cdot R_{X},\phi_j)\trans}{1+\lambda\rho_j}\rho_j& \sum_{j}\frac{\{{V}(R_Z,\phi_j)\}^{\otimes2}}{1+\lambda\rho_j}\rho_j
\end{matrix}
\right]
\end{align*}
Note that for each $1\leq i,j\leq p+q-1$, by Assumption~\ref{a:finite},
\begin{align*}
|e_i\trans\Sigma_{\beta_0,\lambda}e_j|^2\leq \sum_{j\geq1}\{\|V(g_0\cdot R_X,\phi_j)\|_2^2+\|V(R_Z,\phi_j)\|_2^2\}<+\infty.
\end{align*}
Therefore, $\lim_{\lambda\to0}\Sigma_{\beta_0,\lambda}=0_{p+q-1,p+q-1}$ follows from the dominated convergence theorem, which concludes the proof of Proposition~\ref{prop:ru}.

\end{proof}


\subsubsection*{Proof regarding $S_{x,z}$ in Proposition~\ref{cor:sxz}}

We deduce from Proposition~\ref{prop:ru} that for $(x,z)\in\X\times\Z$ and $(\tilde g,\tilde\beta)\in\H_m\times\B$, it holds that
\begin{align*}
g(x\trans\tilde\beta)+\tilde g\,'(x\trans\tilde\beta)x\trans  Q_{\beta_0}\theta+z\trans\gamma&=f_{(g,\theta,\gamma)}\big\{x\trans\tilde\beta,\tilde g\,'(x\trans\tilde\beta)Q_{\beta_0}\trans x,z\big\}\\
&=\big\l S^\dagger_{x\trans\tilde\beta,\tilde g'(x\trans\tilde\beta)Q_{\beta_0}\trans x,z},(g,\theta,\gamma)\big\r
=\l S_{x,z}(\tilde g,\tilde\beta),(g,\theta,\gamma)\r.
\end{align*}
This completes the proof of Proposition~\ref{cor:sxz}.

\subsubsection{A useful corollary (Corollary~\ref{n}) and its proof}\label{app:prep3}

Proposition~\ref{prop:ru} enables us to obtain the relation between $|f_{(g,\theta,\gamma)}|$ and $\|(g,\theta,\gamma)\|$, which is stated in the following lemma.

\begin{corollary}\label{n}
It is true that the $f$ in \eqref{f} and $S^\dagger_{s,u,v}$ in Proposition~\ref{prop:ru} satisfies
\begin{align*}
& \sup_{s\in\I}|f_{(g,\theta,\gamma)}(s,u,v)|\leq\|(g,\theta,\gamma)\|\sup_{s\in\I}\|S^\dagger_{s,u,v}\|\leq c'\lambda ^{-1/(4m)}(1+\|u\|_2+\|v\|_2)\|(g,\theta,\gamma)\|,
\end{align*}
where $c'>0$ is a generic constant independent of $u$ and $v$.
\end{corollary}

%

\begin{proof}[\underline{Proof of Corollary~\ref{n}}]

By Proposition~\ref{prop:ru}, it holds that
\begin{align*}
\|S^\dagger_{s,u,v}\|^2&=\l S^\dagger_{s,u,v},S^\dagger_{s,u,v}\r=H^\dagger_{s,u,v}(s)+v\trans N^\dagger_{s,u,v} +v\trans T^\dagger_{s,u,v}
\\
&=K(s,s)+\Biggl[\begin{matrix}
u-A_X(s)\\
v-A_Z(s)
\end{matrix}\Biggl]\trans(\Omega_{\beta_0}+\Sigma_{\beta_0,\lambda})^{-1}\Biggl[\begin{matrix}
u-A_X(s)\\
v-A_Z(s)
\end{matrix}\Biggl].
\end{align*}
This implies that
\begin{align*}
\|S^\dagger_{s,u,v}\|^2&\leq K(s,s)+\big\|(\Omega_{\beta_0}+\Sigma_{\beta_0,\lambda})^{-1}\big\|^2_2\times\big\{\|v\|_2^2+\|z\|_2^2+\|A_X(s)\|_2^2+\|A_Z(s)\|_2^2\big\}.
\end{align*}
In addition, it is true in view of Proposition~\ref{prop:eigen} that
\begin{align}\label{kss}
\sup_{s\in\I}K(s,s)=\sup_{s\in\I}\sum_{j\geq1}\frac{\phi_j^2(s)}{1+\lambda\rho_j}\leq \sum_{j\geq1}\frac{\|\phi_j\|_\infty^2}{1+\lambda\rho_j}\leq c\lambda^{-1/(2m)}\,.
\end{align}
By Assumption~\ref{a:finite} and Proposition~\ref{prop:eigen} and the Cauchy-Schwarz inequality,
\begin{align*}
\sup_{s\in\I}\|A_X(s)\|_2&=\sup_{s\in\I}\bigg\|\sum_{j\geq1}\frac{{V}(g_0'\cdot R_X,\phi_j)}{1+\lambda\rho_j}\,\phi_j(s)\bigg\|_2\leq \bigg\|\sum_{j\geq1}\frac{{V}(g_0'\cdot R_X,\phi_j)}{1+\lambda\rho_j}\bigg\|_2\\
&\leq c\bigg\{\sum_{j\geq1}\|{V}(g_0'\cdot R_X,\phi_j)\|_2^2\bigg\}^{1/2}\times\bigg\{\sum_{j\geq1}\frac{1}{(1+\lambda\rho_j)^2}\bigg\}^{1/2}\leq c\lambda^{-1/(2m)}.
\end{align*}
Similarly, we have $\sup_{s\in\I}\|A_Z(s)\|_2\leq c\lambda^{-1/(2m)}$. In addition, note that $\lambda_{\min}(\Omega_{\beta_0})\geq c$ due to Assumption~\ref{a:rank}; see the proof of Proposition~\ref{prop:2.1}. Moreover, we obtain from Theorem~\ref{prop:ru} that $\lim_{\lambda\to0}\Sigma_{\beta_0,\lambda}=0$, so that $\|(\Omega_{\beta_0}+\Sigma_{\beta_0,\lambda})^{-1}\|\leq c$, for $n$ large enough. Since $\X$ and $\Z$ are compact, we obtain from the above derivation that
\begin{align*}
\sup_{s\in\I}\|S^\dagger_{s,u,v}\|\lesssim \lambda^{-1/(4m)}(1+\|u\|_2+\|v\|_2).
\end{align*}

Furthermore, by the Cauchy-Schwarz inequality, we obtain that
\begin{align*}
\sup_{s\in\I}|f_{(g,\theta,\gamma)}(s,u,v)|&=\sup_{s\in\I}\big|\big\l S^\dagger_{s,u,v},(g,\theta,\gamma)\big\r\big|\leq\|(g,\theta,\gamma)\|\times\sup_{s\in\I}\|S^\dagger_{s,u,v}\|
\\
&\lesssim \lambda^{-1/(4m)}(1+\|u\|_2+\|v\|_2)\times\|(g,\theta,\gamma)\|.
\end{align*}
The proof of Corollary~\ref{n} is therefore complete.

\end{proof}

\subsubsection{Proof regarding $S_\lambda^{\circ}$ in Proposition~\ref{cor:sxz}}\label{app:proof:plambda}

For $\tilde w=(\tilde g,\tilde\theta,\tilde\gamma)\in\Theta$, by \eqref{ru1}, we find
\begin{align}\label{ru222}
\big\l S_\lambda^\circ (g),\tilde w\big\r&=\big\l(H^\circ_{\lambda}(g),N^\circ_{\lambda}(g),T^\circ_{\lambda}(g)),(\tilde g,Q_{\beta_0}\trans\tilde\beta,\tilde\gamma)\big\r\notag\\
%
%
&=\Biggl[\begin{matrix}
\tilde\theta\,\,\\
\tilde\gamma\,\,
\end{matrix}
\Biggl]\trans \E\big\{\sigma^2(X,Z)H_\lambda^\circ(g)(X\trans\beta_0){G}+\sigma^2(X,Z){G}{G}\trans\big\}\Biggl[\begin{matrix}
N_\lambda^\circ(g)\\
T_\lambda^\circ(g)
\end{matrix}
\Biggl]\notag\\
&\quad+\l H_\lambda^\circ(g)+A_X\trans Q_{\beta_0} N_\lambda^\circ(g)+A_Z\trans T_\lambda^\circ(g),\tilde g\r_K,
\end{align}
where ${G}$ is defined in \eqref{gu}, and where, in view of \eqref{eab}, we used the fact that
\begin{align*}
\E\big\{\sigma^2(X,Z)\tilde g(X\trans\beta_0){G}\big\}=\Biggl[
\begin{matrix}
\l A_X,\tilde g\r_K\\
\l A_Z,\tilde g\r_K
\end{matrix}
\Biggl]\,.
\end{align*}
Note that in view of $M_\lambda$ defined in \eqref{mlambda}, we have
\begin{align}\label{ru22}
\lambda J(g,\tilde g)=\l M_\lambda(g),\tilde g\r_K\,.
\end{align}
Therefore, by comparing \eqref{ru222} and \eqref{ru22}, we find
\begin{align*}
\left\{
\begin{array}{l}
0=\E\big\{\sigma^2(X,Z)H_\lambda^\circ(g)(X\trans\beta_0)\,{G}\big\}+\E\{\sigma^2(X,Z){G} {G}\trans\}\Biggl[
\begin{matrix}
N_\lambda^\circ(g)\\
T_\lambda^\circ(g)
\end{matrix}
\Biggl]\,,\\
M_\lambda(g)=H_\lambda^\circ(g)+A_X\trans Q_{\beta_0} N_\lambda^\circ(g)+A_Z\trans T_\lambda^\circ(g)\,.
\end{array}
\right.
\end{align*}
Substituting $H_\lambda^\circ(g)$ by $M_\lambda(g)-A_X\trans Q_{\beta_0} N_\lambda^\circ(g)-A_Z\trans T_\lambda^\circ(g)$, we obtain from the first equation in the above equation system that
\begin{align*}
0&=\E\{\sigma^2(X,Z){G} {G}\trans\}\Biggl[
\begin{matrix}
N_\lambda^\circ(g)\\
T_\lambda^\circ(g)
\end{matrix}
\Biggl]+\E\Big[\sigma^2(X,Z)\big\{M_\lambda(g)(X\trans\beta_0)-A_X(X\trans\beta_0)\trans Q_{\beta_0} N_\lambda^\circ(g)-A_Z(X\trans\beta_0)\trans T_\lambda^\circ(g)\big\}{G}\Big]\notag\\
&=\E\big\{\sigma^2(X,Z)M_\lambda(g)(X\trans\beta_0)\,{G}\big\}+\left[\E\{\sigma^2(X,Z){G} {G}'\}-\E\Bigg\{\sigma^2(X,Z){G}\Biggl[
\begin{matrix}
Q_{\beta_0}\trans A_X(X\trans\beta_0)\\
A_Z(X\trans\beta_0)
\end{matrix}
\Biggl]\trans \Bigg\}\,\right]\Biggl[
\begin{matrix}
N_\lambda^\circ(g)\\
T_\lambda^\circ(g)
\end{matrix}
\Biggl]\notag\\
&=\E\big\{\sigma^2(X,Z)M_\lambda(g)(X\trans\beta_0)\,{G}\big\}+(\Omega_{\beta_0}+\Sigma_{\beta_0,\lambda})\Biggl[
\begin{matrix}
N_\lambda^\circ(g)\\
T_\lambda^\circ(g)
\end{matrix}
\Biggl]\,.
\end{align*}
Therefore, we deduce from the above equation that
\begin{align*}
\Biggl[
\begin{matrix}
N_\lambda^\circ(g)\\
T_\lambda^\circ(g)
\end{matrix}
\Biggl]&=-(\Omega_{\beta_0}+\Sigma_{\beta_0,\lambda})^{-1}\E\big\{\sigma^2(X,Z)M_\lambda(g)(X\trans\beta_0)\,{G}\big\}\notag\\
&=-(\Omega_{\beta_0}+\Sigma_{\beta_0,\lambda})^{-1}\E\Bigg\{\sigma^2(X,Z)M_\lambda(g)(X\trans\beta_0)\Biggl[
\begin{matrix}
g_0'(X\trans\beta_0)\,Q_{\beta_0}\trans X\\
Z
\end{matrix}
\Biggl]\Bigg\}\notag\\
&=-(\Omega_{\beta_0}+\Sigma_{\beta_0,\lambda})^{-1}\E\Bigg\{\sigma_0^2(X\trans\beta_0)M_\lambda(g)(X\trans\beta_0)\Biggl[
\begin{matrix}
g_0'\cdot R_X(X\trans\beta_0)\\
R_Z(X\trans\beta_0)
\end{matrix}
\Biggl]\Bigg\}\\
&=-(\Omega_{\beta_0}+\Sigma_{\beta_0,\lambda})^{-1}
\Bigg[\begin{matrix}
Q_{\beta_0}\trans{V}(g_0'\cdot R_X,M_\lambda g)\\
{V}(R_Z,M_\lambda g)
\end{matrix}\Bigg],
\end{align*}
which concludes the proof.

\subsection{Proof of Proposition~\ref{thm:score}}\label{app:thm:score}

From Lemma~\ref{prop:jointee}, we obtain that, for any $(g,\theta,\gamma)\in\H_m\times\mathbb R^{p-1}\times\mathbb R^q$, it is true that
\begin{align*}
0&=F_n(g,\theta,\gamma)\\
&=-\frac{1}{n}\sum_{i=1}^n\big\{Y_i-\hat g_n(X_i\trans \hat\beta_n)-Z_i\trans\hat\gamma_n\big\}\Big\{ g(X_i\trans \hat\beta_n)+\hat g_n'(X_i\trans \hat\beta_n)X_i\trans Q_{\beta_0}\theta+Z_i\trans\gamma\Big\}+\lambda J(\hat g_n, g)\\
&\quad-\frac{1}{n}\sum_{i=1}^n\big\{Y_i-\hat g_n(X_i\trans \hat\beta_n)-Z_i\trans\hat\gamma_n\big\}\hat g_n'(X_i\trans \hat\beta_n)X_i\trans (Q_{\hat\beta_n}-Q_{\beta_0})\theta.
\end{align*}
Therefore, by applying Proposition~\ref{cor:sxz}, we obtain that, for any $w=(g,Q_{\beta_0}\trans\beta,\gamma)\in\Theta$, it holds that
\begin{align*}
&\bigg\l-\frac{1}{n}\sum_{i=1}^n\{Y_i-\hat g_n(X_i\trans \hat\beta_n)-Z_i\trans\hat\gamma_n\} S_{X_i,Z_i}(\hat g_n,\hat\beta_n)+S_\lambda^\circ(\hat g_n),w\bigg\r\\
&\qquad=\frac{1}{n}\sum_{i=1}^n\{Y_i-\hat g_n(X_i\trans \hat\beta_n)-Z_i\trans\hat\gamma_n\}\hat g_n'(X_i\trans \hat\beta_n)X_i\trans (Q_{\hat\beta_n}-Q_{\beta_0})\theta.
\end{align*}
Note that Proposition~\ref{prop:2.1} implies that $\|\theta\|_2\leq\{\lambda_{\min}(\Omega_{\beta_0})\}^{-1/2}\|w\|$ for $w=(g,\theta,\gamma)\in\Theta$. We therefore deduce from the above equation that
\begin{align}\label{temp}
&\|\mathcal R_n(\hat g_n,\hat\beta_n,\hat\gamma_n)\|\notag\\
&=\bigg\|-\frac{1}{n}\sum_{i=1}^n\{Y_i-\hat g_n(X_i\trans \hat\beta_n)-Z_i\trans\hat\gamma_n\} S_{X_i,Z_i}(\hat g_n,\hat\beta_n)+S_\lambda^\circ(\hat g_n)\bigg\|\notag\\
&=\sup_{\|w\|=1}\bigg|\bigg\l-\frac{1}{n}\sum_{i=1}^n\{Y_i-\hat g_n(X_i\trans \hat\beta_n)-Z_i\trans\hat\gamma_n\} S_{X_i,Z_i}(\hat g_n,\hat\beta_n)+S_\lambda^\circ(\hat g_n),w\bigg\r\bigg|\notag\\
&=\sup_{\|w\|=1}\bigg|\frac{1}{n}\sum_{i=1}^n\{Y_i-\hat g_n(X_i\trans \hat\beta_n)-Z_i\trans\hat\gamma_n\}\hat g_n'(X_i\trans \hat\beta_n)X_i\trans (Q_{\hat\beta_n}-Q_{\beta_0})\theta\bigg|\notag\\
&\leq\{\lambda_{\min}(\Omega_{\beta_0})\}^{-1/2} \bigg\|\frac{1}{n}\sum_{i=1}^n\{Y_i-\hat g_n(X_i\trans \hat\beta_n)-Z_i\trans\hat\gamma_n\}\hat g_n'(X_i\trans \hat\beta_n)X_i\bigg\|_2\|Q_{\hat\beta_n}-Q_{\beta_0}\|_2,
\end{align}
By Theorem~\ref{thm:rate} and Proposition~\ref{prop:2.1}, it holds that, for $r_n=\lambda^{1/2}+n^{-1/2}\lambda^{-1/(4m)}$,
\begin{align*}
&\|Q_{\hat\beta_n}-Q_{\beta_0}\|_2\leq \|\hat\beta_n-\beta_0\|_2=\|Q_{\beta_0}Q_{\beta_0}\trans(\hat\beta_n-\beta_0)\|_2+\|\beta_0\beta_0\trans(\hat\beta_n-\beta_0)\|_2\leq c\|\hat w_n-w_0\|=O_p(r_n).
\end{align*}
Therefore, we obtain from \eqref{temp} that
\begin{align}\label{new09}
&\|\mathcal R_n(\hat g_n,\hat\beta_n,\hat\gamma_n)\|\leq O_p(r_n)\times\big(\|\mathbb W_{1,n}\|_2+\|\mathbb W_{2,n}\|_2\big),
\end{align}
where
\begin{align}\label{v12}
&\mathbb W_{1,n}=\frac{1}{n}\sum_{i=1}^n\e_i\hat g_n'(X_i\trans \hat\beta_n)X_i,\notag\\
&\mathbb W_{2,n}=\frac{1}{n}\sum_{i=1}^n\big\{g_0(X_i\trans\beta_0)+Z_i\trans\gamma_0-\hat g_n(X_i\trans \hat\beta_n)-Z_i\trans\hat\gamma_n\big\}\hat g_n'(X_i\trans \hat\beta_n)X_i.
\end{align}

For the first term $\mathbb W_{1,n}$ in \eqref{v12}, define
\begin{align}\label{ccn}
\mathcal C_{c_1,c_2}(n)&=\big\{(g,\beta)\in\mathcal H_m\times\mathcal B :\|(g,Q_{\beta_0}\trans\beta,0)-(g_0,0,0)\|\leq c_1r_n,\|g'\|_\infty\leq c_2\big\},\notag\\
\mathcal V_{c_1,c_2}(n)&=\big\{f_{g,\beta}:f_{g,\beta}(x)=g(x\trans\beta)x;x\in\X;(g,\beta)\in \mathcal C_{c_1,c_2}(n)\big\}.
\end{align}
Recalling the definition of $\|\cdot\|$-norm in \eqref{norm}, it follows from Theorem~\ref{thm:rate} and Lemma~\ref{prop:2.1} that
\begin{align}\label{new5}
&\|(\hat g_n,Q_{\beta_0}\trans\hat\beta_n,0)-(g_0,0,0)\|^2\notag
\\
&\leq2\|(\hat g_n,Q_{\beta_0}\trans\hat\beta_n,\hat\gamma_n)-(g_0,0,\gamma_0)\|^2+2\|\E\{\sigma^2(X,Z)ZZ\trans\}\|_2\times\|\hat\gamma_n-\gamma_0\|_2^2=O_p(r_n^2).
\end{align}
In addition, by Theorem~\ref{prop:rate}, it holds that $\|\hat g_n'\|_\infty=O_p(1)$.
This implies that for any $e>0$, there exists $c_{1},c_{2}>0$ (depending on $e$) such that $\P\{(\hat g_n,\hat\beta_n)\notin\mathcal C_{c_{1},c_{2}}(n)\}<e$ for $n$ large enough. Since $\X$ is compact, we obtain that
\begin{align*}
&\sup_{x\in\X}\|g_1(x\trans\beta_1)x-g_2(x\trans\beta_2)x\|_2\leq c\sup_{x\in\X}|g_1(x\trans\beta_1)-g_2(x\trans\beta_2)|\\
&\qquad=c\sup_{x\in\X}|g_1(x\trans\beta_1)-g_1(x\trans\beta_2)|+c\sup_{x\in\X}|g_1(x\trans\beta_2)-g_2(x\trans\beta_2)|\\
&\qquad\leq c\|g_1'\|_\infty\|\beta_1-\beta_2\|_2+c\|g_1-g_2\|_\infty\leq cc_2\|\beta_1-\beta_2\|_2+c\|g_1-g_2\|_\infty.
\end{align*}
Following similar argument as in the proof of Lemma~\ref{rate8} in Section~\ref{app:rate8}, we obtain that
\begin{align*}
{\log N_{}(\eta,{\mathcal V}_{c_1,c_2}(n),\|\cdot\|_\infty)}\leq c\eta^{-1/m},
\end{align*}
where $c>0$ depends only on the constants $c_1,c_2$. Furthermore, it holds that, for any $f_{g,\beta}\in\mathcal V_{c_1,c_2}(n)$, $\|f_{g,\beta}\|_\infty\leq c\|g\|_\infty+c\lambda^{-1/(4m)}r_n\leq c_0$, for some constant $c_0>0$. By the maximal inequality in Corollary~2.2.5 in \cite{vaart1996}, it holds that, for any constant $e>0$ and $\delta>0$,
\begin{align*}
&\P\big(\|\mathbb W_{1,n}\|_2\geq\delta n^{-1/2}r_n^{-1}\big)\leq\P\big\{\|\mathbb W_{1,n}\|_2\geq\delta n^{-1/2}r_n^{-1};(\hat g_n,\hat\beta_n)\in\mathcal C_{c_1,c_2}(n)\big\}+\P\{(\hat g_n,\hat\beta_n)\notin\mathcal C_{c_1,c_2}(n)\}\\
&\leq\P\left(\sup_{f\in{\mathcal V}_{c_1,c_2}(n)}\bigg\|\frac{1}{\sqrt n}\sum_{i=1}^n\e_if(X_i)\bigg\|_2>\delta r_n^{-1}\right)+e/2\leq\delta^{-1}r_n\E\left(\sup_{f\in{\mathcal V}_{c_1,c_2}(n)}\bigg\|\frac{1}{\sqrt n}\sum_{i=1}^n\e_if(X_i)\bigg\|_2\right)+e/2\\
&\leq \delta^{-1}r_n\int_0^{c_0}\sqrt{\log N_{}(\eta,{\mathcal V}_{c_1,c_2}(n),\|\cdot\|_\infty)}\,\d\eta+e/2\leq c r_n+e/2\leq e,
\end{align*}
for $n$ large enough. As a consequence, the above equation implies that
\begin{align}\label{new10}
\|\mathbb W_{1,n}\|_2=o_p(n^{-1/2}r_n^{-1}).
\end{align}

For the second term $\mathbb W_{2,n}$ in \eqref{v12}, note that
\begin{align*}
&\sup_{x\in\X,z\in\Z}|g_0(x\trans\beta_0)+v\trans\gamma_0-\hat g_n(x\trans \hat\beta_n)-v\trans\hat\gamma_n|
\\
&\leq\sup_{x\in\X,z\in\Z}|\hat g_n(x\trans\hat\beta_n)-\{\hat g_n(x\trans\beta_0)+g_0'(x\trans\beta_0)x\trans Q_{\beta_0}Q_{\beta_0}\trans\hat\beta_n\}|\\
&\quad+\sup_{x\in\X,z\in\Z}|\{\hat g_n(x\trans\beta_0)+g_0'(x\trans\beta_0)x\trans Q_{\beta_0}Q_{\beta_0}\trans\hat\beta_n+v\trans\hat\gamma_n\}-\{g_0(x\trans\beta_0)+v\trans\gamma_0\}|
\\
&=\sup_{x\in\X}|r_x(\hat g_n,\hat\beta_n)|+\sup_{x\in\X,z\in\Z}\big|\l\hat w_n-w_0,S_{x,z}(g_0,\beta_0)\r\big|.
\end{align*}
By Proposition~\ref{cor:sxz} and Corollary~\ref{n},
\begin{align*}
\sup_{x\in\X,z\in\Z}\|S_{x,z}(g_0,\beta_0)\|&=\sup_{x\in\X,z\in\Z}\|S^\dagger_{x\trans\beta_0,g_0'(x\trans\beta_0)Q_{\beta_0}\trans x,z}\|\\
&\leq c\lambda^{-1/(4m)}\sup_{x\in\X,z\in\Z}\big\{1+\|g_0'(x\trans\beta_0)Q_{\beta_0}\trans x\|_2^2+\|z\|_2^2\big\}\leq c\lambda^{-1/(4m)}.
\end{align*}
In addition, we have $\|\hat g_n\|_\infty=O_p(1)$ in view of Theorem~\ref{prop:rate}. Therefore, we deduce from the above equation, Lemma~\ref{lem:unionbound2} and Theorem~\ref{thm:rate} that
\begin{align*}
\|\mathbb W_{2,n}\|_2&\leq c\|\hat g_n'\|_\infty\times\sup_{x\in\X,z\in\Z}|g_0(x\trans\beta_0)+v\trans\gamma_0-\hat g_n(x\trans \hat\beta_n)-v\trans\hat\gamma_n|
\\
&\leq O_p(1)\times\Big[\|\hat w_n-w_0\|^{3/2}+\lambda^{-1/(4m)}\|\hat w_n-w_0\|^2+\lambda^{-1/(4m)}\|\hat w_n-w_0\|\Big]
\\
&=O_p(r_n^{3/2}+\lambda^{-1/(4m)}r_n^2+\lambda^{-1/(4m-4)}r_n)=o_p(n^{-1/2}r_n^{-1}),
\end{align*}
where we used Assumption~\ref{a:reg} in the last step. The proof is therefore complete by combining the above equation with \eqref{new09} and \eqref{new10}.

\section{Proof of theoretical results in Sections~\ref{sec:jointbahadur} and \ref{sec:marginalbahadur}}\label{app:C}

This section is organized as follows. First, in Section~\ref{app:c1} we prove the joint Bahadur representation in Theorem~\ref{thm:bahadur}. Its proof is based on several useful lemmas (Lemmas~\ref{lem:n-or}--\ref{lem:bounds}), which are proved in Sections~\ref{app:c2}--\ref{app:thm:bahadur}, respectively.
Next, in Section~\ref{app:thm:asymp} we give the proof of the joint weak convergence in Theorem~\ref{thm:asymp}.

\subsection{Proof of joint Bahadur representation in Theorem~\ref{thm:bahadur}
}\label{app:c1}

Define
\begin{align}\label{rwx}
r_x(g,\beta)=g(x\trans\beta)-\big\{g(x\trans\beta_0)-g_0'(x\trans\beta_0)x\trans Q_{\beta_0}Q_{\beta_0}\trans\beta\big\},\qquad x\in\X,(g,\beta)\in\H_m\times\B.
\end{align}
Note that $r_x(g_0,\beta_0)\equiv 0$ for any $x\in\X$. The following lemma decomposes $\hat w_n-w_0+\tilde{\mathcal R}_n(w_0)$ into several terms, and we will bound each of them separately. Its proof is given in Section~\ref{app:lem:n-or} of the supplement.

\begin{lemma}\label{lem:n-or}
For the $r_x$ in \eqref{rwx} and $S_{x,z}$ in Proposition~\ref{cor:sxz}, it is true that
\begin{align}\label{diff}
&\hat w_n-w_0+\tilde{\mathcal R}_n(w_0)=\mathcal R_n(\hat g_n,\hat\beta_n,\hat\gamma_n)+\mathbb V_{1,n}+\mathbb V_{2,n}+\mathbb V_{3,n}-n^{-1/2}\mathbb M_n(\hat w_n-w_0),
\end{align}
where for $w\in\Theta$,
\begin{align*}
&\mathbb V_{1,n}=\frac{1}{n}\sum_{i=1}^n\e_i\{S_{X_i,Z_i}(\hat g_n,\hat\beta_n)-S_{X_i,Z_i}(g_0,\beta_0)\},\notag\\
&\mathbb V_{2,n}=-\frac{1}{n}\sum_{i=1}^nr_{X_i}(\hat g_n,\hat\beta_n) S_{X_i,Z_i}(\hat g_n,\hat\beta_n),\notag\\
&\mathbb V_{3,n}=\frac{1}{n}\sum_{i=1}^n\l\hat w_n-w_0,S_{X_i,Z_i}(g_0,\beta_0)\r\{S_{X_i,Z_i}(g_0,\beta_0)-S_{X_i,Z_i}(\hat g_n,\hat\beta_n)\},\\
&\M_n(w)=\frac{1}{\sqrt n}\sum_{i=1}^n\Big[\l w,S_{X_i,Z_i}(g_0,\beta_0)\r S_{X_i,Z_i}(g_0,\beta_0)-\E\big\{\l w,S_{X,Z}(g_0,\beta_0)\r S_{X,Z}(g_0,\beta_0)\big\}\Big].
\end{align*}
\end{lemma}
The following lemma states the union bounds for $r_x$ and $S_{x,z}$ and is proved in Section~\ref{app:lem:suprdiff} of the supplement. In particular, the bound on $r_x$ characterizes the approximation error in the expansion \eqref{ad}.

\begin{lemma}\label{lem:unionbound2}
Suppose for $i=1,2$, $\|g_i\|_\infty^2+J(g_i,g_i)<\infty$ and $\beta_i\trans\beta_0>0$. Then, it holds that 
\begin{gather*}
\sup_{x\in\X}|r_x(g_1,\beta_1)|\lesssim\|\beta_1-\beta_0\|_2^{3/2}+\|g'_1-g_0'\|_\infty\times\|\beta_1-\beta_0\|_2,
\\
\sup_{x\in\X,z\in\mathcal Z}\|S_{x,z}(g_1,\beta_1)-S_{x,z}(g_2,\beta_2)\|\lesssim\lambda^{-1/(4m-4)}\|\beta_1-\beta_2\|_2+\|\beta_1-\beta_2\|_2^{1/2}+\|g_1'-g_2'\|_\infty.
\end{gather*}
\end{lemma}
Proposition~\ref{thm:score} together with Lemma~\ref{lem:n-or} and \ref{lem:unionbound2} enable us to derive bounds for each term on the right-hand side of \eqref{diff} w.r.t.~the $\|\cdot\|$-norm. The bound on $\|\mathcal R_n(\hat g_n,\hat\beta_n,\hat\gamma_n)\|$ is derived in Proposition~\ref{thm:score}; the bound on $\|\mathbb M_n(\hat w_n-w_0)\|$ and $\|\mathbb V_{1,n}\|$ follows from the empirical process theory; $\|\mathbb V_{2,n}\|$ and $\|\mathbb V_{3,n}\|$ are bounded via the union bound in Lemma~\ref{lem:unionbound2}. These bounds are formally stated in the following lemma, which is proved in Section~\ref{app:thm:bahadur} of the supplement.

\begin{lemma}\label{lem:bounds}
Under Assumptions~\ref{a:rank}--\ref{a:finite}, it holds that $\|\mathbb V_{1,n}\|=o_p(n^{-1/2})$, and
\begin{align*}
&\|\mathbb V_{2,n}\|+\|\mathbb V_{3,n}\|=O_p\big\{\lambda^{-1/(4m)} (\lambda^{-1/(4m-4)}r_n^2+r_n^{3/2})\big\},\qquad \|\mathbb M_n(\hat w_n-w_0)\|=O_p(a_n),
\end{align*}
where $a_n=\lambda^{-(4m-1)/(8m^2)}r_n\sqrt{\log\log(n)}$ and $r_n$ is defined in Theorem~\ref{thm:rate}.
\end{lemma}

Now, by applying Assumption~\ref{a:reg}, the proof of the joint Bahadur representation is complete by combining Lemmas~\ref{lem:n-or}--\ref{lem:bounds}.

Next, we will prove Lemmas~\ref{lem:n-or}--\ref{lem:bounds} in Sections~\ref{app:c2}--\ref{app:thm:bahadur}, respectively.

\subsubsection{Proof of Lemma~\ref{lem:n-or}}\label{app:lem:n-or}\label{app:c2}

Recall from \eqref{proxy} that, for $w=(g,\theta,\gamma)\in\H_m\times\mathbb R^{p-1}\times\mathbb R^q$,
\begin{align}
\tilde {\mathcal R}_n(w):=-\frac{1}{n}\sum_{i=1}^n\big\{Y_i-\big\l w,S_{X_i,Z_i}(g_0,\beta_0)\big\r\big\}S_{X_i,Z_i}(g_0,\beta_0)+S_\lambda^\circ(g),
\end{align}
where $S_{x,z}(g,\beta)$ and $S_\lambda^\circ(g)$ are defined in Proposition~\ref{cor:sxz}.
Let $\hat w_{\oracle}=(\hat g_\oracle,\hat\theta_\oracle,\hat\gamma_\oracle)\in\Theta$ be the solution to the joint estimating equation $\tilde{\mathcal R}_n(w)=0$, that is,
\begin{align}\label{oracle}
\tilde{\mathcal R}_n(\hat w_\oracle)=-\frac{1}{n}\sum_{i=1}^n\big\{Y_i-\l \hat w_\oracle,S_{X_i,Z_i}(g_0,\beta_0)\r\big\}S_{X_i,Z_i}(g_0,\beta_0)+S_\lambda^\circ(\hat g_\oracle)=0.
\end{align}
Next, we investigate the asymptotic properties of the oracle estimator $\hat w_\oracle$ defined in \eqref{oracle}.
Computing the Fr\'echet derivative of the functional $\tilde{\mathcal R}_n(w)$ defined in \eqref{proxy} yields
\begin{align*}
\mathcal D \tilde{\mathcal R}_n(w)w_1=\frac{1}{n}\sum_{i=1}^n\l w_1,S_{X_i,Z_i}(g_0,\beta_0)\r S_{X_i,Z_i}(g_0,\beta_0)+S_\lambda^\circ(g_1),\qquad\quad w,w_1\in\Theta.
\end{align*}
Note that $\E\{\l\D \tilde{\mathcal R}_n(w)w_1,w_2\r\}=\l w_1,w_2\r$ for any $w_1,w_2\in\Theta$, which implies $\E\{\D \tilde{\mathcal R}_n(w)\}\equiv id$, for any $w\in\Theta$, where $id$ denotes the identity operator. Furthermore, the second order Fr\'echet derivative of $\tilde{\mathcal R}_n$ vanishes: $\D^2 \tilde{\mathcal R}_n\equiv0$. As a consequence, in view of the fact that $\tilde{\mathcal R}_n(\hat w_\oracle)=0$, we obtain
\begin{align*}
\mathcal D \tilde{\mathcal R}_n(w_0)(\hat w_\oracle-w_0)=\tilde{\mathcal R}_n(\hat w_\oracle)-\tilde{\mathcal R}_n(w_0)=-\tilde{\mathcal R}_n(w_0).
\end{align*}
The above equation and \eqref{oracle} imply
\begin{align}\label{wordiff}
\hat w_\oracle-w_0&=\D \tilde{\mathcal R}_n(w_0)(\hat w_\oracle-w_0)-\big[\D \tilde{\mathcal R}_n(w_0)(\hat w_\oracle-w_0)-\E\{\D \tilde{\mathcal R}_n(w_0)\}(\hat w_\oracle-w_0)\big]\notag
\\
&
=-\tilde{\mathcal R}_n(w_0)-n^{-1/2}\mathbb M_n(\hat w_\oracle-w_0),
\end{align}
where, for $w\in\Theta$, $\mathbb M_n(w)$ is defined in Lemma~\ref{lem:n-or}.

Next, we compute the difference $\hat w_n-\hat w_\oracle$. Observing the above equation and the fact that $\E \{\mathcal D\tilde{\mathcal R}_n(w_0)\}=id$, we obtain
\begin{align}\label{newnew}
\hat w_n-\hat w_\oracle&=\E \{\mathcal D\tilde{\mathcal R}_n(w_0)\}(\hat w_n-\hat w_\oracle)\\
&=n^{-1/2}\mathbb M_n(\hat w_\oracle-\hat w_n)+\frac{1}{ n}\sum_{i=1}^n\l \hat w_n-\hat w_\oracle,S_{X_i,Z_i}(g_0,\beta_0)\r S_{X_i,Z_i}(g_0,\beta_0)+S_\lambda^\circ(\hat g_n)-S_\lambda^\circ(\hat g_\oracle).\notag
\end{align}
It follows from Theorems~\ref{thm:score} and equation \eqref{oracle} that
\begin{align*}
&-\frac{1}{n}\sum_{i=1}^n\{Y_i-\hat g_n(X_i\trans\hat\beta_n)-Z_i\trans\hat\gamma_n\}S_{X_i,Z_i}(\hat g_n,\hat\beta_n)+S_\lambda^\circ(\hat g_n)=\mathcal R_n(\hat g_n,\hat\beta_n,\hat\gamma_n),\\
&-\frac{1}{n}\sum_{i=1}^n\{Y_i-\l \hat w_\oracle,S_{X_i,Z_i}(g_0,\beta_0)\r\}S_{X_i,Z_i}(g_0,\beta_0)+S_\lambda^\circ(\hat g_\oracle)=0,
\end{align*}
where $\mathcal R_n(\hat g_n,\hat\beta_n,\hat\gamma_n)$ is defined in \eqref{scoren2}.
We therefore obtain from the above equation that
\begin{align}\label{wdiff}
S_\lambda^\circ(\hat g_n)-S_\lambda^\circ(\hat g_\oracle)
&=\frac{1}{n}\sum_{i=1}^n\{Y_i-\hat g_n(X_i\trans\hat\beta_n)-Z_i\trans\hat\gamma_n\}S_{X_i,Z_i}(\hat g_n,\hat\beta_n)\notag\\
&\ -\frac{1}{n}\sum_{i=1}^n\{Y_i-\l \hat w_\oracle,S_{X_i,Z_i}(g_0,\beta_0)\r\}S_{X_i,Z_i}(g_0,\beta_0)+\mathcal R_n(\hat g_n,\hat\beta_n,\hat\gamma_n).
\end{align}
Recall that $\E \{\mathcal D\tilde{\mathcal R}_n(w_0)\}=id$. Therefore,
\begin{align*}
\hat w_n-\hat w_\oracle&=\E \{\mathcal D\tilde{\mathcal R}_n(w_0)\}(\hat w_n-\hat w_\oracle)\\
&=n^{-1/2}\mathbb M_n(\hat w_\oracle-\hat w_n)+\frac{1}{ n}\sum_{i=1}^n\l \hat w_n-\hat w_\oracle,S_{X_i,Z_i}(g_0,\beta_0)\r S_{X_i,Z_i}(g_0,\beta_0)+S_\lambda^\circ(\hat g_n)-S_\lambda^\circ(\hat g_\oracle)
\\
&=n^{-1/2}\mathbb M_n(\hat w_\oracle-\hat w_n)+\frac{1}{ n}\sum_{i=1}^n\l \hat w_n-\hat w_\oracle,S_{X_i,Z_i}(g_0,\beta_0)\r S_{X_i,Z_i}(g_0,\beta_0)\\
&\quad+\frac{1}{n}\sum_{i=1}^n\{Y_i-\hat g_n(X_i\trans\hat\beta_n)-Z_i\trans\hat\gamma_n\}S_{X_i,Z_i}(\hat g_n,\hat\beta_n)\notag\\
&\quad-\frac{1}{n}\sum_{i=1}^n\{Y_i-\l \hat w_\oracle,S_{X_i,Z_i}(g_0,\beta_0)\r\}S_{X_i,Z_i}(g_0,\beta_0)+\mathcal R_n(\hat g_n,\hat\beta_n,\hat\gamma_n).
\end{align*}
Rearranging the above equation yields 
\begin{align*}
\hat w_n-\hat w_\oracle
&=n^{-1/2}\mathbb M_n(\hat w_\oracle-\hat w_n)+\frac{1}{n}\sum_{i=1}^n\e_i\{S_{X_i,Z_i}(\hat g_n,\hat\beta_n)-S_{X_i,Z_i}(g_0,\beta_0)\}\\
&\quad+\frac{1}{ n}\sum_{i=1}^n\l \hat w_n-w_0,S_{X_i,Z_i}(g_0,\beta_0)\r \{S_{X_i,Z_i}(g_0,\beta_0)-S_{X_i,Z_i}(\hat g_n,\hat\beta_n)\}\\
&\quad+\frac{1}{n}\sum_{i=1}^n\{\l\hat w_n,S_{X_i,Z_i}(g_0,\beta_0)\r-\hat g_n(X_i\trans\hat\beta_n)-Z_i\trans\hat\gamma_n\}S_{X_i,Z_i}(\hat g_n,\hat\beta_n)+\mathcal R_n(\hat g_n,\hat\beta_n,\hat\gamma_n).
\end{align*}
Recalling the definition of $r_x(g,\beta)$ in \eqref{rwx}, we obtain from the above equation and \eqref{wordiff} that
\begin{align*}
&\hat w_n-w_0=(\hat w_n-\hat w_\oracle)+(\hat w_\oracle-w_0)
\\
&=n^{-1/2}\mathbb M_n(\hat w_\oracle-\hat w_n)+\mathbb V_{1,n}+\mathbb V_{2,n}+\mathbb V_{3,n}+\mathcal R_n(\hat g_n,\hat\beta_n,\hat\gamma_n)-\tilde{\mathcal R}_n(w_0)\\
&=-\tilde{\mathcal R}_n(w_0)+\mathbb V_{1,n}+\mathbb V_{2,n}+\mathbb V_{3,n}-n^{-1/2}\mathbb M_n(\hat w_n-w_0)+\mathcal R_n(\hat g_n,\hat\beta_n,\hat\gamma_n).
\end{align*}
This completes the proof of Lemma~\ref{lem:n-or}.

\subsubsection{Proof of Lemma~\ref{lem:unionbound2}}\label{app:lem:suprdiff}

First, we prove the union bound on $r_x(g,\beta)$. Recall the definition of $r_x(g,\beta)$ in \eqref{rwx} that
\begin{align*}
r_x(g,\beta)&=g(x\trans \beta)-g(x\trans   \beta_0)-g_0'(x\trans   \beta_0)x\trans Q_{\beta_0}Q_{\beta_0}\trans\beta,\qquad x\in\X.
\end{align*}
We have $r_x(g,\beta)=I_1(x)+I_2(x)+I_3(x)$, where
\begin{align}\label{i1234}
I_1(x)&=g(x\trans \beta)-g(x\trans   \beta_0)-g'(x\trans\beta_0)x\trans (\beta-\beta_0),\notag\\
I_2(x)&=g'(x\trans \beta_0)x\trans(\beta-\beta_0)-g_0'(x\trans   \beta_0)x\trans(\beta-\beta_0),\notag\\
I_3(x)&=g_0'(x\trans   \beta_0)\{x\trans(\beta-\beta_0)-x\trans Q_{\beta_0}Q_{\beta_0}\trans\beta\}.
\end{align}
For the first term $I_1$, by Taylor's theorem, we obtain
\begin{align*}
|I_1(x)|=\bigg|\int_{x\trans \beta_0}^{x\trans \beta} g''(t)(x\trans \beta-t)dt\bigg|\leq\bigg|\int_{\I} [g''(t)]^2dt\bigg|^{1/2}\times\bigg|\int_{x\trans \beta_0}^{x\trans \beta} (x\trans \beta-t)^2dt\bigg|^{1/2}.
\end{align*}
By Lemma~\ref{lem:agmon}, we obtain that for any $e\in(0,1)$,
\begin{align}\label{out}
\int_\I[g''(t)]^2dt\leq c\bigg\{e^{-2}\int_{\I}|g(t)|^2dt+e^{2m-2}\int_\I|g^{(m)}(t)|^2dt\bigg\}\leq c.
\end{align}
In addition, note that by Lemma~\ref{lem:qbeta},
\begin{align*}
\|\beta-\beta_0\|_2\leq\|Q_{\beta_0}Q_{\beta_0}\trans\beta\|_2+\|\beta_0\beta_0\trans(\beta-\beta_0)\|_2\leq\|Q_{\beta_0}\trans\beta\|_2+\|Q_{\beta_0}\trans\beta\|_2^2\leq 2\|Q_{\beta_0}\trans\beta\|_2.
\end{align*}
Furthermore, by Proposition~\ref{prop:2.1}, we obtain from the above equation that
\begin{align}\label{temp1}
\|\beta-\beta_0\|_2\leq2\|Q_{\beta_0}\trans\beta\|_2=2\|Q_{\beta_0}\trans(\beta-\beta_0)\|_2\leq c\|w-w_0\|.
\end{align}
Therefore, by Proposition~\ref{prop:2.1} we obtain that
\begin{align*}
\sup_{x\in\X}|I_1(x)|\leq c\|\beta-\beta_0\|_2^{3/2}\leq c\|w-w_0\|^{3/2}.
\end{align*}
For the second term $I_2$ in \eqref{i1234}, applying Lemma~\ref{lem:norms} and \eqref{temp1} yields
\begin{align*}
\sup_{x\in\X}|I_2(x)|\leq c\|g'-g_0'\|_\infty\times\|\beta-\beta_0\|_2\leq c\lambda^{-1/(4m-4)}\|w-w_0\|^2.
\end{align*}
For the third term $I_3$ in \eqref{i1234}, by Lemma~\ref{lem:qbeta} and \eqref{temp1},
\begin{align*}
\sup_{x\in\X}|I_3(x)|=\sup_{x\in\X}|g_0'(x\trans\beta_0)x\trans\beta_0\beta_0\trans(\beta-\beta_0)|\leq c\|g_0'\|_\infty\times|\beta\trans\beta_0-1|\leq c\|Q_{\beta_0}\trans\beta\|_2^2\leq c\|\beta-\beta_0\|^2.
\end{align*}
Combining the bounds for $I_1$--$I_3$ yields
\begin{align*}
\sup_{x\in\X}|r_x(g,\beta)|&\lesssim\|\beta-\beta_0\|_2^{3/2}+\|g'-g_0'\|_\infty\times\|\beta-\beta_0\|_2+\|\beta-\beta_0\|_2^2,
\end{align*}
which concludes the proof of the union bound on $r_x(g,\beta)$.

Next, we prove the union bound regarding $S_{x,z}(g,\beta)$. We first deduce from Proposition~\ref{cor:sxz} that, for $(g_1,\beta_1),(g_2,\beta_2)\in\H_m\times\B$,
\begin{align*}
&\l S_{x,z}(g_1,\beta_1),S_{x,z}(g_2,\beta_2)\r
\\
&=\big\l\big(H_{x,z}(g_1,\beta_1),N_{x,z}(g_1,\beta_1),T_{x,z}(g_1,\beta_1)\big),\big(H_{x,z}(g_2,\beta_2),N_{x,z}(g_2,\beta_2),T_{x,z}(g_2,\beta_2)\big)\big\r\\
&=H_{x,z}(g_2,\beta_2)(x\trans\beta_1)+g_1'(x\trans\beta_1)x\trans Q_{\beta_0}N_{x,z}(g_2,\beta_2)+z\trans T_{x,z}(g_2,\beta_2)\\
&=K(x\trans\beta_1,x\trans\beta_2)-A_X(x\trans\beta_1)\trans Q_{\beta_0} N_{x,z}(g_1,\beta_1)-A_Z(x\trans\beta_1)\trans T_{x,z}(g_1,\beta_1)\\
&\quad+g_1'(x\trans\beta_1)x\trans  Q_{\beta_0}N_{x,z}(g_2,\beta_2)+z\trans T_{x,z}(g_2,\beta_2)\\
&=\l K_{x\trans\beta_1},K_{x\trans\beta_2}\r_K+\Biggl[\begin{matrix}
Q_{\beta_0}\trans  \{g_1'(x\trans\beta_1)x-A_X(x\trans\beta_1)\}\\
z-A_Z(x\trans\beta_1)
\end{matrix}\Biggl]\trans(\Omega_{\beta_0}+\Sigma_{\beta_0,\lambda})^{-1}\Biggl[\begin{matrix}
Q_{\beta_0}\trans  \{g_2'(x\trans\beta_2)x-A_X(x\trans\beta_2)\}\\
z-A_Z(x\trans\beta_2)
\end{matrix}\Biggl].
\end{align*}
Therefore, we obtain 
\begin{align}\label{sxzdiff}
&\|S_{x,z}(g_1,\beta_1)-S_{x,z}(g_2,\beta_2)\|^2\notag\\
&=\|S_{x,z}(g_1,\beta_1)\|^2+\|S_{x,z}(g_2,\beta_2)\|^2-2\l S_{x,z}(g_1,\beta_1),S_{x,z}(g_2,\beta_2)\r\notag\\
&=\| K_{x\trans\beta_1}-K_{x\trans\beta_2}\|_K^2+\Biggl[\begin{matrix}
Q_{\beta_0}\trans\big[\{g_1'(x\trans\beta_1)-g_2'(x\trans\beta_2)\} x-\{A_X(x\trans\beta_1)-A_X(x\trans\beta_2)\}\big]\\
-\{A_Z(x\trans\beta_1)-A_Z(x\trans\beta_2)\}
\end{matrix}\Biggl]\trans\notag\\
&\quad\quad\times(\Omega_{\beta_0}+\Sigma_{\beta_0,\lambda})^{-1}\Biggl[\begin{matrix}
Q_{\beta_0}\trans\big[\{g_1'(x\trans\beta_1)-g_2'(x\trans\beta_2)\}  x-\{A_X(x\trans\beta_1)-A_X(x\trans\beta_2)\}\big]\\
-\{A_Z(x\trans\beta_1)-A_Z(x\trans\beta_2)\}
\end{matrix}\Biggl].
\end{align}
By Proposition~\ref{prop:ru}, $\lim_{\lambda\to0}\Sigma_{\beta_0,\lambda}=0$, so $\|(\Omega_{\beta_0}+\Sigma_{\beta_0,\lambda})^{-1}\|\leq c$ for $n$ large enough. 
Therefore, we obtain from the above equation and the Cauchy-Schwarz inequality that
\begin{align*}
&\|S_{x,z}(g_1,\beta_1)-S_{x,z}(g_2,\beta_2)\|^2\leq \| K_{x\trans\beta_1}-K_{x\trans\beta_2}\|_K^2+c\|\{g_1'(x\trans\beta_1)-g_0'(x\trans\beta_2)\}Q_{\beta_0}\trans x\|_2^2
\\
&\qquad\qquad+c\|A_X(x\trans\beta_1)-A_X(x\trans\beta_2)\|_2^2+c\|A_Z(x\trans\beta_1)-A_Z(x\trans\beta_2)\|_2^2
\\
&\quad\leq \| K_{x\trans\beta_1}-K_{x\trans\beta_2}\|_K^2
+c|g_1'(x\trans\beta_1)-g_2'(x\trans\beta_2)|\\
&\qquad\qquad+\|K_{x\trans\beta_1}-K_{x\trans\beta_2}\|^2\times\sum_{j\geq1}\{\|V(g_0'\cdot R_X,\phi_j)\|_2^2+\|V(R_Z,\phi_j)\|_2^2\}
\\
&\quad\leq c\|K_{x\trans\beta_1}-K_{x\trans\beta_2}\|_K^2+c|g_1'(x\trans\beta_1)-g_2'(x\trans\beta_2)|,
\end{align*}
where we applied Assumption~\ref{a:finite} in the last step.
Therefore, we obtain
\begin{align*}
\sup_{x\in\X,z\in\Z}\|S_{x,z}(g_1,\beta_1)-S_{x,z}(g_2,\beta_2)\|^2\leq c(I_4+I_5),
\end{align*}
where
\begin{align*}
I_4=\sup_{x\in\X}\|K_{x\trans\beta_1}-K_{x\trans\beta_2}\|_K^2,\qquad I_5=\sup_{x\in\X}|g_1'(x\trans\beta_1)-g_2'(x\trans\beta_2)|^2.
\end{align*}
For the first term $I_4$, by the mean-value theorem, there exists a constant $\eta$ in between $x\trans\beta_1$ and $x\trans\beta_2$ such that
\[
\|K_{x\trans\beta_1}-K_{x\trans\beta_2}\|_K^2=K(x\trans\beta_1,x\trans\beta_1)-2K(x\trans\beta_1,x\trans\beta_2)+K(x\trans\beta_2,x\trans\beta_2)=\partial_1\partial_2 K(\eta,\eta)(x\trans\beta_1-x\trans\beta_2)^2.\label{newnew}
\]
It follows from Theorem~6.3 of \cite{ferreira2013positive} that $\tilde K:=\partial_1\partial_2K$ is the reproducing kernel of the Sobolev space of order $(m-1)$ on $\I$, denoted by $\H_{m-1}$, equipped with the inner product
$$
\l h_1,h_2\r_{\tilde K}=V(h_1,h_2)+\lambda\int_{\I}h_1^{(m-1)}(s)\,h_2^{(m-1)}(s)\,\d s,\qquad h_1,h_2\in\H_{m-1}.
$$
By the same derivation as in the proof of Lemma~\ref{lem:norms}, we obtain that
$\|\tilde K\|_\infty\leq c\lambda^{-1/(2m-2)}$, which implies that $$I_4\leq c\lambda^{-1/(2m-2)}\|\beta_1-\beta_2\|_2^2.$$
For the second term $I_5$, we have
\begin{align*}
I_5&\leq2\sup_{x\in\X}|g_1'(x\trans\beta_1)-g_1'(x\trans\beta_2)|^2+2\sup_{x\in\X}|g_1'(x\trans\beta_2)-g_2'(x\trans\beta_2)|^2.
\end{align*}
In view of \eqref{out}, we obtain by the Cauchy-Schwarz inequality that
\begin{align*}
\sup_{x\in\X}|g_1'(x\trans\beta_1)-g_1'(x\trans\beta_2)|&=\sup_{x\in\X}\bigg|\int_{x\trans\beta_1}^{x\trans\beta_2}g''(\xi)\d \xi\bigg|\leq c \bigg[\int_\I \{g''(\xi)\}^2\d\xi\bigg]^{1/2}\|\beta_1-\beta_2\|_2^{1/2}\leq c\|\beta_1-\beta_2\|_2^{1/2}.
\end{align*}
Therefore, we deduce from the above two equations that
$$
I_5\leq c\|\beta_1-\beta_2\|_2+c\|g_1'-g_2'\|_\infty^2.
$$
In conclusion, we obtain
\begin{align*}
\sup_{x\in\X,z\in\Z}\|S_{x,z}(g_1,\beta_1)-S_{x,z}(g_2,\beta_2)\|^2&\lesssim \lambda^{-1/(2m-2)}\|\beta_1-\beta_2\|_2^2+\|\beta_1-\beta_2\|_2+\|g_1'-g_2'\|_\infty^2.
\end{align*}
This completes the proof of the union bounds in Lemma~\ref{lem:unionbound2}. 

Next, for the claim regarding $(\hat g_n,\hat\beta_n)$, in view of Theorem~\ref{prop:rate}, we obtain that $$(\hat g_n,\hat\beta_n)\in\{g\in\H_m:\|g\|_\infty^2+J(g,g)<\infty\}\times\{\beta\in\B:\beta\trans\beta_0>0\},$$ on an event with probability tending to 1 as $n\to\infty$. Therefore, we apply the union bounds for $r_n$ and $S_{x,z}$ to obtain from Lemmas~\ref{lem:norms} and \ref{lem:normK} that, for $w_1=(g_1,Q_{\beta_0}\trans\beta_1,\gamma_1)$ and $w_2=(g_2,Q_{\beta_0}\trans\beta_2,\gamma_2)$, 
\[
\|g_1'-g_2'\|_\infty\lesssim \lambda^{-1/(4m-4)}\|w_1-w_2\|.
\]
In addition, Proposition~\ref{prop:2.1} and Lemmas~\ref{lem:qbeta} implies that
\begin{align*}
\|\beta_1-\beta_2\|_2\leq \|Q_{\beta_0}\trans(\beta_1-\beta_2)\|_2+\|\beta_0\beta_0\trans\beta_1\|_2+\|\beta_0\beta_0\trans\beta_2\|_2\leq 2\|Q_{\beta_0}\trans(\beta_1-\beta_2)\|_2\lesssim \|w_1-w_2\|.
\end{align*}
Therefore, the claim regarding $r_x(\hat g_n,\hat\beta_n)$ and $S_{x,z}(\hat g_n,\hat\beta_n)$ follows combining the above union bounds, together with the above two equations and Theorem~\ref{thm:rate}.

\subsubsection{Proof of Lemma~\ref{lem:bounds}}\label{app:thm:bahadur}

For the first term $\mathbb V_{1,n}$ in \eqref{diff}, define
\begin{align*}
&\mathcal G_{c_1,c_2,c_3}=\{g\in\H_m:\|g\|_\infty<c_1,\|g'\|_\infty<c_2,J(g,g)\leq c_3\},
\\
&\mathcal E_{c_1,c_2,c_3,c_4}(n)=\big\{(g,\beta)\in\mathcal G_{c_1,c_2,c_3}\times\mathcal B :\|(g,Q_{\beta_0}\trans\beta,0)-(g_0,0,0)\|\leq c_4r_n\big\},
\\
&\tilde{\mathcal W}_{c_1,c_2,c_3}=\big\{f_{g,\beta}(x,z)=S_{x,z}(g,\beta)-S_{x,z}(g_0,\beta_0):S_{x,z}(g,\beta)\text{ defined in Proposition~\ref{cor:sxz}},\\
&\hspace{8cm}x\in\X;z\in\Z;(g,\beta)\in\mathcal G_{c_1,c_2,c_3}\big\},
\\
&\mathcal W_{c_1,c_2,c_3,c_4}(n)=\big\{f_{g,\beta}(x,z)=S_{x,z}(g,\beta)-S_{x,z}(g_0,\beta_0):S_{x,z}(g,\beta)\text{ defined in Proposition~\ref{cor:sxz}},\\
&\hspace{8cm}x\in\X;z\in\Z;(g,\beta)\in\mathcal E_{c_1,c_2,c_3,c_4}(n)\big\}.
\end{align*}
Note that $\mathcal W_{c_1,c_2,c_3,c_4}(n)\subset \tilde{\mathcal W}_{c_1,c_2,c_3}$.
By equation~\eqref{new5} and Theorem~\ref{prop:rate}, it holds that $\|(\hat g_n,Q_{\beta_0}\trans\hat\beta_n,0)-(g_0,0,0)\|=O_p(r_n)$, $\|\hat g_n\|_\infty=O_p(1)$, and $J(\hat g_n,\hat g_n)=O_p(1)$. Therefore, for any $e>0$, there exist $n_e\in\mathbb Z_+$ and constants $c_1,c_2,c_3,c_4>0$\  (depending on $e$) such that $\P\{(\hat g_n,\hat\beta_n)\notin\mathcal E_{c_1,c_2,c_3,c_4}(n)\}<e/2$ for $n>n_e$.

It follows from  Lemma~\ref{lem:unionbound2} that, for any $f_{g,\beta}\in\mathcal W_{c_1,c_2,c_3,c_4}(n)$, and for the constant $c>0$ in Lemma~\ref{lem:unionbound2},
\begin{align*}
&\sup_{x\in\X,z\in\Z}\|f_{g,\beta}(x,z)\|=\sup_{x\in\X,z\in\Z}\|S_{x,z}( g,\beta)-S_{x,z}(g_0,\beta_0)\|\\
&\qquad\leq c\lambda^{-1/(4m-4)}\|(w-w_0,Q_{\beta_0}\trans\beta,0)\|+c\|(w-w_0,Q_{\beta_0}\trans\beta,0)\|^{1/2}\leq cc_1(\lambda^{-1/(4m-4)}r_n+r_n^{1/2}).
\end{align*}
The above equation implies that $m_n:=cc_1(\lambda^{-1/(4m-4)}r_n+r_n^{1/2})$ is an envelop of $\mathcal W_{c_1,c_2,c_3,c_4}(n)$.

Next, we show that
\begin{align}\label{san2}
\log N\big(\eta,\mathcal W_{c_1,c_2,c_3,c_4}(n),\|\cdot\|_\infty\big)\lesssim\eta^{-1/(m-1)}.
\end{align}
Note that $\lambda^{-1}=O(n^{2m/(2m+1)})=o(n^{2m/3})$ by Assumption~\ref{a:5.0}, so that $r_n=\lambda^{1/2}+n^{-1/2}\lambda^{-1/(4m)}=O(\lambda^{1/(2m-2)})$.
It follows from Lemma~\ref{lem:unionbound2} that, for $(g_1,\beta_1),(g_2,\beta_2)\in\mathcal E_{c_1,c_2,c_3,c_4}(n)$,
\begin{align}\label{san}
&\sup_{x\in\X,z\in\Z}\|f_{g_1,\beta_1}(x,z)-f_{g_2,\beta_2}(x,z)\|=\sup_{x\in\X,z\in\Z}\|S_{x,z}(g_1,\beta_1)-S_{x,z}(g_2,\beta_2)\|\notag
\\
&\qquad\leq c\lambda^{-1/(2m-2)}\|\beta_1-\beta_2\|_2^2+c\|\beta_1-\beta_2\|_2+c\|g_1'-g_2'\|_\infty\notag
\\
&\qquad\leq  c\|\beta_1-\beta_2\|_2^{1/2}(\lambda^{-1/(4m-4)}r_n^{1/2}+1)+c\|g'_1-g_2'\|_\infty\notag
\\
&\qquad\leq c\|\beta_1-\beta_2\|_2^{1/2}+c\|g'_1-g_2'\|_\infty.
\end{align}
Observe that $g'\in\H_{m-1}$ when $g\in\H_m$.
By Theorem~2.4 in \cite{vandegeer2000} and Lemma~4.1 of \cite{pollard1990empirical}, we have, for some absolute constants $c,c'>0$, it holds that
\begin{align*}
&\log N\big(\eta,\{g':g\in\mathcal G_{c_1,c_2,c_3}\},\|\cdot\|_\infty\big)\leq c\exp(c'\eta^{-1/{(m-1)}}),
\\
&N(\eta,\mathcal B,\|\cdot\|_2)\leq c'\eta^{-p+1}.
\end{align*}
Let $\{g_1, g_2, \ldots, g_{r_1}\}$ be a $\eta$-net for $\{g':g\in\mathcal G_{c_1,c_2,c_3}\}$ w.r.t.~the $\|\cdot\|_\infty$-norm, let $\{\beta_1, \beta_2, \ldots, \beta_{r_2}\}$ be a $\eta^2$-net of $\B$ w.r.t.~$\|\cdot\|_2$-norm, where $r_1\leq c\exp(c\eta^{-1/(m-1)})$ and $r_2\leq c\eta^{-2p+2}$. Then, it holds that the set
\begin{align*}
\tilde{\mathcal{W}}:=\Big\{\tilde f_{g_i,\beta_j}:\tilde f_{g_i,\beta_j}(x,z)=S_{x,z}(g_i,\beta_j);1\leq i\leq r_1;1\leq j\leq r_2\Big\}\subset\tilde{\mathcal W}_{c_1,c_2,c_3}
\end{align*}
forms a $\eta$-bracket of $\tilde{\mathcal W}_{c_1,c_2,c_3}$. To see this, by applying Lemma~\ref{lem:g'infinity}, we obtain that, for any $f_{g,\beta,\gamma}\in\tilde{\mathcal W}_{c_1,c_2,c_3}$, there exists $\tilde f_{g_i,\beta_j}\in \tilde{\mathcal W}$, for some $i$, $j$, defined by $\tilde f_{g_i,\beta_j}(x,z)=S_{x,z}(g_i,\beta_j)$, such that by equation \eqref{san},
\begin{align*}
&\sup_{x\in\X,z\in\Z}\|f_{g,\beta}(x,z)-\tilde f_{g_i,\beta_j}(x,z)\|\leq c\|\beta_1-\beta_2\|_2^{1/2}+c\|g'_1-g_2'\|_\infty\leq c\eta ,
\end{align*}
which implies that $\log N_{}\big(\eta, \mathcal{W}_{c_1, c_2,c_3,c_4}(n),\|\cdot\|_{\infty}\big)\leq \log(r_1r_2) \lesssim \eta^{-1/(m-1)}$  and proves \eqref{san2}.

Now, note that $\E\{\e f(X,Z)\}=0$ for $f\in\mathcal W_{c_1,c_2,c_3,c_4}$, so that
\begin{align*}
\|\mathbb V_{1,n}\|&\leq\sup_{(g,\beta)\in\mathcal E_{c_1,c_2,c_3,c_4}(n)}\bigg\|\frac{1}{n}\sum_{i=1}^n\e_i\{S_{X_i,Z_i}(g_0,\beta_0)-S_{X_i,Z_i}(g,\beta)\}\bigg\|
\\
&=n^{-1/2}\sup_{f\in{\mathcal W}_{c_1,c_2,c_3,c_4}(n)}\bigg\|\frac{1}{\sqrt n}\sum_{i=1}^n\e_if(X_i,Z_i)\bigg\|.
\end{align*}
By Markov's inequality and the maximal inequality in Corollary~2.2.5 in \cite{vaart1996}, it holds that, for any constants $\delta>0$, $e>0$, and $n$ large enough
\begin{align*}
\P\big(\|\mathbb V_{1,n}\|\geq\delta n^{-1/2}\big)&\leq\P\big\{\|\mathbb V_{1,n}\|\geq\delta n^{-1/2};(\hat g_n,\hat\beta_n)\in\mathcal E_{c_1,c_2,c_3,c_4}(n)\big\}+\P\big\{(\hat g_n,\hat\beta_n)\notin\mathcal E_{c_1,c_2,c_3,c_4}(n)\big\}\\
&\leq\P\left(\sup_{f\in{\mathcal W}_{c_1,c_2,c_3,c_4}(n)}\bigg\|\frac{1}{\sqrt n}\sum_{i=1}^n\e_if(X_i,Z_i)\bigg\|>\delta\right)+e/2
\\
&\leq\delta^{-1}\E\left(\sup_{f\in{\mathcal W}_{c_1,c_2,c_3,c_4}(n)}\bigg\|\frac{1}{\sqrt n}\sum_{i=1}^n\e_if(X_i,Z_i)\bigg\|\right)+e/2\\
&\leq \delta^{-1}\int_0^{m_n}\sqrt{\log N_{}(\eta,{\mathcal W}_{c_1,c_2,c_3,c_4}(n),\|\cdot\|)}\,\d\eta+e/2
\\
&\leq c \delta^{-1}(\lambda^{-1/(4m-4)}r_n+r_n^{1/2})^{(2m-3)/(2m-2)}+e/2<e.
\end{align*}
The above equation implies that
\begin{align}\label{v2n}
\|\mathbb V_{1,n}\|=o_p(n^{-1/2}).
\end{align}

For the second term $\mathbb V_{2,n}$ in \eqref{diff}, note that by Lemma~\ref{lem:unionbound2} and the fact that $\|\hat w_n-w_0\|=O_p(r_n)$ in view of Theorem~\ref{thm:rate}, it holds that
\begin{align*}
\sup_{x\in\X}|r_x(\hat g_n,\hat\beta_n)|\lesssim\|\hat\beta_n-\beta_0\|_2^{3/2}+\|\hat g_n'-g_0'\|_\infty\times\|\hat\beta_n-\beta_0\|_2.
\end{align*}
By Lemma~\ref{lem:qbeta} and Proposition~\ref{prop:2.1}, we obtain
\begin{align*}
\|\hat\beta_n-\beta_0\|_2\leq\|Q_{\beta_0}\trans\hat\beta_n\|_2+\|\beta_0\beta_0\trans(\hat\beta_n-\beta_0)\|_2\leq 2\|Q_{\beta_0}\trans\hat\beta_n\|_2\leq c\|\hat w_n-w_0\|.
\end{align*}
In addition, in view of Lemmas~\ref{lem:normK} and \ref{lem:norms}, we deduce that
\begin{align*}
\|\hat g_n-g_0\|_\infty\leq c\lambda^{-1/(4m-4)}\|\hat g_n-g_0\|_K\leq c\lambda^{-1/(4m-4)}\|\hat w_n-w_0\|.
\end{align*}
Combining the above three equations yields
\begin{align*}
\sup_{x\in\X}|r_x(\hat g_n,\hat\beta_n)|\leq c\|\hat w_n-w_0\|^{3/2}+c\lambda^{-1/(4m-4)}\|\hat w_n-w_0\|^2=O_p(r_n^{3/2}+\lambda^{-1/(4m-4)}r_n^{2}).
\end{align*}
Recall from Proposition~\ref{cor:sxz} and Corollary~\ref{n} that
\begin{align*}
\|S_{x,z}(\hat g_n,\hat\beta_n)\|=\|S^\dagger_{x\trans\hat\beta_n,\hat g_n'(x\trans\hat\beta_n)Q_{\beta_0}\trans x,z}\|\leq c\lambda^{-1/(4m)}\big\{1+\|\hat g_n'(x\trans\hat\beta_n)Q_{\beta_0}\trans x\|_2^2+\|z\|_2^2\big\}.
\end{align*}
Note that in view of Theorem~\ref{prop:rate} it holds that $\|\hat g_n'\|_\infty=O_p(1)$, so that the above equation together with the fact that $\X,\Z$ are compact implies
\begin{align}\label{new}
\sup_{x\in\X,z\in\Z}\|S_{x,z}(\hat g_n,\hat\beta_n)\|=O_p(\lambda^{-1/(4m)}).
\end{align}
Combining \eqref{new} with Corollary~\ref{n} and Theorem~\ref{thm:rate} yields
\begin{align}\label{v3n}
\|\mathbb V_{2,n}&\|\leq\sup_{x\in\X}|r_x(\hat g_n,\hat\beta_n)|\times\sup_{x\in\X,z\in\Z}\|S_{x,z}(\hat g_n,\hat\beta_n)\|=O_p\{\lambda^{-1/(4m)}(r_n^{3/2}+\lambda^{-1/(4m-4)}r_n^{2})\}.
\end{align}

For the third term $\mathbb V_{3,n}$ in \eqref{diff}, we obtain from Lemmas~\ref{lem:unionbound2} and Theorem~\ref{thm:rate} that
\begin{align*}
\sup_{x\in\X,z\in\Z}\|S_{x,z}(g_0,\beta_0)-S_{x,z}(\hat g_n,\hat\beta_n)\|&\lesssim\lambda^{-1/(4m-4)}\|\hat w_n-w_0\|+\|\hat w_n-w_0\|^{1/2}
\\
&=O_p(\lambda^{-1/(4m-4)}r_n+r_n^{1/2}).
\end{align*}
In view of \eqref{new} and Theorem~\ref{thm:rate}, we deduce that
\begin{align}\label{v4n}
\|\mathbb V_{3,n}\|&\leq\|\hat w_n-w_0\|\times\sup_{x\in\X,z\in\Z}\|S_{x,z}(g_0,\beta_0)\|\times\sup_{x\in\X,z\in\Z}\|S_{x,z}(g_0,\beta_0)-S_{x,z}(\hat g_n,\hat\beta_n)\|\notag\\
&=O_p\{\lambda^{-1/(4m)}r_n(\lambda^{-1/(4m-4)}r_n+r_n^{1/2})\}.
\end{align}

Finally, we bound the fourth term $n^{-1/2}\|\mathbb M_n(\hat w_n-w_0)\|$. We apply the following lemma, which a modified version of Lemma~S.3 in \cite{cheng2015}.
\begin{lemma}\label{lem:s3}
For the constant $c_0>1$ in Lemma~\ref{lem:norms}, define
\begin{align*}
\mathcal G=\big\{f=f_1+f_2:f_1(x)=x\trans\theta,f_2\in\H_m,\|f_1\|_\infty\leq1;\|f_2\|_\infty\leq1;\lambda J(f_2,f_2)\leq c_0\lambda^{1/(2m)};x,\theta\in\mathbb R^{p+q}\big\}.
\end{align*}
Suppose that $\psi$ satisfies
\begin{align*}
|\psi_{X,Z}(f)-\psi_{X,Z}( \tilde{f})| \leq c \lambda^{1/(4m)}\|f-\tilde{f}\|_{\infty}, \qquad f, \tilde{f} \in \mathcal{G},
\end{align*}
almost surely, for some constant $c>0$. Define
\begin{align*}
\mathbb Z_n(f)=\frac{1}{\sqrt{n}} \sum_{i=1}^n\big[\psi_{X_i,Z_i}(f) S_{X_i,Z_i}(g_0,\beta_0)-\E\{\psi_{X,Z}(f) S_{X,Z}(g_0,\beta_0)\}\big],\qquad f\in\mathcal G.
\end{align*}
Then, it holds that
\begin{align*}
\lim_{n\to\infty}\P\left(\sup_{f\in\mathcal G}\frac{\|\mathbb Z_n(f)\|}{\lambda^{-(2m-1)/(8m^2)}\|f\|_\infty^{1-1/(2m)}+n^{-1/2}}\leq \{5\log\log(n)\}^{1/2}\right)=1.
\end{align*}

\end{lemma}

To apply Lemma~\ref{lem:s3}, for the constant $c_0>1$ in Lemma~\ref{lem:norms}, define
\begin{align*}
&\mathcal C_{c}(n)=\big\{(g,\beta,\gamma)\in\mathcal H_m\times\mathcal B\times\mathbb R^q :\|(g,Q_{\beta_0}\trans\beta,\gamma)-(g_0,0,\gamma_0)\|\leq cr_n\big\},
\\
&\mathcal G_c(n)=\Big\{f_w:f_{w}(s,x,z)=c_0^{-1}\lambda^{1/(4m)}\|w-w_0\|^{-1}\big\{g(s)-g_0(s)+x\trans Q_{\beta_0} Q_{\beta_0}\trans(\beta-\beta_0)+z\trans(\gamma-\gamma_0)\big\};\\
&\qquad\qquad\qquad \text{for } w=(g,Q_{\beta_0}\trans\beta,\gamma)\neq w_0; \ f_w\equiv0 \text{ for }w=w_0;\  (s,x,z)\in\I\times\X\times\Z,(g,\beta,\gamma)\in\mathcal C_c(n)\Big\},
\end{align*}
and for the $S_{x,z}(g,\beta)$ in Proposition~\ref{cor:sxz}, and for $f_w\in\mathcal G_n$, define
\begin{align*}
&\psi_{X,Z}(f_w)=\lambda^{1/(4m)}f_w\{X\trans\beta_0,g_0'(X\trans\beta_0) X,Z\},\\
&\mathbb Z_n(f_w)=n^{-1/2}\sum_{i=1}^n\Big(f_w\{X_i\trans\beta_0,g_0'(X_i\trans\beta_0)X_i,Z_i\} S_{X_i,Z_i}(g_0,\beta_0)-\E\big[f_w\{X\trans\beta_0,g_0'(X\trans\beta_0)X,Z\} S_{X,Z}(g_0,\beta_0)\big]\Big).
\end{align*}
Note that by definition, $\mathbb M_n(\hat w_n-w_0)=c_0\lambda^{-1/(4m)}\|\hat w_n-w_0\|\mathbb Z_n(f_{\hat w_n})$. By Theorem~\ref{thm:rate}, for $c'>0$ large enough, $\P\{(\hat g_n,\hat\beta_n,\hat\gamma_n)\notin\mathcal C_{c'}(n)\}=o(1)$ as $n\to\infty$.
It follows from Proposition~\ref{prop:2.1} and Lemma~\ref{lem:norms} that, for $n$ large enough,
\begin{align*}
&c_0^{-1}\lambda^{1/(4m)}\|w-w_0\|^{-1}\|g-g_0\|_\infty\leq 1\\
&c_0^{-1}\lambda^{1/(4m)}\|w-w_0\|^{-1}\sup_{x\in\X,z\in\Z}|x\trans Q_{\beta_0} Q_{\beta_0}\trans(\beta-\beta_0)+z\trans(\gamma-\gamma_0)|\leq c\lambda^{1/(4m)}\\
&c_0^{-2}\lambda^{1+1/(2m)}\|w-w_0\|^{-2} J(g-g_0,g-g_0)\leq c_0^{-2}\lambda^{1/(2m)}\leq c_0\lambda^{1/(2m)}.
\end{align*}
In addition, it holds that
\begin{align*}
|\psi_{X,Z}(f_{w})-\psi_{X,Z}(\tilde f_{w})|\leq \lambda^{1/(4m)}\|f_{w}-\tilde f_{w}\|_{\infty}, \qquad a.s.
\end{align*}
The above calculations show that $\mathcal G_n$ and $\psi$ satisfy the conditions in Lemma~\ref{lem:s3}. Furthermore, it follows from Corollary~\ref{n} that, for any $f_w\in\mathcal G_{c'}(n)$,
\begin{align*}
\|f_w\|_\infty=\lambda^{1/(4m)}\|w-w_0\|^{-1}\sup_{x\in\X,z\in\Z}\l S^\dagger_{s,Q_{\beta_0}\trans x,z},w-w_0\r\leq \lambda^{1/(4m)}\sup_{x\in\X,z\in\Z}\|S_{s,u,v}^\dagger\|\leq c,
\end{align*}
where $S^\dagger_{s,u,v}$ is defined in Proposition~\ref{prop:ru}. Therefore, we apply Lemma~\ref{lem:s3} and deduce from the above equation that
\begin{align*}
&\P\left(\|\mathbb Z_n(f_{\hat w_n})\|\leq c(\lambda^{-(2m-1)/(8m^2)}+n^{-1/2})\sqrt{\log\log(n)}\right)
\\
&\leq \P\{(\hat g_n,\hat\beta_n,\hat\gamma_n)\notin\mathcal C_{c'}(n)\}+\P\left(\|\mathbb Z_n(f_{\hat w_n})\|\leq c(\lambda^{-(2m-1)/(8m^2)}+n^{-1/2})\sqrt{\log\log(n)};(\hat g_n,\hat\beta_n,\hat\gamma_n)\in\mathcal C_{c'}(n)\right)
\\
&\leq\P\{(\hat g_n,\hat\beta_n,\hat\gamma_n)\notin\mathcal C_{c'}(n)\}+\P\left(\frac{\|\mathbb Z_n(f_{\hat w_n})\|}{\lambda^{-(2m-1)/(8m^2)}\|f_{\hat w_n}\|_\infty^{1-1/(2m)}+n^{-1/2}}\leq c\sqrt{\log\log(n)}\right)
\\
&\leq \P\{(\hat g_n,\hat\beta_n,\hat\gamma_n)\notin\mathcal C_{c'}(n)\}+\P\left(\sup_{f_w\in\mathcal G_{c'}(n)}\frac{\|\mathbb Z_n(f_{w})\|}{\lambda^{-(2m-1)/(8m^2)}\|f_{w}\|_\infty^{1-1/(2m)}+n^{-1/2}}\leq c\sqrt{\log\log(n)}\right)
\\
&=o(1).
\end{align*}
The above equation implies that
\begin{align*}
\|\mathbb Z_n(f_{\hat w_n})\|=O_p\big(\lambda^{-(2m-1)/(8m^2)}\sqrt{\log\log(n)}\big).
\end{align*}
As a consequence, recalling the definition of $\mathbb M_n$ in Lemma~\ref{lem:n-or}, we obtain from the above equation and Theorem~\ref{thm:rate} that, for the $a_n$ defined in Assumption~\ref{a:reg},
\begin{align*}
&n^{-1/2}\|\mathbb M_n(\hat w_n-w_0)\|=n^{-1/2}\lambda^{-1/(4m)}\|\hat w_n-w_0\|\times\|\mathbb Z_n(f_{\hat w_n})\|\\
&\leq O_p(n^{-1/2}\lambda^{-1/(4m)}r_n)\times O_p\big(\lambda^{-(2m-1)/(8m^2)}\sqrt{\log\log(n)}\big)
=O_p(n^{-1/2}a_n)=o_p(n^{-1/2}).
\end{align*}

The proof of Theorem~\ref{thm:bahadur} is complete by combining the above equation with  \eqref{v2n}, \eqref{v3n} and \eqref{v4n}.

\subsection{Proof of joint weak convergence in Theorem~\ref{thm:asymp}}\label{app:thm:asymp}

The proof of Theorem~\ref{thm:asymp} is organized as follows. First, in Section~\ref{app:1}, we state a preparatory theorem, namely Theorem~\ref{thm:asymp:prep}. Note that in Theorem~\ref{thm:asymp:prep}, we do not require Condition~\eqref{newcond}. Section, we derive a useful lemma regarding $A_X,A_Z$ and $R_X,R_Z$ in Section~\ref{app:lem:cross}. The proof of Theorem~\ref{thm:asymp} is given in Section~\ref{app:thm:asymp2}. Finally, in Section~\ref{app:psudo}, we verify that $\Omega^+$ defined in \eqref{Omegastar} is the Moore–Penrose inverse of $\Omega$ in \eqref{Omega}.

\subsubsection{A preparatory theorem (Theorem~\ref{thm:asymp:prep}) and its proof}\label{app:1}

\begin{theorem}\label{thm:asymp:prep}
Assume Assumptions~\ref{a:rank}--\ref{a:reg} are valid. For $s\in\I$, suppose the asymptotic variance $\sigma_{(s)}^2$ is defined in \eqref{sigmas2}. In addition, suppose
\begin{align}\label{cd}
c(s)=\lim_{\lambda\to0}\lambda^{1/(4m)}\Bigg[\begin{matrix}
M_\lambda A_X(s)\\
M_\lambda A_Z(s)
\end{matrix}\Bigg],\qquad d(s)
=\lim_{\lambda\to0}\lambda^{1/(4m)}\Bigg[\begin{matrix}
A_X(s)\\
A_Z(s)
\end{matrix}\Bigg],\qquad s\in\I.
\end{align}
Then, it holds that, for the $\big(H^\circ_{\lambda}(g),N^\circ_{\lambda}(g),T^\circ_{\lambda}(g)\big)$ defined in Proposition~\ref{cor:sxz},
\begin{align}\label{main}
&\left[\begin{matrix}
\sqrt{n}\lambda^{1/(4m)}\{\hat g_n(s)-g_0(s)+H_\lambda^\circ g_0(s)\}\\
\sqrt{n}(\hat\beta_n-\beta_0+N_\lambda^\circ g_0)\\
\sqrt{n}(\hat\gamma_n-\gamma_0+T_\lambda^\circ g_0)\\
\end{matrix}\right]\notag
\\[0.5cm]
&\hspace{1cm}\converged N\left(\Bigg[\begin{matrix}
0 \\ 0_{p+q}
\end{matrix}\Bigg],\Bigg[\begin{matrix}
\sigma_{(s)}^2-2d(s)\trans\Omega^+ c(s)+d(s)\trans\Omega^+d(s) & \{c(s)-d(s)\}\trans\Omega^+\\
\Omega^+\{c(s)-d(s)\} & \Omega^+
\end{matrix}\Bigg]\right)\,,
\end{align}
where $M_\lambda$ is defined in \eqref{mlambda} and $\Omega^+$ is the Moore–Penrose inverse of $\Omega$ in \eqref{Omega} defined in \eqref{Omegastar}.

\end{theorem}

\begin{proof}[\underline{Proof of Theorem~\ref{thm:asymp:prep}}]

We first recall and introduce some notations.
In the sequel, for notational convenience, in view of $S_{x,z}$ defined in Proposition~\ref{cor:sxz} in Proposition~\ref{cor:sxz}, we write $U_i=(X_i,Z_i)$ and $S_{U_i}=S_{X_i,Z_i}(g_0,\beta_0)$ for brevity, that is, $S_{U_i}=(H_{U_i},N_{U_i},T_{U_i})$, where
\begin{align}\label{stemp}
\Biggl[\begin{matrix}
N_{U_i}\\
T_{U_i}
\end{matrix}\Biggl]&=(\Omega_{\beta_0}+\Sigma_{\beta_0,\lambda})^{-1}\Biggl[\begin{matrix}
Q_{\beta_0}\trans  \{g_0'(X_i\trans\beta_0)X_i-A_X(X_i\trans\beta_0)\}\\
Z_i-A_Z(X_i\trans\beta_0)
\end{matrix}\Biggl]\,,\notag\\
H_{U_i}&=K_{X_i\trans\beta_0}-A_X\trans Q_{\beta_0} N_{U_i}-A_Z\trans T_{U_i}.
\end{align}
Suppose $(s,\zeta,v)\in\I\times\mathbb R^p\times\mathbb R^q$. In addition, let $\xi_0:=(s,Q_{\beta_0}\trans \zeta,v)\in\I\times\mathbb R^{p-1}\times\mathbb R^q$ and let $S^\dagger_{\xi_0}$ be defined as in Proposition~\ref{prop:ru} by $S_{\xi_0}^\dagger=(H_{\xi_0}^\dagger,N_{\xi_0}^\dagger,T_{\xi_0}^\dagger)\in\H_m\times\mathbb R^{p-1}\times\mathbb R^q$, where
\begin{align*}
\Biggl[\begin{matrix}
N_{\xi_0}^\dagger\\
T_{\xi_0}^\dagger
\end{matrix}\Biggl]&=(\Omega_{\beta_0}+\Sigma_{\beta_0,\lambda})^{-1}\Biggl[\begin{matrix}
Q_{\beta_0}\trans \{\zeta-A_X(s)\}\\
z-A_Z(s)
\end{matrix}\Biggl],\\
 H_{\xi_0}^\dagger&=aK_{s}-A_X\trans Q_{\beta_0} N_{\xi_0}^\dagger-A_Z\trans T_{\xi_0}^\dagger\,,
\end{align*}
where $A_X$ and $A_Z$ are defined in \eqref{abx}, so that
\begin{align}\label{sxi0}
\l S_{\xi_0}^\dagger,(g,\theta,\gamma)\r=ag(s)+\zeta\trans Q_{\beta_0}\theta+v\trans\gamma,\qquad (g,\theta,\gamma)\in\Theta.
\end{align} 

To apply the Cramér-Wold theorem, for arbitrary (fixed) $a\in\mathbb R$, $\zeta\in\mathbb R^{p}$ and $v\in\mathbb R^q$, we aim to show the weak convergence of
\begin{align}\label{Jn}
\hspace{-0.5em}\mathbb J_n&:=\sqrt{n}\Big[\lambda^{1/(4m)}a\{\hat g_n(s)-g_0(s)+H_\lambda^\circ g_0(s)\}+\zeta\trans(\hat\beta_n-\beta_0+N_\lambda^\circ g_0)+v\trans(\hat\gamma_n-\gamma_0+T_\lambda^\circ g_0)\Big].
\end{align}
The following decomposition follows directly from the definition of $\mathbb J_n$ in \eqref{Jn}:
\begin{align*}
&\mathbb J_n
=I_1+I_2+I_3+I_4,
\end{align*}
where, recalling the definition of $S_\lambda^\circ(g_0)=\big(H^\circ_{\lambda}(g_0),N^\circ_{\lambda}(g_0),T^\circ_{\lambda}(g_0)\big)$ in Proposition~\ref{cor:sxz},
\begin{align}\label{i123}
I_1&=\sqrt n\lambda^{1/(4m)}\Bigg[a\bigg\{\hat g_n(s)-g_0(s)-\frac{1}{n}\sum_{i=1}^n\e_i H_{U_i}(s)+H_\lambda^\circ(g_0)\bigg\}\notag\\
&\hspace{0.5cm}+\zeta\trans Q_{\beta_0}\bigg\{Q_{\beta_0}\trans(\hat\beta_n-\beta_0)-\frac{1}{n}\sum_{i=1}^n\e_i N_{U_i}+N_\lambda^\circ (g_0)\bigg\}+v\trans\bigg\{\hat\gamma_n-\gamma_0-\frac{1}{n}\sum_{i=1}^n\e_iT_{U_i}+T_\lambda^\circ (g_0)\bigg\}\Bigg];\notag
\\
I_{2}&=\frac{1}{\sqrt n}\sum_{i=1}^n\e_i\big\{a\lambda^{1/(4m)}H_{U_i}(s)+\zeta\trans Q_{\beta_0}N_{U_i}+v\trans T_{U_i}\big\};\notag
\\
I_3&=\sqrt n(1-\lambda^{1/(4m)})\bigg[\zeta\trans Q_{\beta_0}\bigg\{Q_{\beta_0}\trans(\hat\beta_n-\beta_0)-\frac{1}{n}\sum_{i=1}^n\e_i N_{U_i}+N_\lambda^\circ (g_0)\bigg\}\notag\\
&\hspace{3.5cm}+v\trans\bigg\{\hat\gamma_n-\gamma_0-\frac{1}{n}\sum_{i=1}^n\e_iT_{U_i}+T_\lambda^\circ (g_0)\bigg\}\bigg];\notag\\
I_4&=\sqrt{n}\zeta\trans\beta_0\beta_0\trans(\hat\beta_n-\beta_0).
\end{align}
Next, we investigate the asymptotic properties for all the terms $I_1$--$I_4$ in \eqref{i123}, where we will show that the second term $I_2$ is the leading term.

For the first term $I_1$, using the representation in equation~\eqref{sxi0}, we have
\begin{align*}
I_1=\sqrt n\lambda^{1/(4m)}\bigg\l S_{\xi_0}^\dagger,\hat w_n-w_0-\frac{1}{n}\sum_{i=1}^n\e_iS_{U_i}+S_\lambda^\circ(g_0)\bigg\r.
\end{align*}
By the joint Bahadur representation in Theorem~\ref{thm:bahadur} and Corollary~\ref{n}, we deduce that
\begin{align}\label{i1112}
|I_{1}|&=\sqrt n\lambda^{1/(4m)}\bigg|\bigg\l S_{\xi_0}^\dagger,\hat w_n-w_0-\frac{1}{n}\sum_{i=1}^n\e_iS_{U_i}-S_\lambda^\circ(g_0)\bigg\r\bigg|\notag\\
&\leq \sqrt n\lambda^{1/(4m)}\|S_{\xi_0}^\dagger\|\times\bigg\|\hat w_n-w_0-\frac{1}{n}\sum_{i=1}^n\e_iS_{U_i}-S_\lambda^\circ(g_0)\bigg\|=O_p(\sqrt na_n)=o_p(1).
\end{align}

For the second term $I_{2}$ in \eqref{i123}, for $1\leq i\leq n$, define
\begin{align}\label{etai}
\eta_i&:=\e_i\big\{a\lambda^{1/(4m)}H_{U_i}(s)+\zeta\trans Q_{\beta_0}N_{U_i}+v\trans T_{U_i}\big\}\notag
\\ &=\e_i\big[a\lambda^{1/(4m)}\{K(X_i\trans\beta_0,s)-A_X(s)\trans Q_{\beta_0}N_{U_i}-A_Z(s)\trans T_{U_i}\}+\zeta\trans Q_{\beta_0} N_{U_i}+v\trans T_{U_i}\big].
\end{align}
In view of \eqref{stemp},
\begin{align}\label{etai2}
&\eta_i= \e_i\left\{a\lambda^{1/(4m)}K(X_i\trans\beta_0,s)+\Biggl[\begin{matrix}
Q_{\beta_0}\trans \{\zeta-a\lambda^{1/(4m)}A_X(s)\}\\
v-a\lambda^{1/(4m)}A_Z(s)
\end{matrix}\Biggl]\trans(\Omega_{\beta_0}+\Sigma_{\beta_0,\lambda})^{-1}\Biggl[\begin{matrix}
N_{U_i}\\
T_{U_i}
\end{matrix}\Biggl]\right\}\notag
\\
&=\e_i\left\{a\lambda^{1/(4m)}K_{s}(X_i\trans\beta_0)+\Biggl[\begin{matrix}
Q_{\beta_0}\trans \{\zeta-a\lambda^{1/(4m)}A_X(s)\}\\
v-a\lambda^{1/(4m)}A_Z(s)
\end{matrix}\Biggl]\trans(\Omega_{\beta_0}+\Sigma_{\beta_0,\lambda})^{-1}\Biggl[\begin{matrix}
Q_{\beta_0}\trans \{g_0'(X_i\trans\beta_0)X_i-A_X(X_i\trans\beta_0)\}\\
Z_i-A_Z(X_i\trans\beta_0)
\end{matrix}\Biggl]\right\}.
\end{align}
Observe that $\E(\eta_i)=0$ and recall the definition of $A_X$ and $A_Z$ in \eqref{abx}. Let $\sigma_\eta^2=\var(\eta_i)$.
We obtain from the above equation that
\begin{align}\label{sigmaeta}
\sigma_\eta^2&=a^2\lambda^{1/(2m)}\E\{\sigma_0^2(X\trans\beta_0)K_{s}^2(X\trans\beta_0)\}\notag\\
&\qquad+2a\lambda^{1/(4m)}\Biggl[\begin{matrix}
Q_{\beta_0}\trans \{\zeta-a\lambda^{1/(4m)}A_X(s)\}\\
v-a\lambda^{1/(4m)}A_Z(s)
\end{matrix}\Biggl]\trans(\Omega_{\beta_0}+\Sigma_{\beta_0,\lambda})^{-1}\notag\\
&\qquad\times\E\left\{\sigma^2(X,Z)K_{s}(X\trans\beta_0)\Biggl[\begin{matrix}
Q_{\beta_0}\trans \{g_0'(X\trans\beta_0)X-A_X(X\trans\beta_0)\}\\
Z-A_Z(X\trans\beta_0)
\end{matrix}\Biggl]\right\}\notag\\
&\qquad+\Biggl[\begin{matrix}
Q_{\beta_0}\trans \{\zeta-a\lambda^{1/(4m)}A_X(s)\}\\
v-a\lambda^{1/(4m)}A_Z(s)
\end{matrix}\Biggl]\trans(\Omega_{\beta_0}+\Sigma_{\beta_0,\lambda})^{-1}\Biggl[\begin{matrix}
Q_{\beta_0}\trans \{\zeta-a\lambda^{1/(4m)}A_X(s)\}\\
v-a\lambda^{1/(4m)}A_Z(s)
\end{matrix}\Biggl]\notag\\
&=a^2\lambda^{1/(2m)}{V}(K_s,K_s)\notag\\
&\qquad+2a\lambda^{1/(4m)}\Biggl[\begin{matrix}
Q_{\beta_0}\trans\{ \zeta-a\lambda^{1/(4m)}A_X(s)\}\\
v-a\lambda^{1/(4m)}A_Z(s)
\end{matrix}\Biggl]\trans(\Omega_{\beta_0}+\Sigma_{\beta_0,\lambda})^{-1}\Biggl[\begin{matrix}
Q_{\beta_0}\trans\{A_X(s)-{V}(A_X,K_s)\}\\
A_Z(s)-{V}(A_Z,K_s)
\end{matrix}\Biggl]\notag\\
&\qquad+\Biggl[\begin{matrix}
Q_{\beta_0}\trans\{ \zeta-a\lambda^{1/(4m)}A_X(s)\}\\
v-a\lambda^{1/(4m)}A_Z(s)
\end{matrix}\Biggl]\trans(\Omega_{\beta_0}+\Sigma_{\beta_0,\lambda})^{-1}\Biggl[\begin{matrix}
Q_{\beta_0}\trans\{ \zeta-a\lambda^{1/(4m)}A_X(s)\}\\
v-a\lambda^{1/(4m)}A_Z(s)
\end{matrix}\Biggl].
\end{align}

Observe equation \eqref{to} and the fact that
\begin{align*}
A_X(s)-V(K_s,A_X)=(M_\lambda A_X)(s),
\\
A_Z(s)-V(K_s,A_Z)=(M_\lambda A_Z)(s).
\end{align*}
If we denote
\begin{align*}
&c_{X}(s)=\lim_{\lambda\to0}\lambda^{1/(4m)}(M_\lambda A_X)(s),\qquad c_{Z}(s)=\lim_{\lambda\to0}\lambda^{1/(4m)}(M_\lambda A_Z)(s),\\
&d_{X}(s)=\lim_{\lambda\to0}\lambda^{1/(4m)} A_X(s),\qquad ~~~~~~d_{Z}(s)=\lim_{\lambda\to0}\lambda^{1/(4m)}A_Z(s),
\end{align*}
then in view of the assumption in \eqref{cd}, we have
\begin{align*}
c(s)=\Bigg[\begin{matrix}
c_X(s)\\c_Z(s)
\end{matrix}\Bigg]=\lim_{\lambda\to0}\lambda^{1/(4m)}\Bigg[\begin{matrix}
M_\lambda A_X(s)\\
M_\lambda A_Z(s)
\end{matrix}\Bigg],\qquad d(s)=\Bigg[\begin{matrix}
d_X(s)\\d_Z(s)
\end{matrix}\Bigg]
=\lim_{\lambda\to0}\lambda^{1/(4m)}\Bigg[\begin{matrix}
A_X(s)\\
A_Z(s)
\end{matrix}\Bigg].
\end{align*}
Therefore, \eqref{sigmaeta} together with the fact that $\lim_{\lambda\to0}\Sigma_{\beta_0,\lambda}=0$ yields that, as $\lambda\to0$,
\begin{align}\label{uareta}
\sigma_\eta^2&\longrightarrow a^2\sigma_{(s)}^2+2a\Biggl[\begin{matrix}
Q_{\beta_0}\trans \{\zeta-ad_X(s)\}\\
v-ad_Z(s)
\end{matrix}\Biggl]\trans\Omega_{\beta_0}^{-1}\Biggl[\begin{matrix}
Q_{\beta_0}\trans c_{X}(s)\\
c_Z(s)
\end{matrix}\Biggl]+\Biggl[\begin{matrix}
Q_{\beta_0}\trans \{\zeta-ad_X(s)\}\\
v-ad_Z(s)
\end{matrix}\Biggl]\trans\Omega_{\beta_0}^{-1}\Biggl[\begin{matrix}
Q_{\beta_0}\trans \{\zeta-ad_X(s)\}\\
v-ad_Z(s)
\end{matrix}\Biggl]\notag\\
&=\left[\begin{array}{c}
a\\
\zeta\\
v
\end{array}\right]\trans\Bigg[\begin{matrix}
\sigma_{(s)}^2-2d(s)\trans\Omega^+ c(s)+d(s)\trans\Omega^+d(s) & \{c(s)-d(s)\}\trans\Omega^+\\
\Omega^+\{c(s)-d(s)\} & \Omega^+
\end{matrix}\Bigg]\left[\begin{array}{c}
a\\
\zeta\\
v
\end{array}\right],
\end{align}
where $\Omega^+$ is defined in \eqref{Omegastar}.

Next, we verify the Lindeberg's condition of the triangular array of the $\eta_i$'s in \eqref{etai}. In view of \eqref{abx}, by the Cauchy-Schwarz inequality and Assumption~\ref{a:finite}, we obtain
\begin{align}\label{ox}
\sup_{s\in\I}\|A_X(s)\|_2&\leq\sum_{j\geq1}\frac{\|{V}(g_0'\cdot R_{X},\phi_j)\|_2}{1+\lambda\rho_j}\|\phi_j\|_\infty\notag\\
&\leq\sup_{j\geq1}\|\phi_j\|_\infty\times \bigg\{\sum_{j\geq1}\|{V}(g_0'\cdot R_{X},\phi_j)\|_2^2\bigg\}^{1/2}\times\bigg\{\sum_{j\geq1}(1+\lambda\rho_j)^{-2}\bigg\}^{1/2}\notag\\
&\leq c\bigg\{\sum_{j\geq1}(1+\lambda\rho_j)^{-2}\bigg\}^{1/2}\leq c\lambda^{-1/(4m)}.
\end{align}
Similarly, we have $\sup_{s\in\I}\|A_Z(s)\|_2\leq c\lambda^{-1/(4m)}$.
Hence, by Corollary~\ref{n} and Proposition~\ref{prop:ru} and the fact that $\|K_s\|_\infty\leq c\lambda^{-1/(4m)}$, in view of \eqref{etai2}, we obtain that, for $n$ large enough,
\begin{align*}
|\eta_i|&\leq a|\e_i|\times\lambda^{1/(4m)}|K(X_i\trans\beta_0,s)|\\
&\qquad+|\e_i|\times\left|\Bigg[\begin{matrix}
Q_{\beta_0}\trans\{\zeta-a\lambda^{1/(4m)}A_X(s)\}\\ v-a\lambda^{1/(4m)}A_Z(s)
\end{matrix}\Bigg](\Omega_{\beta_0}+\Sigma_{\beta_0,\lambda})^{-1}\Biggl[\begin{matrix}
Q_{\beta_0}\trans  \{g_0'(X_i\trans\beta_0)X_i-A_X(X_i\trans\beta_0)\}\\
Z_i-A_Z(X_i\trans\beta_0)
\end{matrix}\Biggl]\right|
\\
&\leq a\lambda^{-1/(4m)}|\e_i|+c|\e_i|\times\|(\Omega_{\beta_0}+\Sigma_{\beta_0,\lambda})^{-1}\|_2\times \Big\{\|\zeta-a\lambda^{1/(4m)}A_X(s)\|_2+\|v-a\lambda^{1/(4m)}A_Z(s)\|_2\Big\}\\
&\hspace{3cm}\times\Big\{\|g_0'(X_i\trans\beta_0)X_i-A_X(X_i\trans\beta_0)\|_2+\|Z_i-A_Z(X_i\trans\beta_0)\|_2\Big\}
\\
&\leq c\lambda^{-1/(4m)}|\e_i|,
\end{align*}
where we used \eqref{ox} and the fact that $\X$ and $\Z$ are compact. By the Cauchy-Schwarz inequality and Markov's inequality, under the assumption that $n\lambda^{1/m}\to\infty$, for any $e>0$, we have
\begin{align*}
\sigma_\eta^{-2}\mathbb {E} \left[\eta^{2}\cdot \mathbf {1} \{|\eta|>e\sqrt n \sigma_\eta\sigma_\eta\}\right]
&\leq c\sigma_\eta^{-2}\lambda^{-1/(2m)}\mathbb {E} \left[\e^2\cdot \mathbf {1} \{|\e|>c^{-1}e \sqrt n\lambda^{1/(4m)}\}\right]\\
&\leq c\sigma_\eta^{-2}\lambda^{-1/(2m)}\{\mathbb {E}(\e^4)\}^{1/2}\times\big[\P \{|\e|>c^{-1}e \sqrt n\lambda^{1/(4m)}\}\big]^{1/2}\\
&\leq c\sigma_\eta^{-2}\lambda^{-1/(2m)}(c^2e^{-2}n^{-1}\lambda^{-1/(2m)})\E(\e^4)=o(1),\qquad n\to0.
\end{align*}
Therefore, by Lindeberg's CLT, we obtain from the above equation and \eqref{uareta} that
\begin{align}\label{i11}
&I_{2}=\frac{1}{n}\sum_{i=1}^n\eta_i
\\
&\quad\converged N\left(0_{p+q},[a;\zeta\trans;v\trans]\Bigg[\begin{matrix}
\sigma_{(s)}^2-2d(s)\trans\Omega^+ c(s)+d(s)\trans\Omega^+d(s) & \{c(s)-d(s)\}\trans\Omega^+\\
\Omega^+\{c(s)-d(s)\} & \Omega^+
\end{matrix}\Bigg][a;\zeta\trans;v\trans]\trans\right).\notag
\end{align}

For the third term $I_3$ in \eqref{i123}, by the joint Bahadur representation in Theorem~\ref{thm:bahadur} and Proposition~\ref{prop:2.1}, we have
\begin{align}\label{i3}
|I_3|&\leq \sqrt n\bigg|\zeta\trans Q_{\beta_0}\bigg\{Q_{\beta_0}\trans(\hat\beta_n-\beta_0)-\frac{1}{n}\sum_{i=1}^n\e_i N_{U_i}+N_\lambda^\circ (g_0)\bigg\}+v\trans\bigg\{\hat\gamma_n-\gamma_0-\frac{1}{n}\sum_{i=1}^n\e_iT_{U_i}+T_\lambda^\circ (g_0)\bigg\}\bigg|\notag\\
&\leq c\sqrt n\max\{\|\zeta\|_2,\|v\|_2\}\times\bigg\|\hat w_n-w_0-\frac{1}{n}\sum_{i=1}^n\e_iS_{U_i}-S_\lambda^\circ(g_0)\bigg\|=O_p(\sqrt na_n)=o_p(1).
\end{align}
For the fourth term $I_4$ in \eqref{i123}, by Lemma~\ref{lem:qbeta} and Theorem~\ref{thm:rate}, we have
\begin{align*}
|I_4|\leq c \sqrt{n}|\beta_0\trans\hat\beta_n-1|\leq c\sqrt{n}\|Q_{\beta_0}\trans\hat\beta_n\|^2\leq c\sqrt{n}\|\hat w_n-w_0\|^2=O_p(\sqrt n r_n^2)=o_p(1).
\end{align*}

In conclusion, combining the above equation with \eqref{i11}--\eqref{i3} yields
\begin{align*}
\mathbb J_n&=\sqrt{n}\Big[\lambda^{1/(4m)}a\{\hat g_n(s)-g_0(s)+H_\lambda^\circ g_0(s)\}+\zeta\trans(\hat\beta_n-\beta_0+N_\lambda^\circ g_0)+v\trans(\hat\gamma_n-\gamma_0+T_\lambda^\circ g_0)\Big]
\\[0.2cm]
&\converged N\left(0_{p+q},[a;\zeta\trans;v\trans]\Bigg[\begin{matrix}
\sigma_{(s)}^2-2d(s)\trans\Omega^+ c(s)+d(s)\trans\Omega^+d(s) & \{c(s)-d(s)\}\trans\Omega^+\\
\Omega^+\{c(s)-d(s)\} & \Omega^+
\end{matrix}\Bigg][a;\zeta\trans;v\trans]\trans\right).
\end{align*}
The proof of Theorem~\ref{thm:asymp:prep} is therefore complete by applying the Cramér-Wold theorem.

\end{proof}

\subsubsection{A lemma (Lemma~\ref{lem:cross}) under Condition~\eqref{newcond} and its proof}\label{app:lem:cross}

\begin{lemma}\label{lem:cross}
Assume Condition~\eqref{newcond} is valid. First, it holds that $$\sup_{s\in\I}\{\|A_X(s)\|_2+\|A_Z(s)\|_2\}\leq c.$$ 
Second, we have
\begin{align*}
c(s)=\lim_{\lambda\to0}\lambda^{1/(4m)}\Bigg[\begin{matrix}
M_\lambda A_X(s)\\
M_\lambda A_Z(s)
\end{matrix}\Bigg]=0,\qquad d(s)=\lim_{\lambda\to0}\lambda^{1/(4m)}\Bigg[\begin{matrix}
A_X(s)\\
A_Z(s)
\end{matrix}\Bigg]=0.
\end{align*}
Third, it is true that $\|{V}(g_0'\cdot R_X,M_\lambda g_0)\|_2^2+\|{V}(R_Z,M_\lambda g_0)\|_2^2=o(n^{-1})$ as $n\to\infty$.

\end{lemma}

\begin{proof}[\underline{Proof of Lemma~\ref{lem:cross}}]

We shall prove the claims regarding $A_X$ and $R_Z$; the claims regarding $A_Z$ and $R_Z$ follow from the same argument. In view of equations~\eqref{abx} and \eqref{mlambda}, we obtain
\begin{align*}
&\lambda^{1/(4m)}A_X=\lambda^{1/(4m)}\sum_{j\geq1}\frac{V(g_0'\cdot R_{X},\phi_j)}{1+\lambda\rho_j}\phi_j\,;
\\
&\lambda^{1/(4m)}(M_\lambda A_X)=\lambda^{1+1/(4m)}\sum_{j\geq1}V(A_X,\phi_j)\frac{\rho_j\phi_j}{1+\lambda\rho_j}=\lambda^{1+1/(4m)}\sum_{j\geq1}\frac{V(g_0'\cdot R_{X},\phi_j)}{(1+\lambda\rho_j)^2}\rho_j\phi_j.
\end{align*}
First, by the Cauchy-Schwarz inequality and condition~\ref{newcond}, it follows that
\begin{align}\label{ox2}
\sup_{s\in\I}\|A_X(s)\|_2&\leq\sum_{j\geq1}\frac{\|{V}(g_0'\cdot R_{X},\phi_j)\|_2}{1+\lambda\rho_j}\|\phi_j\|_\infty\notag\\
&\leq c\bigg\{\sum_{j\geq1}\|{V}(g_0'\cdot R_{X},\phi_j)\|_2^2(1+\rho_j)^\mu\bigg\}^{1/2}\times\bigg\{\sum_{j\geq1}(1+\lambda\rho_j)^{-2}(1+\rho_j)^{-\mu}\bigg\}^{1/2}\notag\\
&\leq c\bigg\{\sum_{j\geq1}(1+\rho_j)^{-\mu}\bigg\}^{1/2}\leq c.
\end{align}
Second, it holds that 
\begin{align*}
&\lambda^{1/(4m)}\sup_{s\in\I}\|(M_\lambda A_X)(s)\|_2\leq \lambda^{1+1/(4m)}\sum_{j\geq1}\frac{\|{V}(g_0'\cdot R_{X},\phi_j)\|_2}{(1+\lambda\rho_j)^2}\rho_j\|\phi_j\|_\infty\\
&\leq c\lambda^{1+1/(4m)}\bigg[\sum_{j\geq1}\|{V}(g_0'\cdot R_X,\phi_j)\|_2^2\rho_j^{\mu}\bigg]^{1/2}\times\bigg[\sum_{j\geq1}\frac{\rho_j^{2-\mu}}{(1+\lambda\rho_j)^4}\bigg]^{1/2}\leq c\lambda^{\mu/2}=o(1).
\end{align*}
Third, in view of the equation below \eqref{mlambda}, we deduce that
\begin{align*}
{V}(g_0'\cdot R_X,M_\lambda g_0)=\sum_{j\geq1}{V}(g_0'\cdot R_X,\phi_j){V}(M_\lambda g_0,\phi_j)=\lambda\sum_{j\geq1}{V}(g_0'\cdot R_X,\phi_j){V}(g_0,\phi_j)\frac{\rho_j}{1+\lambda\rho_j}.
\end{align*}
Hence, by the Cauchy-Schwarz inequality and the fact that $J(g_0,g_0)=\sum_{j\geq1}|{V}(g_0,\phi_j)|^2\rho_j<\infty$, we obtain
\begin{align}\label{nn2}
\|{V}(g_0'\cdot R_X,M_\lambda g_0)\|_2^2&\leq\lambda^{1+\mu}\sum_{j\geq1}\|{V}(g_0'\cdot R_X,\phi_j)\|_2^2\rho_j^{\mu}\frac{(\lambda\rho_j)^{1-\mu}}{1+\lambda\rho_j}\times\sum_{j\geq1}|{V}(g_0,\phi_j)|^2\frac{\rho_j}{1+\lambda\rho_j}\notag\\
&\leq\lambda^{1+\mu}\sum_{j\geq1}\|{V}(g_0'\cdot R_X,\phi_j)\|_2^2\rho_j^{\mu}\times J(g_0,g_0)\leq c\lambda^{1+\mu}=o(n^{-1}),
\end{align}
where we apply the condition that $\sqrt n \lambda^{(\mu+1)/2}=o(1)$.
This completes the proof of Lemma~\ref{lem:cross}.

\subsubsection{Proof of Theorem~\ref{thm:asymp}}\label{app:thm:asymp2}

Now, we prove Theorem~\ref{thm:asymp} using the preparatory theorem Theorem~\ref{thm:asymp:prep} and Lemma~\ref{lem:cross} on the cross terms. Similar to the proof of Theorem~\ref{thm:asymp:prep}, to apply the Cramér-Wold theorem, for arbitrary (fixed) $a\in\mathbb R$, $\zeta\in\mathbb R^{p}$ and $v\in\mathbb R^q$, we aim to show the weak convergence of
\begin{align*}
\sqrt{n}\Big[a\lambda^{1/(4m)}\{\hat g_n(s)-g_0(s)+M_\lambda g_0(s)\}+\zeta\trans(\hat\beta_n-\beta_0)+v\trans(\hat\gamma_n-\gamma_0)\Big].
\end{align*}
Recalling the definition of $\mathbb J_n$ in \eqref{Jn}, we obtain the following decomposition:
\begin{align*}
\sqrt{n}\Big[a\lambda^{1/(4m)}\{\hat g_n(s)-g_0(s)+M_\lambda g_0(s)\}+\zeta\trans(\hat\beta_n-\beta_0)+v\trans(\hat\gamma_n-\gamma_0)\Big]=a\mathbb J_n+\mathbb I_n,
\end{align*}
where recalling the definition of $S_\lambda^\circ(g_0)=\big(H^\circ_{\lambda}(g_0),N^\circ_{\lambda}(g_0),T^\circ_{\lambda}(g_0)\big)$ in Proposition~\ref{cor:sxz},
\begin{align}\label{i1233}
\mathbb I_n&=\sqrt n\Big[a\lambda^{1/(4m)}\{M_\lambda g_0(s)-H_\lambda^\circ (g_0)\}-\zeta\trans Q_{\beta_0}N_\lambda^\circ (g_0)-v\trans T_\lambda^\circ (g_0)\Big].
\end{align}
First, applying Theorem~\ref{thm:asymp:prep} and Lemma~\ref{lem:cross}, we obtain
\begin{align}\label{nn1}
\mathbb J_n\converged N\left(0_{p+q},[1;\zeta\trans;v\trans]\Bigg[\begin{matrix}
\sigma_{(s)}^2 & 0_{p+q}\trans\\
0_{p+q} & \Omega^+
\end{matrix}\Bigg][1;\zeta\trans;v\trans]\trans\right).
\end{align}
Second, for the term $\mathbb I_n$ in \eqref{i1233}, in view of Proposition~\ref{cor:sxz}, we have
\begin{align*}
\mathbb I_n=-\sqrt n\Bigg[\begin{matrix}
Q_{\beta_0}\trans{V}(g_0'\cdot R_X,M_\lambda g_0)\\
{V}(R_Z,M_\lambda g_0)
\end{matrix}\Bigg]\trans(\Omega_{\beta_0}+\Sigma_{\beta_0,\lambda})^{-1}\Bigg[
\begin{matrix}
Q_{\beta_0}\trans\{\zeta-\lambda^{1/(4m)}A_X(s)\}\\
v-\lambda^{1/(4m)}A_Z(s)
\end{matrix}
\Bigg],
\end{align*}
which in view of Lemma~\ref{lem:cross} implies that for $n$ large enough,
\begin{align*}
|\mathbb I_n|&\leq \sqrt n\|(\Omega_{\beta_0}+\Sigma_{\beta_0,\lambda})^{-1}\|_2\times\Big\{\|{V}(g_0'\cdot R_X,M_\lambda g_0)\|_2+\|{V}(R_Z,M_\lambda g_0)\|_2\Big\}\\
&\hspace{2cm}\times\Big\{\|\zeta-\lambda^{1/(4m)}A_X(s)\|_2+\|v-\lambda^{1/(4m)}A_Z(s)\|_2\Big\}
\\
&\leq c\sqrt n\Big\{\|{V}(g_0'\cdot R_X,M_\lambda g_0)\|_2+\|{V}(R_Z,M_\lambda g_0)\|_2\Big\}\Big[\|\zeta\|_2+\|z\|_2+\lambda^{1/(4m)}\{\|A_X(s)\|_2+\|A_Z(s)\|_2\}\Big]
\\
&\leq c\sqrt n\Big\{\|{V}(g_0'\cdot R_X,M_\lambda g_0)\|_2+\|{V}(R_Z,M_\lambda g_0)\|_2\Big\}=o(1).
\end{align*}
Combining the above equation with \eqref{nn1} yields
\begin{align*}
&\sqrt{n}\Big[a\lambda^{1/(4m)}\{\hat g_n(s)-g_0(s)+M_\lambda g_0(s)\}+\zeta\trans(\hat\beta_n-\beta_0)+v\trans(\hat\gamma_n-\gamma_0)\Big]\\
&\hspace{3cm}\converged N\left(0_{p+q},[1;\zeta\trans;v\trans]\Bigg[\begin{matrix}
a^2\sigma_{(s)}^2 & 0_{p+q}\trans\\
0_{p+q} & \Omega^+
\end{matrix}\Bigg][1;\zeta\trans;v\trans]\trans\right).
\end{align*}
The proof of Theorem~\ref{thm:asymp} is therefore complete by applying the Cramér-Wold theorem.

\end{proof}

\subsubsection{Verification of the pseudo inverse of $\Omega$}\label{app:psudo}

To show that $\Omega^+$ defined in \eqref{Omegastar} is the Morre-Penrose inverse of $\Omega$ in \eqref{Omega}, let $ D=\Big[
\begin{smallmatrix}
Q_{\beta_0}&0_{p,q}\\
0_{q,p-1}&I_q
\end{smallmatrix}\Big]$. We have $\Omega^+= D( D\trans\Omega D)^{-1} D\trans$.
Observing the definition of $\Omega_1$ and $\Omega_3$ in \eqref{Omega123}, we obtain from \eqref{beta0null} that $\Omega_1\beta_0=0_p$ and $\Omega_3\beta_0=0_q$. Therefore, we obtain that $\Omega D D\trans=\Omega$. Direct calculations imply that $\Omega^+$ is the Moore-Penrose inverse of $\Omega$:
\begin{align*}
&\Omega\Omega^+\Omega=\Omega D( D\trans\Omega D)^{-1} D\trans\Omega= D D\trans\Omega D( D\trans\Omega D)^{-1} D\trans\Omega D D\trans= D D\trans\Omega D D\trans=\Omega,\\
&\Omega^+\Omega\Omega^+= D( D\trans\Omega D)^{-1} D\trans\Omega D( D\trans\Omega D)^{-1} D\trans= D( D\trans\Omega D)^{-1} D\trans=\Omega^+,\\
&(\Omega\Omega^+)\trans= D( D\trans\Omega D)^{-1} D\trans\Omega= D( D\trans\Omega D)^{-1} D\trans\Omega D D\trans= D D\trans= D D\trans\Omega D( D\trans\Omega D)^{-1} D\trans=\Omega\Omega^+,\\
&(\Omega^+\Omega)\trans=\Omega D( D\trans\Omega D)^{-1} D\trans= D D\trans\Omega D( D\trans\Omega D)^{-1} D\trans= D D\trans= D( D\trans\Omega D)^{-1} D\trans\Omega D D\trans=\Omega^+\Omega.
\end{align*}
The proof is complete.

\subsection{Proof of Theorem~\ref{thm:marginal} and Proposition~\ref{prop:rootn}}\label{app:marginal}

This section is organized as follows. We first state and prove in Section~\ref{app:new} a preparatory theorem under Assumptions~\ref{a:rank}--\ref{a:reg}. Next, we prove Theorem~\ref{thm:marginal} in Section~\ref{app:11}. The proof of Proposition~\ref{prop:rootn} is given in Section~\ref{app:prop:rootn}.

\subsubsection{A preparatory theorem (Theorem~\ref{thm:marginal2}) and its proof}\label{app:new}

\begin{theorem}\label{thm:marginal2}
Suppose Assumptions~\ref{a:rank}--\ref{a:reg} are valid. Then, it holds that 
\begin{gather}
\sup_{s\in\I}\bigg|\hat g_n(s)-g_0(s)+H^\circ_\lambda  g_0(s)- \frac{1}{n}\sum_{i=1}^n\e_i H_{X_i,Z_i}(g_0,\beta_0)(s)\bigg|=o_p(n^{-1/2}\lambda^{-1/(4m)}),\notag
\\
\left\|\Bigg[\begin{matrix}
\hat\beta_n-\beta_0\\
\hat\gamma_n-\gamma_0
\end{matrix}\Bigg]+\Bigg[\begin{matrix}
Q_{\beta_0}N_\lambda^\circ(g_0)\\
T_\lambda^\circ(g_0)
\end{matrix}\Bigg]-\frac{1}{n}\sum_{i=1}^n\e_i\Bigg[\begin{matrix} Q_{\beta_0}N_{X_i,Z_i}(g_0,\beta_0)\\
T_{X_i,Z_i}(g_0,\beta_0)
\end{matrix}
\Bigg]\right\|_2=o_p(n^{-1/2}),\label{m11}
\end{gather}
where $(H_{x,z},N_{x,z},T_{x,z})$ and $(H_\lambda^\circ,N_\lambda^\circ,T_\lambda^\circ)$ are defined in Proposition~\ref{cor:sxz}, respectively.

\end{theorem}

\begin{proof}[\underline{Proof of Theorem~\ref{thm:marginal2}}]
We first show the marginal Bahadur representation of $\hat g_n(s)$ in \eqref{m11}. It follows from Theorem~\ref{thm:bahadur} and Lemma~\ref{lem:normK} that
\begin{align*}
&\bigg\|\hat g_n-g_0+H^\circ_\lambda  g_0- \frac{1}{n}\sum_{i=1}^n\e_i H_{X_i,Z_i}(g_0,\beta_0)\bigg\|_K
\\
&\leq\bigg\|\hat w_n-w_0-\frac{1}{n}\sum_{i=1}^n\epsilon_iS_{X_i,Z_i}(g_0,\beta_0)+S_\lambda^\circ(g_0)\bigg\|=o_p(n^{-1/2}).
\end{align*}
Observe that by Proposition~\ref{prop:eigen}, it holds that
\begin{align}\label{ksnorm}
\sup_{s\in\I}\|K_s\|_K^2=\sup_{s\in\I}\sum_{j\geq1}\frac{\phi_j^2(s)}{1+\lambda\rho_j}\leq \sum_{j\geq1}\frac{\|\phi_j\|_\infty^2}{1+\lambda\rho_j}\leq c\lambda^{-1/(2m)}.
\end{align}
We therefore deduce from the above two equations and the Cauchy-Schwarz inequality that
\begin{align*}
&\sup_{s\in\I}\bigg|\hat g_n(s)-g_0(s)+H^\circ_\lambda  g_0(s)- \frac{1}{n}\sum_{i=1}^n\e_i H_{X_i,Z_i}(g_0,\beta_0)(s)\bigg|
\\
&=\sup_{s\in\I}\bigg|\bigg\l K_s,\hat g_n-g_0+H^\circ_\lambda  g_0- \frac{1}{n}\sum_{i=1}^n\e_i H_{X_i,Z_i}(g_0,\beta_0) \bigg\r_K\bigg|
\\
&\leq\sup_{s\in\I}\|K_s\|_K^2\times\bigg\|\hat g_n-g_0+H^\circ_\lambda  g_0- \frac{1}{n}\sum_{i=1}^n\e_i H_{X_i,Z_i}(g_0,\beta_0)\bigg\|_K=o_p(n^{-1/2}\lambda^{-1/(4m)}).
\end{align*}
This proves the claim regarding $\hat g_n$ in \eqref{m11}.

Next, we prove the result in \eqref{m11} regarding $[\hat\beta_n\trans,\hat\gamma_n\trans]\trans$. It follows from the joint Bahadur representation in Theorem~\ref{thm:bahadur} and Proposition~\ref{prop:2.1} that
\begin{align*}
\left\|\Bigg[\begin{matrix}
Q_{\beta_0}\trans(\hat\beta_n-\beta_0)\\
\hat\gamma_n-\gamma_0
\end{matrix}\Bigg]+\Bigg[\begin{matrix}
N_\lambda^\circ (g_0)\\
T_\lambda^\circ (g_0)
\end{matrix}\Bigg]-\frac{1}{n}\sum_{i=1}^n\e_i\Bigg[\begin{matrix} N_{X_i,Z_i}(g_0,\beta_0)\\
T_{X_i,Z_i}(g_0,\beta_0)
\end{matrix}
\Bigg]\right\|_2=o_p(n^{-1/2}),
\end{align*}
The above equation implies that
\begin{align*}
\left\|\Bigg[\begin{matrix}
Q_{\beta_0}Q_{\beta_0}\trans(\hat\beta_n-\beta_0)\\
\hat\gamma_n-\gamma_0
\end{matrix}\Bigg]+\Bigg[\begin{matrix}
Q_{\beta_0}N_\lambda^\circ (g_0)\\
T_\lambda^\circ (g_0)
\end{matrix}\Bigg]-\frac{1}{n}\sum_{i=1}^n\e_i\Bigg[\begin{matrix} Q_{\beta_0}N_{X_i,Z_i}(g_0,\beta_0)\\
T_{X_i,Z_i}(g_0,\beta_0)
\end{matrix}
\Bigg]\right\|_2=o_p(n^{-1/2}).
\end{align*}
In view of Lemma~\ref{lem:qbeta}, we have $\|\beta_0\beta_0\trans\hat\beta_n\|_2=O_p(r_n^2)=o_p(n^{-1/2})$, so that
\begin{align*}
\left\|\Bigg[\begin{matrix}
\hat\beta_n-\beta_0\\
\hat\gamma_n-\gamma_0
\end{matrix}\Bigg]+\Bigg[\begin{matrix}
Q_{\beta_0}N_\lambda^\circ (g_0)\\
T_\lambda^\circ (g_0)
\end{matrix}\Bigg]-\frac{1}{n}\sum_{i=1}^n\e_i\Bigg[\begin{matrix} Q_{\beta_0}N_{X_i,Z_i}(g_0,\beta_0)\\
T_{X_i,Z_i}(g_0,\beta_0)
\end{matrix}
\Bigg]\right\|_2=o_p(n^{-1/2})+\|\beta_0\beta_0\trans\hat\beta_n\|_2=o_p(n^{-1/2}).
\end{align*}
This completes the proof of Theorem~\ref{thm:marginal2}.

\end{proof}

\subsubsection{Proof of Theorem~\ref{thm:marginal}}\label{app:11}

Assume \eqref{newcond} is valid. Recall from Proposition~\ref{cor:sxz} and equation \eqref{stemp} that
\begin{align*}
&H_{U_i}(s)=K(X_i\trans\beta_0,s)-A_X(s)\trans Q_{\beta_0} N_{U_i}-A_Z(s)\trans T_{U_i},\\
&H_\lambda^\circ g_0(s)=(M_\lambda g_0)(s)-A_X(s)\trans Q_{\beta_0} N^\circ_{\lambda}(g_0)-A_Z(s)\trans T^\circ_{\lambda}(g_0).
\end{align*}
We have the following decomposition
\begin{align}\label{gn}
\sqrt{n}\lambda^{1/(4m)}\bigg\{\hat g_n(s)-g_0(s)+M_\lambda  g_0(s)- \frac{1}{n}\sum_{i=1}^n\e_i K(X_i\trans\beta_0,s)\bigg\}=I_{1,n}(s)+I_{2,n}(s)+I_{3,n}(s),
\end{align}
where
\begin{align}\label{in1in2}
&I_{1,n}(s)=\sqrt n\lambda^{1/(4m)}\bigg\{\hat g _{n}(s)- g _0(s)+H_\lambda^\circ g_0(s)-\frac{1}{n}\sum_{i=1}^n\e_i H_{U_i}(s)\bigg\},\notag\\
&I_{2,n}(s)=-\sqrt{n}\lambda^{1/(4m)}\{H_\lambda^\circ g_0(s)-M_\lambda g_0(s)\},\notag\\
&I_{3,n}(s)=n^{-1/2}\lambda^{1/(4m)}\sum_{i=1}^n\e_i\{H_{U_i}(s)-K(X_i\trans\beta_0,s)\}.
\end{align}

For the first term $I_{1,n}(s)$ in \eqref{in1in2}, by the Bahadur representation in Theorem~\ref{thm:bahadur} and Lemma~\ref{lem:normK}, we obtain that, for the $S_{U_i}=(H_{U_i},N_{U_i},T_{U_i})$ defined in \eqref{stemp},
\begin{align}\label{temp2}
\bigg\|\hat g_n-g_0+H_\lambda^\circ g_0-\frac{1}{n}\sum_{i=1}^n\e_iH_{U_i}\bigg\|_K\leq c\bigg\|\hat w_n-w_0+S_\lambda^\circ(g_0)-\frac{1}{n}\sum_{i=1}^n\e_i S_{U_i}\bigg\|=O_p(a_n).
\end{align}
Therefore, we obtain from the above equation, \eqref{ksnorm} and \eqref{temp2} that
\begin{align}\label{i1n}
&\sup_{s\in\I}|I_{1,n}(s)|\leq\sqrt n\lambda^{1/(4m)}\sup_{s\in\I}\bigg|\bigg\l K_s,\hat g_n-g_0+H_\lambda^\circ g_0-\frac{1}{n}\sum_{i=1}^n\e_i S_{U_i}\bigg\r_K\bigg|\notag\\
&\leq \sqrt n\lambda^{1/(4m)}\sup_{s\in\I}\|K_s\| _K\times\bigg\|\hat g_n-g_0+H_\lambda^\circ g_0-\frac{1}{n}\sum_{i=1}^n\e_i S_{U_i}\bigg\|_K=O_p(\sqrt n a_n)=o_p(1).
\end{align}

For the second term $I_{2,n}(s)$ in \eqref{in1in2}, from Lemma~\ref{lem:cross}, we obtain that
\begin{align*}
&\sup_{s\in\I}\{\|A_X(s)\|_2+\|A_Z(s)\|_2\}\leq c,
\\
&\|{V}(g_0'\cdot R_X,M_\lambda g_0)\|_2+\|{V}(R_Z,M_\lambda g_0)\|_2=o(n^{-1/2}).
\end{align*}
Therefore, recalling the definition of $N^\circ_{\lambda}(g_0)$ and $T^\circ_{\lambda}(g_0)$ in Proposition~\ref{cor:sxz}, we obtain that
\begin{align}\label{i2n}
&\sup_{s\in\I}|I_{2,n}(s)|=\sqrt{n}\lambda^{1/(4m)}\sup_{s\in\I}\big|A_X(s)\trans Q_{\beta_0} N^\circ_{\lambda}(g_0)+A_Z(s)\trans T^\circ_{\lambda}(g_0)\big|\notag\\
&\leq c\sqrt{n}\lambda^{1/(4m)}\sup_{s\in\I}\big\{\|A_X(s)\|_2\|{V}(g_0'\cdot R_X,M_\lambda g_0)\|_2+\|A_Z(s)\|_2\|{V}(R_Z,M_\lambda g_0)\|_2\big\}=o(1).
\end{align}

For the third term $I_{3,n}(s)$ in \eqref{in1in2}, in view of \eqref{stemp}, we have
\begin{align*}
I_{3,n}(s)&=-n^{-1/2}\lambda^{1/(4m)}\sum_{i=1}^n\e_i\{A_X(s)\trans Q_{\beta_0} N_{U_i}+A_Z(s)\trans T_{U_i}\}\\
&=-n^{-1/2}\lambda^{1/(4m)}\Bigg[\begin{matrix}
Q_{\beta_0}\trans A_X(s)\\
A_Z(s)
\end{matrix}\Bigg]\trans(\Omega_{\beta_0}+\Sigma_{\beta_0,\lambda})^{-1}\sum_{i=1}^n\e_i\Biggl[\begin{matrix}
Q_{\beta_0}\trans  \{g_0'(X_i\trans\beta_0)X_i-A_X(X_i\trans\beta_0)\}\\
Z_i-A_Z(X_i\trans\beta_0)
\end{matrix}\Biggl].
\end{align*}
Therefore, since $\sup_{s\in\I}\{\|A_X(s)\|_2+\|A_Z(s)\|_2\}\leq c$, we obtain that
\begin{align*}
\sup_{s\in\I}|I_{3,n}(s)|&\leq c n^{-1/2}\lambda^{1/(4m)}\Bigg[\bigg\|\sum_{i=1}^n\e_iQ_{\beta_0}\trans  \{g_0'(X_i\trans\beta_0)X_i-A_X(X_i\trans\beta_0)\}\bigg\|_2+\bigg\|\sum_{i=1}^n\e_i\{Z_i-A_Z(X_i\trans\beta_0)\}\bigg\|_2\Bigg].
\end{align*}
Note that $\E[\e_iQ_{\beta_0}\trans  \{g_0'(X_i\trans\beta_0)X_i-A_X(X_i\trans\beta_0)\}]=\E[\e_i\{Z_i-A_Z(X_i\trans\beta_0)\}]=0$ and
\begin{align}\label{new66}
\bigg\|\sum_{i=1}^n\e_iQ_{\beta_0}\trans  g_0'(X_i\trans\beta_0)X_i\bigg\|_2+\bigg\|\sum_{i=1}^n\e_iZ_i\bigg\|_2=O_p(\sqrt n).
\end{align}
In addition, in view of the definition of $A_X$ in \eqref{abx}, we obtain that, for $1\leq k\leq p$, by the Cauchy-Schwarz inequality,
\begin{align*}
&\var[\e_i e_k\trans Q_{\beta_0}\trans A_X(X_i\trans\beta_0)]=\var\bigg[\e_ie_k\trans Q_{\beta_0}\trans  \sum_{j\geq1}\frac{{V}(g_0'\cdot R_{X},\phi_j)}{1+\lambda\rho_j}\phi_j(X_i\trans\beta_0)\bigg]\\
&=\sum_{j\geq1}\frac{e_k\trans Q_{\beta_0}\trans{V}(g_0'\cdot R_X,\phi_j){V}(g_0'\cdot R_X,\phi_j)\trans Q_{\beta_0}e_k}{(1+\lambda\rho_j)^2}\leq\sum_{j\geq1}\frac{\|{V}(g_0'\cdot R_X,\phi_j)\|_2^2}{(1+\lambda\rho_j)^2}
\\
&\leq\bigg[\sum_{j\geq1}\|{V}(g_0'\cdot R_X,\phi_j)\|_2^2(1+\rho_j)^\mu\bigg]^{1/2}\times\bigg[\sum_{j\geq1}\|{V}(g_0'\cdot R_X,\phi_j)\|_2^2(1+\lambda\rho_j)^{-4}(1+\rho_j)^{-\mu}\bigg]^{1/2}
\\
&\leq\bigg[\sum_{j\geq1}\|{V}(g_0'\cdot R_X,\phi_j)\|_2^2(1+\rho_j)^\mu\bigg]^{1/2}\times\bigg[\sum_{j\geq1}\|{V}(g_0'\cdot R_X,\phi_j)\|_2^2\bigg]^{1/2}\times\bigg[\sum_{j\geq1}(1+\lambda\rho_j)^{-4}(1+\rho_j)^{-\mu}\bigg]^{1/2} 
\\
&=O(1),
\end{align*}
where we applied Assumptions~\ref{a:finite} and \eqref{newcond} and the fact that $\sum_{j\geq1}(1+\rho_j)^{-\mu}<+\infty$ in the last step. Similarly, we have $\var[\e_i A_Z(X_i\trans\beta_0)]=O(1)$.
Therefore, we deduce that
\begin{align*}
\bigg\|\sum_{i=1}^n\e_iQ_{\beta_0}\trans  A_X(X_i\trans\beta_0)\bigg\|_2+\bigg\|\sum_{i=1}^n\e_iA_Z(X_i\trans\beta_0)\bigg\|_2=O_p(\sqrt n).
\end{align*}
Combining the above equation and \eqref{new66} yields
\begin{align}\label{i3n}
\sup_{s\in\I}|I_{3,n}(s)|\leq cn^{-1/2}\lambda^{1/(4m)}\times O_p(\sqrt n)=o_p(1).
\end{align}
The proof of the result on $\hat g_n$ in \eqref{m22} is complete by combining \eqref{i1n}--\eqref{i3n}.

Next, suppose that condition \eqref{newcond} holds and we prove the result regarding $[\hat g_n\trans,\hat\beta_n\trans]\trans$ in Theorem~\ref{thm:marginal}. We first have the following decomposition:
\begin{align}\label{f1}
&\Bigg[\begin{matrix}
\hat\beta_n-\beta_0\\
\hat\gamma_n-\gamma_0
\end{matrix}\Bigg]-\frac{1}{n}\sum_{i=1}^n\e_i\Omega^+\Bigg[\begin{matrix} g_0'(X_i\trans\beta_0)\{X_i-R_X(X_i\trans\beta_0)\}\\
Z_i-R_Z(X_i\trans\beta_0)
\end{matrix}
\Bigg]=I_{1,n}+I_{2,n}+I_{3,n}, 
\end{align}
where
\begin{align*}
I_{1,n}&=\Bigg[\begin{matrix}
\hat\beta_n-\beta_0\\
\hat\gamma_n-\gamma_0
\end{matrix}\Bigg]+\Bigg[\begin{matrix}
Q_{\beta_0}N_\lambda^\circ (g_0)\\
T_\lambda^\circ (g_0)
\end{matrix}\Bigg]-\frac{1}{n}\sum_{i=1}^n\e_i\Bigg[\begin{matrix} Q_{\beta_0}N_{X_i,Z_i}(g_0,\beta_0)\\
T_{X_i,Z_i}(g_0,\beta_0)
\end{matrix}
\Bigg],
\\
I_{2,n}&=-\Bigg[\begin{matrix}
Q_{\beta_0}N_\lambda^\circ (g_0)\\
T_\lambda^\circ (g_0)
\end{matrix}\Bigg],
\\
I_{3,n}&=\frac{1}{n}\sum_{i=1}^n\e_i\left\{\Bigg[\begin{matrix} Q_{\beta_0}N_{X_i,Z_i}(g_0,\beta_0)\\
T_{X_i,Z_i}(g_0,\beta_0)
\end{matrix}
\Bigg]-\Omega^+\Bigg[\begin{matrix} g_0'(X_i\trans\beta_0)\{X_i-R_X(X_i\trans\beta_0)\}\\
Z_i-R_Z(X_i\trans\beta_0)
\end{matrix}
\Bigg]\right\}
\end{align*}

For the first term $I_{1,n}$, it follows from \eqref{m11} that $\|I_{1,n}\|_2=o_p(n^{-1/2})$. For the second term $I_{2,n}$, denote $D=\Big[
\begin{smallmatrix}
Q_{\beta_0}&0_{p,q}\\
0_{q,p-1}&I_q
\end{smallmatrix}\Big]$, Recall from the definition of $N^\circ_{\lambda}(g_0)$ and $T^\circ_{\lambda}(g_0)$ in Proposition~\ref{cor:sxz} that
\begin{align*}
\Bigg[\begin{matrix}
Q_{\beta_0}N^\circ_{\lambda}(g_0)\\
T^\circ_{\lambda}(g_0)
\end{matrix}\Bigg]&=-D(\Omega_{\beta_0}+\Sigma_{\beta_0,\lambda})^{-1}D
\Bigg[\begin{matrix}
{V}(g_0'\cdot R_X,M_\lambda g_0)\\
{V}(R_Z,M_\lambda g_0)
\end{matrix}\Bigg],
\end{align*}
which in view of the fact that $\|D\|_2\leq 1$ implies that
\begin{align*}
\|I_{2,n}\|_2\leq\|(\Omega_{\beta_0}+\Sigma_{\beta_0,\lambda})^{-1}\|_2\times\Big\{\|{V}(g_0'\cdot R_X,M_\lambda g_0)\|_2+\|{V}(R_Z,M_\lambda g_0)\|_2\Big\}.
\end{align*}
In addition, we obtain from Lemma~\ref{lem:cross} that
\begin{align*}
&\|{V}(g_0'\cdot R_X,M_\lambda g_0)\|_2+\|{V}(R_Z,M_\lambda g_0)\|_2=o(n^{-1/2}).
\end{align*}
We therefore deduce that
\begin{align}\label{f2}
\|I_{2,n}\|_2\leq c\Big\{\|{V}(g_0'\cdot R_X,M_\lambda g_0)\|_2+\|{V}(R_Z,M_\lambda g_0)\|_2\Big\}=o(n^{-1/2}).
\end{align}

For the third term $I_{3,n}$, recalling the definition of $N_{x,z}(g,\beta)$ and $T_{x,z}(g,\beta)$ in Proposition~\ref{cor:sxz}, we have
\begin{align}\label{s}
I_{3,n}&=\Biggl[
\begin{matrix}
Q_{\beta_0}& 0_{p,q}\\
 0_{q,p-1}&I_q
\end{matrix}
\Biggl]\times\frac{1}{n}\sum_{i=1}^n\e_i\left\{(\Omega_{\beta_0}+\Sigma_{\beta_0,\lambda})^{-1}\Biggl[\begin{matrix}
Q_{\beta_0}\trans  \{g_0'(X_i\trans\beta_0)X_i-A_X(X_i\trans\beta_0)\}\\
Z_i-A_Z(X_i\trans\beta_0)
\end{matrix}\Biggl]\right.\notag
\\
&\left.\hspace{5cm}-\Omega_{\beta_0}^{-1}\Bigg[\begin{matrix} g_0'(X_i\trans\beta_0)Q_{\beta_0}\trans\{X_i-R_X(X_i\trans\beta_0)\}\\
Z_i-R_Z(X_i\trans\beta_0)
\end{matrix}
\Bigg]\right\}\notag
\\
&=D(I_{3,1,n}+I_{3,2,n}),
\end{align}
where
\begin{align}\label{b}
I_{3,1,n}&=(\Omega_{\beta_0}+\Sigma_{\beta_0,\lambda})^{-1}D\times\frac{1}{n}\sum_{i=1}^n\e_i\Biggl[\begin{matrix}
g_0'(X_i\trans\beta_0)R_X(X_i\trans\beta_0)-A_X(X_i\trans\beta_0)\\
R_Z(X_i\trans\beta_0)-A_Z(X_i\trans\beta_0)
\end{matrix}\Biggl],\notag
\\
I_{3,2,n}&=\big\{(\Omega_{\beta_0}+\Sigma_{\beta_0,\lambda})^{-1}-\Omega_{\beta_0}^{-1}\big\}D\times\frac{1}{n}\sum_{i=1}^n\e_i\Bigg[\begin{matrix} g_0'(X_i\trans\beta_0)\{X_i-R_X(X_i\trans\beta_0)\}\\
Z_i-R_Z(X_i\trans\beta_0)
\end{matrix}
\Bigg].
\end{align}
For $I_{3,1,n}$, observe from the definition of $A_X$ in \eqref{abx} that
\begin{align*}
g_0'\cdot R_X-A_X=\sum_{j\geq1}{V}(g_0'\cdot R_{X},\phi_j)\phi_j-\sum_{j\geq1}\frac{{V}(g_0'\cdot R_{X},\phi_j)}{1+\lambda\rho_j}\phi_j=\lambda\sum_{j\geq1}\frac{{V}(g_0'\cdot R_{X},\phi_j)\rho_j}{1+\lambda\rho_j}\phi_j.
\end{align*}
We therefore obtain from the Cauchy-Schwarz inequality and condition~\ref{newcond} that, for $1\leq k\leq p$,
\begin{align*}
&\var\big[\e_i e_k\trans\{g_0'(X_i\trans\beta_0)R_X(X_i\trans\beta_0)-A_X(X_i\trans\beta_0)\}\big]=\lambda^2\E\Bigg[\e_i^2e_k\trans\bigg\{\sum_{j\geq1}\frac{{V}(g_0'\cdot R_{X},\phi_j)\rho_j}{1+\lambda\rho_j}\phi_j(X_i\trans\beta_0)\bigg\}^{\otimes2}e_k\Bigg]
\\
&=\sum_{j\geq1}\frac{e_k\trans{V}(g_0'\cdot R_{X},\phi_j)^{\otimes2}e_k\rho_j^2}{(1+\lambda\rho_j)^2}\leq\lambda^2\sum_{j\geq1}\frac{\|{V}(g_0'\cdot R_{X},\phi_j)\|_2^2\rho_j^2}{(1+\lambda\rho_j)^2} 
\\
&\leq\bigg[\sum_{j\geq1}\|{V}(g_0'\cdot R_X,\phi_j)\|_2^2(1+\rho_j)^\mu\bigg]^{1/2}\times\bigg[\lambda^4\sum_{j\geq1}\|V(g_0'\cdot R_X,\phi_j)\|_2^2(1+\lambda\rho_j)^{-4}(1+\rho_j)^{-\mu}\rho_j^4\bigg]^{1/2}
\\
&\leq\bigg[\sum_{j\geq1}\|{V}(g_0'\cdot R_X,\phi_j)\|_2^2(1+\rho_j)^\mu\bigg]^{1/2}\times\bigg[\sum_{j\geq1}\|V(g_0'\cdot R_X,\phi_j)\|_2^2\bigg]^{1/2}\times\bigg[\sum_{j\geq1}(1+\rho_j)^{-\mu}\bigg]^{1/2}.
\end{align*}
It therefore follows from the dominated convergence theorem that $$\var\big[\e_i e_k\trans\{g_0'(X_i\trans\beta_0)R_X(X_i\trans\beta_0)-A_X(X_i\trans\beta_0)\}\big]=o(1),\qquad 1\leq k\leq p.$$
Similarly, we have $\var\big[\e_i e_\ell\trans\{R_Z(X_i\trans\beta_0)-A_Z(X_i\trans\beta_0)\}\big]=o(1)$, for $1\leq\ell\leq q$. Therefore, observing that $\|D\|_2\leq1$, we deduce that
\begin{align}\label{r}
\|I_{3,1,n}\|_2&\leq \|(\Omega_{\beta_0}+\Sigma_{\beta_0,\lambda})^{-1}\|_2\times\left\|\frac{1}{n}\sum_{i=1}^n\e_i\Biggl[\begin{matrix}
g_0'(X_i\trans\beta_0)R_X(X_i\trans\beta_0)-A_X(X_i\trans\beta_0)\\
R_Z(X_i\trans\beta_0)-A_Z(X_i\trans\beta_0)
\end{matrix}\Biggl]\right\|_2=o_p(n^{-1/2}).
\end{align}

For the second term $I_{3,2,n}$ in \eqref{b} observe that $\big\|\big\{(\Omega_{\beta_0}+\Sigma_{\beta_0,\lambda})^{-1}-\Omega_{\beta_0}^{-1}\big\}\big\|_2=o(1)$ due to Proposition~\ref{prop:ru}, and observe that
\begin{align*}
\left\|\frac{1}{n}\sum_{i=1}^n\e_i\Bigg[\begin{matrix} g_0'(X_i\trans\beta_0)\{X_i-R_X(X_i\trans\beta_0)\}\\
Z_i-R_Z(X_i\trans\beta_0)
\end{matrix}
\Bigg]\right\|_2=O_p(n^{-1/2}).
\end{align*}
Therefore, we deduce that $\|I_{3,2,n}\|_2=o_p(n^{-1/2})$. Combining this with \eqref{s} and \eqref{r} yields
\begin{align*}
\|I_{3,n}\|_2\leq \|D\|_2\times(\|I_{3,1,n}\|_2+\|I_{3,2,n}\|_2)=o_p(n^{-1/2}).
\end{align*}
Combining the above equation with \eqref{f1} and \eqref{f2} concludes the proof of \eqref{m22}.

\subsubsection{Proof of Proposition~\ref{prop:rootn}}\label{app:prop:rootn}

Since $\X$ and $\Z$ are compact, we deduce that
\begin{align*}
\left\|\frac{1}{n}\sum_{i=1}^n\e_i\Omega^+\Bigg[\begin{matrix} g_0'(X_i\trans\beta_0)\{X_i-R_X(X_i\trans\beta_0)\}\\
Z_i-R_Z(X_i\trans\beta_0)
\end{matrix}
\Bigg]\right\|_2=O_p(n^{-1/2}).
\end{align*}
Therefore, it follows from \eqref{m22} in Theorem~\ref{thm:marginal} that
\begin{align*}
\left\|\Bigg[\begin{matrix}
\hat\beta_n-\beta_0\\
\hat\gamma_n-\gamma_0
\end{matrix}
\Bigg]\right\|_2&\leq\left\|\Bigg[\begin{matrix}
\hat\beta_n-\beta_0\\
\hat\gamma_n-\gamma_0
\end{matrix}\Bigg]-\frac{1}{n}\sum_{i=1}^n\e_i\Omega^+\Bigg[\begin{matrix} g_0'(X_i\trans\beta_0)\{X_i-R_X(X_i\trans\beta_0)\}\\
Z_i-R_Z(X_i\trans\beta_0)
\end{matrix}
\Bigg]\right\|_2
\\
&\quad+\left\|\frac{1}{n}\sum_{i=1}^n\e_i\Omega^+\Bigg[\begin{matrix} g_0'(X_i\trans\beta_0)\{X_i-R_X(X_i\trans\beta_0)\}\\
Z_i-R_Z(X_i\trans\beta_0)
\end{matrix}
\Bigg]\right\|_2=O_p(n^{-1/2}),
\end{align*}
which implies that $\|\hat\beta_n-\beta_0\|_2+\|\hat\gamma_n-\gamma_0\|_2=O_p(n^{-1/2})$ and completes the proof of Proposition~\ref{prop:rootn}.

\section{Proof of theoretical results in Section~\ref{sec:stat}}\label{app:sec:stat}

\subsection{Proof of Theorem~\ref{thm:minimax}}\label{app:thm:minimax}

For the upper bound in (i), by taking $\lambda\asymp n^{-2m/(2m+1)}$, the upper bound in (i) follows from Theorem~\ref{thm:rate}, lemma~\ref{lem:normK}, and the fact that
\begin{align*}
\mathfrak R^2(\hat g_n,\hat\beta_n,\hat\gamma_n)&\leq c_1^{-1}{V}(\hat g_n-g_0,\hat g_n-g_0)+c\|\hat w_n-w_0\|^2+c\|\beta_0\beta_0\trans\hat\beta_n\|^2
\\
&\leq c_1^{-1}\|\hat w_n-w_0\|^2+c\|\hat w_n-w_0\|^4 =O_p(n^{-2m/(2m+1)}),
\end{align*}
where we applied Lemma~\ref{lem:qbeta} and Proposition~\ref{prop:2.1}, and the assumption that $\sigma_0^2(s)f_{X\trans\beta_0}(s)\geq c_1>0$ for any $s\in\I$.

For the lower bound in (ii), it suffices to prove the lower bound for $V(\tilde g_n-g_0,\tilde g_n-g_0)$, since
\begin{align*}
\mathfrak R^2(\tilde g_n,\tilde\beta_n,\tilde\gamma_n)\geq\int_{\I}\{\tilde g_n(s)-g_0(s)\}^2\d s\geq c_2^{-1}V(\tilde g_n-g_0,\tilde g_n-g_0),
\end{align*}
where we applied the assumption that $\|\sigma_0^2f_{X\trans\beta_0}\|_\infty\leq c_2$ for some $c_2>0$ in Assumption~\ref{a:rank}. Let $\mathcal F_0\subset\mathcal F$ be the set of the joint distribution $F_{X,Y,Z}$ according to model \eqref{model} where $\gamma_0=0_q$ and $\e|X,Z\sim N(0,1)$ a.s. It therefore suffices to show the lower bound for $F_{X,Y,Z}\in\mathcal F_0$.
By Theorem~2.5 in \cite{tsybakov2009nonparametric}, in order to show the lower bound, we need to show that, for $M\geq2$  to be specified later, $\H_m$ contains elements $ g _1,\ldots, g _{M+1}$ that satisfy the following two conditions:
\begin{enumerate}[label=(C\arabic*),noitemsep,series=conditionC]
\item\label{c1} $V(g _{j}- g _{j'},g _{j}- g _{j'})\geq c_0n^{-2m/(2m+1)}$, for $1\leq j<j'\leq M+1$;

\item\label{c2} $M^{-1}\sum_{j=2}^{M+1}\mathcal{K}(P_j,P_1)\leq\alpha\log M$, where $0<\alpha<1/8$, $\mathcal{K}$ is the Kullback-Leibler divergence, and $P_j$ denotes joint distribution of $\{(X_{i,j},Y_{i,j},Z_{i,j})\}_{i=1}^n$, where $Y_{i,j}=g_j(X_{i,j}\trans\beta_0)+\e_{i,j}$, for $1\leq i\leq n$ and $1\leq j\leq M+1$.
\end{enumerate}

Define $\nu_n=\lfloor n^{1/(2m+1)}\rfloor$. Recall the eigenfunctions $\{\phi_j\}_{j\geq1}$ in Proposition~\ref{prop:eigen}. For any $\omega=(\omega_{\nu_n+1},\ldots,\omega_{2\nu_n})\in\{0,1\}^{\nu_n}$, let
\begin{align}\label{betaomega}
 g _{\omega}=c_1n^{-1/2}\sum_{k=\nu_n+1}^{2\nu_n}\omega_{k}\,\phi_{k}\,,
\end{align}
where $c_1>0$ is a constant independent of $n$ to be specified later. We first verify that the $ g _\omega$'s are elements in $\H$. Since, by Proposition~\ref{prop:eigen}, $\{\phi_{j}\}_{j\geq 1}$ diagonalizes the operator $J$ defined in \eqref{J}, we have 
\begin{align*}
\Vert g _\omega\Vert_K^2&=c_1^2\,n^{-1}\Bigg\l\sum_{k=\nu_n+1}^{2\nu_n}\omega_{k}\,\phi_{k}\,,\sum_{k'=\nu_n+1}^{2\nu_n}\omega_{k'}\,\phi_{k'}\Bigg\r_K=c_1^2\,n^{-1}\sum_{k=\nu_n+1}^{2\nu_n}\omega_{k}^2\,\Vert\phi_{k}\Vert_K^2\\
&\leq c\,n^{-1}\sum_{k=\nu_n+1}^{2\nu_n}(1+\lambda k^{2m})\,.
\end{align*}
Note that the above equation holds for any $\lambda>0$.
Therefore, taking $\lambda=1$ yields
\begin{align*}
\Vert g _\omega\Vert_K^2\leq c\,n^{-1}\sum_{k=\nu_n+1}^{2\nu_n}k^{2m}\leq c\,n^{-1}\nu_n^{1+2m}\leq c\,,
\end{align*}
which shows that, for any $\omega\in\{0,1\}^{\nu_n}$,  $ g _\omega$ defined in \eqref{betaomega} is an element of $\H_m$.

By the Varshamov-Gilbert bound (Lemma~2.9 in \citealp{tsybakov2009nonparametric}), for $\nu_n\geq8$, there exists a subset $\Omega=\{\omega^{(1)},\ldots,\omega^{(M+1)}\}\in\{0,1\}^{\nu_n}$ with $M\geq 2^{\nu_n/8}$ such that, $\omega^{(1)}=(0,\ldots,0)$ and for any $1\leq j< j'\leq M+1$,
\begin{align*}
{\rm H}\big(\omega^{(j_1)},\omega^{(j_2)}\big)\geq\frac{\nu_n}{8}\,,
\end{align*}
where ${\rm H}(\cdot,\cdot)$ is the Hamming distance. For $1\leq j\leq M+1$, let $\omega^{(j)}=(\omega^{(j)}_{\nu_n+1},\ldots,\omega^{(j)}_{2\nu_n})$. Let $ g _1,\ldots, g _{M+1}$ denote the functions defined as in \eqref{betaomega} that corresponds to $\omega^{(1)},\ldots,\omega^{(M+1)}\in\Omega$. For $1\leq j<j'\leq M+1$, in view of
\eqref{betaomega},
\begin{align}\label{betadiff}
 g _{j}- g _{j'}=c_1\,n^{-1/2}\sum_{k=\nu_n+1}^{2\nu_n}\big(\omega^{(j)}_{k}-\omega^{(j')}_{k}\big)\,\phi_{k}\,.
\end{align}
By Proposition~\ref{prop:eigen}, since $\{\phi_{j}\}_{j\geq 1}$ diagonalizes the operator $V$ defined in \eqref{V}, we deduce from \eqref{betadiff} that
\begin{align*}
{V}(g _{j}-g _{j'},g _{j}-g _{j'})&=c_1^2\,n^{-1}\sum_{k=\nu_n+1}^{2\nu_n}\big(\omega^{(j)}_{k}-\omega^{(j')}_{k}\big)^2\,{V}(\phi_{k},\phi_{k})\\
&=c_1^2\,n^{-1}\sum_{k=\nu_n+1}^{2\nu_n}\one\big\{\omega^{(j)}_{k}\neq\omega^{(j')}_{k}\big\}\\
&=c_1^2\,n^{-1}{\rm H}\big(\omega^{(j)},\omega^{(j')}\big)\geq c_1^2\,8^{-1}n^{-1}\nu_n\geq c_1^2\,8^{-1}n^{-2m/(2m+1)}\,.
\end{align*}
By taking $c_0=c_1^2/16$, the above equation implies that Condition~\ref{c1} is valid.

For any $1\leq j<j'\leq M+1$, in view of \eqref{betadiff}, applying the formula for the KL-divergence between Gaussian distributions, we obtain that
\begin{align*}
\mathcal{K}(P_{j},P_{j'})&=\frac{n}{2}\,\E\big[g_j(X
\trans\beta_0)-g_{j'}(X\trans\beta_0)\big]^2=\frac{n}{2}{V}(g _{j}-g _{j'},g _{j}-g _{j'})\\
&=c_1^2\,\sum_{k=\nu_n+1}^{2\nu_n}\big(\omega^{(j)}_{k}-\omega^{(j')}_{k}\big)^2\leq c_1^2\,\nu_n\,.
\end{align*}
Therefore, for any $0<\alpha<1/8$, by taking $0<c_1<\sqrt{\alpha\log 2/8}$ in \eqref{betaomega}, we have
\begin{align*}
\frac{1}{M}\sum_{j=2}^{M+1}\mathcal{K}(P_j,P_1)\leq c_1^2\,\nu_n\leq\frac{\alpha\nu_n\log2}{8}\leq\alpha\log M\,,
\end{align*}
which verifies Condition~\ref{c2} and completes the proof.

\subsection{Proof of Lemma~\ref{prop:bias}}\label{app:prop:bias}

We start by stating and proving the following useful lemma regarding $M_\lambda$.

\begin{lemma}\label{lem:Mlambda}
For any $g\in\H_m$, let $M_\lambda^{1/2}$ be the square-root operator of $M_\lambda$ in \eqref{mlambda} defined by
\begin{align*}
M_\lambda^{1/2}g=\sqrt\lambda\sum_{j\geq1}{V}(g,\phi_j)\sqrt{\frac{\rho_{j}}{1+\lambda\rho_{j}}}\phi_{j},\qquad g\in\H_m.
\end{align*}
Then, it holds that $\|M_\lambda^{1/2} (g)\|_K^2=\lambda J(g,g)$. Furthermore, it is true that $$\|M_\lambda (g)\|_K^2\leq\lambda J(g,g)\leq\|g\|_K^2.$$
\end{lemma}

\begin{proof}[\underline{Proof of Lemma~\ref{lem:Mlambda}}]

First, it is clear that $M_\lambda^{1/2}$ is self-adjoint. Then, for $g\in\H_m$, it holds that
\begin{align*}
\|M_\lambda^{1/2}g\|_K^2=\l M_\lambda^{1/2}g,M^{1/2}_\lambda g\r_K=\l M_\lambda g,g\r_K=\lambda J(g,g)\leq\|g\|_K^2.
\end{align*}
In addition, it is true that
\begin{align*}
\|M_\lambda g\|_K^2=\|M_\lambda^{1/2}(M_\lambda^{1/2}g)\|_K^2=\lambda J(M_\lambda^{1/2} g,M_\lambda^{1/2} g)\leq \|M_\lambda^{1/2}g\|_K^2.
\end{align*}
Combining the above two equations yields that $\|M_\lambda g\|_K^2\leq\lambda J(g,g)$ and $\|M_\lambda g\|_K^2\leq\|g\|_K^2$, which completes the proof of Lemma~\ref{lem:Mlambda}.

\end{proof}

Now, we prove Lemma~\ref{prop:bias}.
We start by introducing the notation
\begin{align*}
\delta_n=\hat g_n-g_0+M_\lambda g_0-\frac{1}{n}\sum_{i=1}^n\e_i K_{X_i\trans\beta_0}.
\end{align*}
It follows from Theorem~\ref{thm:marginal} that
\begin{align*}
\sup_{s\in\I}\sqrt n\lambda^{1/(4m)}|\delta_n(s)|=\sup_{s\in\I}\sqrt n\lambda^{1/(4m)}\bigg|\hat g_n(s)-g_0(s)+M_\lambda g_0(s)-\frac{1}{n}\sum_{i=1}^n\e_i K(X_i\trans\beta_0,s)\bigg|=o_p(1).
\end{align*}
We obtain from the definition of $\delta_n$ the following decomposition:
\begin{align*}
M_\lambda\hat g_n(s)-M_\lambda g_0(s)=M_\lambda(\delta_n)(s)-M_\lambda^2g_0(s)+\frac{1}{n}\sum_{i=1}^n\e_iM_\lambda(K_{X_i\trans\beta_0})(s),\qquad s\in\I.
\end{align*}
First, applying lemma~\ref{lem:Mlambda} yields
\begin{align*}
\|M_\lambda(\delta_n)\|_K\leq \|\delta_n\|_K=o_p(n^{-1/2}).
\end{align*}
Therefore,
\begin{align*}
\sup_{s\in\I}\sqrt n\lambda^{1/(4m)}|M_\lambda(\delta_n)(s)|&\leq\sqrt n\lambda^{1/(4m)}\sup_{s\in\I}|\l M_\lambda(\delta_n),K_s\r_K|
\\
&\leq \sqrt n\lambda^{1/(4m)}\|M_\lambda(\delta_n)\|_K\times\sup_{s\in\I}\|K_s\|_K=o_p(1).
\end{align*}

Second, since $M_\lambda$ is self-adjoint, we have
\begin{align*}
&\sup_{s\in\I}\sqrt n\lambda^{1/(4m)}|M_\lambda^2 g_0(s)|=\sqrt n\lambda^{1/(4m)}\sup_{s\in\I}|\l M_\lambda^2 g_0,K_s\r_K|=\sqrt n\lambda^{1/(4m)}\sup_{s\in\I}|\l M_\lambda g_0,M_\lambda K_s\r_K|
\\
&\qquad\leq\sqrt n\lambda^{1/(4m)}\|M_\lambda g_0\|_K\times\sup_{s\in\I}\|M_\lambda (K_s)\|_K\leq  O(\sqrt n\lambda^{(2m+1)/(4m)})\times o(1)=o(1).
\end{align*}

Third, observe that $\E\{\e^2\phi_j(X\trans\beta_0)\phi_{j'}(X\trans\beta_0)\}=V(\phi_j,\phi_{j'})=\delta_{jj'}$. Hence,
\begin{align*}
\E\bigg\|\frac{1}{n}\sum_{i=1}^n\e_iK_{X_i\trans\beta_0}\bigg\|_K^2&=\frac{1}{n^2}\sum_{i_1,i_2=1}^n\E\Big(\e_{i_1}\e_{i_2}\l K_{X_{i_1}\trans\beta_0},K_{X_{i_2}\trans\beta_0}\r_K\Big)=\frac{1}{n^2}\sum_{i=1}^n\E\Big\{\e_{i}^2K(X_i\trans\beta_0,X_i\trans\beta_0)\Big\}
\\
&=n^{-1}\E\bigg\{\e^2\sum_{j\geq1}\frac{\phi_j(X\trans\beta_0)\phi_j(X\trans\beta_0)}{1+\lambda\rho_j}\bigg\}=n^{-1}\sum_{j\geq1}\frac{1}{1+\lambda\rho_j}=O(n^{-1}\lambda^{-1/(2m)}).
\end{align*}
By the Cauchy-Schwarz inequality, we deduce from the above equation that
\begin{align*}
&\sup_{s\in\I}\sqrt n\lambda^{1/(4m)}\bigg|\frac{1}{n}\sum_{i=1}^n\e_iM_\lambda(K_{X_i\trans\beta_0})(s)\bigg|=\sup_{s\in\I}\sqrt n\lambda^{1/(4m)}\bigg|\bigg\l  \frac{1}{n}\sum_{i=1}^n\e_iK_{X_i\trans\beta_0},M_\lambda K_s\bigg\r_K\bigg|
\\
&\leq\sqrt n\lambda^{1/(4m)}\bigg\|\frac{1}{n}\sum_{i=1}^n\e_iK_{X_i\trans\beta_0}\bigg\|_K\times\sup_{s\in\I}\|M_\lambda K_s\|_K=O_p(1)\times\sup_{s\in\I}\|M_\lambda K_s\|_K=o_p(1).
\end{align*}
The proof is therefore complete.

\subsection{Proof of Theorem~\ref{thm:cb:boot}}\label{app:thm:bootstrap}

The following proposition is useful for proving Theorem~\ref{thm:cb:boot}. Its proof leverages the marginal Bahadur representation established in Theorem~\ref{thm:marginal}.

\begin{proposition}\label{thm:process}
If Assumptions~\ref{a:rank}--\ref{a:b} and \eqref{newcond} hold, then 
\begin{align}\label{t}
T_n:=\sqrt{n}\lambda^{1/(4m)}\sup_{s\in\I}|\hat g_{n}(s)+M_\lambda \hat g_n(s)-g_0(s)|\converged T:=\sup_{s\in\I}|\mathbb G(s)|,
\end{align}
where $\{\mathbb G(s)\}_{s\in\I}$ is a mean-zero Gaussian process with covariance kernel
\begin{align*}
C_{\mathbb G}(s,t)=\lim_{\lambda\to0}\lambda^{1/(2m)}\sum_{j\geq1}\frac{\phi_j(s)\phi_j(t)}{(1+\lambda\rho_j)^2},\qquad\quad s,t\in\I.
\end{align*}

\end{proposition}

\begin{proof}[\underline{Proof of Proposition~\ref{thm:process}}]

Recalling the definition of the process $\G_n$ in Proposition~\ref{thm:process}, we have the following decomposition
\begin{align}\label{gn2}
\G_n(s)=\sqrt{n}\lambda^{1/(4m)}\{\hat g_n(s)-g_0(s)+M_\lambda \hat g_n(s)\}=I_{1,n}(s)+I_{2,n}(s)+I_{3,n}(s).
\end{align}
where
\begin{align}\label{in1in22}
&I_{1,n}(s)=\sqrt{n}\lambda^{1/(4m)}\bigg\{\hat g_n(s)-g_0(s)+M_\lambda g_0(s)-\frac{1}{n}\sum_{i=1}^n\e_i K(X_i\trans\beta_0,s)\bigg\},\notag\\
&I_{2,n}(s)=\sqrt{n}\lambda^{1/(4m)}\{M_\lambda\hat g_n(s)-M_\lambda g_0(s)\},\notag\\
&I_{3,n}(s)=n^{-1/2}\lambda^{1/(4m)}\sum_{i=1}^n\e_i K(X_i\trans\beta_0,s).
\end{align}
It follows from Propositions~\ref{thm:marginal} and \ref{prop:bias}, respectively, that the first two terms $I_{1,n}(s)$ and $I_{2,n}(s)$ in \eqref{in1in22} satisfy
\begin{align}\label{new6}
\sup_{s\in\I}\big\{|I_{1,n}(s)|+|I_{2,n}(s)|\big\}=o_p(1).
\end{align}
We therefore focus on the third term $I_{3,n}(s)$ in \eqref{in1in22}. Since $\E\big\{\e K(X\trans\beta_0,\cdot)\big\}\equiv0$, we compute
\begin{align}\label{var}
\cov\big\{I_{3,n}(s_1),I_{3,n}(s_2)\big\}&=\lambda^{1/(2m)}\cov\big\{\e K(X\trans\beta_0,s_1),\e K(X\trans\beta_0,s_2)\big\}\notag\\
&=\lambda^{1/(2m)}\,\E\Bigg[\e^2\bigg\{\sum_{k\geq1}\frac{\phi_k(X\trans\beta_0)}{1+\lambda\rho_{k}}\phi_{k}(s_1)\bigg\}\bigg\{\sum_{k'\geq1}\frac{\phi_{k'}(X\trans\beta_0)}{1+\lambda\rho_{k'}}\phi_{k'}(s_2)\bigg\}\Bigg]\notag\\
&=\lambda^{1/(2m)}\sum_{k,k'\geq1}\frac{\phi_{k}(s_1)\phi_{k'}(s_2){V}(\phi_k,\phi_{k'})}{(1+\lambda\rho_{k})(1+\lambda\rho_{k'})}\notag\\
&=\lambda^{1/(2m)}\sum_{k\geq1}\frac{\phi_{k}(s_1)\phi_{k}(s_2)}{(1+\lambda\rho_{k})^2},\qquad s_1,s_2\in\I.
\end{align}

We first prove the weak convergence of the finite-dimensional marginal distributions of $\G_n$. By the Cram\'er-Wold theorem, we shall show that, for any $r\in\mathbb{N}$, $(c_1,\ldots,c_r)\trans\in\mathbb{R}^r$ and $s_1,\ldots,s_r\in\I$,
\begin{align}\label{valid}
\sum_{j=1}^rc_j\G_n(s_j)\overset{d.}{\longrightarrow}\sum_{j=1}^rc_j \mathbb G(s_j)\,.
\end{align}
For $1\leq i\leq n$, let $\mathfrak U_{i,r}=\lambda^{1/(4m)}\e_i\sum_{j=1}^rc_j K(X_i\trans\beta_0,s_j)$. In view of \eqref{new6}, we deduce that
\begin{align*}
\sum_{j=1}^rc_j\G_n(s_j)&=\sum_{j=1}^rc_jI_{3,n}(s_j)+\sum_{j=1}^rc_j\{I_{1,n}(s_j)+I_{2,n}(s_j)\}=\frac{1}{\sqrt n}\sum_{i=1}^n\mathfrak U_{i,r}+o_p(1)\,.
\end{align*}
By \eqref{var}, we obtain
\begin{align*}
\sigma_{\mathfrak U}^2&:=\var(\mathfrak U_{i,r})=\lambda^{1/(2m)}\sum_{j_1,j_2=1}^rc_{j_1}c_{j_2}\,\cov\big\{\e_iK(X_i\trans\beta_0,s_{j_1}),\e_iK(X_i\trans\beta_0,s_{j_2})\big\}\\
&=\lambda^{1/(2m)}\sum_{j_1,j_2=1}^rc_{j_1}c_{j_2}\sum_{k\geq1}\frac{\phi_{k}(s_1)\phi_{k}(s_2)}{(1+\lambda\rho_{k})^2}=n^{-1}\lambda^{1/(2m)}\sum_{j_1,j_2=1}^rc_{j_1}c_{j_2}C_{\mathbb G}(s_{j_1},s_{j_2})+o(1)\\
&=\lambda^{1/(2m)}\sum_{j_1,j_2=1}^rc_{j_1}c_{j_2}\cov\big\{I_{3,n}(s_{j_1}),I_{3,n}(s_{j_2})\big\}+o(1)\,.
\end{align*}
When $\sum_{j_1,j_2=1}^rc_{j_1}c_{j_2}C_{\mathbb G}(s_{j_1},s_{j_2})=0$, we have that $\sum_{j=1}^rc_j \mathbb G(s_j)$ has a degenerate distribution with a point mass at zero, so that \eqref{valid} is followed by the Markov's inequality. When $\sum_{j_1,j_2=1}^rc_{j_1}c_{j_2}C_{\mathbb G}(s_{j_1},s_{j_2})\neq0$, to prove \eqref{valid}, we shall check that $\{\mathfrak U_{i,r}\}_{i=1}^n$ satisfy Lindeberg's condition. It follows from Proposition~\ref{prop:eigen} that
\begin{align*}
\sup_{s_1,s_2\in\I}|K(s_1,s_2)|=\sup_{s_1,s_2\in\I}\sum_{j\geq1}\frac{|\phi_j(s_1)\phi_j(s_2)|}{1+\lambda\rho_j}\leq \sum_{j\geq1}\frac{\|\phi_j\|_\infty^2}{1+\lambda\rho_j}\leq c\lambda^{-1/(2m)}\,.
\end{align*}
Therefore, it holds that
\begin{align*}
|\mathfrak U_{i,r}|=\lambda^{1/(4m)}|\e_i|\sum_{j=1}^rc_j |K(X_i\trans\beta_0,s_j)|\leq c\lambda^{-1/(4m)}|\e_i|.
\end{align*}
By the Cauchy-Schwarz inequality and Markov's inequality, under the assumption that $n\lambda^{1/m}\to\infty$, for any $e>0$, we have
\begin{align*}
&\sigma_{\mathfrak U}^{-2}\mathbb {E} \left[\mathfrak U_{i,r}^{2}\cdot \mathbf {1} \{|\mathfrak U_{i,r}|>e\sqrt n \sigma_{\mathfrak U}\}\right]\\
&\leq c\lambda^{-1/(2m)}\sigma_{\mathfrak U}^{-2}\mathbb {E} \left[\e_i^2\cdot \mathbf {1} \{|\e_i|>c^{-1}e\sqrt n \lambda^{1/(4m)}\sigma_{\mathfrak U}\}\right]\\
&\leq c\lambda^{-1/(2m)}\sigma_{\mathfrak U}^{-2}\times\{\E(\e_i^4)\}^{1/2}\times\big\{\P(|\e_i|>c^{-1}e\sqrt n\lambda^{1/(4m)}\sigma_{\mathfrak U})\big\}^{1/2}\\
&\leq c\lambda^{-1/(2m)}\sigma_{\mathfrak U}^{-2}(c^2e^{-2}n^{-1}\lambda^{-1/(2m)}\sigma_{\mathfrak U}^{-2})\E(\e_i^4)=O(n^{-1}\lambda^{-1/m})=o(1),\qquad n\to\infty.
\end{align*}
Therefore, by Lindeberg's CLT, we deduce that
\begin{align*}
\sum_{j=1}^rc_j\G_n(s_j)& =\frac{1}{\sqrt n}\sum_{i=1}^n\mathfrak U_{i,r}+o_p(1)\converged \sum_{j=1}^rc_j \mathbb G(s_j) \sim 
{\cal N} \bigg (0,\sum_{j_1,j_2=1}^rc_{j_1}c_{j_2}C_{\mathbb G}(s_{j_1},s_{j_2})\bigg),
\end{align*}
which proves \eqref{valid}.


Next, we show the equicontinuity of $\{\G_n(s)\}_{s\in\I}$ under the assumptions of  Proposition~\ref{thm:process}. We first focus on the leading term $I_{3,n}$ in \eqref{in1in2}, where we recall that
\begin{align}\label{hnst}
I_{3,n}(s)=\sum_{i=1}^n\mathfrak U_i(s)=n^{-1/2}\lambda^{1/(4m)}\sum_{i=1}^n\e_iK(X_i\trans\beta_0,s),\qquad s\in\I\,.
\end{align}
Let $\Psi(s)=x^2$ and let $\Vert U\Vert_\Psi=\inf\{c>0:\E\{\Psi(|U|/c)\}\leq 1\}$ denote the Orlicz norm for a real-valued random variable $U$. For some metric $d$ on $\I$, let $\D(w,d)$ denote the $w$-packing number of the metric space $(\I,d)$, where $d$ is an appropriate metric specified below. Note that $\E\{\e_iK(X_i\trans\beta_0,s)\}=0$ for any $s\in\I$, and $\E\{\e^2\phi_j(X\trans\beta_0)\phi_{j'}(X\trans\beta_0)\}=\E\{\sigma_0^2(X\trans\beta_0)\phi_j(X\trans\beta_0)\phi_{j'}(X\trans\beta_0)\}={V}(\phi_j,\phi_{j'})=\delta_{jj'}$. In view of the definition of $V$ in \eqref{V} and by applying Assumption~\ref{holder}, we deduce that
\begin{align}\label{m1}
\E|I_{3,n}(s_1)-I_{3,n}(s_2)|^2&=  \lambda^{1/(2m)}\,\E\big|\e\{K_{s_1}(X\trans\beta_0)-K_{s_2}(X\trans\beta_0)\}\big|^2\notag\\
&=\lambda^{1/(2m)}\E\big[\sigma_0^2(X\trans\beta_0)\{K_{s_1}(X\trans\beta_0)-K_{s_2}(X\trans\beta_0)\}^2\big]\notag\\
&=\lambda^{1/(2m)}\|K_{s_1}-K_{s_2}\|_V^2\leq c\lambda^{1/(2m)-\nu}|s_1-s_2|^{2\vartheta}.
\end{align}
We therefore deduce from \eqref{m1} that 
\begin{align}\label{m2}
\Vert I_{3,n}(s_1)-I_{3,n}(s_2)\Vert_\Psi\leq c\,\lambda^{1/(2m)-\nu}\,|s_1-s_2|^{\vartheta}.
\end{align}

Next, we shall show that, there exists a metric $d$ on $\I$ such that, for any $e>0$, 
\begin{align}\label{p}
\lim_{\delta\to0}\,\limsup_{n\to\infty}\,\P\bigg\{\sup_{d(s_1,s_2)\leq\delta}|I_{3,n}(s_1)-I_{3,n}(s_2)|>e\bigg\}=0\,,
\end{align}
where we distinguish the following two cases: $\vartheta>1$ and $0\leq\vartheta\leq1$.

In the first case where $\vartheta>1$, let $d_1(s_1,s_2)=|s_1-s_2|^{\vartheta}$. In view of \eqref{m1}, we have $\Vert I_{3,n}(s_1)-I_{3,n}(s_2)\Vert_\Psi\leq c\,d_1(s_1,s_2)$. Note that the packing number of $\I$ with respect to the metric $d_1$ satisfies $\D(\zeta,d_1)\lesssim \zeta^{-1/\vartheta}$. By Theorem~2.2.4 in \cite{vaart1996}, for any $e,\eta>0$,
\begin{align*}
&\P\bigg\{\sup_{d_1(s_1,s_2)\leq\delta}|I_{3,n}(s_1)-I_{3,n}(s_2)|>e\bigg\}\\
&\leq c\,\bigg\Vert\sup_{d_1(s_1,s_2)\leq\delta}|I_{3,n}(s_1)-I_{3,n}(s_2)|\bigg\Vert_\Psi\leq c\int_0^{\eta}\sqrt{\D(\zeta,d_1)}d\zeta+\delta\,\D(\eta,d_1)\\
&\leq c\int_0^{\eta}\zeta^{-1/(2\vartheta)}d\zeta+\delta\,\eta^{-1/\vartheta}=c\,\eta^{(2\vartheta-1)/(2\vartheta)}+\delta\,\eta^{-1/\vartheta}\,.
\end{align*}
Using $\eta=\sqrt{\delta}$ and $\vartheta>1/2$, it follows that \eqref{p} holds by taking the metric $d=d_1$.

In the second case where $0\leq\vartheta\leq1$, in order to show \eqref{p}, we shall use Lemma~\ref{lem:kley} in Section~\ref{app:aux}, which is a modified version of Lemma~A.1 in \cite{kley2016quantile}. Let 
\begin{align*}
d_2(s_1,s_2)=|s_1-s_2|^{2}
\end{align*}
and let
$\overline\eta= \lambda^{(1/m-2\nu)/(2-\vartheta)}$. In view of \eqref{m2}, we have, when $d_2(s_1,s_2)\geq\overline\eta/2>0$,
\begin{align*}
\Vert I_{3,n}(s_1)-I_{3,n}(s_2)\Vert_\Psi\leq c\, \lambda^{1/(2m)-\nu}\big\{d_2(s_1,s_2)\big\}^{\vartheta/2}\leq c\, d_2(s_1,s_2)\,.
\end{align*}
By Assumption~\ref{a:02} and Markov's inequality, by taking $c'>c_0^{-1}$, for the constant $c_0>0$ in Assumption~\ref{a:02}, it holds that
\begin{align*}
\P\Big(\max_{1\leq i\leq n}|\e_i|\geq c'\log n\Big)\leq n^{1-c'c_0}\E\{\exp(c_0|\e|)\}=n^{1-c'c_0}\E\big[\E\{\exp(c_0|\e|)\}|X,Z\big]=o(1),
\end{align*}
which implies that it suffices to restrict the proof on the event where $\max_{1\leq i\leq n}|\e_i|< c\log n$. Hence, by Proposition~\ref{prop:eigen},
\begin{align}\label{bound}
\sup_{s\in\I}|\mathfrak U_i(s)|&\leq n^{-1/2}\lambda^{1/(4m)}|\e_i|\sup_{s\in\I}K(X_i\trans\beta_0,s)=n^{-1/2}\lambda^{1/(4m)}|\e_i|\sup_{s\in\I}\sum_{j\geq1}\frac{\phi_j(X_i\trans\beta_0)\phi_j(s)}{1+\lambda\rho_j}\notag\\
&\leq cn^{-1/2}\lambda^{1/(4m)}\log(n)\times\sup_{j\geq1}\|\phi_j\|^2_\infty\times\sum_{j\geq1}\frac{1}{1+\lambda\rho_j}\leq c\,n^{-1/2}\lambda^{-1/(4m)}\log(n).
\end{align}
In addition, by equation \eqref{m1},
\begin{align*}
\sup_{s_1,s_2\in\I}\E|I_{3,n}(s_1)-I_{3,n}(s_2)|^2\leq c\,\lambda^{1/(2m)-\nu}.
\end{align*}
By Bernstein's inequality, combining the above equation with \eqref{bound}, we deduce that, for $n$ large enough, for any $s_1,s_2\in\I$ and for any $e>0$,
\begin{align}\label{bern}
&\P\Big\{ |I_{3,n}(s_1)-I_{3,n}(s_2)|>e/4\Big\}\leq 2\exp\bigg\{-\frac{e^2/16}{c\lambda^{1/(2m)-\nu}+en^{-1/2}\lambda^{-1/(4m)}\log (n)/6}\bigg\}\,.
\end{align}

Now, note that $\D(\zeta,d_2)\leq c\,\zeta^{-1}$ and recall that in the case of $0\leq\vartheta\leq1$ we have assumed $b<a$. By Lemma~\ref{lem:kley} and \eqref{bern}, there exists a set $\tilde\I\subset \I$ that contains at most $\D(\overline\eta,d_2)$ points, such that, for any $\delta,e>0$ and $\eta>\overline\eta$, as $n\to\infty$,
\begin{align*}
&\P\bigg\{\sup_{d_2(s_1,s_2)\leq\delta}|I_{3,n}(s_1)-I_{3,n}(s_2)|>e\bigg\}\notag\\
&\leq c\,\bigg\{\int_{\overline\eta/2}^\eta\sqrt{\D(\zeta,d_2)}d\zeta+(\delta+2\overline\eta)\,\D(\eta,d_2)\Bigg\}^2+\P\Bigg\{\sup_{\substack{d_2(s_1,s_2)\leq\overline\eta\\s_1,s_2\in\tilde\I}}|I_{3,n}(s_1)-I_{3,n}(s_2)|>e/4\Bigg\}\notag\\
&\leq c\,\bigg\{\int_{\overline\eta/2}^\eta\zeta^{-1/2}d\zeta+(\delta+2\lambda^{(1/m-2\nu)/(2-\vartheta)})\eta^{-1}\bigg\}^2\notag\\
&\qquad+\D(\overline\eta,d_2)\times\sup_{s_1,s_2\in\tilde\I}\P\Big\{ |I_{3,n}(s_1)-I_{3,n}(s_2)|>e/4\Big\}\notag\\
&\leq c\,(\eta+\delta^2\eta^{-2})+c\,\lambda^{-(1/m-2\nu)/(2-\vartheta)}\exp\bigg\{-\frac{e^2/16}{c\lambda^{1/(2m)-\nu}+en^{-1/2}\lambda^{-1/(4m)}\log (n)/6}\bigg\}\notag\\
&\leq c\,(\eta+\delta^2\eta^{-2})+o(1)\,,
\end{align*}
where in the last step we used $\lambda^{-1}\lesssim n^{2m/(2m+1)}$ and $\lambda\lesssim n^{-2m/(2m+1)}$ by Assumptions~\ref{a:rate} and \ref{a:b}, so that $\lambda^{1/(2m)-\nu}=O(n^{(2mb-1)/(2m+1)})$ so that $n^{-1/2}\lambda^{-1/(4m)}\lesssim n^{-m/(2m+1)}$. Therefore, by taking $\eta=\sqrt\delta$, we deduce from the above equation that, when $0\leq\vartheta\leq 1$, for any $e>0$, \eqref{p} holds by taking the metric $d=d_2$. 

As for the remaining terms $I_{1,n}$ and $I_{2,n}$ in \eqref{gn2}, in view of \eqref{new6}, for any $e>0$ and for the metric $d$,
\begin{align*}
&\lim_{\delta\to0}\,\limsup_{n\to\infty}\,\P\bigg\{\sup_{d(s_1,s_2)\leq\delta}\big|\{I_{1,n}(s_1)+I_{1,n}(s_1)\}-\{I_{2,n}(s_2)+I_{2,n}(s_2)\}\big|>e\bigg\}\\
&\leq \limsup_{n\to\infty}\,\P\bigg\{\sup_{s\in\I}|I_{1,n}(s)|+\sup_{s\in\I}|I_{2,n}(s)|>e/2\bigg\}=0\,,
\end{align*}
Combining this result  with \eqref{p} proves that $\G_n$ is asymptotic uniformly equicontinuous w.r.t.~the metric $d$ in \eqref{p} (that is, $d=d_1$ when $\vartheta>1$, and $d=d_2$ when $0\leq\vartheta\leq 1$), which entails the asymptotic tightness of $\G_n$. The proof of the theorem is therefore complete by applying Theorems~1.5.4 and 1.5.7 in \cite{vaart1996}, together with the continuous mapping theorem.

\end{proof}

\subsubsection*{Proof of Theorem~\ref{thm:cb:boot}}

Let ${\rm BL}_1\{\ell^\infty(\I)\}$ denote the collection of all functionals $h:\ell^\infty(\I)\to[-1,1]$ such that $h$ is uniformly Lipschitz: for any $g_1,g_2\in \ell^\infty(\I)$, $|h(g_1)-h(g_2)|\leq \Vert g_1-g_2\Vert_\infty=\sup_{s\in\I}|g_1(s)-g_2(s)|$. We shall show that conditionally on the data $\{(X_i,Y_i,Z_i)\}_{i=1}^n$, the bootstrap estimators $\G_{n,b}^*$ converges to the same limit as $\G_n$ in Proposition~\ref{thm:process}. To achieve this, by Theorem~23.7 in \cite{van2000asymptotic}, we shall prove that, for the Gaussian process $\mathbb G$ in Proposition~\ref{thm:process}, as $n\to\infty$,
\begin{align*}
\sup_{h\in {\rm BL}_1\{\ell^\infty(\I)\}}|\E_M\{h(\G_{n,1}^*)\}-\E\{h(\mathbb G)\}|=o_p(1)\,,
\end{align*}
where $\E_M$ denote the conditional expectation given the data $\{(X_i,Y_i,Z_i)\}_{i=1}^n$.
Let $$\mathbb G_n(s)=\sqrt n\lambda^{1/(4m)}\{\hat g_n(s)+M_\lambda\hat g_n(s)-g_0(s)\},\qquad s\in\I.$$
We apply Lemma~3.1 in \cite{bucher2019note} by proving that, for any fixed $B\geq2$,
\begin{align}\label{gnq}
(\G_n,\G_{n,1}^*,\ldots,\G_{n,B}^*)\weakconverge (\mathbb G,\mathbb G_1,\ldots,\mathbb G_B)\quad\text{ in }\{\ell^\infty(\I)\}^{B+1}\,,\qquad n\to\infty,
\end{align}
where $\mathbb G_1,\ldots,\mathbb G_B$ are i.i.d.~copies of the Gaussian process $\mathbb G$ in Proposition~\ref{thm:process}.

In view of $\ell_n$ in \eqref{elln}, define the multiplier bootstrap sum of squared error
\begin{align*}
\ell_{n,b}^*(g,\beta,\gamma)=\frac{1}{n}\sum_{i=1}^nW_{i,b}\big\{Y_i-g(X_i\trans \beta)-Z_i\trans\gamma\big\}^2,\qquad g\in\H_m,\beta\in\mathcal B,\gamma\in\mathbb R^q.
\end{align*}
Note the fact that $|W_{i,b}|\leq1+\sqrt 2$ almost surely, $\E(W_{i,b})=\var(W_{i,b})=1$.
Following the proof of Theorem~\ref{thm:rate}, we obtain that $\hat w_{n,b}^*=(\hat g_{n,b}^*,Q_{\beta_0}\trans\hat\beta_{n,b}^*,\hat\gamma_{n,b}^*)$ is such that $\| \hat w_{n,b}^*-w_0\|=O_p(r_n)$, for the $r_n$ defined in Theorem~\ref{thm:rate}.

For $1\leq b\leq B$, it follows from similar derivation in Lemma~\ref{lem:n-or} that
\begin{align}
&\hat w_{n,b}^*-w_0-\frac{1}{n}\sum_{i=1}^n\epsilon_iW_{i,b}S_{X_i,Z_i}(g_0,\beta_0)+S_\lambda^\circ(g_0)\notag\\
&\hspace{3cm}=\mathbb V_{1,n,b}^*+\mathbb V_{2,n,b}^*+\mathbb V_{3,n,b}^*-n^{-1/2}\mathbb M_{n,b}^*(\hat w_{n,b}^*-w_0)+\mathcal R_n(\hat g_{n,b}^*,\hat\beta_{n,b}^*,\hat\gamma_{n,b}^*),
\end{align}
where for $S_{x,z}$ and $S_\lambda^\circ$ in \eqref{cor:sxz},
\begin{align*}
&\mathbb V_{1,n,b}^*=\frac{1}{n}\sum_{i=1}^n\e_i W_{i,b}\{S_{X_i,Z_i}(\hat g_{n,b}^*,\hat\beta_{n,b}^*)-S_{X_i,Z_i}(g_0,\beta_0)\},\notag\\
&\mathbb V_{2,n,b}^*=-\frac{1}{n}\sum_{i=1}^nr_{X_i}(\hat g_{n,b}^*,\hat\beta_{n,b}^*) S_{X_i,Z_i}(\hat g_{n,b}^*,\hat\beta_{n,b}^*),\notag\\
&\mathbb V_{3,n,b}^*=\frac{1}{n}\sum_{i=1}^n\l\hat w_{n,b}^*-w_0,S_{X_i,Z_i}(g_0,\beta_0)\r\{S_{X_i,Z_i}(g_0,\beta_0)-S_{X_i,Z_i}(\hat g_{n,b}^*,\hat\beta_{n,b}^*)\},
\\
&\mathcal R_{n,b}^*( \hat g_{n,b}^*,\hat\beta_{n,b}^*,\hat\gamma_{n,b}^*)=\frac{1}{n}\sum_{i=1}^nW_{i,b}\big\{Y_i- \hat g_{n,b}^*(X_i\trans \hat\beta_{n,b}^*)-Z_i\trans\hat\gamma_{n,b}^*\big\}S_{X_i,Z_i}( \hat g_{n,b}^*,\hat\beta_{n,b}^*)+S_\lambda^\circ(\hat g_{n,b}^*).
\end{align*}
Following the proof of Theorem~\ref{thm:bahadur}, we obtain the following joint Bahadur representation for the multiplier bootstrap estimator
\begin{align*}
\bigg\|\hat w_{n,q}^*-w_0-\frac{1}{n}\sum_{i=1}^n\epsilon_i W_{i,b}S_{X_i,Z_i}(g_0,\beta_0)+S_\lambda^\circ(g_0)\bigg\|=o_p(n^{-1/2}).
\end{align*}
As a consequence, following the proof of Theorem~\ref{thm:marginal} and Proposition~\ref{prop:bias}, for the $H_{x,z}$ and $H_\lambda^\circ$ defined in Proposition~\ref{cor:sxz}, respectively, it holds that
\begin{align*}
\sup_{s\in\I}\bigg|\hat g_{n,b}^*(s)-g_0(s)+M_\lambda  \hat g_n(s)- \frac{1}{n}\sum_{i=1}^n\e_iW_{i,b} K(X_i\trans\beta_0,s)\bigg|=o_p(n^{-1/2}\lambda^{-1/(4m)}).
\end{align*}
Now, we have the following decomposition for the $\mathbb G_{n,b}^*$ defined in \eqref{gnqstar}
\begin{align}
\G_{n,b}^*(s)&=\sqrt{n}\lambda^{1/(4m)}\{\hat g_{n,b}^*(s)+M_\lambda \hat g_{n,b}^*(s)-\hat g_n(s)-M_\lambda\hat g_n(s)\}\notag\\
&=I_{1,n,b}^*(s)+I_{2,n,b}^*(s)+I_{3,n,b}^*(s)+I_{4,n,b}^*(s),
\end{align}
where
\begin{align}
&I_{1,n,b}^*(s)=\sqrt{n}\lambda^{1/(4m)}\bigg\{\hat g_{n,b}^*(s)-g_0(s)+M_\lambda \hat g_{n,b}^*(s)-\frac{1}{n}\sum_{i=1}^n\e_iW_{i,b} K(X_i\trans\beta_0,s)\bigg\},\notag\\
&I_{2,n,b}^*(s)=-\sqrt{n}\lambda^{1/(4m)}\Big[\{M_\lambda\hat g_{n,b}^*(s)-M_\lambda g_0(s)\}-\{M_\lambda\hat g_{n}^*(s)-M_\lambda g_0(s)\}\Big],\notag\\
&I_{3,n,b}^*(s)=-\sqrt{n}\lambda^{1/(4m)}\bigg\{\hat g_{n}(s)-g_0(s)+M_\lambda\hat g_n(s)-\frac{1}{n}\sum_{i=1}^n\e_i K(X_i\trans\beta_0,s)\bigg\},\notag\\
&I_{4,n,b}^*(s)=n^{-1/2}\lambda^{1/(4m)}\sum_{i=1}^n\e_i(W_{i,b}-1) K(X_i\trans\beta_0,s).
\end{align}
Following the proof of Proposition~\ref{thm:process}, we deduce by Theorem~\ref{thm:marginal} that
\begin{align*}
\sup_{s\in\I}\{|I_{1,n,b}^*(s)|+|I_{2,n,b}^*(s)|+|I_{3,n,b}^*(s)|\}=o_p(1),\qquad n\to\infty.
\end{align*}
Therefore, in view of the proof of Proposition~\ref{thm:process}, in order to prove \eqref{gnq}, we shall show that, for the $\{I_{4,n}(s)\}_{s\in\I}$ in \eqref{in1in22}, it is true that
\begin{align*}
\mathbb I_n:=(I_{4,n},I_{4,n,1}^*,\ldots,I_{4,n,B}^*)\weakconverge (\mathbb G,\mathbb G_1,\ldots,\mathbb G_B),\qquad\mbox{ in }\{\ell^\infty(\I)\}^{B+1}. 
\end{align*}
Applying Theorems~1.5.4 and 1.5.7 in \cite{vaart1996}, the proof of the above weak convergence therefore relies on the finite-dimensional weak convergence of $\mathbb I_n$ and the asymptotic tightness of  $\mathbb I_n$. Note that the finite-dimensional weak convergence of $\mathbb I_n$ follows similar argument as in the proof of finite-dimensional weak convergence of $I_{4,n}$ in Proposition~\ref{thm:process}. In addition, the asymptotic tightness of $\mathbb I_n$ follows from the proof of the asymptotic tightness of $\{I_{4,n}(s)\}_{s\in\I}$, where we equip the space $\{\ell^\infty(\I)\}^{B+1}$ with the metric $d'\{(s_1,\ldots,s_{B+1}),(t_1,\ldots,t_{B+1})\}=\max_{1\leq b\leq B+1}\{d(s_b,t_b)\}$, for $(s_1,\ldots,s_{B+1}),(t_1,\ldots,t_{B+1})\in\I^{B+1}$, and for $d$ chosen as in the proof of Proposition~\ref{thm:process}. 
Then, the bootstrap statistic $\hat T_{n,b}^*$ defined by \eqref{gnqstar} in Algorithm~\ref{algo:band} satisfies $\hat T_{n,b}^*\converged T$ as $n\to\infty$, conditionally on the data $\{(X_i,Y_i,Z_i)\}_{i=1}^n$, where the random variable  $T$ is defined in \eqref{t}.
Furthermore, \eqref{consetboot} follows by applying Lemma~4.2 in \cite{bucher2019note}. The proof of Theorem~\ref{thm:cb:boot} is therefore complete.

\subsection{Proof of Proposition~\ref{thm:rt:boot}}\label{app:thm:extreme}

We first state and prove the following useful lemma. For $\mathbb G$ defined in Proposition~\ref{thm:process}, define
\begin{align*}
T_{\EE}:= \max\Big \{\sup_{s\in\EE^+}\mathbb G(s),\sup_{s\in\EE^-}\{-\mathbb G(s)\}\Big \}.
\end{align*}

\begin{lemma} \label{cordet}

Suppose Assumptions~\ref{a:rank}--\ref{a:b} and \eqref{newcond} are valid. Then,
\begin{align}\label{eq:rele}
\lim_{n\to\infty}\P\Big\{\hat d_\infty>\Delta+\frac{\mathcal Q_{1-\alpha}(T_{\EE})}{\sqrt{n}\lambda^{1/(4m)}}\Big\}=
\left\{\begin{array}{ll}
\vspace{-0.5em}0 & \quad \text{if}\ d_\infty<\Delta\,\\
\vspace{-0.5em}\alpha &\quad \text{if}\ d_\infty=\Delta\,\\
1 &\quad \text{if}\ d_\infty>\Delta\,
\end{array}\right. ,
%
\end{align}
where $\mathcal Q_{1-\alpha}(T_{\mathcal E})$ denotes the $(1-\alpha)$-quantile of $T_{\mathcal E}$.

\end{lemma}

\begin{proof}[\underline{Proof of Lemma~\ref{cordet}}]

First, by Assumption~\ref{a:5.0}, it is true that $n^{-1/2}\lambda^{-1/(4m)}\log(n\lambda^{1/(2m)})=o(1)$. Then, it follows from Theorem~B.1 in \cite{dette2021bio} and Proposition~\ref{thm:process} that
\begin{align}\label{te}
\sqrt{n}\lambda^{1/(4m)}(\hat d_\infty-d_\infty)\converged T_{\EE}= \max\Big \{\sup_{s\in\EE^+}\mathbb G(s),\sup_{s\in\EE^-}\{-\mathbb G(s)\}\Big \}\, . 
\end{align}
Observe that that when $d_\infty<\Delta$, it follows from \eqref{te} that
\begin{align*}
&\quad\lim_{n\to\infty}\P\bigg\{\hat d_\infty>\Delta+\frac{\mathcal Q_{1-\alpha}(T_\EE)}{\sqrt{n}\lambda^{1/(4m)}}\bigg\}\\
&=\lim_{n\to\infty}\P\Big\{\sqrt{n}\lambda^{1/(4m)}(\hat d_\infty-d_\infty)>\sqrt{n}\lambda^{1/(4m)}(\Delta-d_\infty)+\mathcal Q_{1-\alpha}(T_\EE)\Big\}=0\,,
\end{align*}
since $\sqrt{n}\lambda^{1/(4m)}\to\infty$ as $n\to\infty$, where $T_\EE$ is defined in \eqref{te}. If $d_\infty=\Delta$,
\begin{align*}
&\lim_{n\to\infty}\P\bigg\{\hat d_\infty>\Delta+\frac{\mathcal Q_{1-\alpha}(T_\EE)}{\sqrt{n}\lambda^{1/(4m)}}\bigg\}=\lim_{n\to\infty}\P\Big\{\sqrt{n}\lambda^{1/(4m)}(\hat d_\infty-d_\infty)>\mathcal Q_{1-\alpha}(T_\EE)\Big\}=\alpha\,.
\end{align*}
In the case of $d_\infty>\Delta$,
\begin{align*}
&\quad\lim_{n\to\infty}\P\bigg\{\hat d_\infty>\Delta+\frac{\mathcal Q_{1-\alpha}(T_{\EE})}{\sqrt{n}\lambda^{1/(4m)}}\bigg\}\\
&=\lim_{n\to\infty}\P\Big\{\sqrt{n}\lambda^{1/(4m)}(\hat d_\infty-d_\infty)>\sqrt{n}\lambda^{1/(4m)}(\Delta-d_\infty)+\mathcal Q_{1-\alpha}(T_\EE)\Big\}=1\,.
\end{align*}
The proof of Lemma~\ref{cordet} is therefore complete.

\end{proof}

For the proof of Proposition~\ref{thm:rt:boot},
observe that, by the continuous mapping theorem, Proposition~\ref{thm:process} and Lemma~B.3 in \cite{dette2020functional}, conditional on the data $\{(X_i,Y_i,Z_i)\}_{i=1}^n$, it is true that $\hat T_{\EE,n,b}^*\converged T_\EE$ in \eqref{te}, the same limit as $\sqrt{n}\lambda^{1/(4m)}(\hat d_\infty-d_\infty)$. Hence, the assertion in Proposition~\ref{thm:rt:boot} follows from arguments similar to the ones in the proof of \eqref{eq:rele} in Section~\ref{app:thm:extreme}.

\subsection{Proof of Proposition~\ref{prop:tn}}\label{app:prop:tn}

Under the null hypothesis in \eqref{hp}, it is true that $g_0(x\trans\beta_0)+z\trans\gamma_0-y=0$, which implies that
\begin{align*}
\hat T_n&=\sqrt n\lambda^{1/(4m)}\Big[\big\{\hat g_n(x\trans\beta_0)+M_\lambda g_0(x\trans\beta_0)-g_0(x\trans\beta_0)\big\}+g_0'(x\trans\beta_0)x\trans Q_{\beta_0}Q_{\beta_0}\trans(\hat\beta_n-\beta_0)+z\trans(\hat\gamma_n-\gamma_0)\Big]\\
&\quad+\sqrt n\lambda^{1/(4m)}\Big[\hat g_n(x\trans\hat\beta_n)-\big\{\hat g_n(x\trans\beta_0)+g_0'(x\trans\beta_0)x\trans Q_{\beta_0}Q_{\beta_0}\trans\hat\beta_n\big\}\Big]\\
&\quad+\sqrt n\lambda^{1/(4m)}\Big[M_\lambda\hat g_n(x\trans\hat\beta_n)-M_\lambda g_0(x\trans\beta_0)\Big].
\end{align*}
Define
\begin{align*}
\sigma_n:=\sqrt{\frac{n}{\lambda^{-1/(2m)}\sigma_{(x\trans\beta_0)}^2+v_{x,z}\trans\Omega_{\beta_0}^{-1}v_{x,z}}},
\end{align*}
where $v_{x,z}=\big[g_0'(x\trans\beta_0)x\trans Q_{\beta_0};z\trans\big]\trans\in\mathbb R^{p+q-1}$.
It holds that $\sigma_n=O(\sqrt n\lambda^{1/(4m)})$.
Therefore, we deduce that
\begin{align*}
(\sigma_{(x\trans\beta_0)}^2+\lambda^{1/(2m)}v_{x,z}\trans\Omega_{\beta_0}^{-1}v_{x,z})^{-1/2}\hat T_n=I_{n,1}+I_{n,2}+I_{n,3}+I_{n,4},
\end{align*}
where
\begin{align*}
I_{n,1}&=\sigma_n\Big[\big\{\hat g_n(x\trans\beta_0)+M_\lambda g_0(x\trans\beta_0)-g_0(x\trans\beta_0)\big\}+g_0'(x\trans\beta_0)x\trans Q_{\beta_0}Q_{\beta_0}\trans(\hat\beta_n-\beta_0)+z\trans(\hat\gamma_n-\gamma_0)\Big],
\\
I_{n,2}&=\sigma_n\Big[\hat g_n(x\trans\hat\beta_n)-\big\{\hat g_n(x\trans\beta_0)+g_0'(x\trans\beta_0)x\trans Q_{\beta_0}Q_{\beta_0}\trans\hat\beta_n\big\}\Big],
\\
I_{n,3}&=\sigma_n\{M_\lambda\hat g_n(x\trans\hat\beta_n)-M_\lambda  g_0(x\trans\hat\beta_n)\},
\\
I_{n,4}&=\sigma_n\{M_\lambda g_0(x\trans\hat\beta_n)-M_\lambda g_0(x\trans\beta_0)\}.
\end{align*}

For the first term $I_{n,1}$, by Theorem~\ref{thm:asymp} and the continuous mapping theorem, we obtain that
\begin{align*}
\sigma_n I_{1,n}=\sqrt{\frac{n}{\lambda^{-1/(2m)}\sigma_{(x\trans\beta_0)}^2+v_{x,z}\trans\Omega_{\beta_0}^{-1}v_{x,z}}}I_{1,n}\converged N(0,1).
\end{align*}
For the second term $I_{n,2}$, applying Assumption~\ref{a:reg}, Lemma~\ref{lem:unionbound2} and Theorem~\ref{thm:rate}, we obtain
\begin{align*}
\sup_{x\in\X}|r_x(\hat g_n,\hat\beta_n)|&\leq c\|\hat w_n-w_0\|^{3/2}+c\lambda^{-1/(4m-4)}\|\hat w_n-w_0\|^2
\\
&=O_p(r_n^{3/2}+\lambda^{-1/(4m-4)}r_n^2)=o_p(n^{-1/2}\lambda^{1/(4m)}).
\end{align*}
The above equation implies that
\begin{align*}
|I_{n,2}|\leq O(\sqrt n\lambda^{1/(4m)})\times\sup_{x\in\X}|r_x(\hat g_n,\hat\beta_n)|=o_p(1).
\end{align*}
For the third term $I_{n,3}$, we deduce from Lemma~\ref{prop:bias} that
\begin{align*}
|I_{n,3}|\leq O(\sqrt n\lambda^{1/(4m)})\times\sup_{s\in\I}|M_\lambda\hat g_n(s)-M_\lambda g_0(s)|=o_p(1).
\end{align*}
For the fourth term $I_{n,4}$, it follows from Assumption~\ref{holder} that,  for $s_1,s_2\in\I$,
\begin{align*}
\lambda^{\nu}\sum_{j\geq1}\frac{|\phi_{j}(s_1)-\phi_{j}(s_2)|^2}{(1+\lambda\rho_{j})^2}\leq c|s_1-s_2|^{2\vartheta}.
\end{align*}
Therefore, by the Cauchy-Schwarz inequality, Theorem~\ref{thm:rate} and Proposition~\ref{prop:2.1}, we obtain from the above equation that
\begin{align*}
&|M_\lambda g_0(x\trans\hat\beta_n)-M_\lambda g_0(x\trans\beta_0)|=\lambda\bigg|\sum_{j\geq1}V(g_0,\phi_j)\frac{\rho_j}{(1+\lambda\rho_j)^2}\{\phi_j(x\trans\hat\beta_n)-\phi_j(x\trans\beta_0)\}\bigg|
\\
&\leq\bigg\{\lambda^{2-\nu}\sum_{j\geq1}V^2(g_0,\phi_j)\frac{\rho_j^2}{(1+\lambda\rho_j)^2}\bigg\}^{1/2}\times \bigg\{\lambda^\nu\sum_{j\geq1}\frac{\{\phi_j(x\trans\hat\beta_n)-\phi_j(x\trans\beta_0)\}^2}{(1+\lambda\rho_j)^2}\bigg\}^{1/2}
\\
&\leq\lambda^{(1-\nu)/2}\bigg\{\sum_{j\geq1}V^2(g_0,\phi_j)\rho_j\bigg\}^{1/2} |x\trans(\hat\beta_n-\beta_0)|^\vartheta
\\
&=c\lambda^{(1-\nu)/2}\times \sqrt{J(g_0,g_0)}\times\|\hat\beta_n-\beta_0\|_2^{\vartheta}\leq c\lambda^{(1-\nu)/2}r_n^\vartheta.
\end{align*}
Therefore, we have
\begin{align*}
|I_{n,4}|=O(\sqrt n\lambda^{1/(4m)})\times|M_\lambda g_0(x\trans\hat\beta_n)-M_\lambda g_0(x\trans\beta_0)|\lesssim \sqrt n\lambda^{1/(4m)+(1-\nu)/2}r_n^\vartheta.
\end{align*}
Note that $n\lambda^{1/(4m)+(1-\nu)/2}r_n^\vartheta\lesssim n\lambda^{1/(2m)+1-\nu+\vartheta}\lesssim\lambda^{\vartheta-\nu}=o(1)$ in view of Assumptions~\ref{a:bias} and \ref{holder}, which is guaranteed by the condition that $\nu=1/(2m)$, $\vartheta>1$. This shows that $|I_{n,4}|=o_p(1)$. The proof is therefore complete by applying Slutsky's lemma.

\section{Auxiliary lemmas and their proofs}\label{app:aux}

The following lemma states that  
the mapping $ \beta\mapsto  \theta_\beta:= Q_{\beta_0}\trans\beta$ in \eqref{ad} is one-to-one
 in a neighborhood of $\beta_0$ 
 and  that the orthogonal projection onto ${\rm span}\{\beta_0\}^\perp$ dominates the difference $\beta-\beta_0$.
\begin{lemma}\label{lem:qbeta}
Let $Q_{\beta_0}\in\mathbb R^{p\times(p-1)}$ be such that $[\beta_0,Q_{\beta_0}]\in\mathbb R^{p\times p}$ is an orthogonal matrix. Define $\mathcal B_{p-1}(\beta_0)=\{Q_{\beta_0}\trans\beta:\beta\in\mathcal B,\beta_0\trans\beta>0\}\subset\mathbb R^{p-1}$.
\begin{enumerate}[label={\rm(\arabic*)},series=Qbeta,nolistsep]
%
\item\label{Qbeta1} For any $\beta\in\mathcal B$ such that $\beta_0\trans\beta>0$, it holds that $\beta=Q_{\beta_0}Q_{\beta_0}\trans\beta+\beta_0\sqrt{1-\|Q_{\beta_0}\trans\beta\|_2^2}$. As a consequence, the mapping defined by $\beta\mapsto Q_{\beta_0}\trans\beta$ is a bijection between $\mathcal B$ and $\mathcal B_{p-1}(\beta_0)$.

\item\label{Qbeta2} We have the orthogonal decomposition $\beta-\beta_0=Q_{\beta_0}Q_{\beta_0}\trans\beta+\beta_0(\beta_0\trans\beta-1)$, and
\begin{equation*}
\|\beta_0(\beta_0\trans\beta-1)\|_2=|\beta_0\trans\beta-1|\leq \|Q_{\beta_0}\trans\beta\|_2^2\,.
\end{equation*}
\end{enumerate}

\end{lemma}

\begin{proof}[\underline{Proof of Lemma~\ref{lem:qbeta}}]

For part \ref{Qbeta1},
observe that
\begin{align}\label{id0}
\beta_0\trans \beta_0=1\,,\quad Q_{\beta_0}\trans Q_{\beta_0}=I_{p-1}\,,\quad Q_{\beta_0}Q_{\beta_0}\trans+\beta_0\beta_0\trans =I_p.
\end{align}
For any $\beta\in\mathbb{R}^p$, it holds that $\beta=\beta_0\beta_0\trans\beta+Q_{\beta_0}Q_{\beta_0}\trans\beta$, where $\beta\mapsto \beta_0\beta_0\trans \beta$ and $\beta\mapsto Q_{\beta_0}Q_{\beta_0}\trans\beta$ define the orthogonal projection of $\beta$ onto ${\rm span}\{\beta_0\}$ and ${\rm span}\{\beta_0\}^\perp$, respectively. Hence,
\begin{align*}
1=\|\beta\|_2^2=\|Q_{\beta_0}Q_{\beta_0}\trans\beta\|_2^2+\|\beta_0\beta_0\trans \beta\|_2^2=\|Q_{\beta_0}\trans\beta\|_2^2+(\beta_0\trans\beta)^2.
\end{align*}
When $\beta\in\B$ and $\beta_0\trans\beta>0$, it holds that $\beta_0\trans\beta=\sqrt{1-\|Q_{\beta_0}\trans\beta\|_2^2}.$ This implies that $$\beta=\beta_0(\beta_0\trans\beta)+Q_{\beta_0}Q_{\beta_0}\trans\beta=Q_{\beta_0}Q_{\beta_0}\trans\beta+\beta_0\sqrt{1-\|Q_{\beta_0}\trans\beta\|_2^2}.$$


For part \ref{Qbeta2}, recall that $Q_{\beta_0}\trans\beta_0=0$. For any $\beta\in\B$, we have $\|\beta\|^2_2=\beta\trans\beta=1$, so that in view of \eqref{id0}, we deduce that
\begin{align}\label{betatemp}
\beta-\beta_0&=(Q_{\beta_0}Q_{\beta_0}\trans+\beta_0\beta_0\trans )(\beta-\beta_0)=Q_{\beta_0}Q_{\beta_0}\trans\beta+\beta_0\beta_0\trans (\beta-\beta_0)=Q_{\beta_0}Q_{\beta_0}\trans\beta+\beta_0(\beta_0\trans \beta-1)\notag\\
&=Q_{\beta_0}Q_{\beta_0}\trans\beta+\big\{(\beta\trans\beta_0\beta_0\trans \beta)^{1/2}-1\big\}\beta_0=Q_{\beta_0}Q_{\beta_0}\trans\beta+\big[\big\{\beta\trans(I_p-Q_{\beta_0}Q_{\beta_0}\trans)\beta\big\}^{1/2}-1\big]\beta_0\notag\\
&=Q_{\beta_0}Q_{\beta_0}\trans\beta+\big\{\big(\beta\trans\beta-\beta\trans Q_{\beta_0}Q_{\beta_0}\trans\beta\big)^{1/2}-1\big\}\beta_0\notag\\
&=Q_{\beta_0}Q_{\beta_0}\trans\beta+\big\{\big(1-\|Q_{\beta_0}\trans\beta\|_2^2\big)^{1/2}-1\big\}\beta_0\,.
\end{align}
Observing the fact that $|(1-a)^{1/2}-1|\leq |a|$ for any $a\in\mathbb R$ such that $|a|\leq1$, we deduce that
\begin{align*}
\|\beta_0(\beta_0\trans \beta-1)\|_2&=\|(\beta-\beta_0)-Q_{\beta_0}Q_{\beta_0}\trans\beta\|_2=\big\|\big\{\big(1-\|Q_{\beta_0}\trans\beta\|_2^2\big)^{1/2}-1\big\}\beta_0\big\|_2\\
&=\|\beta_0\|_2\times\big|\big(1-\|Q_{\beta_0}\trans\beta\|_2^2\big)^{1/2}-1\big|=\big|\big(1-\|Q_{\beta_0}\trans\beta\|_2^2\big)^{1/2}-1\big|\leq \|Q_{\beta_0}\trans\beta\|_2^2\,,
\end{align*}
which concludes the proof.

\end{proof}

The following lemma establishes the relation between the $\|\cdot\|_K$-norm and the corresponding $\|\cdot\|$-norm in \eqref{norm}.

\begin{lemma}\label{lem:normK}
Let $(g,\theta,\gamma)\in\Theta$. Then, it holds that $\|g\|_K\leq c\|w\|$, where $c>0$ depends only on $\|g_0'\|_\infty,\|\sigma_0^2\|_\infty,\X,\Z$ and $\Omega_{\beta_0}$ in \eqref{Omega}.
\end{lemma}

\begin{proof}[\underline{Proof of Lemma~\ref{lem:normK}}]
Recall the definition of $\|\cdot\|_K$ in \eqref{V}. Let $D_X=\sup_{x\in\X}\|x\|_2$ and $D_Z=\sup_{z\in\Z}\|z\|_2$. We deduce that
\begin{align*}
\|g\|_K^2&=\E[\sigma_0^2(X\trans\beta_0)\{g(X\trans\beta_0)\}^2]+\lambda J(g,g)\\
&\leq 2\E\big[\sigma_0^2(X\trans\beta_0)\{g(X\trans\beta_0)+g_0'(X\trans\beta_0)X\trans Q_{\beta_0}Q_{\beta_0}\trans\beta+Z\trans\gamma\}^2\big]+\lambda J(g,g)\\
&\quad+2\E\big[\sigma_0^2(X\trans\beta_0)\{g_0'(X\trans\beta_0)X\trans Q_{\beta_0}Q_{\beta_0}\trans\beta+Z\trans\gamma\}^2\big]\\
&\leq 2\|w\|^2+2\|\sigma_0^2\|_\infty\big(\|g_0'\|_\infty D_X^2\|Q_{\beta_0}\trans\beta\|_2^2+D_Z^2\|\gamma\|_2^2\big)\leq c\|w\|^2,
\end{align*}
where we applied Proposition~\ref{prop:2.1} and Assumption~\ref{a:rank2}.

\end{proof}

The following lemma is useful for proving Theorem~\ref{thm:rate}.

\begin{lemma}\label{lem:g'infinity}
Suppose $g\in\H_m$. Then, it holds that
\begin{align*}
&|g'(s_1)-g'(s_2)|\leq c|s_1-s_2|^{1/2}\{\|g\|_\infty^2+J(g,g)\}^{1/2},\qquad s_1,s_2\in\I.
\\
&\|g'\|_\infty^2\leq c\{\|g\|_\infty^2+J(g,g)\},
\end{align*}
where $c>0$ is an absolute constant.
\end{lemma}

\begin{proof}[\underline{Proof of Lemma~\ref{lem:g'infinity}}]

For the first claim, by the Cauchy-Schwarz inequality and the Sobolev interpolation lemma in Lemma~\ref{lem:agmon}, we obtain
\begin{align*}
|g'(s)-g'(s_0)|^2&=\bigg|\int_{s_0}^sg''(t)dt\bigg|^2\leq |s-s_0|\times\int_\I\{g''(t)\}^2dt
\\
&\leq c |s-s_0|\times\bigg[\int_\I\{g(t)\}^2dt+\int_\I\{g^{(m)}(t)\}^2dt\bigg]\leq c|s-s_0|\{\|g\|_\infty^2+J(g,g)\}.
\end{align*}

For the second claim, by Taylor's theorem and the Cauchy-Schwarz inequality, we obtain that for any $s,s_0\in\I$,
\begin{align*}
&\big\{g(s)-g(s_0)-g'(s_0)(s-s_0)\big\}^2=\frac{1}{4}\bigg|\int_{s_0}^s g''(t)(s-t)dt\bigg|^2\leq \frac{1}{4}\bigg|\int_{s_0}^s \{g''(t)\}^2dt\bigg|\times\bigg|\int_{s_0}^s|s-t|^2dt\bigg|
\\
&\qquad\leq \frac{1}{12}\bigg|\int_{\I} \{g''(t)\}^2dt\bigg|\times|s-s_0|^3\leq c\bigg|\int_{\I} \{g''(t)\}^2dt\bigg|\leq c\{\|g\|_\infty^2+J(g,g)\},
\end{align*}
where in the last step we applied the Sobolev interpolation lemma in Lemma~\ref{lem:agmon}. By choosing $s\in\I$ such that $|s-s_0|=|\mathcal I|/2$, where $|\I|$ denotes the length of $\I$, then the above equation yields
\begin{align*}
\{g'(s_0)\}^2\leq |s-s_0|^{-2}\times\Big[4\|g\|_\infty^2+c\{\|g\|_\infty^2+J(g,g)\}\Big]\leq c\{\|g\|_\infty^2+J(g,g)\},
\end{align*}
which completes the proof.

\end{proof}

\begin{lemma}\label{lem:norms}
For $j=0$ or 1, there exists a constant $c_0>1$ such that, for $w=(g,Q_{\beta_0}\trans\beta,\gamma)\in\Theta$, it holds that
\begin{align*}
\sup_{s\in\I}|g^{(j)}(s)|\leq c_0\lambda^{-1/(4m-4j)}\|g\|_K.
\end{align*}
As a consequence, in view of Proposition~\ref{lem:normK}, for $w_1=(g_1,Q_{\beta_0}\trans\beta_1,\gamma_1)$ and $w_2=(g_2,Q_{\beta_0}\trans\beta_2,\gamma_2)$, it holds that, for $j=0$ or 1,
\begin{align}\label{'diff}
\sup_{x\in\X}\big|g_1^{(j)}(x\trans\beta_0)-g_2^{(j)}(x\trans\beta_0)\big|\leq c\,\lambda^{-1/(4m-4j)}\Vert w_1-w_2\Vert.
\end{align}
\end{lemma}

\begin{proof}[\underline{Proof of Lemma~\ref{lem:norms}}]

First, for $j=0$, by the Cauchy-Schwarz inequality and Proposition~\ref{prop:eigen},
\begin{align}\label{suptemp}
\sup_{s\in\I}|g(s)|^2&=\sup_{s\in\I}\big|\l g,K_s\r_K\big|^2\leq \|g\|_K^2\times\sup_{s\in\I}\|K_s\|_K^2\notag
\\
&=\|g\|_K^2\times\sup_{s\in\I}\sum_{j\geq1}\frac{\phi_j^2(s)}{1+\lambda\rho_j}\leq\|g\|_K^2\times\sum_{j\geq1}\frac{\|\phi_j\|_\infty^2}{1+\lambda\rho_j}\leq c\lambda^{-1/(2m)}\|g\|_K^2.
\end{align}

Next, for $j=1$, since $\sigma_0^2(s)f_{X\trans\beta_0}(s)\leq c$ for any $s\in\I$ and for some $c>0$. Hence, applying Lemma~\ref{lem:agmon} and recalling the definition of the operator $V$ in \eqref{V}, for any $e\in(0,1)$, we obtain that
\begin{align*}
&\E\big[\sigma_0^2(X\trans\beta_0)\{g'(X\trans\beta_0)\}^2\big]+\lambda J(g,g)
\\
&= \int_{\I} |g'(s)|^2\sigma_0^2(s) f_{X\trans\beta_0}(s)ds+\lambda \int_\I|g^{(m)}(s)|^2ds
\\
&\leq c\bigg\{e^{-1}\int_{\I} |g(s)|^2\sigma_0^2(s) f_{X\trans\beta_0}(s)ds+e^{2m-1}\int_{\I} |g^{(m)}(s)|^2\sigma_0^2(s) f_{X\trans\beta_0}(s)ds\bigg\}+\lambda \int_\I|g^{(m)}(s)|^2ds\\
&= ce^{-1}{V}(g,g)+ce^{2m-1}\int_{\I} |g^{(m)}(s)|^2\sigma_0^2(s) f_{X\trans\beta_0}(s)ds+\lambda \int_\I|g^{(m)}(s)|^2ds\\
&\leq ce^{-1}{V}(g,g)+ce^{2m-1}\int_{\I}|g^{(m)}(s)|^2ds+\lambda \int_\I|g^{(m)}(s)|^2ds\\
&\leq ce^{-1}{V}(g,g)+ce^{2m-1}J(g,g)+\lambda J(g,g)\leq c(e^{-1}+e^{2m-1}\lambda^{-1}+\lambda)\|g\|_K^2.
\end{align*}
Therefore, by taking $e=\lambda^{1/(2m)}$, we deduce from the above equation that 
\begin{align}\label{new1}
\E\big[\sigma^2_0(X\trans\beta_0)\{g'(X\trans\beta_0)\}^2\big]+\lambda J(g,g)\leq c\lambda^{-1/(2m)}\|g\|_K^2.
\end{align}
Now, let
$\H_{m-1}=\big\{g:\I\to\mathbb{R}\,\big|\,g^{(k)}\text{ is absolutely continuous},\,0\leq k\leq m-1\,;g^{(m)}\in L^2(\I)\big\}$ denote the Sobolev space on $\I$ of order $m-1>1/2$, and denote
\begin{align*}
J_{m-1}(h_1,h_2)=\int_{\I}h_1^{(m-1)}(s)\,h_2^{(m-1)}(s)\, ds.
\end{align*}
Following the construction of the eigen-system in Proposition~\ref{prop:eigen}, let $\{(\eta_j,\psi_j)\}_{j\geq1}$ be the eigenvalue-eigenfunction pairs such that $\sup_{j\geq1}\|\psi_j\|_{\infty}<\infty$, $\eta_j\asymp j^{2m-2}$, and
\begin{align*}
V(h_j,h_{j'})=\delta_{jj'},\qquad J_{m-1}(h_1,h_2)=\eta_j\,\delta_{jj'},
\end{align*}
which correspond to the noramlized solutions of the following ordinary differential equation (ODE):
\begin{align*}
\begin{cases}
(-1)^{m}\psi_{j}^{(2m-2)}(s)=\eta_{j} \,\psi_{j}(s)f_{X\trans\beta_0}(s)\,\sigma_0^2(s),\\
\psi_{j}^{(\theta)}(s)=0\,,\qquad\qquad s\in\partial\I \,,\,m-1\leq\theta\leq 2m-3\,.
\end{cases}
\end{align*}
Then, following Section~\ref{sec:rkhs}, $\H_{m-1}$ is a reproducing kernel Hilbert space equipped with the reproducing kernel $\tilde K$ defined by
\begin{align*}
\tilde K(s,t )=\sum_{j\geq1}\frac{\psi_j(s)\psi_j(t)}{1+\lambda\eta_j},
\end{align*}
with inner product defined by
\begin{equation*}
\l h_1,h_2\r_{\tilde K}=V(h_1,h_2)+\lambda J_{m-1}(h_1,h_2),\qquad h_1,h_2\in\H_{m-1}.
\end{equation*}
Moreover, any $h\in\H_{m-1}$ allows for the series expansion
\begin{align*}
h=\sum_{j\geq1}\frac{\l h,\psi_j\r_{\tilde K}}{1+\lambda\eta_j}\psi_j
\end{align*}
with convergence in $\H_{m-1}$ w.r.t.~the $\|\cdot\|_{\tilde K}$-norm. Therefore, by the same calculations in equation \eqref{suptemp}, we obtain that for any $h\in\H_{m-1}$, it holds that
\begin{align*}
\sup_{s\in\I}|h(s)|^2\leq c\lambda^{-1/(2m-2)}\|h\|_{\tilde K}^2=c\lambda^{-1/(2m-2)}\{V(h,h)+\lambda J_{m-1}(h,h)\}.
\end{align*}
Now, taking $h=g'\in\H_{m-1}$, we obtain from the above equation and equation \eqref{new1} that
\begin{align*}
\sup_{s\in\I}|g'(s)|^2\leq c\lambda^{-1/(2m-2)}\{V(g',g')+\lambda J_{m-1}(g',g')\}= c\lambda^{-1/(2m-2)}\|g\|_K^2,
\end{align*}
which concludes the proof.

\end{proof}

The following lemma states the Sobolev embedding theorem; see, for example, Theorem~3.13 of \cite{oden2012}.
\begin{lemma}[Sobolev embedding theorem]\label{lem:embedding}
Let $\I$ be a bounded domain of $\mathbb R$ and suppose $g\in\H_m(\I)$. Then, there exists a polynomial $g_1(s)=\sum_{\ell=1}^mb_js^{j-1}$ with degree $(m-1)$ such that, for the $J$ defined in \eqref{J}, it holds that
\begin{align*}
\sup_{s\in\I}|g(s)-g_1(s)|^2\leq cJ(g,g),
\end{align*}
where the constant $c>0$ only depends on $\I$.
\end{lemma}

The following lemma is a modified version of the interpolation inequality in Sobolev space; see, for example, Corollary~3.1 in \cite{agmon2010}.

\begin{lemma}[Sobolev interpolation lemma]\label{lem:agmon}
Suppose $g\in\H_m(\I)$ and $\varpi$ is a function on $\I$ such that $c_0\leq \varpi(s)\leq c_0^{-1}$ for some $c_0>0$ and for any $s\in\I$. Then, for any $0<e<1$ and $0\leq k\leq m$, it holds that
\begin{align*}
\int_{\I}|g^{(k)}(s)|^2\varpi(s)ds\leq c\,\bigg\{e^{-k}\int_{\I}|g(s)|^2\varpi(s)ds+e^{2m-k}\int_{\I}|g^{(m)}(s)|^2\varpi(s)ds\bigg\}\,,
\end{align*}
where $c>0$ is a constant that does not depend on $e$.

\end{lemma}

The following lemma is applied to prove Theorem~\ref{thm:rate}; see Lemma~8.4 of \cite{vandegeer2000}.

\begin{lemma} \label{lem:geer}
Let $\mathcal{G}$ be function class and $x_1,\ldots,x_n\in\X$. For $f\in\mathcal G$, define $\|f\|_n^2=n^{-1}\sum_{i=1}^nf(x_i)$. Suppose $\log N_{}(\delta, \mathcal{G},\|\cdot\|_{\infty}) \leq c_1 \delta^{-\alpha}$, $\sup _{f \in \mathcal{G}}\|f\|_n \leq c_2$, and $\epsilon$ is sub-Gaussian satisfying Assumption~\ref{a:02}, for some constants $0<\alpha<2$ and $c_1,c_2>0$. Then for some constant $c>0$, it holds that for all $T \geq c$,
$$
\mathbb{P}\Bigg\{\sup _{f \in \mathcal{G}} \frac{|n^{-1/2} \sum_{i=1}^n \epsilon_i f(x_i)|}{\|f\|_n^{1-\alpha/2}} \geq T\Bigg\} \leq c \exp (-T^2/c^2) .
$$
\end{lemma}

The following lemma is useful for proving Theorem~\ref{thm:rate}; see Lemma 5.16 of \cite{vandegeer2000}.
\begin{lemma}\label{lem:geer2}
 Suppose $\mathcal{G}$ is a class of uniformly bounded functions and for some $0<v<2$,
\begin{align*}
\sup _{\delta>0}\, \delta^v \log N_{[\,]}(\delta, \mathcal{G},\|\cdot\|_{\infty})<\infty .
\end{align*}
Then, for every $c_0>0$, there exists a constant $c>0$ such that
\begin{align*}
\limsup _{n\to\infty}\, \mathbb{P}\Bigg(\sup _{f \in \mathcal{G},\|f\|_*>cn^{-1 /(2+v)}}\bigg|\frac{\|f\|_n}{\|f\|_*}-1\bigg|>c_0\Bigg)=0.
\end{align*}
\end{lemma}

The following lemma is a modified version of Lemma~A.1 in \cite{kley2016quantile}, which we use to prove Proposition~\ref{thm:process}.
\begin{lemma}\label{lem:kley}
Let $\I$ be a compactly supported subset of $\mathbb R$.
For any non-decreasing, convex function $\Psi:\mathbb{R}^+\to\mathbb{R}^+$ with $\Psi(0)=0$ and for any real-valued random variable $Z $, let $\Vert Z\Vert_\Psi=\inf\{c>0:\E\{\Psi(|Z|/c)\}\leq 1\}$ denote the Orlicz norm. Let $\{H(s)\}_{s\in\I}$ be a separable stochastic process with $\Vert H(s_1)-H(s_2)\Vert_\Psi\leq c\,d(s_1,s_2)$ for any $s_1,s_2\in\I$ with $d(s_1,s_2)\geq \overline\eta/2\geq0$ and for some constant $c>0$. Let $\D(\cdot,d)$ denote the packing number of the metric space $(\I,d)$. Then, for any $\delta>0$, $\eta>\overline\eta$, there exists a random variable $S$ and a constant $K>0$ such that
\begin{align*}
&\sup_{d(s_1,s_2)\leq\delta}|H(s_1)-H(s_2)|\leq S+2\sup_{\substack{d(s_1,s_2)\leq\overline\eta\\s_1,s_2\in\I}}|H(s_1)-H(s_2)|,\\
\text{and }~~~~~&\Vert S\Vert_\Psi\leq K\bigg[\int_{\overline\eta/2}^\eta\Psi^{(-1)}\{\D(\e,d)\}d\e+(\delta+2\overline\eta)\,\Psi^{(-1)}\{\D^2(\eta,d)\}\bigg]\,,
\end{align*}
where $\Psi^{(-1)}$ is the inverse function of $\Psi$, and the set $\I$ contains at most $\D(\overline\eta,d)$ points. 

\end{lemma}

\end{document}